\documentclass{amsart}
\usepackage[a4paper, margin=1.5cm]{geometry}
\usepackage{amsmath}
\usepackage{amsfonts}
\usepackage{amssymb}
\usepackage{dsfont}
\usepackage{graphicx}
\usepackage{subcaption}
\usepackage{cite}
\usepackage[hidelinks]{hyperref}
\usepackage{changepage}
\usepackage{amsthm}
\usepackage{xcolor}

\usepackage[inline,shortlabels]{enumitem}
\usepackage{stmaryrd}
\usepackage{bbold}

\newcommand{\ind}{\mathbf{1}}
\newcommand{\ca}{{\mathcal A}}
\newcommand{\cb}{{\mathcal B}}
\newcommand{\cc}{{\mathcal C}}
\newcommand{\cd}{{\mathcal D}}

\newcommand{\ce}{{\mathcal E}}
\newcommand{\cf}{{\mathcal F}}
\newcommand{\cg}{{\mathcal G}}
\newcommand{\ch}{{\mathcal H}}

\newcommand{\cn}{{\mathcal N}}

\newcommand{\cp}{{\mathcal P}}

\newcommand{\E}{{\mathbb E}}
\newcommand{\F}{{\mathbb F}}

\newcommand{\N}{{\mathbb N}}
\renewcommand{\P}{{\mathbb P}}

\newcommand{\R}{{\mathbb R}}

\newcommand{\hk}{\mathfrak{h}}

\newcommand{\dom}{{\mathrm {dom}}\;}
\newcommand{\supp}{{\mathrm {supp}}\;}

\newcommand{\norm}[1]{\mathop{\parallel\! #1 \! \parallel}\nolimits}

\newcommand{\inv}[1]{\mathop{\frac{1}{ #1}}\nolimits}
\newcommand{\expp}[1]{\mathop {\mathrm{e}^{ #1}}}

\newcommand{\Var}{{\rm Var}\;}
\newcommand{\geqcv}{\,\geq _\mathrm{cvx}\,}
\newcommand{\leqcv}{\,\leq _\mathrm{cvx}\,}
\newcommand{\lscv}{\, <_\mathrm{cvx}\,}

\newcommand{\rd}{\mathrm{d}}
\newcommand{\rg}{\mathrm{g}}
\newcommand{\m}{\mathrm{m}}

\newcommand{\cstb}{c_\mathrm{b}}
\newcommand{\cstc}{c_\mathrm{c}}
\newcommand{\cstd}{c_\mathrm{d}}

\theoremstyle{plain}
\newtheorem{theo}{Theorem}[section]
\newtheorem*{theo*}{Theorem}

\newtheorem{cor}[theo]{Corollary}
\newtheorem*{cor*}{Corollary}
\newtheorem{prop}[theo]{Proposition}
\newtheorem{lem}[theo]{Lemma}
\newtheorem{defi}[theo]{Definition}

\newtheorem{theorem}[theo]{Theorem}
\theoremstyle{remark}
\newtheorem{rem}[theo]{Remark}
\newtheorem{exple}[theo]{Example}

\AtEndEnvironment{rem}{\hfill$\Diamond$}
\AtEndEnvironment{exple}{\hfill$\Diamond$}

\usepackage[draft]{fixme}
\FXRegisterAuthor{jf}{ajf}{JFD}

\fxusetheme{colorsig}

\definecolor{orangeTableau}{rgb}{0.75,0.35,0.0}

\usepackage{bbm}
\usepackage{pgfplots}
\usepackage{pgfplotstable}
\pgfplotsset{compat=1.18}
\usepackage{tikz}
\usepgfplotslibrary{patchplots}

\usetikzlibrary{patterns.meta}

\usetikzlibrary{external}
\newcommand{\cov}{{\mathrm{Cov}}}
\newcommand{\diag}{{\mathrm{Diag}}}
\newcommand{\corr}{{\mathrm{corr}}}

\newcommand{\leqlorenz}{{\leq _\mathrm{Lorenz}}}

\newcommand{\cl}{{\mathcal L}}

\usetikzlibrary{calc}

\usepackage{etoolbox}

\usepackage{verbatim}
\newcommand{\indep}{\mathrel{\perp \!\!\! \perp}}

\newcommand{\indicatrice}[1]{\mathop{\mathds{1}_{\{#1\}}}\nolimits}

\newcommand{\im}{\mathrm{im}}

\newcommand{\transpose}{\top}
\newcommand{\rank}{\mathrm{rank}}
\newcommand{\id}{\mathrm{Id}}
\newcommand{\idmatrix}{\mathrm{I}}
\newcommand{\nsamples}{60} 
\newcommand{\Cv}{\mathrm{Cv}}

\newcommand{\csgeq}{\; \succcurlyeq_\mathrm{Cz} \;}
\newcommand{\csleq}{\; \preccurlyeq_\mathrm{Cz} \;}

\newcommand{\csls}{\; \prec_\mathrm{Cz} \;}

\newcommand{\cseq}{\; =_\mathrm{Cz} \;}

\newcommand{\sgn}{\mathrm{sgn}}
\newcommand{\ddom}{{\mathrm{D}}}
\newcommand{\rmf}{{\mathrm{f}}}

\usepackage{xcolor}
\definecolor{BlueBernoulli}{RGB}{0,0,255}

\definecolor{GreenBernoulli}{RGB}{0,255,0}

\definecolor{RedBernoulli}{RGB}{255,0,0}

\definecolor{BackgroundBernoulli}{RGB}{0,0,0}

\newcommand{\OpacityPlot}{0.5}

\usepackage{xparse}
\newcommand{\ColorCm}{green}
\newcommand{\ColorOut}{red}       
\newcommand{\ColorCp}{blue}     
\newcommand{\ArrowStyle}{->}

\NewDocumentEnvironment{RegionPlotBernoulli}{m m m +m +m}{%
  \pgfmathsetmacro{\Pp}{#1}%
  \pgfmathsetmacro{\Rp}{#2}%
  \pgfmathsetmacro{\Qval}{#3}%
  \pgfmathsetmacro{\hPp}{1-\Pp}%
  \pgfmathsetmacro{\hRp}{\Qval-\Rp}%
  \pgfmathsetmacro{\Qbar}{1-\Qval}%
  \pgfmathparse{abs(\Qval-0.5)<0.0001 ? 1 : 0}%
  \ifnum\pgfmathresult=1
    \edef\XTicks{\Qval,1}%
    \def\XTickLabels{$\vphantom{1} q$,$\vphantom{q} 1$}%
  \else
    \edef\XTicks{\Qval,\Qbar,1}%
    \def\XTickLabels{$\vphantom{1} q$,$1-\vphantom{1} q$,$\vphantom{q} 1$}
  \fi
  \begin{tikzpicture}
    \begin{axis}[
        xmin=0, xmax=1+0.05,
        ymin=0, ymax=\Qval+0.05,
        width=6cm,
        height=(\Qval+0.05)/(1+0.05)*6cm,
        scale only axis,
        clip=false,
        xtick/.expand once=\XTicks,
        xticklabels/.expand once=\XTickLabels,
        xticklabel style={anchor=north},
        ytick={\Qval}, yticklabels={$q$},
        axis lines=left,
        axis line style={\ArrowStyle},
        xlabel={$p$},
        ylabel={$r$},
        xlabel style={at={(axis description cs:1,0)},anchor=west},
        ylabel style={at={(axis description cs:0,1)},anchor=south,rotate=270},
    ]
    
    #5 % contenu à symétriser, 1ere copie

    \begin{scope}[shift={(1,\Qval)}, xscale=-1, yscale=-1]
    #5 % 2eme copie, miroir
    \end{scope}

    \fill (1,\Qval) circle (0.01) node[above] {$h(A_0)$};
    \fill (0,0) circle (0.01) node[below left] {$A_0$};

    \addplot[black, dashdotted, forget plot] coordinates {(0,0) (1,\Qval)};
    \addplot[black, dashed, forget plot] coordinates {(0,\Qval) (1,\Qval)};
    \addplot[black, dashed, forget plot] coordinates {(1,0) (1,\Qval)};

    \pgfmathparse{abs(\Qval-0.5)<0.0001 ? 1 : 0}%
    \ifnum\pgfmathresult=1
      \addplot[black, dashed, forget plot] coordinates {(\Qval,0) (\Qval,\Qval)};
    \else
      \addplot[black, dashed, forget plot] coordinates {(\Qval,0) (\Qval,\Qval)};
      \addplot[black, dashed, forget plot] coordinates {(\Qbar,0) (\Qbar,\Qval)};
    \fi

    #4 

  \end{axis}
  \end{tikzpicture}
}{}

\usepackage{makecell}

\title[The Csisz\'ar order]{The Csisz\'ar order for ranking the importance of input random variables}

 \author[C. Butucea]{Cristina Butucea}
 \address{
 Cristina Butucea,
 CREST, ENSAE, 91120 Palaiseau, France.}
 \email{cristina.butucea@ensae.fr}

 \author[J.F. Delmas]{Jean-Fran\c cois Delmas}
 \address{
 Jean-Fran\c cois Delmas,
CERMICS, ENPC, Institut Polytechnique de Paris, CNRS, Marne-la-vall\'ee,
 France.}
 \email{jean-francois.delmas@enpc.fr}

 \author[A. Dutfoy]{Anne Dutfoy}
 \address{
 Anne Dutfoy,
 EDF Research \& Development, Industrial Risk Management Department, 91120 Palaiseau, France.}
 \email{anne.dutfoy@edf.fr}

 \author[A. Schoonaert]{Antoine Schoonaert}
 \address{
 Antoine Schoonaert,
 EDF Research \& Development, Industrial Risk Management Department, 91120 Palaiseau, France.}
 \email{antoine.schoonaert@edf.fr}

\date{\today}

\begin{document}

\begin{abstract}
Csiszár indices $S_f(X,Y)$, defined as the $f$-divergence between the joint distribution of a general random couple $(X,Y)$ and the product of its marginal distributions, are widely used  to quantify dependence between the input variable $X$ and the  output variable $Y$. 
We show that the ranking of two input variables $X$ and $X'$ in terms of their importance w.r.t.\  the output variable, given by the comparison of the  Csiszár indices $S_f(X,Y)$ and $S_f(X',Y)$ depends on the choice of $f$.
We introduce the Csiszár order of input variables as the uniform ordering of Csiszár indices over all proper  convex functions $f$. We give a sufficient condition to prove Csisz\'ar order of two input multivariate random variables.  
We fully characterize it in two cases: when $(X,X',Y)$ is  a Gaussian vector  (in the one-dimensional case this amount to compare the square of their correlation); and when $X$, $X'$ and $Y$ are  Bernoulli random variables. 
Finally, motivated by sensitivity analysis,  we show that, when $(X,X',Y)$ is a three-dimensional Gaussian random vector, thresholding $X$, $X'$, and $Y$ either preserves the Csiszár order or makes the input variables incomparable; in particular  thresholding  never reverses the ranking.

\end{abstract}



\thanks{The research is partially supported by the ANRT fellowship No. 2024/0845. The authors warmly thank Michaël Baudin for numerous discussions and the suggestion to consider the thresholded Gaussian case.}

\maketitle

\section{Introduction}

A fundamental problem is to compare the influence of input variables on a response variable.
In the literature, two main perspectives are commonly considered: variable selection and feature importance in statistical learning on the one hand, and global sensitivity analysis for uncertainty quantification on the other hand. 
Variable selection and feature importance aim at identifying the input variable, say $X$ or $X'$, that best allows one to predict the output variable, say $Y$, from the corresponding input. In contrast, global sensitivity analysis seeks to quantify the variability of the output $Y$ induced by variations in the input $X$.
Despite their different motivations and methodological frameworks, both approaches aim at quantifying the relationship between the inputs and the response.

Global sensitivity analysis typically quantifies the contribution of an input variable to the variability of the output through variance-based decompositions. A prominent example is provided by Sobol indices \cite{sobol1993sensitivity}. In the one-dimensional setting, the Sobol index of $X$ with respect to the output $Y$ is defined by:
\[
\frac{\Var(\E[Y \mid X])}{\Var(Y)} \cdot
\]
This index provides an interpretable decomposition of the total variance. These methods are generally developed under the assumption of independent inputs and focus on second-order moments. More general measures of the strength of dependence can be found, such as the Chaterjee's rank correlation in \cite{chatterjee2021new} and the conditional convex order \cite{ansari2025ordering}.

A different line of work quantifies the relationship between an input and the output through statistical measures of dependence. This includes entropy-based indices \cite{auder2008global}, %kala2021global
optimal transport approaches \cite{ansari2026dependence}. 
Unlike variance-based methods, these approaches exploit the joint distribution of the input and the output and therefore capture more general forms of dependence.

In parallel, machine learning methods define variable importance through algorithmic criteria. Random forests provide importance scores based on impurity reduction or permutation schemes \cite{breiman2001random}, some of which have been shown to converge asymptotically to Sobol-type quantities \cite{benard2022mean}. Other approaches, such as the Lasso, perform variable selection implicitly through an $\ell_1$ regularization. Finally, linear correlation remains one of the most widely used measures of association in practice because of its simplicity and interpretability, despite being restricted to linear relationships.

Among distributional dependence measures, Csisz\'ar indices \cite{da2015global,rahman2016f} provide a unified framework based on $f$-divergences. Given two random variables $X$ and $Y$, their Csisz\'ar index is defined by:
\[
S_f(X,Y)=D_f(P_{(X,Y)} \Vert P_X \otimes P_Y),
\]
where $D_f$ denotes an $f$-divergence. 
The general definition of $D_f$ and $S_f$ are recalled in Section~\ref{sec:Definitions}. These indices use the joint distribution of $(X,Y)$ and not only the second order moments, apply to arbitrary dimensions and do not require to have independent inputs. These indices can also be written in terms of the conditional distribution of $Y$ given $X$ and the marginal distribution of $Y$, see Remark~\ref{rem:conditional}.

The Csisz\'ar indices characterize independence: if $f$ is strictly convex at $1$ with $f(1) = 0$, then, according to~\cite[Corollary~7]{butucea2026csisz}:
\begin{equation}
    \label{eq:Sf-indep}
S_f(X,Y)=0
%\quad\text{if and only if}\quad
\quad \Longleftrightarrow\quad
%X 
\text{$X$ and $Y$ are independent.}
%\indep Y,
\end{equation}
This makes $f$-divergence indices particularly attractive for independence testing: 
for the mutual information 
Csisz\'ar index $S_{\mathrm{KL}}$ 
 with $f_\mathrm{KL}(t) = t \log t$ 
see~\cite{berrett2019nonparametric, gretton2010consistent, ai2024testing}; 
for the Neyman's $\chi^2$ Csisz\'ar index $S_{\mathrm{N}}$ with $f_\mathrm{N}(t) = (1-t^2)/t$ see~\cite{berrett2021usp} in a discrete setting; 
for the Hellinger Csisz\'ar index $S_{\mathrm{H}}$ with $f_\mathrm{H}(t) = (\sqrt{t}-1)^2$ see~\cite{su2008nonparametric}; 
for the Total variation Csisz\'ar index $S_{\mathrm{TV}}$ with $f_\mathrm{TV}(t) = \vert t -1\vert$ see~\cite{gretton2010consistent, diakonikolas2016new}. 
See also~\cite{butucea2026csisz} %bouzebda2009new
 for an approach based on copulas and~\cite{neykov2021minimax} for testing conditional independence. 

While these measures are equivalent for testing independence, they do not quantify dependence in the same way for various choices of $f$. For ranking or sensitivity analysis based on 
Csisz\'ar indices,  
see:
\cite{borgonovo2007new} for the total variation distance;
\cite{suzuki2009mutual} for the Neyman $\chi^2$ divergence, 
\cite{cover2006element} for the Kullback-Leibler divergence, and~\cite{geenens2022hellinger} for the Hellinger distance. 
However, the influence of the choice of $f$ on the particular ranking is typically not discussed.

\bigskip 

In this paper, we analyse the ranking of multivariate input random variables $X$ and $X'$ with respect to (w.r.t.)\ a multivariate  output random variable $Y$ using Csisz\'ar indices. Very surprisingly, the next example shows that these indices are not always ordered in the same way as the correlations.

\begin{exple}
[Correlation vs.\ Csisz\'ar indices]
\label{exple:OrdreUniformeVsCorr}
Let $X$ be uniformly distributed on $\{1,2,3\}$ and $X'$ be uniformly distributed on $\{0,4\}$. Consider the additive model $Y=X+X'$, so that $Y$ is uniformly distributed on $\{1,2,3,5,6,7\}$. Then, we have:
\[
\corr(X,Y)^2 < \corr(X',Y)^2,
\qquad
\Var(X) < \Var(X')
\quad
\text{while}
\quad
S_f(X,Y) \geq S_f(X',Y)
\quad \text{for all $f$},
\]
where $f$ is a proper convex function, and with the last inequality being strict if $f$ is strictly convex at some point in $(0,3)$ and $f(0)+f(3)<\infty$.
Thus, the ordering induced by all Csisz\'ar indices may be opposite for all choices of $f$ to the rankings suggested by variance and correlation. 
\end{exple}

More importantly, even within the class of $f$-divergences, the ranking of variables may depend on the particular choice of the convex function $f$.

\begin{exple}
[Csisz\'ar indices may not give the same ranking]
\label{exple:ContreExempleBernoulli}
Let $X$, $X'$, and $Y$ be Bernoulli random variables satisfying $\P(X=1)=\P(Y=1)=3/5$, $\P(X'=1)=1/2$, $\P(X=1,Y=1)=7/25$ and $\P(X'=1,Y=1)=11/50$. Then, we have:
\[
S_{\mathrm{TV}}(X,Y)<S_{\mathrm{TV}}(X',Y)
\quad
\text{while}
\quad
S_{\mathrm{H}}(X,Y)>S_{\mathrm{H}}(X',Y),
\]
where $S_{\mathrm{TV}}$ and $S_{\mathrm H}$ denote the total variation and Hellinger Csisz\'ar indices, respectively. Therefore, two Csisz\'ar indices may lead to opposite rankings of the same pair of variables. 
\end{exple}

This difference in rankings shows that, in the absence of a guarantee of uniform ordering, the choice of the function $f$ is crucial. 
The large variety of available $f$-divergences naturally raises the following question: 

\textbf{Given two input variables, can one determine which one  is the most informative w.r.t.\ the output variable, simultaneously for all Csisz\'ar indices?} 

\medskip
In order to answer this question, we introduce  the  Csisz\'ar order  as  the  ordering of the Csisz\'ar indices uniformly over all proper convex functions $f$. 
We analyse this preorder, and give sufficient conditions to ensure the comparison of multivariate input variables; and then we investigate thoroughly  the Gaussian and the Bernoulli models.

\subsection{Definitions}
\label{sec:Definitions}

We first recall the definition of the $f$-divergence. 
Let $\F$ be the collection of convex functions $f$ on $\R$ taking values in $(-\infty, +\infty]$ such that its domain $\dom(f)=\{f < +\infty\}$ is in $\R_+$, $f(1)=0$ (so that $f$ is a proper convex function), $1$ belongs to the interior of $\dom(f)$ and $f$ is right continuous at 0 in $(-\infty, +\infty]$. 
For $f\in \F$, we call $f^*(t) = t f\left( 1/t \right) $ for $t >0$ the convex inversion; it  belongs also to $\F$.  

If $P$ and $Q$ are probability measures on the same  measurable space $(\Omega, \mathcal{F})$, we define $p/q$ $\lambda$-a.e.\ with $\lambda = P + Q$. The ratio $p/q$ is understood as the Radon-Nikodym derivative $p = \rd P/ \rd \lambda$ divided by the Radon-Nikodym derivative $q = \rd Q/\rd \lambda$, see~\cite{butucea2026csisz}  for more details. For $f\in \F$, we define the $f$-divergence between $P$ and $Q$ by:
\begin{equation}
  \label{eq:definition_df}
  D_f(P\Vert Q) = \int_\Omega f\left(\frac{p}{q} \right) \, q \, \rd \lambda + f^*(0) \cdot P \left(q = 0 \right),
\end{equation}
where the last term above is 0 when $P$ is dominated by $Q$. 
Notice that:
\begin{equation}
   \label{eq:decom_Df}
  D_f(P\Vert Q)
    = \int_{0\leq p<q} f\left(\frac{p}{q}\right) \, q  \, \rd \lambda + \int_{0\leq  q<p} f^*\left(\frac{q}{p}\right)\, p \, \rd \lambda
    =D_{f^*}(Q\Vert P).
\end{equation}

Eventually, for a random variable $Z$ defined on $\Omega$, we denote by $P_Z$ its probability distribution defined on the measurable space $(\Omega_Z, \cf_Z)$, where $Z$ takes its values. If $P_Z$ is dominated by a measure, say $\lambda$, which is denoted by $P_Z\ll \lambda $, 
we simply write  $p_Z$ for its corresponding  density: $P_Z(\rd z)= p_Z(z)\, \lambda(\rd z)$.

Let $X$ and $Y$ be two random variables defined on $\Omega$ taking values in possibly different measurable spaces. Following \cite{liese2006divergences,polyanskiy2025information}, we define the Csisz\'ar index between $X$ and $Y$ as the $f$-divergence between the joint distribution $P_{(X,Y)}$ and the product distribution $P_X \otimes P_Y$:
\[
S_f(X,Y) = D_f(P_{(X,Y)} \Vert P_X \otimes P_Y).
\]
Let us mention that \cite{da2015global, rahman2016f} use another convention for the Csisz\'ar index given therein by $D_f(P_X \otimes P_Y \Vert P_{(X,Y)})$ which is equal to $D_{f^*}(P_{(X,Y)} \Vert P_X \otimes P_Y)$ thanks to~\eqref{eq:decom_Df}. Notice  the convention for the Csisz\'ar index
 has no impact on results which hold for all $f \in \F$ which is equivalent to $f^* \in \F$.

\begin{rem}[Using conditional distribution]   \label{rem:conditional}
    Let $P_{X|Y=y}$ denote the conditional distribution of $X$ w.r.t.\ $Y$ when it exists. In this case, the function 
    $D_f(P_{X\mid Y} \Vert P_X)$ is a random variable  and, according to~\cite{butucea2026csisz}, we have:
    \[
      S_f(X,Y) =  \E[D_f(P_{X\mid Y} \Vert P_X)] = \int D_f(P_{X\mid Y=y} \Vert P_X) \, P_Y(\rd y).
    \]
Similarly, we can define the 
Csisz\'ar index between $X$ and a sub-$\sigma$-field $\cg$, provided the  conditional distribution $P_{X\mid \cg}$ of $X$ w.r.t.\ $\cg$ exists, by $
 S_f(X,\cg) =   \E[D_f(P_{X\mid \cg} \Vert P_X)]$.
  Since by definition $P_{X\mid Y}=P_{X\mid \sigma(Y)}$, with $\sigma(Y)$ the $\sigma$-field generated by $Y$, we have $ S_f(X,Y) =S_f(X,\sigma(Y)) $. Then, we immediately deduce that, provided $P_{X|Y}$ exists, for any random variable $Z$ generating the same $\sigma$-algebra as $Y$, then the Csisz\'ar indices are equal, that is, $\sigma(Y)=\sigma(Z)$ implies $S_f(X,Y)=S_f(X, Z)$ for all $f\in \F$.
\end{rem}

The Csisz\'ar index $S_f$ has several useful properties,  see~\cite{butucea2026csisz}: it is symmetric, nonnegative, bounded from above by $f(0)+ f^*(0)$, and  vanishes whenever $X$ and $Y$ are independent. Moreover, if $f$ is strictly convex at $1$, then $S_f(X,Y)=0$ if and only if $X$ and $Y$ are independent.

\begin{defi}[Csisz\'ar order]
    \label{def:Cs order}
For three random variables $X$, $X'$, and $Y$, we say that $X$ is smaller or equal to $X'$ in the Csisz\'ar order w.r.t.\ $Y$ if $S_f(X,Y) \leq S_f(X',Y)$ for all $f \in \F$, and denote this by: 
\[
X \csleq X' \quad \text{w.r.t.}\quad Y.
\]
We say that $X$ and $X'$ are equal in the Csisz\'ar order, denoted by $X \cseq X'$ w.r.t.\ $Y$, if both $X \csleq X'$ and $X' \csleq X$ w.r.t.\ $Y$. 
\end{defi}

Notice that  the random variables $X$ and $X'$ do not take values in the same measurable space \emph{a priori}.
We  deduce from Remark~\ref{rem:conditional}, that the  Csisz\'ar  order defined w.r.t.\ $Y$ coincides with the  Csisz\'ar  order defined w.r.t.\ $Y'$ if $\sigma(Y)=\sigma(Y')$. Recall a function is bi-measurable if it is measurable and if the image of any measurable set is a measurable set. 
According to Remark~\ref{rem:conditional}, if $\varphi$ is a bi-measurable bijection, then we have $X \cseq \varphi(X)$. For example, we have $U \cseq \log(U)$ w.r.t.\ $Y$, for $U$ a uniform random variable on $[0,1]$.

\medskip

To simplify, we shall assume from now on that the random variables $X$ and $X'$  take values in Polish spaces, so that  the conditional distributions  $P_{X \mid Y}$ and $P_{X' \mid Y}$ are well defined. 
Since the relation $ \csleq$ is reflexive and transitive on the conditional distributions w.r.t.\ $Y$, it is a preorder; however it is not a partial order   because, if $X \csleq X'$ and $X' \csleq X$ w.r.t.\ $Y$ it does not necessarily imply that $P_{X \mid Y} =P_{X' \mid Y}$. 
It may also happen for a triplet $(X,X',Y)$ that neither $X \csleq X'$ nor $X \csgeq X'$ w.r.t.\ $Y$ holds; in this case, we say that $X$ and $X'$ are not ordered in the Csisz\'ar order w.r.t.\ $Y$. 
Finally, we say that $X$ is smaller than $X'$ in the Csisz\'ar order (denoted by $X \csls X'$ w.r.t.\ $Y$) if $X \csleq X'$ w.r.t.\ $Y$ and $X$ and $X'$ are not equal in the Csisz\'ar order w.r.t.\ $Y$.

\begin{defi}[Csisz\'ar sets]
    \label{def:Csiszar sets}
For a fixed pair $(X',Y)$, we define the upper Csisz\'ar set of $X'$ w.r.t.\ $Y$, denoted $\cp_Y^+(X') $,  as the set of conditional distributions $P_{X \mid Y}$ such that $X$ is larger than or equal to $X'$ in the Csisz\'ar order w.r.t.\ $Y$, and analogously for  the lower Csisz\'ar set denoted $\cp_Y^-(X') $:
\begin{equation}
\label{eq:defP+P-}
\cp_Y^+(X') = \left\{ P_{X \mid Y} \, \colon \, X \csgeq X'  \text{ w.r.t.\ } Y\right\} 
\quad \text{and} \quad 
\cp_Y^-(X') = \left\{ P_{X \mid Y} \, \colon \, X \csleq X'  \text{ w.r.t.\ } Y\right\}.
\end{equation}
We also define the set of variables that are equal to $X'$ in the Csisz\'ar order w.r.t.\ $Y$:
\begin{equation}
\label{eq:defP=}
\cp_Y(X') = \left\{ P_{X \mid Y} \, \colon \, X \cseq X'  \text{ w.r.t.\ } Y\right\} = \cp_Y^+(X') \cap \cp_Y^-(X').
\end{equation}
\end{defi}
Our goal is to characterize these sets for specific families of probabilistic models. 
Note that these sets are nonempty as the Csisz\'ar order is reflexive. 
Moreover, the union $\cp_Y^+(X') \cup \cp_Y^-(X')$ does not exhaust the set of all possible conditional distributions $P_{X\mid Y}$ in general. Indeed, as illustrated by Example~\ref{exple:OrdreUniformeVsCorr} and Example~\ref{exple:ContreExempleBernoulli}, there may exist a conditional distribution $P_{X\mid Y}$ such that $S_f(X,Y)<S_f(X',Y)$ and $S_{f'}(X,Y)>S_{f'}(X',Y)$ for two distinct functions $f,f'\in\F$, that is, $P_{X\mid Y}$ belongs to neither $\cp_Y^+(X')$ nor $\cp_Y^-(X')$.
Among the conditional distributions $P_{X\mid Y}$, the minimal element corresponds to $X$ being independent of $Y$, that is,  $P_X \in \cp_Y^-(X')$.

\subsection{Content and structure of the paper}
The remainder of the paper is organized as follows. In Section~\ref{sec:characterization}, we first establish a general characterization of the Csisz\'ar order in terms of the convex order of the likelihood ratio between the joint distribution and the product of its marginals. This result provides a general criterion for identifying models that admit Csisz\'ar order. Next, we discuss more general properties that a global dependence order is often expected to satisfy.

In Section~\ref{sec:Gaussian}, we study the Gaussian setting. We show that the Csisz\'ar index of a Gaussian vector $(X,Y)$ is a non-decreasing function of the eigenvalues of the correlation matrix of $(X,Y)$. 
This result yields sufficient conditions on the Gaussian vectors $(X,Y)$ and $(X',Y)$ in order to rank them using the  Csisz\'ar order. When $Y$ is one-dimensional or when $X$ and $X'$ are one-dimensional, these conditions are shown to be necessary, leading to a complete characterization of the  Csisz\'ar order.

In Section~\ref{sec:AdditiveModel}, we  consider a more general model. 
In particular, for the additive model $Y = X + X'$, where $X$ and $X'$ are independent Cauchy random variables with scale parameters $\alpha$ and $\alpha'$, respectively, we have $X \csleq X'$ w.r.t.\ $Y$ if and only if $\alpha \leq \alpha'$. This generalization extends to all Lévy processes.

In Section~\ref{sec:Bernoulli}, we turn to discrete models and provide a complete characterization of the Csisz\'ar order for  Bernoulli random variables in terms of the model parameters. In particular,  we  show that given $X'$ and $Y$, one can compare all $X$ to $X'$ w.r.t.\ $Y$ if and only if  $(X',Y)$ is degenerate. Intuitively, this means that being able to compare given  input variables is not a straightforward  result.  

In sensitivity analysis it may be important to analyse how the large values of the inputs $X$ or $X'$ induce changes in the high values of $Y$. In Section~\ref{sec:ThresholdedGaussian}, 
we investigate thresholded Gaussian models: for  $(X,X',Y)$  a three-dimensional Gaussian  vector, we consider the Bernoulli variables obtained by thresholding the components of $(X,X',Y)$. We prove that the choice of the threshold cannot alter the ranking of $X$ and $X'$ given by the Csisz\'ar order.  More precisely,  the thresholded variables are ranked in the same order as the initial Gaussian variables or they are not comparable:  thresholding Gaussian random variables never reverses the ranking.

\section{Characterization and properties of the Csisz\'ar order} \label{sec:characterization}

In this section we first give a characterization of the Csisz\'ar order in terms of the density ratio. Next we discuss the properties that are generally required for a global dependence ordering as sumed up by \cite{ansari2025ordering}.

\subsection{Characterization of the Csisz\'ar order}
\label{sec:ConnectionConvexOrder}

Let $Z$ and $Z'$ be two real  random variables. We say that $Z$ is  smaller than $Z'$ in the convex order, denoted by $Z \leqcv Z'$, if for all proper convex measurable functions $\varphi$ defined on $\R$, such that $ \varphi(Z)$ and $\varphi(Z')$ are quasi-integrable, we have:
\begin{equation}
    \label{eq:ineg-cvx}
    \E\left[ \varphi(Z)\right] \leq \E\left[ \varphi(Z')\right].
\end{equation}
We also
say that $Z$ is strictly smaller that $Z'$ in the convex order, denoted by  $Z \lscv Z'$, if 
$Z \leqcv Z'$ and  the inequality~\eqref{eq:ineg-cvx} is strict 
 for at least one proper convex function $\varphi$. 
Notice that if $Z$ and $Z'$ are quasi-integrable, then $Z \leqcv Z'$ and  $Z' \leqcv Z$ is equivalent to  $Z$ and $Z'$ having the same  distribution (see \cite[Theorem~3.A.42]{shaked2007stochastic}). We recall that a real valued random variable $Z$ is quasi-integrable if either its positive part $Z_+$ or its negative part $Z_-$ are integrable.

Let $(X,Y)$ be a random variable defined on $\Omega$ taking values in a product space $\Omega_X \times \Omega_Y$
endowed with the product $\sigma$-field $\cf_X\otimes \cf_Y$. 
According to the elementary Lemma~\ref{lem:abs-cont-XY},  if $P_{(X,Y)}$ is dominated by  a  product  of $\sigma$-finite measures, say $\lambda=\lambda_X\otimes \lambda_Y$,  with density $p_{(X,Y)}$, if and only if it is also dominated by  $P_X\otimes P_Y$ with corresponding Radon-Nikodym derivative given by the density ratio: 
\begin{equation}
    \label{eq:R}
R_{(X,Y)}(x,y)=\frac{p_{(X,Y)}(x,y)}{p_X(x)\, p_Y(y)},
\end{equation}
where $p_{(X,Y)}= \rd P_{(X,Y)}/ \rd (\lambda_X \otimes \lambda_Y)$, $p_X=\rd P_X/\rd \lambda_X$ and $p_Y=\rd P_Y/\rd \lambda_Y$.

Furthermore, in this case the expression of $S_f(X,Y)$ % in~\eqref{eq:def_S_f} 
reduces to:
%the term with $f^*(0)$ in~\eqref{eq:definition_df} (with $f$ replaced by $f^*$) which appears also in~\eqref{eq:def_S_f} has no contribution, and thus: 
\begin{equation}
  \label{eq:Sf_et_Z}
  S_f(X,Y) = \E \left[ f\left( Z \right) \right], 
  \quad \text{with} \quad Z = R_{(X,Y)}(X,Y)  \quad \text{and} \quad (X,Y) \sim P_X \otimes P_Y.
\end{equation}
Let $X'$ be another random variable defined on $\Omega$ taking values in a measurable space $(\Omega_{X'}, \cf_{X'})$. 
We assume similarly that  $P_{(X',Y)} $ is dominated by a product of $\sigma$-finite measures, and we define the random variable $Z'$ as $Z$ in~\eqref{eq:Sf_et_Z} with $X$ replaced by $X'$. The following proposition relates the convex order between $Z$ and $Z'$ to the Csisz\'ar order between $X$ and $X'$ w.r.t.\ $Y$;  it will be used in particular in  Section~\ref{sec:Bernoulli} on Bernoulli random variables.  The proof is given in Section~\ref{sec:ProofsConnectionConvexOrder}.

\begin{prop}[Characterization of the Csisz\'ar order]
\label{prop:equivalence_ConvexOrder_CsiszarOrder}
Let $(X,X',Y)$ be a multivariate random variable taking values in a product of measurable spaces  such that each probability measure $P_{(X,Y)}$ and $P_{(X',Y)}$ is 
dominated by  the product of two $\sigma$-finite measures. 
Then, we have: 
\begin{equation}
  \label{eq:equivalence_ConvexOrder_CsiszarOrder}
 Z \leqcv Z'
% \quad  \text{if and only if} \quad 
 \quad  \Longleftrightarrow \quad 
 X \csleq X' \quad \text{w.r.t.}\quad Y  .
\end{equation}
\end{prop}

It is obvious that any input variable $X$ independent of $Y$ is minimal for the Csiszár order; the existence of a maximal input variable $X'$ with a given distribution has been considered in~\cite{scarsini1990ordering} and solved in the Bernoulli case, see also Remark~\ref{rem:ScarsiniBernoulli}.

The conditional convex order introduced by~\cite{ansari2025ordering}  
for one-dimensional output random variables is different in general than the Csiszár order, see 
Remarks~\ref{rem:Connection_Ccvx_Bernoulli} for the Bernoulli model; however when the output variables in~\cite{ansari2025ordering}  are equal these two notions coincide in the one-dimensional Gaussian model, see Remark~\ref{rem:Connection_Ccvx_Gaussian}.

\subsection{Csisz\'ar order as a global dependence ordering}

The following axioms are considered desirable properties for a global dependence order between two pairs of input and output random variables, as formalized by \cite{ansari2025ordering}. We analyze them for the particular case of the  Csiszár order introduced here. 

We shall say that the  input/output $(X,Y)$ has  less dependence than the input/output $(X', Y')$ if $S_f(X,Y) \le S_f(X', Y')$ for all $f\in \F$, and then simply write $(X,Y) \csleq (X',Y')$. This is consistent with the  Csiszár order in the sense that 
$ X \cseq X'$ w.r.t.\ $Y$ is equivalent to $(X,Y) \csleq (X',Y)$, where the outputs are the same. For the sake of comparison with~\cite{ansari2025ordering}, we consider only the  Csiszár order on the couple input/output instead of the order on the input w.r.t. a given output, the results on the latter are similar and left to the reader. 

The  Csiszár order on the couple input/output (or simply on the input w.r.t. a given output) satisfies the axioms of~\cite{ansari2025ordering} except the characterization of perfect dependence (iv). However, our order is valid for random variables of arbitrary dimensions whereas \cite{ansari2025ordering} considers that the output $Y$ is real-valued and the input is $\R^d$-valued.  Indeed, we have:

\begin{enumerate}[(i)]
  \item \textbf{Law invariance:} \label{axiom:LawInvariance}
  If $(X,Y) \stackrel{\cl}{=} (X',Y')$, then $ (X,Y) \cseq (X',Y')$.  
  \item \textbf{Quasi-ordering:} The relation $\csleq$ is reflexive and transitive.
  \item \textbf{Characterization of independence:}  $(X,Y) \csleq (X',Y')$  for all $(X',Y')$ if and only if $X$ is independent of $Y$ (see~\eqref{eq:Sf-indep}).
  \item \textbf{Characterization of perfect dependence} fails as  $(X',Y') \csleq (X,Y)$ for all $(X',Y')$   
  does not imply that  $Y$ is perfectly dependent of $X$, that is, $Y$ is $\sigma(X)$-measurable. 

Thanks to \cite[Lemma~1~(iv)]{butucea2026csisz},
we get that $(X,Y)$ is maximal, that is, $(X',Y') \csleq (X,Y)$ for all $(X',Y')$, if and only if the distributions $P_{(X,Y)}$ and $P_X \otimes P_Y$ are mutually singular. 
In particular if $Y$ is perfectly dependent of $X$, then $(X,Y)$ is maximal provided that the distribution of $X$ or of $Y$ has no atom.  

  The couple 
  $(X,Y)$ can be maximal without perfect dependence between $X$ and $Y$: consider  $Y = X + Z$ with $X$ a standard Gaussian random variable and $Z$ a non-degenerate Bernoulli random variable independent of $X$, the distributions $P_{(X,Y)}$ and $P_X \otimes P_Y$ are singular even though 
   $Y$ is not perfectly dependent of $X$.
   
  \item \textbf{Information monotonicity:} \label{axiom:InformationMonotonicity}  $(X,Y) \csleq ((X,Z), Y)$ (that is, $X \csleq (X,Z)$ w.r.t.\ $Y$) for any random variable $Z$ (see~\cite[Corollary~8]{butucea2026csisz}).
  \item \textbf{Characterization of conditional independence:} \label{axiom:CharacterizationConditionalIndependence} $(X,Y) \cseq ((X,Z), Y)$ (that is, $X \cseq (X,Z)$ w.r.t.\ $Y$)  if and only if $Y$ and $Z$ are conditionally independent given $X$ (see~\cite[Theorem~4]{butucea2026csisz}).
  \item \textbf{Transformation invariance:} 
  \label{axiom:transformation_invariance}
   $ (X,Y) \cseq (h(X), g(Y)) $ for all bi-measurable injections $h$ and $g$ (see Proposition~\ref{prop:invariance_csiszar_index}).
  Notice that $(h,g)$ can be replaced by any measurable function $(h,g)$ such that 
$(\rd P_{(X,Y)}/\rd P_X \otimes \rd P_Y)(x,y)$ is a measurable function of $(h(x),g(y))$. This result is stronger than in~\cite{ansari2025ordering} where $Y$ was assumed to be real and $g$  increasing. 
\end{enumerate}
To introduce the next axiom, we shall denote, for $Z=(Z_i)_{1\le i\le d}$ a random variable in $\R^d$, its distribution function  $F_Z=(F_{Z_i})
_{1\le i\le d}$, where $F_{Z_i}$ is the cdf of $Z_i$, and set $F^{-1}_Z=(F^{-1}_{Z_i})_{1\le i\le d}$, where $F^{-1}_{Z_i}$ is the generalized right-continuous inverse of $F_{Z_i}$, as well as 
  $F^{-1}_Z\circ F_Z=(F^{-1}_{Z_i}\circ F_{Z_i})_{1\le i\le d}$.
  
\begin{enumerate}[(i), resume]
  \item \textbf{Distributional invariance:} 
  For $(X,Y)$ taking values in  $\R^{d_X}\times \R^{d_Y}$, we have 
$ (X,Y) \cseq (F_X (X), Y)\cseq (X, F_Y(Y))\cseq (F_X (X), F_Y(Y))
$.  (Indeed, we have $X=F_X^{-1} \circ F_X(X)$ a.s.\ and thus using twice the contraction of the Csisz\'ar index by functional transformation, see~\cite[Corollary~8]{butucea2026csisz}, 
  we get that $S_f(X, Y) = S_f(F_X^{-1}\circ F_X(X), Y) \leq S_f(F_X(X), Y) \leq S_f(X, Y)$ and thus $S_f(F_X(X), Y)=S_f(X,Y)$
  for all $f\in \F$; the other equalities are proven similarly.)
\end{enumerate}

The conditional convex order (ccx) introduced in~\cite{ansari2025ordering} for inputs in $\R^d$ for some $d\ge 1$ and real valued outputs
satisfies the 8 axioms (under some restrictions), see 
 Theorem~1.4 therein. 

%%%%%%%%%%%%%%%%%%%%%%%%%%%%%%%
\section{Csisz\'ar index for Gaussian random vectors}
\label{sec:Gaussian}
%%%%%%%%%%%%%%%%%%%%%%%%%%%%%%%

Let $\cn(\mu, \Sigma)$ denote the probability distribution of a $d$ dimensional Gaussian distribution with mean $\mu$ and covariance matrix $\Sigma$. The image of $\Sigma$ and the support of $\cn(\mu, \Sigma)$ are defined as follows:
\[
\im(\Sigma) = \{ \Sigma x \, \colon \, x \in \R^d\} \quad \text{and} \quad \supp(\cn(\mu, \Sigma)) = \im(\Sigma) + \mu.
\]
For a dimension $d \in \N^*$, we denote by $\idmatrix_d$ the $d \times d$ identity matrix.
For  $X$ and $Y$ two square integrable  random vectors,  the  correlation
matrix $C$  between $X$ and $Y$
is:
\begin{equation}
   \label{eq:def-corr}
   C:=\corr(X,Y)  = \Sigma_X^{-1/2}  \cov(X,Y)  \Sigma_Y^{-1/2},
\end{equation}
where $\Sigma_X$ and  $\Sigma_Y$ are the covariance matrices  of $X$ and
$Y$, and with the convention that $C$ is the null matrix if $X$ or $Y$
is a.s.\ constant (that is, if $\min(\rank(\Sigma_X), \rank(\Sigma_Y))=0$).  See the definition of powers of covariance matrices in Section~\ref{sec:MatrixDefinitions}. Notice that all eigenvalues of $C$ belong to [-1,1]. 

\subsection{\texorpdfstring{$f$}{f}-divergence between Gaussian distributions}
Let  $P_1$ and $P_2$ be two $d$-dimensional Gaussian probability distributions.
The next
lemma  provides  a formula  for  the  corresponding $f$-divergence
if their  supports are equal  and not reduced  to a singleton.   We
consider   possibly  degenerate   distributions   and   show  that   the
$f$-divergence  between two  Gaussian measures  can be  expressed as  an
$f$-divergence between  Gaussian distributions with  diagonal covariance matrices; this  extends \cite[Corollary~2]{nielsen2021information}. In the sequel, for a $d\times d$ covariance matrix $M$ of rank $r$ we denote by $R_M$ a $d\times r$ matrix of eigenvectors associated to the $r$ positive eigenvalues of $M$ and by $\Lambda_M^+$ the $r\times r$ diagonal matrix with the positive eigenvalues of $M$. 
The proof     of     Lemma~\ref{lem:f_div_Gaussian}      is     given     in
Section~\ref{sec:proof-G}.

\begin{lem}[$f$-divergence between two Gaussians]
    \label{lem:f_div_Gaussian}
 Let $d\geq 1$.    For $i=1, 2$, let $P_i=\cn(\mu_i, \Sigma_i)$ be a
 $d$-dimensional  Gaussian probability distribution. If $P_1$ and $P_2$
 have the same support (that is, $\mu_1-\mu_2 \in \im(\Sigma_1)=\im (\Sigma_2)$) and  $r=\dim(\im(\Sigma_i))>0$, then
 for all $f \in \F$, we have:
    \begin{equation*}
        D_f(P_1 \,\Vert\, P_2) = D_f(\cn(\mu, \Lambda_M^+) \,\Vert\,  \cn(0,\idmatrix_{r})),
    \end{equation*}
    where the matrix $M = \Sigma_2^{-1/2} \Sigma_1 \Sigma_2^{-1/2}$ has rank $r$ and $\mu = R_M^\transpose \Sigma_2^{-1/2} (\mu_1 -\mu_2) $. If furthermore $\mu_1 = \mu_2$, we get:
    \begin{equation*}
            D_f(P_1 \,\Vert\, P_2)
            = D_f(\cn(0,\Lambda_M^+) \,\Vert\, \cn(0,\idmatrix_{r})).
      \end{equation*}
\end{lem}

\begin{rem}[The degenerated cases]
    \label{rem:P1P2-degenerated-cases}
Let $f\in \F$.  If the supports of $P_1$ and $P_2$ are not equal, the distributions are mutually singular, and $
D_f(P_1 \Vert P_2) = f(0) + f^*(0)$.
On the other hand, if their support are equal and reduced to a singleton (that is, $r=\dim(\im(\Sigma_i))=0$ in the above lemma), then the
probability measures $P_1$ and $P_2$ coincide with the same Dirac mass
and $D_f(P_1 \Vert P_2) =0$.
\end{rem}

\begin{rem}[Closed-form expressions and bounds]
Lemma~\ref{lem:f_div_Gaussian} reduces the computation of $D_f(P_1\Vert P_2)$ between two Gaussians with the same support to the case of a centered, diagonalized pair $\cn(\mu,\Lambda_M^+)$ and $\cn(0,\idmatrix_r)$. 
When $\mu=0$, closed-form expressions for $D_f(P_1\Vert P_2)$ exist for specific choices of $f$, in particular for the Kullback-Leibler divergence (and its reverse) and for the $\alpha$-divergences (which correspond to the $f$-divergences with $f_\alpha(t) = \big(t^\alpha - \alpha t - (1-\alpha)\big)/(\alpha(1-\alpha))$ for $\alpha \in \R \setminus \{0,1\}$).
Note that the $\alpha$-divergences include the squared Hellinger distance ($\alpha = 1/2$), the Pearson $\chi^2$-divergence ($\alpha = 2$) and the Neyman $\chi^2$-divergence ($\alpha = -1$), and that letting $\alpha$ tend to $1$ (resp.\ $0$) gives the Kullback-Leibler divergence (resp.\ the reverse Kullback-Leibler divergence).
For a general $f\in\F$, however, no such closed form is available.
In some cases, explicit bounds are available; see \cite{arbas2023polynomial} and \cite{devroye2018total} on the total variation distance between two Gaussians.

\end{rem}

%%%%%%%%%%%%%%%%%%%%%%%%%%%%%%%%%%%%%%%%%%%%%%%%%
\subsection{Csisz\'ar index for Gaussian vectors}
%%%%%%%%%%%%%%%%%%%%%%%%%%%%%%%%%%%%%%%%%%%%%%%%%

In this section we show that the Csisz\'ar index $S_f(X,Y)$ for $d$-dimensional Gaussian vectors $X$ and $Y$ can be written as some function $\kappa_f^{(d)}$ of the singular values of the correlation matrix $C= \corr(X,Y)$ defined in \eqref{eq:def-corr}. We first introduce the function $\kappa_f^{(d)}$ together with its properties. For some choices of $f$, explicit expressions of the Csisz\'ar index can be found in Table~\ref{tab:expression_gf1}. Then, we show the explicit relation between the Csisz\'ar index and the function $\kappa_f^{(d)}$ and give particular examples in more details. Finally, we study a Gaussian vector $(X,X',Y)$ and discuss Csisz\'ar order, \emph{i.e.\ }find conditions on the vector's distribution such that:
\[
X \csleq X' \quad \text{w.r.t.} \quad Y.
\]

\begin{rem}[The assumption that $(X,X',Y)$ is a Gaussian vector is slightly too strong]
\label{rem:XX'Y_GV_toostrong}
In this section, we consider a Gaussian random vector $(X,X',Y)$. Since we are only interested in Csiszár indices between $X$ and $Y$, $X'$ and $Y$, the dependence between $X$ and $X'$ can be arbitrary. So the results of this section hold under the slightly weaker assumption that $(X,Y)$ and $(X',Y)$ are two Gaussian random vectors.
\end{rem}

The study of the Csisz\'ar index for Gaussian vectors
relies on the following function. 
For  $f \in \F$, we set $\kappa_f^{(0)}=0$ and for $d\geq 1$ we define the function:
\[
\kappa_f^{(d)} : [0,1]^d \to [0, + \infty] \quad \text{by}\quad\kappa_f^{(d)} (x) =  D_{f}(\cn(0,\Gamma_x)) \, \Vert \, \cn(0,\idmatrix_{2d})) ,
\]
where $\Gamma_x$ is a diagonal matrix of size $2d$, with $x=(x_1, \ldots, x_d)$:
\begin{equation}
   \label{eq:def-Gamma-x}
\Gamma_x= \diag(1-\sqrt{x_1}, \dots, 1-\sqrt{x_d}, 1+\sqrt{x_1}, \dots, 1+\sqrt{x_d}).
\end{equation}
Table~\ref{tab:expression_gf1} presents the expressions of $\kappa_f^{(1)}$ and $\kappa_f^{(d)}$ for some functions $f \in \F$; which can also be deduced from~\cite{micheas2006measuring}.

The proof of the following properties of $\kappa_f^{(d)}$ is provided in Section~\ref{sec:proof-G}. The proof of the  monotonicity of $\kappa_f^{(d)} $  relies on  the theory of convex order from~\cite{shaked2007stochastic}; see also~\cite{ali1966general} for a direct but more technical proof.

\begin{prop}[Properties of $\kappa_f^{(d)}$] 
  \label{prop:gaussian_gfd_properties}
    Let $d\geq 1$ and $x=(x_1, \ldots, x_d)\in [0, 1]^d$. We have: %The function $\kappa_f^{(d)}$ satisfies the following properties:
    \begin{enumerate}[(i)]
    \item \label{enum:gaussian_gfd_properties_permutation} Permutation invariant: $\kappa_f^{(d)}(x_{\sigma(1)}, \ldots, x_{\sigma(d)}) = \kappa_f^{(d)}(x_1, \ldots, x_d)$ for any permutation $\sigma$ of $\{1, \ldots, d\}$.
    \item \label{enum:gaussian_gfd_properties_monoticity} Monotonicity: $\kappa_f^{(d)}$ is coordinate-wise non-decreasing.
     \item \label{enum:gaussian_gfd_properties_min} Minimum at $0$: $\kappa_f^{(d)}(0, \ldots, 0)=0$.
     \item \label{enum:gaussian_gfd_properties_max} Maximum at 1: If $x_i = 1$ for some $i\in \{1, \ldots d\}$, then $\kappa_f^{(d)}(x) = f(0) + f^*(0)$.
       \item \label{enum:gaussian_gfd_properties_dim} Dimension reduction: If  $x_i = 0$ for some $i\in \{1, \dots , d\}$, then:
        \[
        \kappa_f^{(d)}(x) = \kappa_f^{(d-1)}(x_1, \dots, x_{i-1}, x_{i+1}, \dots, x_d).
        \]
    \end{enumerate}
\end{prop}

See Lemma~\ref{lem:cvx-kappaf} for some comment on the convexity of $\kappa_f^{(2)}$.

\begin{table}[htbp]
    \renewcommand{\arraystretch}{1.25}
    \centering
    \begin{tabular}{|w{c}{5cm}|w{c}{4.7cm}|w{c}{5.4cm}|}
        \hline
        \Gape[4pt]{\makecell{$D_f$}}
        & \Gape[4pt]{\makecell{$\kappa_f^{(1)}(x),$ \\ $x\in [0,1)$}}
        & \Gape[4pt]{\makecell{$\kappa_f^{(d)}(x),$ \\ $x=(x_1,\ldots,x_d)\in [0,1]^d$}} \\
        \hline
        \makecell{Kullback-Leibler divergence, \\ $f_{\mathrm{KL}}(t) = t \ln(t)$}
        & $\displaystyle -\frac{1}{2} \log(1-x) \vphantom{\sum_{i=1}^d}$
        & $\displaystyle \sum_{i=1}^d \kappa_f^{(1)}(x_i) \vphantom{\sum_{i=1}^d}$ \\

        \makecell{Reverse Kullback-Leibler \\ divergence, $f_{\mathrm{KL}^*}(t) = - \ln(t)$}
        & $\displaystyle \frac{1}{2} \log(1-x) + \frac{x}{1-x} \vphantom{\sum_{i=1}^d}$
        & $\displaystyle \sum_{i=1}^d \kappa_f^{(1)}(x_i) \vphantom{\sum_{i=1}^d}$ \\

        \makecell{$\alpha$-divergence, $\alpha \in \R \setminus \{0,1\}$, \\ $f_{\alpha}(t) = c_\alpha(\alpha t - t^\alpha + (1-\alpha))$}
        & $\displaystyle c_\alpha \left( 1- \sqrt{\frac{(1-x)^{1-\alpha}}{(1-(1-\alpha)^2 x)_+}}\right) \vphantom{\prod_{i=1}^d}$
        & $\displaystyle c_\alpha \left( 1- \sqrt{\prod_{i=1}^d \frac{(1-x_i)^{1-\alpha}}{(1-(1-\alpha)^2 x_i)_+}}\right) \vphantom{\prod_{i=1}^d}$ \\

        \hline
    \end{tabular}
    \caption{Expressions of $\kappa_f^{(1)}$ and $\kappa_f^{(d)}$ (with $c_\alpha = 1/(\alpha(1-\alpha))$  and the convention $1/0=+\infty$).}
    \label{tab:expression_gf1}
\end{table}

The formula of the Csisz\'ar index $S_f(X,Y)$ between Gaussian vectors can be simplified into the function $\kappa_f^{(k)}$ for a properly chosen $k$. The proof of the next theorem is given in Section~\ref{sec:proof-G}.
We stress that the Gaussian vectors $X$ and $Y$ can be degenerated.

\begin{theorem}[Csisz\'ar  Index for Gaussian Vectors]
\label{thm:csiszar_index_Gaussian}
Let $(X, Y)$  be a Gaussian vector and 
$C  = \corr(X,Y)$   be   the              correlation               matrix
between $X$ and $Y$, see~\eqref{eq:def-corr}.   Let
$\lambda  \in (0,1]^k$  be the  vector  of the  positive eigenvalues  of
$C C^\transpose$  (and of  $C^\transpose C$) with $k\geq  0$ the  rank of
$C$. We have for $f\in \F$:
\begin{equation}
   \label{eq:Sf=gf}
S_f(X,Y) = \kappa_f^{(k)}(\lambda).
\end{equation}
\end{theorem}

\begin{rem}[Connection to canonical correlation analysis] 
\label{rem:CCA}
The ordered values of $\sqrt{\lambda_i}$, $i=1,\dots,k$, are the nonzero singular values of the correlation matrix $C$ and are classically known in the multivariate statistics literature as the canonical correlations between $X$ and $Y$, 
see~\cite{hardle2024canonical}.
It appears already in~\cite{gel1959random}, based on explicit computations,  that the mutual information $S_\mathrm{KL}$ is a function of the canonical correlations for  non-degenerate Gaussian vectors.
Our Theorem~\ref{thm:csiszar_index_Gaussian} above shows that the general Csiszár index $S_f(X,Y)$ between two jointly Gaussian vectors depends on their joint distribution only through their canonical correlations.
\end{rem}

\begin{exple}[Csisz\'ar index for Gaussian random vectors]
    $ $
    \begin{enumerate}[(i)]
        \item Let $(X,Y)$ be a $2$-dimensional Gaussian vector with correlation $\rho = \corr(X,Y)$, so that  $C= \rho$ is a number in [-1,1] and thus $C C^\transpose=\rho^2$. We get that for $f\in \F$:
        \[
        S_f(X,Y) = \kappa_f^{(1)}(\rho^2) = D_{f} (\cn(0, \diag(1-\rho, 1+\rho) ) ) \, \Vert \, \cn(0,\idmatrix_{2})) .
        \]
      \item Let  $X$ and  $Y$ be  independent Gaussian  random variables
        with $\Var(X)>0$,  and consider comparing $(X,Y)$  with $X$. The
        matrix $C^\transpose  C$, associated to  $((X,Y), X)$ is  then a
        $1 \times 1$ matrix and $C^\transpose C=1$. Thus, we have:
        \[
        S_f((X,Y), X) = \kappa_f^{(1)}(1) = f(0) + f^*(0).
        \]
        This result could also been deduced from Remark~\ref{rem:forgetting_independent_variables}
        %Remark~13~$(iii)$ in \cite{butucea2026csisz} 
        which states that $S_f((X,Y),Z) = S_f(X,Z)$ if $Y$ is independent of $Z$ and $Y$ has a conditional distribution given $X$. In this example, we get $S_f((X,Y), X) = S_f(X, X) = f(0) + f^*(0)$.
      \item Let  $X = (X_1, X_2)  \sim \cn(0, \idmatrix_2)$
        and $Y  = (Y_1, Y_2)  \sim \cn(0, \idmatrix_2)$, with  $(X,Y)$ a
        Gaussian vector. Suppose  that $X_2$ is independent  of $Y$. Let
        $\rho = \corr(X_1,Y_1)$ and $\rho' = \corr(X_1,Y_2)$. We get:
        \[
        C = \begin{pmatrix} \rho & \rho'\\ 0 & 0 \end{pmatrix} \quad \text{and} \quad C^\transpose C= \begin{pmatrix} \rho^2 & \rho \rho' \\ \rho \rho'  & \rho'^2 \end{pmatrix}.
        \]
        The eigenvalues of $C^\transpose C$ are $\rho^2 + \rho'^2$ and $0$. We obtain:
        \[
        S_f(X,Y) = \kappa_f^{(2)}(\rho^2 + \rho'^2,0) = \kappa_f^{(1)}(\rho^2 + \rho'^2).
        \]
        We could also use Remark~\ref{rem:forgetting_independent_variables}
        % Remark~13~$(iii)$ in \cite{butucea2026csisz}
        to get that $S_f(X,Y)=S_f(X_1,Y)$.  Since the correlation matrix
        between  $X_1$   and  $Y$  reduces  to   $C=(\rho,  \rho')$  and
        $C C^\transpose$  is a real number equal  to $\rho^2 + \rho'^2$,  we then
        recover that $S_f(X,Y)=S_f(X_1,Y)= \kappa_f^{(1)}(\rho^2 + \rho'^2)$.
    \end{enumerate}
\end{exple}

\begin{rem}[Invariance by sphering transformation]
    \label{rem:invariance_by_sphering_transformation}
    Let $(X,Y)$  be a  Gaussian vector  with $X$  and $Y$  non constant.
    Since  the  correlation  matrix  between $X$  and  $Y$  and  between
    $\Sigma_X^{-1/2}X$ and  $\Sigma_Y^{-1/2}Y$ are  the same,  we deduce
    from Theorem~\ref{thm:csiszar_index_Gaussian} that:
    \[
    S_f(X,Y) = S_f(\Sigma_X^{-1/2}X,\,  \Sigma_Y^{-1/2}Y).
  \]
  When $X$ and $Y$ have full ranks, %that is, when the kernel of $\Sigma_X$ and $\Sigma_Y$ are of dimension 0, 
  then the covariance matrices of $\Sigma_X^{-1/2} X$ and $ \Sigma_Y^{-1/2} Y$ are the identity matrices.
\end{rem}

\subsection{Csisz\'ar ordering for Gaussian vectors}
Let $(X,X',Y)$ be a Gaussian random vector.
In this subsection, we provide a complete characterization of the Csisz\'ar order when either $Y$ is one-dimensional or both $X$ and $X'$ are one-dimensional, together with a sufficient condition in higher dimensions.

\medskip

\subsubsection{Univariate case}
%\textbf{Univariate case:}
We first consider the Csisz\'ar order either in the case where $Y$ is one-dimensional, or in the case where $X$ and $X'$ are one-dimensional.
From Theorem~\ref{thm:csiszar_index_Gaussian}, we observe that, in this setting, the Csisz\'ar index depends on a single quantity only, which we introduce through the function $h$.
Let  $(G_0, G)  $  be a  Gaussian  random vector  with  $G_0\in \R$  and
$G\in \R^d$,  $d\geq 1$.       For
$x=(x_1,       \ldots,       x_{d})\in       \R^{d}$,       we       set
$\norm{x}_2^2= \sum_{i=1}^{d}  x_i^2 $. We  define the real number $h(G_0, G)=h(G, G_0)$ by:
\begin{equation}
   \label{eq:def-h}
h(G, G_0)= \norm{ \corr(G,G_0)  }_2^2.
\end{equation}
Notice that, with $G=(G_1, \ldots, G_d)$ and $\Sigma_G$ the covariance
matrix of $G$:
\[
  \Sigma_G \text{ diagonal} \quad
  \Longrightarrow\quad
  h(G,G_0) = \sum_{i=1}^d \corr(G_i, G_0)^2.
\]
From Theorem~\ref{thm:csiszar_index_Gaussian}, we deduce the following corollary, which characterizes the Csisz\'ar order. Its proof is given in Section~\ref{sec:proof-G}.
We also obtain that whenever $Y$ is one-dimensional or both $X$ and $X'$ are one-dimensional, $X$ and $X'$ are always comparable in the Csisz\'ar order.

\begin{cor}[The univariate case] \label{cor:gaussian_case_dim_1}
Consider a Gaussian vector $(X,X',Y)$ with $X\in \R^{d_X}$, $X'\in
\R^{d_{X'}}$ and $Y\in \R^{d_Y}$.  %We have the following  properties:
\begin{enumerate}[(i)]
    \item \label{enum:gaussian_dim1_case1} If $d_Y=d_X=d_{X'}=1$, we
      have:
      \begin{equation}
        \label{eq:corr-Sf-1}
        X \csleq X' \quad \text{w.r.t.} \quad Y
        \quad \iff \quad
        \corr(X, Y)^2 \leq  \corr(X',Y)^2 .
      \end{equation}
    \item \label{enum:gaussian_dim1_case2} If $d_Y=1$ or if $d_X= d_{X'}=1$, we have
      $S_f(X,Y)= \kappa_f^{(1)}\left(h(X,Y)\right)$ and:
    \begin{equation}
        \label{eq:corr-Sf-X=1}
        X \csleq X' \quad \text{w.r.t.} \quad Y
        \quad \iff \quad 
        h(X, Y) \leq h(X',Y).
      \end{equation}
\end{enumerate}
\end{cor}

\begin{rem}[Connection to the conditional convex order]   \label{rem:Connection_Ccvx_Gaussian}
In the Gaussian model with a  one-dimensional output variable, 
Theorem~3.1 of~\cite{ansari2025ordering}  shows that  the Csiszár order is equivalent to 
the conditional convex order defined therein. 
\end{rem}

In Corollary~\ref{cor:gaussian_case_dim_1}, we observed that the Csisz\'ar order between $X$ and $X'$ w.r.t.\ $Y$, when $(X,X',Y)$ is a three-dimensional Gaussian vector, is entirely characterized by the correlation. However, this structure does not extend in general outside the Gaussian vectors; see for instance Example~\ref{exple:ContreExempleBernoulli}, which shows that Csisz\'ar indices may be ranked differently from correlations, as well as Example~\ref{exple:OrdreUniformeVsCorr}, which shows that the Csisz\'ar order can also disagree with the ordering induced by correlations.

\begin{rem}[Connection to the Linear Regression Model]
Consider the real linear regression model:
\[
Y = \beta^\transpose  X + \varepsilon \quad \text{with } Y \in \R,
\]
where  $X$ and $\varepsilon$  are independent   centered  Gaussian vectors  and $\beta$  is  a coefficient vector. Let $\Sigma_X$ denote the covariance matrix of $X$, which we assume to be invertible. Then, we have $\cov(X,Y) = \Sigma_X \beta$. The coefficient of determination or $R^2$ score of this model is given by:
\[
R^2 %= \frac{\Var(\tilde{Y})}{\Var(Y)} 
= \frac{\cov(Y,X) \Sigma_X^{-1} \cov(X,Y)}{\Var(Y)} = h(X,Y),
\]
where %$\tilde{Y} = \beta^\transpose X$ is the noiseless part of the model and 
$h$ is the function defined in~\eqref{eq:def-h}. Thus the Csisz\'ar index between $X$ and $Y$  is  a function of the  score $R^2$ of the linear regression prediction:
\[
  S_f(X,Y)=\kappa_f^{(1)}(R^2).
\]
The usual interpretation of $R^2$ being high (close to 1) for a response $Y$ strongly dependent of the input $X$ is totally valid for Gaussian vectors and univariate responses. 
% However, as our previous examples have shown, this is no longer obvious for multivariate regression (with multivariate $Y$) even within the restrictive framework of Gaussian vectors.

When, in addition, the entries of $X = (X_1, \dots, X_d)$ are independent, then we have the decomposition $R^2 = \sum_{i=1}^d \corr(X_i,Y)^2$, and each term $\corr(X_i,Y)^2$ coincides with the first-order Sobol index of $X_i$. These same quantities are also connected, in this independent-features setting, to standard random forest importance measures: the sum of the marginal Mean Decrease Impurity (MDI) scores converges to $R^2 \Var(Y)$ \cite{scornet2023trees}, and the Mean Decrease Accuracy (MDA) of variable $X_i$ converges to a quantity related to its total Sobol index \cite{benard2022mean}. Consequently, in this particular case, the order induced by the Csiszár indices on the individual entries of $X$ w.r.t.\ $Y$ coincides with the order induced by the Sobol indices.
\end{rem}

\subsubsection{Multivariate case}
% \textbf{Multivariate case:}
We now consider the case where $(X,X',Y)$ is a Gaussian random vector of arbitrary dimension and we denote by $C = \corr(X,Y)$ and $C' = \corr(X',Y)$ the correlation matrices between $X$ and $Y$, and between $X'$ and $Y$, respectively see~\eqref{eq:def-corr}. Recall that by  convention the correlation between $X$ and $Y$, see~\eqref{eq:def-corr}, is zero when $X$ or  $Y$ are constant.

For two real symmetric square matrices $M$ and $M'$ with the same size, we say that $M$ is smaller than $M'$ in the Loewner order, denoted by $M \preceq M'$,  if $M' - M$ is positive semi-definite. In particular if $\lambda_1 \geq \cdots \geq \lambda_d $ and $\lambda'_1 \geq \cdots \geq \lambda'_d $ are  the eigenvalues of $M$ and $M'$, we recall that $M \preceq M'$ implies that $\lambda_i\leq \lambda'_i$ for all $i=1, \ldots, d$.

We first state a corollary giving a sufficient condition for $X \csleq X'$ w.r.t.\ $Y$. Its proof is immediate from the monotonicity property of $\kappa_f^{(d)}$ in Proposition~\ref{prop:gaussian_gfd_properties}~\ref{enum:gaussian_gfd_properties_monoticity}.

\begin{cor}[Loewner order for the correlation matrices and  the Csisz\'ar order]
\label{cor:loewner-C}
    Let $(X, X',Y)$  be a Gaussian vector and set $C = \corr(X, Y)$ and $C' = \corr(X', Y)$, see~\eqref{eq:def-corr}.  We have:
    \[
    C^\transpose C \preceq (C')^\transpose C' \quad \implies
    \quad 
    X \csleq X' \quad \text{w.r.t.} \quad Y.
    \]
\end{cor}

\begin{rem}[Sufficient conditions for the Csisz\'ar order]
  \label{rem:sufficient_conditions}
  From the properties of $\kappa_f^{(d)}$ given in Proposition~\ref{prop:gaussian_gfd_properties}, we can derive two sufficient conditions for the Csisz\'ar order: 
  \begin{enumerate}[(i)]
      \item \label{it:loewner_rotation} If there exists an orthogonal matrix $U$ such that $C^\transpose C \preceq U^\transpose (C')^\transpose C' U$, then $X \csleq X'$ w.r.t.\ $Y$. This follows from the coordinate-wise monotonicity of $\kappa_f^{(d)}$. It is stronger than Corollary~\ref{cor:loewner-C}, as it allows the correlation matrices to be compared up to an orthogonal rotation. 
      \item \label{it:vp1} If $1$ is an eigenvalue of $C'^\transpose C'$, then the Csisz\'ar index between $X'$ and $Y$ attains its maximal value, and hence $X \csleq X'$ w.r.t.\ $Y$.
  \end{enumerate}
\end{rem}

The converse of Corollary~\ref{cor:loewner-C} does not hold in general, even up to the orthogonal rotation of Remark~\ref{rem:sufficient_conditions}~\ref{it:loewner_rotation}, as shown by the following counterexample.

\begin{exple}[The Csisz\'ar order is not equivalent to the Loewner order]
  Consider a Gaussian vector $(X,X',Y)$ with $X\in \R^2$, $X'\in \R^2$, and $Y\in \R^2$ such that:
  \[
  C = \corr(X,Y) = \begin{pmatrix} \alpha & 0 \\ 0 & \beta \end{pmatrix}
  \quad \text{and} \quad
  C' = \corr(X',Y) = \begin{pmatrix} \alpha' & 0 \\ 0 & \beta' \end{pmatrix},
  \]
  with $\alpha, \alpha', \beta, \beta'$ be four correlation coefficients in $[-1,1]$.
  The eigenvalues of $C^\transpose C$ are $\alpha^2$ and $\beta^2$, and those of $C'^\transpose C'$ are $\alpha'^2$ and $\beta'^2$. Taking $\alpha = 1$, $\alpha', \beta' \in (-1,1) \setminus \{0\}$ and $\beta \in (-1,1)$ such that $\beta^2 < \min(\alpha'^2, \beta'^2)$, we find that $C^\transpose C$ and $C'^\transpose C'$ are not comparable in the Loewner order even up to a rotation, since $\alpha^2 > \max(\alpha'^2,\beta'^2)$ but $\beta^2 < \min(\alpha'^2, \beta'^2)$. However, $X \csgeq X'$ w.r.t.\ $Y$, since the Csisz\'ar index between $X$ and $Y$ attains its maximal value by property~\ref{it:vp1} of Remark~\ref{rem:sufficient_conditions}.
\end{exple}

Moreover, unlike the case where $Y$ is one-dimensional or where $X$ and $X'$ are one-dimensional, the Csiszár order is not total for Gaussian vectors in general as shown in the following example.

\begin{exple}[The Csiszár order is not total in the multivariate Gaussian case]
\label{exple:ranking_depends_on_f_gaussian}
In the multivariate case
% , if the correlation matrices are not ordered in Loewner sense, 
the Csisz\'ar indices can be ordered differently according to the choice of $f$. 
Let $(X,X',Y)$ be a $5$-dimensional Gaussian vector such that $\Var(X)=1$,  $\Var(X') = \Var(Y) = \idmatrix_2$  and:
\[
C = \cov(X,Y) = \begin{pmatrix} \alpha & 0 \end{pmatrix}
\quad \text{and} \quad
C' = \cov(X',Y) = \beta \idmatrix_2 .
\]
Notice %, with the notation of the Corollary~\ref{cor:loewner-C},
that the matrices $C^\transpose C = \diag(\alpha^2, 0)$ and $C'^\transpose C' = \diag(\beta^2, \beta^2)$ of dimension $2\times 2$ are not comparable w.r.t.\ the Loewner order if $\beta^2 < \alpha^2$.
According to Theorem~\ref{thm:csiszar_index_Gaussian}, for every $f \in \F$, the Csisz\'ar indices are given by:
\[
S_f(X,Y) = \kappa_f^{(1)}(\alpha^2)
 \quad \text{and} \quad
 S_f(X',Y) = \kappa_f^{(2)}(\beta^2, \beta^2).
\]
Taking $\alpha = 0.75$ and $\beta = 0.6$ and using the expressions of $\kappa_f^{(2)}$ given in Table~\ref{tab:expression_gf1}, we remark that:
\[
S_{\mathrm{KL}}(X,Y) > S_{\mathrm{KL}}(X',Y)
\quad \text{and} \quad
S_{\mathrm{KL}^*}(X,Y) < S_{\mathrm{KL}^*}(X',Y).
\]
We deduce that the Csisz\'ar order is not total in the multidimensional Gaussian case.
\end{exple}

Let $(X',Y)$ be a Gaussian vector with $Y\in \R^ d$. In order to characterize $\cp_Y^+(X')$ and $\cp_Y^-(X')$ restricted to the conditional Gaussian distributions, that is, $(X,Y)$ is Gaussian, we define the parameter $x \in [0,1]^d$ as the vector of eigenvalues of $C^\transpose C$, with $C= \corr(X,Y)$ and, similarly for $x'$ with $X$ replaced by $X'$. Notice that $x$ does not characterize $P_{X \mid Y}$, but $x = x'$ implies that $X \cseq X'$ w.r.t.\ $Y$. 
For a given $x'\in [0, 1]^d$, the subsets:
\[
\cc_{x'}^+= \{x \in [0, 1]^d \, \colon \, X \csgeq X' \quad\text{w.r.t.} \quad Y\} 
\quad \text{and} \quad
\cc_{x'}^-= \{x \in [0, 1]^d \, \colon \, X \csleq X' \quad\text{w.r.t.} \quad Y\},
\]
contains the parameters $x$ associated to the conditional distribution
$P_{X \mid Y}$ of the sets $\cp_Y^+(X')$ and $\cp_Y^-(X')$, respectively. 
% Recall that $\cp_Y^+(X')$ and $\cp_Y^-(X')$ are defined in Definition~\ref{def:Csiszar sets}. 
As we have seen in Example~\ref{exple:ranking_depends_on_f_gaussian}, $X$ and $X'$ are not necessarily comparable.
Thus, $X \csgeq X'$ corresponds to a vector $x$ in $ \cc_{x'}^+$ and $X \csleq X'$ to  $x $ in $\cc_{x'}^-$. If $X$ and $X'$ are not ordered in Csisz\'ar order w.r.t.\ $Y$ then $x$ belongs to the set $[0, 1]^d \setminus (\cc_{x'}^+\cup \cc_{x'}^-)$. 

\medskip

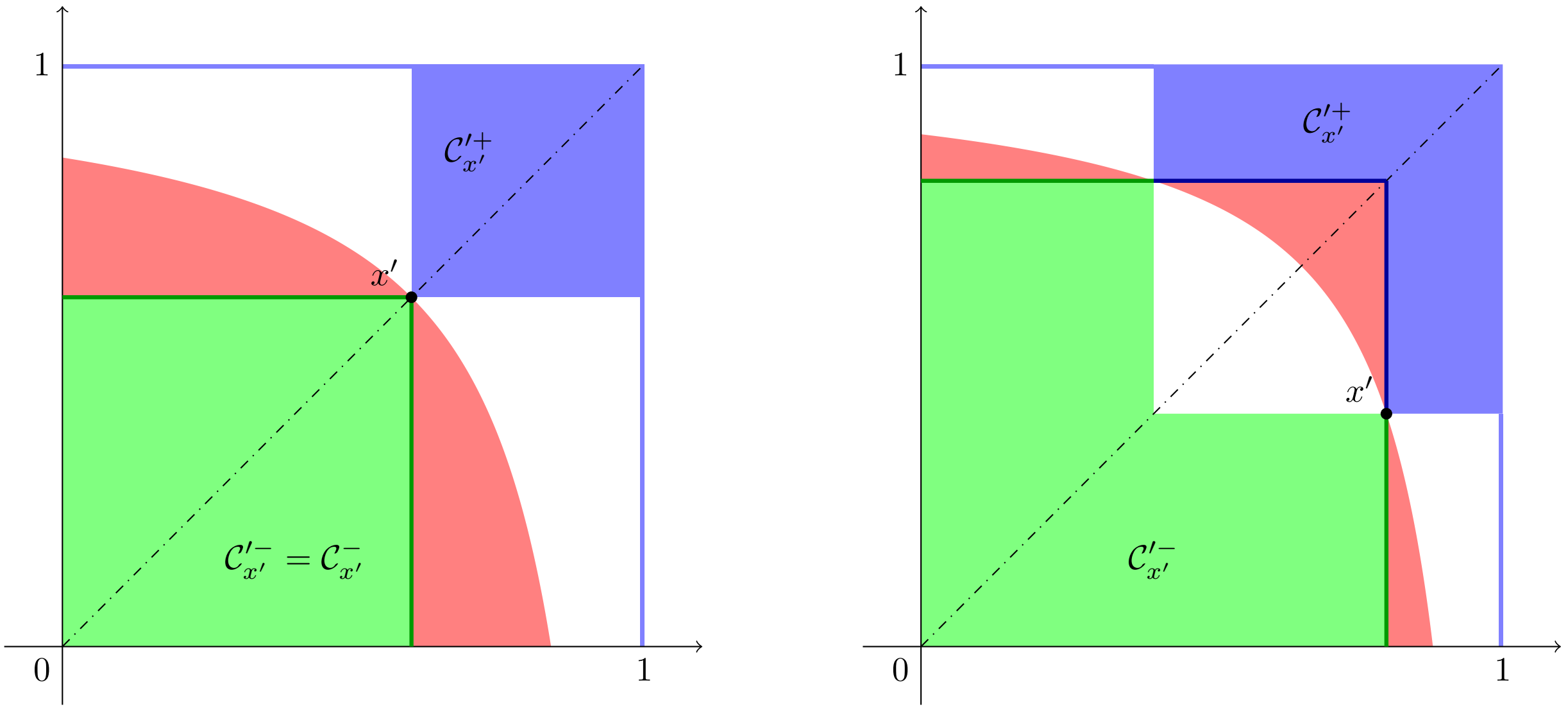
\begin{figure}[htbp]
    \centering

    \begin{subfigure}[t]{0.45\textwidth}
        \centering
        \begin{tikzpicture}[scale=6]

            % Partie basse
            \path[fill=\ColorCm, opacity=\OpacityPlot]
                (0,0) rectangle (0.6,0.6);

            % Partie haute
            \path[fill=\ColorCp, opacity=\OpacityPlot]
                (0.6,0.6) rectangle (1,1);
            \path[fill=\ColorCp, opacity=\OpacityPlot]
                ([xshift=-0.2pt]1,0) rectangle (1,0.6);
            \path[fill=\ColorCp, opacity=\OpacityPlot]
                ([yshift=-0.2pt]0,1) rectangle (0.6,1);

            \path[fill=\ColorOut,opacity=\OpacityPlot]
            plot[domain=0:0.6, smooth, samples = \nsamples](\x,{1-exp(-2*0.916290731874155)/(1-\x)}) -- (0, 0.6) 
            -- cycle; 

            \begin{scope}[cm={0,1,1,0,(0,0)}]
            \path[fill=\ColorOut,opacity=\OpacityPlot]
            plot[domain=0:0.6, smooth, samples = \nsamples](\x,{1-exp(-2*0.916290731874155)/(1-\x)}) -- (0, 0.6) 
            -- cycle;

            \end{scope}

            % Axes et carré
            \draw[dashdotted] (0,0) -- (1,1);
            \draw[\ArrowStyle] (-0.1,0) -- (1.1,0);
            \draw[\ArrowStyle] (0,-0.1) -- (0,1.1);

            \draw[green!60!black, very thick] (0.6,0,0) -- (0.6,0.6);
            \draw[green!60!black, very thick] (0,0.6) -- (0.6,0.6);

            \fill (0.6,0.6) circle (0.01) node[above left] {$x'$};

            \node at (0.4,0.15) {$\cc_{x'}'^- = \cc_{x'}^-$};
            \node at (0.7,0.85) {$\cc_{x'}'^+$};

            \node[below left] at (0,0) {$0$};
            \node[below] at (1,0) {$1$};
            \node[left]  at (0,1) {$1$};

        \end{tikzpicture}
        \caption{$x' = (0.6, 0.6)$}
        \label{fig:A+-_gaussiennne_0.6_0.6}
    \end{subfigure}
    \hfill
    \begin{subfigure}[t]{0.45\textwidth}
        \centering
        \begin{tikzpicture}[scale=6]

            \path[fill=\ColorOut,opacity=\OpacityPlot]
            plot[domain=0.4:0.8, smooth, samples = \nsamples](\x,{1-exp(-2*1.0601317681000455)/(1-\x)}) -- (0.8, 0.8) 
            -- cycle;  

            \path[fill=\ColorOut,opacity=\OpacityPlot]
            plot[domain=0:0.4, smooth, samples = \nsamples](\x,{1-exp(-2*1.0601317681000455)/(1-\x)}) -- (0, 0.8) 
            -- cycle; 

            \begin{scope}[cm={0,1,1,0,(0,0)}]
                \path[fill=\ColorOut,opacity=\OpacityPlot]
                plot[domain=0:0.4,smooth,samples=\nsamples]
                    (\x,{1-exp(-2*1.0601317681000455)/(1-\x)})
                -- (0,0.8)
                -- cycle;

            \end{scope}

            % Partie basse
            \path[fill=\ColorCm, opacity=\OpacityPlot]
                (0,0) rectangle (0.4,0.4)
                (0.4,0) rectangle (0.8,0.4)
                (0,0.4) rectangle (0.4,0.8);

            % Partie haute
            \path[fill=\ColorCp, opacity=\OpacityPlot]
                (0.8,0.8) rectangle (1,1)
                (0.8,0.4) rectangle (1,0.8)
                (0.4,0.8) rectangle (0.8,1);
            
            \path[fill=\ColorCp, opacity=\OpacityPlot]
                ([xshift=-0.2pt]1,0) rectangle (1,0.4);
            \path[fill=\ColorCp, opacity=\OpacityPlot]
                ([yshift=-0.2pt]0,1) rectangle (0.4,1);

            % Axes et carré
            \draw[dashdotted] (0,0) -- (1,1);
            \draw[\ArrowStyle] (-0.1,0) -- (1.1,0);
            \draw[\ArrowStyle] (0,-0.1) -- (0,1.1);

            \draw[blue!60!black, very thick] (0.8,0.8) -- (0.8,0.4);
            \draw[blue!60!black, very thick] (0.8,0.8) -- (0.4,0.8);

            \draw[green!60!black, very thick] (0.8,0) -- (0.8,0.4);
            \draw[green!60!black, very thick] (0,0.8) -- (0.4,0.8);

            \fill (0.8,0.4) circle (0.01) node[above left] {$x'$};

            \node at (0.4,0.15) {$\cc_{x'}'^-$};
            \node at (0.7,0.9) {$\cc_{x'}'^+$};

            \node[below left] at (0,0) {$0$};
            \node[below] at (1,0) {$1$};
            \node[left]  at (0,1) {$1$};
        \end{tikzpicture}
        \caption{$x' = (0.8, 0.4)$}
        \label{fig:A+-_gaussiennne_0.8_0.4}
    \end{subfigure}
    \caption{Two instances of the sets $\cc_{x'}'^-$ in green and $\cc_{x'}'^+$ in blue for two different values of $x'\in [0, 1]^2$. 
    The Csiszár indices are not ordered in the red region, and the comparison is unknown in the white region. 
    }
\label{fig:gf_ordering_dim2}
\end{figure}

To further discuss and illustrate the Csisz\'ar order for Gaussian vectors graphically, we assume that $Y$ is two-dimensional. For $x=(x_1, x_2)$ and $x'=(x'_1, x'_2)$, we set $\tilde x=(x_2, x_1)$ and $x\le x'$ if $x_i\le x'_i$ for $i=1, 2$. The symmetry, monotonicity, and maximum value of $\kappa_f^{(2)}$ (see Proposition~\ref{prop:gaussian_gfd_properties}), together with Corollary~\ref{cor:loewner-C}, imply that if $x \leq x'$ or $x \leq \tilde x'$, then $x \in \cc_{x'}^-$, and similarly if $x \geq x'$ or $x \geq \tilde x'$, then $x \in \cc_{x'}^+$. Moreover, if $x_i = 1$ for some $i$, then $x \in \cc_{x'}^+$. Thus, we get that $\cc_{x'}'^-$ is a subset of $\cc_{x'}^-$ and $\cc_{x'}'^+$ is a subset of $\cc_{x'}^+$, where:
\[
\cc_{x'}'^- = \left\{x \in [0, 1]^2\, \colon\, x  \leq x' \text{ or } x \leq \tilde x' \right\}
\text{ and }
\cc_{x'}'^+ = \left\{x \in [0, 1]^2\, \colon\, x \geq x' \text{ or } x \geq \tilde x' \text{ or } x_i = 1 \text{ for some } i \right\}.
\]
We give an instance of those two sets in Fig.~\ref{fig:gf_ordering_dim2}.
By considering the level set $\{x \in [0,1]^2 \, \colon \, \kappa_\alpha(x) = \kappa_\alpha(x')\}$ and letting $\alpha$ tend first to $1 - 1/\sqrt{\max(x')}$ and then to $1$ (corresponding to the Kullback-Leibler divergence), we obtain the red  subset of $[0,1]^2 \setminus (\cc_{x'}^+ \cup \cc_{x'}^-)$, corresponding to the parameters $x$ associated to the Gaussian vector $X$ not comparable to $X'$ for the Csisz\'ar order.
This gives part of the exact boundaries for $\cc_{x'}^+$ and $\cc_{x'}^-$, shown as dark green and dark blue lines in the figure.
In particular, $\cc_{x'}^-$ is a subset of the square $[0, \max(x')]^2$, with equality if $x_1' = x_2'$.

An analytic characterization of $\cc_{x'}^+$ and $\cc_{x'}^-$ remains an open problem. The white regions in Fig.~\ref{fig:gf_ordering_dim2} are therefore the regions where it is not known whether $x$ belongs to $\cc_{x'}^+$, to $\cc_{x'}^-$, or to neither.

It seems difficult to give other properties on the set $\cc_{x'}^-$, such as its convexity (see Section~\ref{sec:Bernoulli} on Bernoulli random variables where such property is established). For example, the convexity of 
the set $\cc_{x'}^-$ would be a direct consequence of the 
functions $\kappa_f^{(d)}$ being  quasi-convex for all $f\in \F$. The next lemma asserts that the latter assertion is false.

Recall that a function $g:[0, 1]^d\to\R$ is  quasi-convex if, for all $x,y\in [0, 1]^d$ and all $t\in[0,1]$, $g(tx+(1-t)y)\leq \max (g(x),g(y))$, that is,  its sublevel sets $\{g\leq a\}$ are convex for all $a\in\R$.

\begin{lem}[Convexity of $\kappa_f^{(d)}$] 
\label{lem:cvx-kappaf}
Let $d\ge 1$. 
\begin{enumerate}[(i)]
  \item \label{enum:kappa_KL} The function $\kappa_f^{(d)}$ is convex for $f\in \{f_{\mathrm{KL}},f_{\mathrm{KL}^*}\}$.
  \item\label{enum:kappa_alpha} The function $\kappa_\alpha^{(d)}$ is convex for $\alpha \in \R \setminus [0,1]$, and also  for $\alpha \in \R \setminus \{0,1\}$ if  $d =1, 2$. 
  \item \label{enum:kappa_TV} The function $\kappa_{\mathrm{TV}}^{(d)}$ is not quasi-convex for $d \geq 2$. %\AS{and not convex for $d = 1$}.
\end{enumerate}
\end{lem}

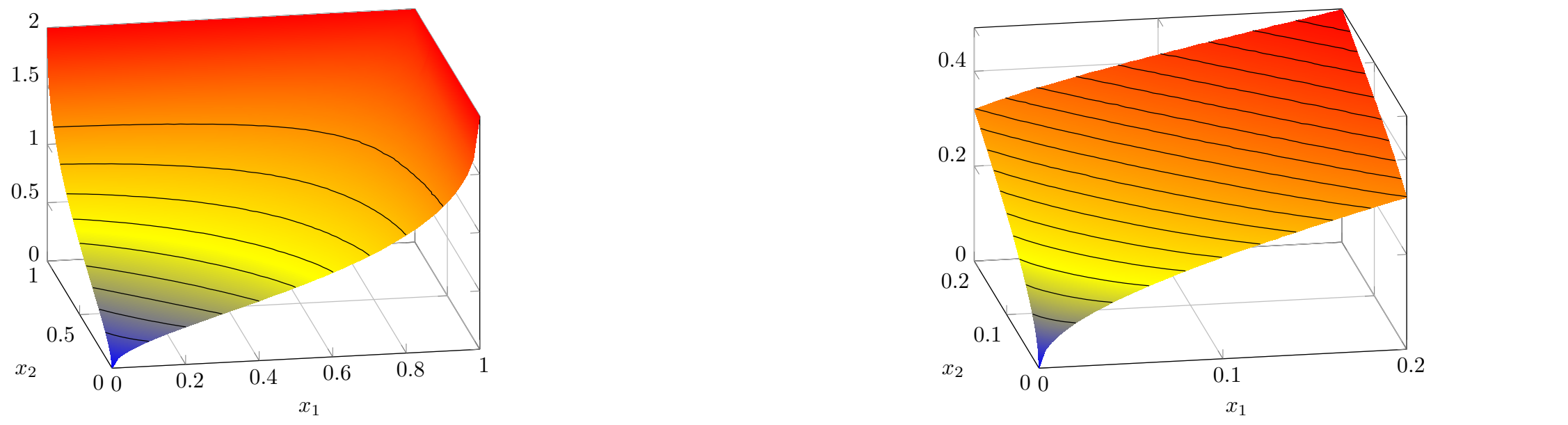
\begin{figure}
    \centering
    \begin{subfigure}[t]{0.45\textwidth}
        \centering
                \begin{tikzpicture}[scale = 0.6]
            \begin{axis}[
                view={-10}{25},
                xlabel={$x_1$},
                ylabel={$x_2$},
                grid=major,
                xmin=0, xmax=1,
                ymin=0, ymax=1,
                domain=0:1,
                y domain=0:1,
                samples=\nsamples,
                enlargelimits=false,
                xtick={0,0.2,0.4,0.6,0.8,1},    
            ]

            \addplot3[
                surf,
                shader=interp,
            ]
            table {Plots/Plots_gf_gaussian/surfaceTV_0_1.dat};

            \addplot3[
                contour prepared={labels=false},
                contour/draw color=black,
            ]
            table {Plots/Plots_gf_gaussian/contoursTV_0_1.dat};

            \end{axis}
        \end{tikzpicture}
        \caption{Contour plot of $\kappa_{\mathrm{TV}}^{(2)}$ on $[0,1]^2$.}
      \end{subfigure}
    \hfill
    \begin{subfigure}[t]{0.45\textwidth}
        \centering

        \begin{tikzpicture}[scale = 0.6]
            \begin{axis}[
                view={-10}{25},
                xlabel={$x_1$},
                ylabel={$x_2$},
                grid=major,
                xmin=0, xmax=0.2,
                ymin=0, ymax=0.2,
                domain=0:0.2,
                y domain=0:0.2,
                samples=\nsamples,
                enlargelimits=false,
                xtick={0,0.1,0.2},    
            ]

            \addplot3[
                surf,
                shader=interp,
            ]
            table {Plots/Plots_gf_gaussian/surfaceTV_0_0.2.dat};

            \addplot3[
                contour prepared={labels=false},
                contour/draw color=black,
            ]
            table {Plots/Plots_gf_gaussian/contoursTV_0_0.2.dat};

            \end{axis}
        \end{tikzpicture} 
        \caption{Zoom on the region $[0,0.2]^2$ to illustrate the non-convexity of the level sets of $\kappa_{\mathrm{TV}}^{(2)}$.}
    \end{subfigure}
    \caption{Contour plot of $\kappa_{\mathrm{TV}}^{(2)}$  estimated by Monte Carlo sampling, showing that the function $\kappa_{\mathrm{TV}}^{(2)}$ is not quasi-convex (see Lemma~\ref{lem:cvx-kappaf}).}
    \label{fig:contour_plot_gf_TV}
\end{figure}

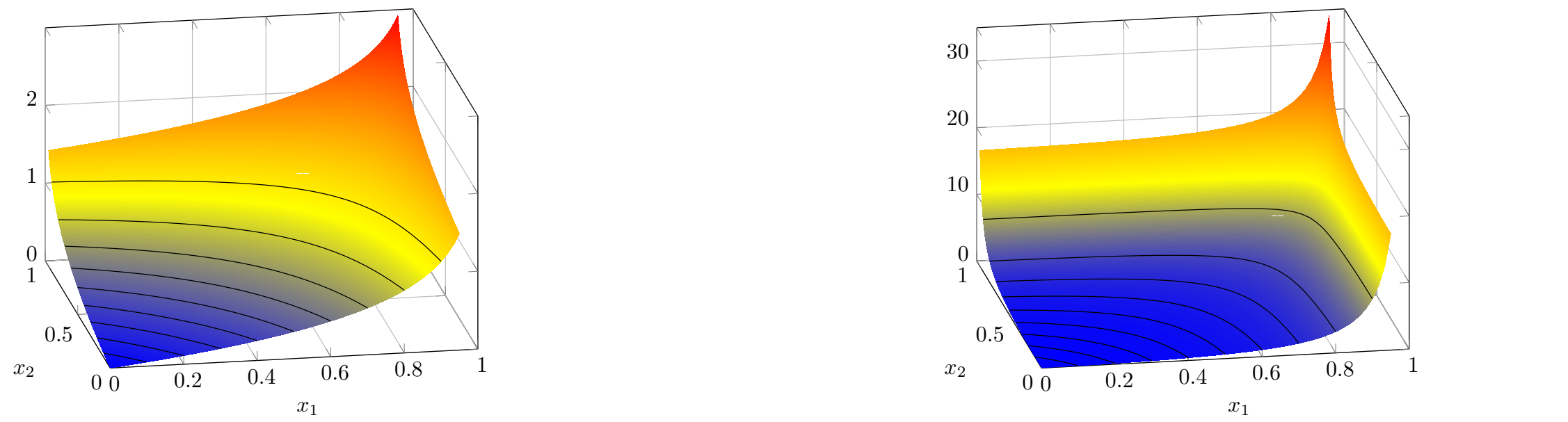
\begin{figure}[htbp]
    \centering

    \begin{subfigure}[t]{0.45\textwidth}
              \centering
        \begin{tikzpicture}[scale = 0.6]
            \begin{axis}[
                view={-10}{25},
                xlabel={$x_1$},
                ylabel={$x_2$},
                grid=major,
                xmin=0, xmax=1,
                ymin=0, ymax=1,
                domain=0:0.95,
                y domain=0:0.95,
                samples=\nsamples,
                enlargelimits=false,
            ]

            \addplot3[
                surf,
                shader=interp,
            ]
            {-0.5*ln(1-x) - 0.5*ln(1-y)};

            \addplot3[
                contour gnuplot={
                    levels={
                        0.05268025782891314,
                        0.11157177565710485,
                        0.17833747196936622,
                        0.25541281188299536,
                        0.34657359027997264,
                        0.45814536593707766,
                        0.6019864021629681,
                        0.8047189562170503,
                        1.151292546497023
                    },
                    draw color=black,
                    labels=false,
                },
                samples=\nsamples,
                domain=0:1,
            ]
            {-0.5*ln(1-x) - 0.5*ln(1-y)};

            \end{axis}
        \end{tikzpicture}
        \caption{Contour plot of $\kappa_{\mathrm{KL}}^{(2)}$.}
        \label{fig:contour_plot_gf_KL_conj}
    \end{subfigure}
    \hfill
    \begin{subfigure}[t]{0.45\textwidth}
        \centering
        \begin{tikzpicture}[scale = 0.6]
            \begin{axis}[
                view={-10}{25},
                xlabel={$x_1$},
                ylabel={$x_2$},
                grid=major,
                xmin=0, xmax=1,
                ymin=0, ymax=1,
                domain=0:0.95,
                y domain=0:0.95,
                samples=\nsamples,
                enlargelimits=false,
            ]

            \addplot3[
                surf,
                shader=interp,
            ]
            {0.5*ln(1-x) + x/(1-x) + 0.5*ln(1-y) + y/(1-y)};

            \addplot3[
                contour gnuplot={
                    levels={
                        0.05843085328219798,
                        0.13842822434289515,
                        0.25023395660206244,
                        0.4112538547836714,
                        0.6534264097200273,
                        1.041854634062923,
                        1.731346931170366,
                        3.1952810437829506,
                        7.848707453502978
                    },
                    draw color=black,
                    labels=false,
                },
                samples=\nsamples,
                domain=0:1,
            ]
            {0.5*ln(1-x) + x/(1-x) + 0.5*ln(1-y) + y/(1-y)};

            \end{axis}
        \end{tikzpicture}
        \caption{Contour plot of $\kappa_{\mathrm{KL}^*}^{(2)}$.}
        \label{fig:contour_plot_gf_KL}
    \end{subfigure}

    \caption{Contour plots of $\kappa_f^{(2)}$ for $f = f_\mathrm{KL}$ and $f = f_{\mathrm{KL}^\star}$. These plots illustrate the convexity of $\kappa_{\mathrm{KL}}^{(2)}$ and $\kappa_{\mathrm{KL}^*}^{(2)}$ (see Lemma~\ref{lem:cvx-kappaf}).}
    \label{fig:contour_plot_gf}
\end{figure}

The proof of this lemma is given in Section~\ref{sec:proof-G}. 
The non-quasi-convexity of $\kappa_{\mathrm{TV}}^{(d)}$ is obtained via a Taylor expansion at $0$ of the function $\kappa_{\mathrm{TV}}^{(d)}$. We illustrate this in Fig.~\ref{fig:contour_plot_gf_TV}, where we provide a contour plot of the function $\kappa_{\mathrm{TV}}^{(2)}$ with a zoom on the region $[0,0.2]^2$.
In Fig.~\ref{fig:contour_plot_gf}, we illustrate that the level sets of $\kappa_{\mathrm{KL}}^{(2)}$ and $\kappa_{\mathrm{KL}^*}^{(2)}$ are convex, which is consistent with Lemma~\ref{lem:cvx-kappaf}.

%%%%%%%%%%%%%%%%%%%%%%%%%%%%%%%%%%%

\section{Beyond the Gaussian model}
\label{sec:AdditiveModel}

In this section we consider the model:
\begin{equation}
  \label{eq:additif}
  Y= m(X+X',\,W)
  \quad\text{with}\quad X' \stackrel{\cl}{=} X + U,
\end{equation} 
 with $X, X', U$ independent $\R^d$-valued random variables,  $W$ another independent random variable on a measurable space $E$, and $m$ a measurable function defined on $\R^d\times E$ taking values in a measurable space.

This model extends the Gaussian model of the vector $(X,X',Y)$ in the previous section to more general distributions such that $X$ and $X'$ are always ordered in Csisz\'ar order w.r.t.\ $Y$. This is stated in the following proposition, whose proof 
relies
on the axioms of the Csisz\'ar order on: law invariance,  informational monotonicity and  conditional independence .

\begin{prop}[Csisz\'ar index for model~\eqref{eq:additif}]   \label{prop:ranking_Sf_using_markov}
  We have in the model~\eqref{eq:additif} that:
  \[
  X \csleq X' \quad \text{w.r.t.} \quad Y.
  \]
\end{prop}

\begin{proof}
 Without loss of generality, we can write $X'=\tilde X +U$, where $X, \tilde X, U$ are independent and $\tilde X$ is distributed as $X$.  
We have for $f\in \F$ that :
\[
S_f(X',Y) = S_f((\tilde X, X'),Y)
\ge S_f(\tilde X,Y) = S_f( X,Y) ,
\]
where we used that  conditionally on $X'$, the random variables $\tilde X$ and $Y$ are independent and then axiom~\ref{axiom:CharacterizationConditionalIndependence} for the first equality, the contraction property from axiom~\ref{axiom:InformationMonotonicity} for the inequality and 
the law invariance from axiom~\ref{axiom:LawInvariance} for the last inequality. 
\end{proof}

\begin{rem}[Global dependence and additive error models]
    \label{rem:+error}
   Consider the additive error models $Y= g(X) + \sigma \varepsilon$ and $Y'= g(X) + \sigma' \varepsilon$ in $\R$ where the noise  $\varepsilon$ is independent of $g(X)$ and $\sigma'>\sigma\ge 0$. In\cite[Theorem~1.8]{ansari2025ordering} the authors state that, provided  the distributions of $Y$ and $Y'$ have no atom, then $(X,Y')$ is less dependent than $(X, Y)$ for the conditional convex order. 
For the Csiszár order on additive error models, we may have either case as shown in the following examples.
\begin{enumerate}
    \item \textbf{Self-decomposable noise.} We recall that $\varepsilon$ is self decomposable if, for all $c\in [0, 1]$ the following equality in distribution holds:   $\varepsilon \stackrel{(d)}{=} c \varepsilon + Z_c$, where $Z_c$ and $\varepsilon$ are independent. For example stable distributions are self decomposable. Assuming the noise is self-decomposable, we get, with $c=\sigma/\sigma'$, that $(X,Y')$ is distributed as $(X, Y'')$ with $Y''= g(X) + \sigma \varepsilon + \sigma' Z_c= Y+ \sigma' Z_c $, where $X$, $\varepsilon$ and $Z_c$ are independent.  
    We get that:     
    $$
    S_f(X, Y')=S_f(X, Y+ \sigma' Z_c) \le 
    S_f(X, (Y, Z_c))=S_f(X, Y)
    $$ for all $f\in \F$, using the law invariance (axiom (i)) for the first equality, the information monotonicity (axiom (v)) for the inequality and 
the     conditional independence (axiom (vi))
    for the second equality. This means that 
   $(X,Y')$ is less dependent than $(X, Y)$ for the Csiszár order: 
       $(X,Y') \csleq (X, Y)$.

\item \textbf{A counter-example.}
    Let $g(X)$ be  uniform on $\{0,2\}$,  $\varepsilon$ uniform on $[0,1] \cup [2,3]$,  and $\sigma'=2> \sigma=1$ (notice that $Y$ and $Y'$ have no atom). Elementary computations yield, as $X$ and $\varepsilon$ are independent,   that $S_f(X,Y') = 2 S_f(X,Y)= (f(0) + f(2))/2$ for $f\in \F$,  that is, $(X,Y')$ is strictly more dependent than $(X, Y)$ for the Csiszár order: $(X,Y) \csls (X, Y')$.

    The intuition behind this counterexample is that larger noise increases the information we can recover on the input from the output. 
     Indeed, if $Y \in [2,3]$, then $g(X)$ is either $0$ or $2$, whereas $g(X) $ is perfectly dependent of $Y'$.
\end{enumerate}
Thus, each model must be carefully studied in order to give sufficient conditions for Csisz\'ar ordering. 
\end{rem}

As Lévy processes are stochastic processes with independent and stationary increments, we have the following direct application.

\begin{cor}[Application to Lévy Processes]
\label{cor:indices_PAIS} Let $X = (X_t)_{ t\geq 0}$ be a Lévy Process in $\R^d$. Set $Y=m(X_{t_0}, W)$ with $W$ a  random variable independent of $X$ and taking values in a 
measurable space $E$, and $m$ a measurable function defined on $\R^d\times E$ taking values in a measurable space. We have:
\[
X_s \csleq X_t 
\quad \text{w.r.t.} \quad Y,
\quad \text{for all}\quad 0\leq s\leq t\leq t_0.
\]
\end{cor}

We illustrate this result in the next example.
\begin{exple}[Application to Gaussian and Cauchy random variables]
\label{exple:levy_process_and_brownian_motion}
We consider the model in $\R^d$:
\begin{equation}
  \label{eq:additif-levy}
Y=X+X' + W \quad \text{with $X$, $X'$ and $W$ independent}.
\end{equation}
\begin{enumerate}[(i)]
    \item If $X$ and $X'$ are real Gaussian random variables with standard deviation $\sigma$ and $\sigma'$, then we have:
    \[
    X \csleq X' \quad \text{w.r.t.} \quad Y 
    \quad \iff \quad 
    \sigma \leq \sigma'
    \quad \iff \quad
    \corr(X,Y)^2  \leq \corr(X',Y)^2.
    \]
   This can also be seen as an extension of Corollary~\ref{cor:gaussian_case_dim_1}~\ref{enum:gaussian_dim1_case1}.

   \item If $X$ and $X'$ are real Cauchy random variables with scaling parameter $\alpha$ and $\alpha'$, then we have:
    \[
    X \csleq X' \quad \text{w.r.t.} \quad Y
    \iff \alpha \leq \alpha'.
    \]
 \end{enumerate}
Notice that we may consider arbitrary location parameters for $X$ and $X'$ and they do not affect the comparison of the Csisz\'ar indices. 
The Csisz\'ar order between the $\alpha$-stable $X$ and $X'$ in the additive model defined in \eqref{eq:additif-levy}, is entirely determined by the scale parameters of the distributions.
\end{exple}

In Corollary~\ref{cor:indices_PAIS}, there exists for the model~\eqref{eq:additif} a Csisz\'ar order for Lévy processes, which is determined by the ordering of variances in the Gaussian case and by the scale parameter in the Cauchy case. However, this structure does not extend in general as shown in the following example.

\begin{exple}[In general, there is no Csisz\'ar order for additive models]
\label{exple:order-depends-onf}
We consider  $X$  uniformly distributed over $\{1,2,3\}$,  $B$  a Bernoulli random variable with parameter $p \in [0,1]$ and $\tilde X$  uniformly distributed over $\{0,4\}$, with $X, B$ and $\tilde X$ independent. We consider the additive model 
\[
Y = X + X'
\]
with $  X' = B \tilde X + 2(1-B)$.
The total variation Csisz\'ar index between $X$ and $Y$, and between $X'$ and $Y$, are given by:
\[
S_{\mathrm{TV}}(X,Y) = \inv{3} \left(p+2 +\frac{2}{3} \left\vert 2p - 1 \right\vert - \frac{\vert 5p-4 \vert}{3} \right) 
\quad \text{and} \quad
S_{\mathrm{TV}}(X',Y) = \frac{p}{3} \left(8 -5p \right),
\]
and the Pearson $\chi^2$ Csisz\'ar index:
\[
S_{\mathrm{P}}(X,Y) = \frac{- 4 p^{2} + 6 p - 2}{p - 2}
\quad \text{and} \quad
S_{\mathrm{P}}(X',Y) = \frac{4}{3} \cdot
\]
Taking $p=0.7$ yields:
\[
S_{\mathrm{TV}}(X,Y) < S_{\mathrm{TV}}(X',Y) \quad \text{and} \quad S_{\mathrm{P}}(X,Y) > S_{\mathrm{P}}(X',Y).
\]
This example shows that, in general, there is no Csisz\'ar order in additive models. 
Obviously the squared correlation of $(X,Y)$ and $(X',Y)$ can always be ranked, however as shown in this example this is not the case for Csisz\'ar order even in the additive model. 
\end{exple}

\begin{exple}[Csisz\'ar order and convex order]
In the model~\eqref{eq:additif-levy}, one compares $S_f(X,Y)$ and $S_f(X',Y)$, where $X'$ is distributed as $X+U$, with $U$ and $X$ independent $\R^d$-valued random variables. Using invariance properties of the Csisz\'ar index, if $X$ and $U$ are integrable, we can, without loss of generality, assume that $\E[U]=0$ for the comparison of the Csisz\'ar index and thus we have $X \leqcv X'$. This implies, that $\Var(X) \leq \Var(X')$ whenever $X$ and $X'$ are square-integrable real-valued random variables, and therefore, provided that $Y$ is square-integrable:
\[
\corr(X,Y)^2= \frac{\Var(X)}{\Var(Y)} \leq \frac{\Var(X')}{\Var(Y)}=  \corr(X',Y)^2.
\]
However, the convex order $X \leqcv X'$ in the model~\eqref{eq:additif-levy}
does not always induce a Csisz\'ar order. In Example~\ref{exple:OrdreUniformeVsCorr}, which is a particular case of Example~\ref{exple:order-depends-onf} with $p=1$ and thus $B=1$, it is easy to verify that $X \lscv X'$ (and hence that $\corr(X,Y)^2 \leq \corr(X',Y)^2$), whereas $X' \csls X$ w.r.t.\ $Y$. Therefore, the convex order does not necessarily imply the Csisz\'ar order in the additive model.
\end{exple}

\section{Csisz\'ar index for Bernoulli random variables}
\label{sec:Bernoulli}

In this section we consider $X,X',Y$ Bernoulli random variables defined on a probability space which is supposed rich enough (for example $\Omega=[0, 1]$ with $\cf$ the Borel sets and $\P$  the Lebesgue measure). 
In this setup, we are able to fully characterize the sets $\cp_Y^+(X')$ and $\cp_Y^-(X')$ (defined in Definition~\ref{def:Csiszar sets}) in terms of the probabilities of the Bernoulli distributions and of their cross-moments $\E[XY] = \P(X=Y=1)$ and $\E[X'Y] = \P(X'=Y=1)$. More precisely, the set of all admissible parameters will be split into blue areas associated to  $\cp_Y^+(X')$, green areas associated to  $\cp_Y^-(X')$ and red areas where we know that $X$ and $X'$ rank differently according to the choice of $f$ in the $f$-divergence.

Unless the distribution of $(X',Y)$ is degenerate, see Proposition~\ref{prop:C+-cup-C-} in Section~\ref{sec:C+-cup-C-}, there exists a Bernoulli random variable $X$ such that $X$ and $X'$ are not comparable w.r.t.\ $Y$, so the Csisz\'ar order is not total.

\subsection{Preliminaries}
\label{sec:Bernoulli-intro}

\begin{figure}
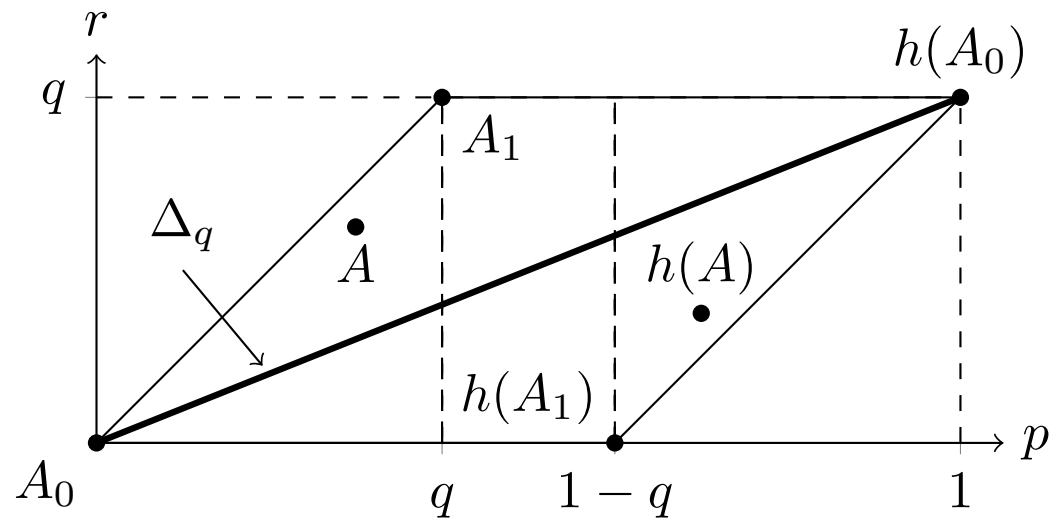

    \centering
    \begin{RegionPlotBernoulli}{0.55}{0.3}{0.4}
      {
      \addplot[solid, very thick] (0,0) -- (1,0.4) -- cycle;
      
      \addplot[dashed] (0.4,0) -- (0.4,0.4) -- cycle;
      \addplot[dashed] (0.6,0) -- (0.6,0.4) -- cycle;

      \draw[\ArrowStyle, thin, shorten >=2pt] (0.1,0.2) -- (0.2,0.08);
      \fill (0.1,0.2) circle (0) node[above] {$\Delta_q$};

      \addplot[solid] (0,0) -- (0.6,0) -- (1,0.4) -- (0.4,0.4) -- cycle;
      \fill (0.4,0.4) circle (0.01) node[below right] {$A_1$};
      \fill (0.6,0) circle (0.01) node[above left] {$h(A_1)$};

      \fill (0.3,0.25) circle (0.01) node[below] {$A$};
      \fill (0.7,0.15) circle (0.01) node[above] {$h(A)$};
      }
      {}
    \end{RegionPlotBernoulli}
    \caption{The domain $D_q$ given by the convex hull of $A_0=(0,0)$, $A_1=(q,q)$, $h(A_0)=(1, q)$ and $h(A_1)=(1-q, 0)$;
    the segment $\Delta_q$ from $A_0$ to $h(A_0)$ (in bold), 
    and an example of the involution $h$ applied to an arbitrary point $A$.} 
    \label{fig:definition_Dq}
\end{figure}

Let $X$, $X'$ and $Y$ be  Bernoulli random variables defined on the same probability space. Let $q=\P(Y=1)\in[0,1]$ be fixed. We define $D_q$ as the set of parameters $p$ and $r$ of the joint distribution of $(X,Y)$ that are compatible with $q$ in our model, where:
\[
p=\P(X=1) \quad \text{and} \quad r=\P(X=Y=1).
\]
We denote such a parameter by $A=(p,r)$, and we can easily see that:
\begin{equation}
    \label{eq:def-Dq}
    D_q=\{A=(p,r)\in [0,1]^2\,:\,(p+q-1)_+ \leq r \leq \min(p,  q) \}.
\end{equation}
The rhombus in Fig.~\ref{fig:definition_Dq} represents the set $D_q$. Consider the involution $h=h_q$ of $D_q$, defined by:
\begin{equation}
    \label{eq:def-hq}
h(p,r) = (1-p, q-r). 
\end{equation}
By the invariance of the Csisz\'ar index with respect to injective transformations of $X$ (or $Y$) (see Proposition~\ref{prop:invariance_csiszar_index}), the Csiszár index of $(X,Y)$ and of $(\tilde X, Y)$, where $\tilde X=1-X$, coincide. 
Notice that $h$ leaves $D_q$ invariant. 

We now introduce several subsets of $D_q$. The diagonal segment of $D_q$: 
\[
\Delta_q=\{A=(p,r)\in D_q\,:\, r=pq\}
\]
is invariant under $h$. It corresponds to the case where $X$ and $Y$ are independent; this includes the cases $p\in\{0,1\}$ and $q\in\{0,1\}$. Notice that the full set $D_q$ reduces to its diagonal, $D_q=\Delta_q$, if and only if $q\in \{0,1\}$. We also define $A_1=(q,q)$ and the set $E_q=\{A_1, h(A_1)\}=\{(q,q),(1-q,0)\}$, corresponding respectively to the degenerate cases $X=Y$ a.s.\ and $X=1-Y$ a.s..

The map $A \mapsto  S_f(X, Y)$  is equal to its minimal value, given by   0, (at least) on the set $\Delta_q$, see \cite[Corollary 7]{butucea2026csisz}; we shall see in Proposition~\ref{prop:bernoulli_cm}
that the map  is equal to its maximal value, given by  $qf^*(q)+(1-q) f^*(1-q)+ 2 q (1-q)f(0)$,  (at least) on $E_q$. We have:
\[
0 =\max_{A\in \Delta_q} S_f(X,Y)=\min_{A\in D_q} S_f(X,Y) \leq \max_{A\in D_q} S_f(X, Y)=\min_{A\in E_q} S_f(X,Y).
\]

As described in Section~\ref{sec:ConnectionConvexOrder}, the Csiszár order relies on the convex order of density ratios.
The density ratio $R_{(X,Y)}$ in \eqref{eq:R} is a discrete random variable supported on at most four atoms, corresponding to the four possible value pairs of $(X,Y)$, namely $(0,0)$, $(1,1)$, $(1,0)$ and $(0,1)$. On the set $\Delta_q$, the variables $X$ and $Y$ are independent, and therefore the density ratio $R_{(X,Y)}$ is identically equal to $1$. 

We define similarly $A'=(p', r') \in D_q$ and $Z'$ with $X$ replaced by $X'$. Thus the model is fully characterized by the parameters $A=(p,r)$, $A'=(p',r')$ and $q$. 

Let us consider $A'=(p',r')\in D_q$ fixed, and consider the regions:
\begin{align*}
C^+_{q,A'}
& = \{A=(p,r)\in D_q\, \colon\, X \csgeq X' \quad\text{w.r.t.}\quad Y\} ,\\
%=\{A=(p,r)\in D_q\, \colon\, Z \, \geqcv \,  Z' \},\\
C^-_{q,A'}
& = \{A=(p,r)\in D_q\, \colon\, X \csleq X' \quad\text{w.r.t.}\quad Y\} .
%=\{A=(p,r)\in D_q\, \colon\, Z \, \leqcv \, Z'\}.
\end{align*}
These sets correspond to the sets $\cp_Y^+(X')$ and $\cp_Y^-(X')$ restricted to Bernoulli random variables, \emph{i.e.\ }to the set of Bernoulli distributions $X$ such that $S_f(X,Y) \geq S_f(X',Y)$ and $S_f(X,Y) \leq S_f(X',Y)$ for all $f\in \F$, respectively. 
By construction, we have that any point $A'$ belongs to both sets: $C^+_{q,A'}\cap C^-_{q,A'}$.
Since the Csiszár index is preserved under $X \mapsto \tilde X=1-X$, the sets $C^+_{q,A'}$ and $C^-_{q,A'}$ are also invariant under $h$, that is:
\begin{equation}
    \label{eq:sym-C}
    C^\pm_{q, A'}= h( C^\pm_{q, A'})=C^\pm_{q, h(A')}.
\end{equation}
In particular, $h(A')\in C^+_{q,A'}\cap C^-_{q,A'}$.

In the following section, we give an explicit description of $C^-_{q, A'}$. Many instances of $C^+_{q, A'}$ (in blue) and $C^-_{q, A'}$ (in green), and of areas (in red) where no Csisz\'ar order exists, are represented in Fig.~\ref{fig:description_cp}, \ref{fig:degenerate_cases}, \ref{fig:description_cp-cas-1demi}, and in the appendix in Fig.~\ref{fig:description_cm} and \ref{fig:A'-on-boundary}.

\subsection{Characterization of the Csiszár Order} 
In this subsection, we characterize the sets $C^-_{q,A'}$ and $C^+_{q,A'}$, which completely determine the Csiszár order between $X$ and $X'$. More precisely, they characterize the regions where $X\csleq X'$, $X\csgeq X'$, and where $X$ and $X'$ are not comparable.

We proceed in two steps. Proposition~\ref{prop:bernoulli_cp} gives an explicit description of $C_{q,A'}^-$ as the convex hull of an explicit set of points. Proposition~\ref{prop:bernoulli_cm} in the appendix provides a description of $C_{q,A'}^+$ whose geometry is more complex.

Notice %$\Delta_q$ is the convex hull of $\{A_0, h_q(A_0)\}$ and 
that the diagonal segment $\Delta_q$ is included in $ C^-_{q,A'}$ as $S_f(X,Y)=0$ is minimal for all $(p,r)\in \Delta_q$. 
Recall also that $\Delta_q\neq D_q$ if and only if $q\in (0, 1)$.
We denote by $\Cv(\ca)$ the convex hull of a set $\ca\subset  \R^2$. The lengthy proof of the next proposition is given in Section~\ref{sec:proof-bernoulli_cp}. It relies on Proposition~\ref{prop:equivalence_ConvexOrder_CsiszarOrder}, giving that the Csisz\'ar order between $X$ and $X'$ holds if and only if the convex order between $Z$ and $Z'$ holds. 
Since $Z$ and $Z'$ are discrete, one can derive necessary and sufficient conditions on the parameters $p$ and $r$ (with $p'$, $r'$, and $q$ fixed) for the convex order of $Z$ and $Z'$, and hence for the Csisz\'ar order of $X$ and $X'$ w.r.t. $Y$, to hold.

\begin{prop}[Description of $C_{q,A'}^-$]
\label{prop:bernoulli_cp}
Let $q\in [0, 1]$ and $A'=(p',r')\in D_q$. 
Set $\beta= \min(q, 1-q)$ and let $I_q=[\beta, 1- \beta]$ be the interval between $q$ and $(1-q)$. We consider the set $\ca_0$ defined by:
\begin{enumerate}[(i)]
\item \label{it:1-D1} $\ca_0=\{A_0\}$ 
if $A'\in \Delta_q$ (that is, $X'$ and $Y$ independent).  
\item \label{it:2-D1} $\ca_0=\{A_0, A'\}$   if $A'\not\in \Delta_q$  and $p' \in I_q$.
\item \label{it:3-D1} $\ca_0=\{A_0, A',A''\}$ otherwise (and thus $p',q\in (0, 1)$), with   $A''=(p'', r'')$ and   $\alpha'= \min(p', 1-p')$:
\begin{equation}
\label{eq:bernoulli_cp''}
p'' = \frac{\alpha'(1-\alpha')}{2\beta - \alpha'} \quad \text{and} \quad 
r'' =qp''  +  \varepsilon   \, p''\, (qp'-r') \frac{\beta}{\alpha' (1- \beta)} \quad\text{where}\quad\varepsilon=\begin{cases}
+1 & \text{if $p'< \beta$,}\\ 
-1 & \text{if $p'>1- \beta$.}
\end{cases}
\end{equation}
Then, we have with $\ca=\ca_0\cup h(\ca_0)$ and $\Cv(\ca)$ is the (closed) convex hull  of $\ca$ that.
\[
C_{q,A'}^-=\Cv(\ca).
\]
\end{enumerate}

\end{prop}

The set $\ca=\ca_0\cup h(\ca_0)$ contains 2 distinct points in case (i), 4 in case (ii) and  6 in case (iii). %; and the points  $\ca$ belongs to $D_q$ and are in convex position.

\begin{figure}[htbp]
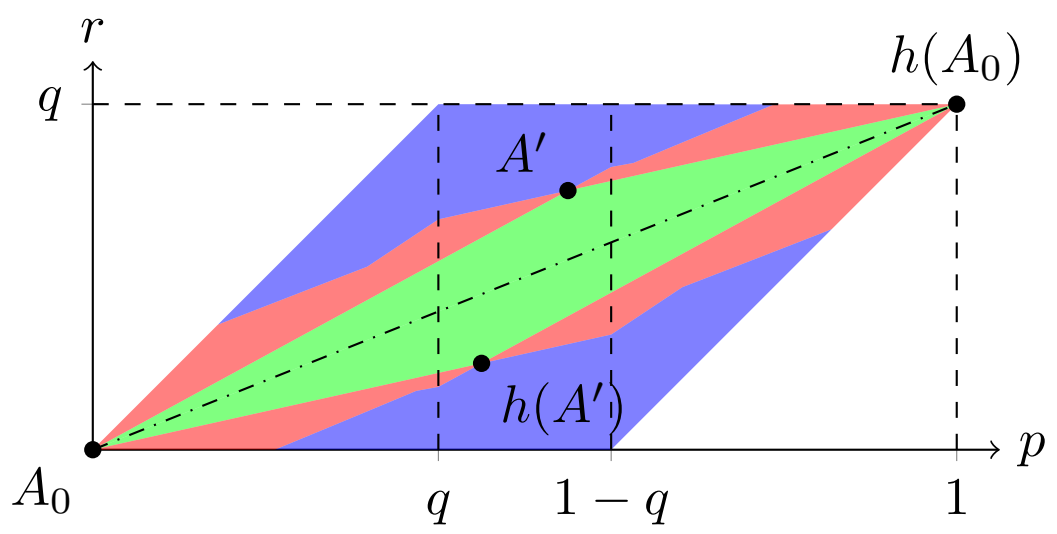
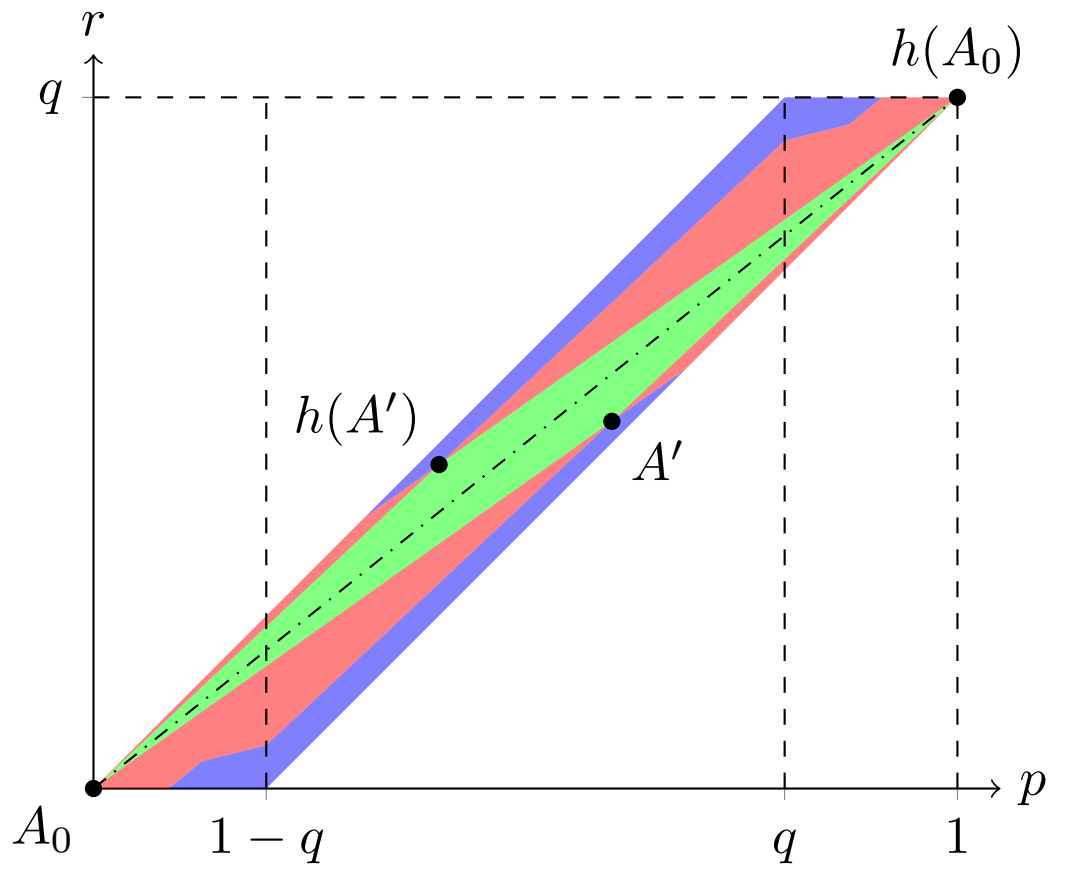
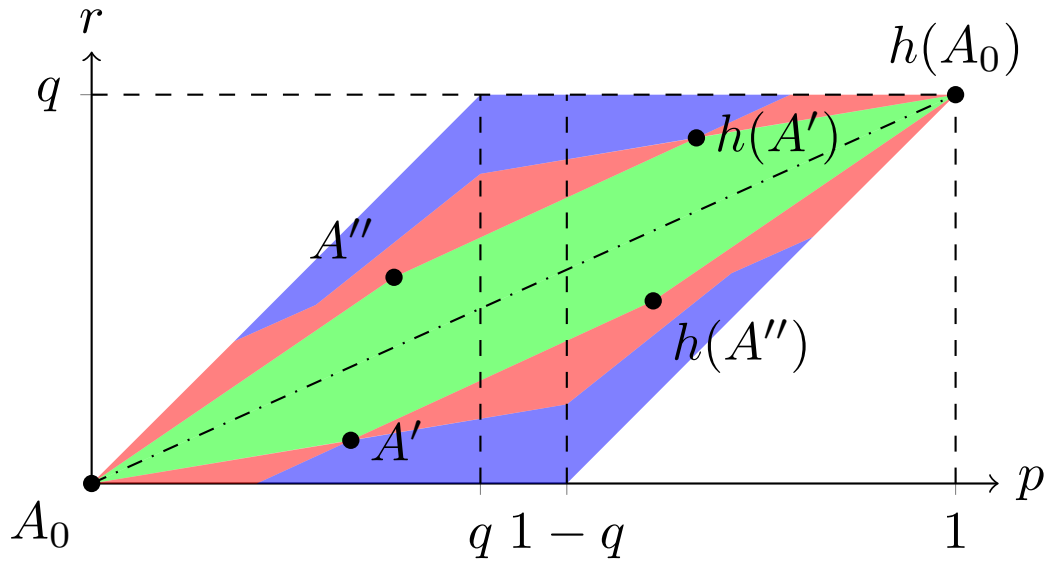
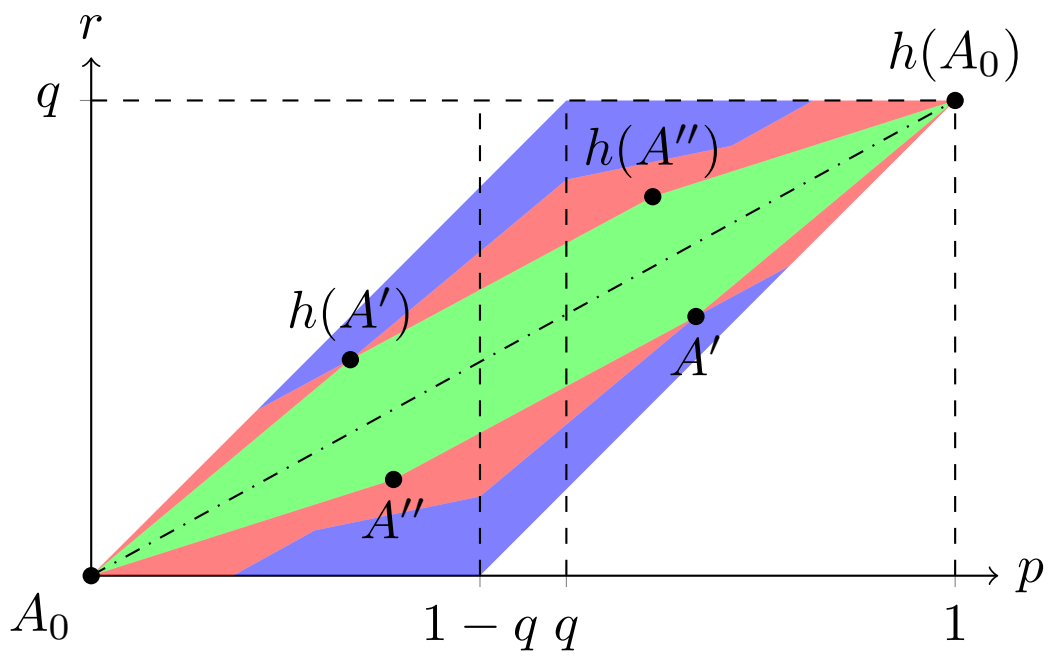

    \centering

\begin{subfigure}[c]{0.49\textwidth}
\centering
\begin{RegionPlotBernoulli}{0.55}{0.3}{0.4}
{
  \fill (\Pp,\Rp) circle (0.01) node[above left] {$A'$};
  \fill (\hPp,\hRp) circle (0.01) node[below right] {$h(A')$};
}
{
  \path[fill=\ColorCm, opacity=\OpacityPlot] (0,0) -- (0.55,0.3) -- (1,0.4) -- cycle;
    
  \path[fill=\ColorOut, opacity=\OpacityPlot]
    (0,0) -- plot[domain=0.144765682192536:0.3180294850975073, smooth, samples = \nsamples](\x,{0.4*\x + 0.0959999999999999*(1 - \x)/(1.09 - \x)}) -- (0.4, 0.2666666666666666) -- (0.55,0.3)  -- cycle;

  \path[fill=\ColorOut, opacity=\OpacityPlot]
    plot[domain=0.6257804885470349:0.789354390525168, smooth, samples = \nsamples](\x,{0.4*\x + 0.0959999999999999*\x/(\x + 0.11)})
    -- (1,0.4) 
    -- (0.55,0.3)
    -- (0.6, 0.3272727272727272)
    -- cycle;
    
  \path[fill=\ColorCp, opacity=\OpacityPlot]
    plot[domain=0.144765682192536:0.3180294850975073, smooth, samples = \nsamples](\x,{0.4*\x + 0.0959999999999999*(1 - \x)/(1.09 - \x)}) -- (0.4, 0.2666666666666666) 
    -- (0.4, 0.2666666666666666)
    -- (0.55,0.3)
    -- (0.6, 0.3272727272727272)
    -- plot[domain=0.6257804885470349:0.789354390525168, smooth, samples = \nsamples](\x,{0.4*\x + 0.0959999999999999*\x/(\x + 0.11)})
    -- (0.4,0.4)
    -- cycle;   
}
\end{RegionPlotBernoulli}

\caption{Case: $p' \in I_q$ ($p' = 0.55$, $r' = 0.3$, $q=0.4$).}
\label{fig:p_prime_Iq_q_leq_1-q}
\end{subfigure}
\hfill
\begin{subfigure}[c]{0.49\textwidth}
\centering
\begin{RegionPlotBernoulli}{0.6}{0.425}{0.8}
  {
    \fill (0.4,0.8-0.425) circle (0.01) node[above left] {$h(A')$};
    \fill (0.6,0.425) circle (0.01) node[below right] {$A'$};
  }
  {
    \path[fill=\ColorCm, opacity=\OpacityPlot]
    (0,0) -- (0.6,0.425) -- (1,0.8) -- cycle;
    
    \path[fill=\ColorOut, opacity=\OpacityPlot]
    (0,0)
    -- plot[domain=0.0871011991847576:0.12554373534619723, smooth, samples = \nsamples](\x,{0.8*\x - 0.088*(1 - \x)/(1.24 - \x)})
    -- (0.2,0.05)
    -- (0.6,0.425)
    -- cycle;

    \path[fill=\ColorOut, opacity=\OpacityPlot]
    (0.6,0.425)
    -- (0.685714285714286, 0.4857142857142861)
    -- (1,0.8)
    -- cycle;
    
    \path[fill=\ColorCp, opacity=\OpacityPlot]
    (0.2,0)
    -- plot[domain=0.0871011991847576:0.12554373534619723, smooth, samples = \nsamples](\x,{0.8*\x - 0.088*(1 - \x)/(1.24 - \x)})
    -- (0.2,0.05)
    -- (0.6,0.425)
    -- (0.685714285714286, 0.4857142857142861)
    -- cycle;   
  }
\end{RegionPlotBernoulli}

\caption{Case: $p' \in I_q$ ($p' = 0.6$, $r' = 0.425$, $q=0.8$).} 
\label{fig:p_prime_Iq_1-q_leq_q}
\end{subfigure}
\vspace{0.5cm}

\begin{subfigure}[c]{0.49\textwidth}
\centering

\begin{RegionPlotBernoulli}{0.3}{0.05}{0.45}
  {
    \fill (0.3,0.05) circle (0.01) node[right] {$A'$};
    \fill (0.65,0.21136363636363642) circle (0.01) node[below right] {$h(A'')$};

    \fill (1-0.3,0.45-0.05) circle (0.01) node[right] {$h(A')$};
    \fill (1-0.65,0.45-0.21136363636363642) circle (0.01) node[above left] {$A''$};
  }
  {
    \path[fill=\ColorCm, opacity=\OpacityPlot]
    (0,0) -- (0.3,0.05) -- (0.65,0.21136363636363642) -- (1,0.45) -- cycle;
    
    \path[fill=\ColorOut, opacity=\OpacityPlot]
    (0,0) -- plot[domain=0.191226924268372:0.3, smooth, samples = \nsamples](\x,{0.45*\x - 0.0935*(1 - \x)/(1.07 - \x)}) -- cycle;
    
    \path[fill=\ColorCp, opacity=\OpacityPlot]
    plot[domain=0.191226924268372:0.3, smooth, samples = \nsamples](\x,{0.45*\x - 0.0935*(1 - \x)/(1.07 - \x)})
    -- (0.55,0.09166666666666669)
    -- plot[domain=0.7405124837953327:0.835889894354067, smooth,samples = \nsamples](\x,{0.45*\x - 0.0935*\x/(\x + 0.03)})
    -- (0.55,0)
    -- cycle;
    
    \path[fill=\ColorOut, opacity=\OpacityPlot]
    (0.3,0.05) -- (0.55,0.09166666666666669) -- plot[domain=0.7405124837953327:0.835889894354067, smooth,samples = \nsamples](\x,{0.45*\x - 0.0935*\x/(\x + 0.03)}) -- (1,0.45) -- (0.65,0.21136363636363642) -- cycle;
  }
  
\end{RegionPlotBernoulli}

  \caption{Case: $p' < \beta$ ($p' = 0.3$, $r' = 0.05$, $q = 0.45$). }
  \label{fig:p_prime_<=_beta}
\end{subfigure}
\hfill
\begin{subfigure}[c]{0.49\textwidth}
\centering

\begin{RegionPlotBernoulli}{0.7}{0.3}{0.55}
  {
    \fill (0.7,0.3) circle (0.01) node[below] {$A'$};
    \fill (0.65,0.43863636363636366) circle (0.01) node[above] {$h(A'')$};

    \fill (0.3,0.25) circle (0.01) node[above] {$h(A')$};
    \fill (0.35,0.1113636363636364) circle (0.01) node[below] {$A''$};

  }
  {
    \path[fill=\ColorCm, opacity=\OpacityPlot]
    (0,0) -- (0.35,0.1113636363636364) -- (0.7,0.3) -- (1,0.55) -- cycle;
    
    \path[fill=\ColorOut, opacity=\OpacityPlot]
    (0,0) 
    -- plot[domain=0.164110105645933:0.25948751620466726, smooth,samples = \nsamples](\x,{0.55*\x - 0.0935*(1 - \x)/(1.03 - \x)})
    -- (0.45,0.09166666666666667)
    -- (0.7,0.3)
    -- (0.35,0.1113636363636364)
    -- cycle;

    \path[fill=\ColorOut, opacity=\OpacityPlot]
    plot[domain=0.7:0.808773075731628, smooth,samples = \nsamples](\x,{0.55*\x - 0.0935*\x/(\x + 0.07)})
    -- (1,0.55) 
    -- cycle;
    
    \path[fill=\ColorCp, opacity=\OpacityPlot]
    plot[domain=0.164110105645933:0.25948751620466726, smooth,samples = \nsamples](\x,{0.55*\x - 0.0935*(1 - \x)/(1.03 - \x)})
    -- (0.45,0.09166666666666667)
    -- (0.7,0.3)
    -- plot[domain=0.7:0.808773075731628, smooth,samples = \nsamples](\x,{0.55*\x - 0.0935*\x/(\x + 0.07)})
    -- (0.45,0)
    -- cycle;
  }
\end{RegionPlotBernoulli}
\caption{Case: $p' > 1-\beta$ ($p' = 0.7$, $r' = 0.3$, $q = 0.55$).}
\label{fig:p_prime_>=_1-beta}
\end{subfigure}
\caption{Description of $C^-_{q, A'}$ (in green) and $C^+_{q, A'}$ (in blue) with $A_0 = (0,0)$, $A'=(p', r')$, $\beta=\min (q, 1-q)$ and $I_q=[\beta, 1-\beta]$.}
\label{fig:description_cp}
\end{figure}

\begin{figure}[htbp]
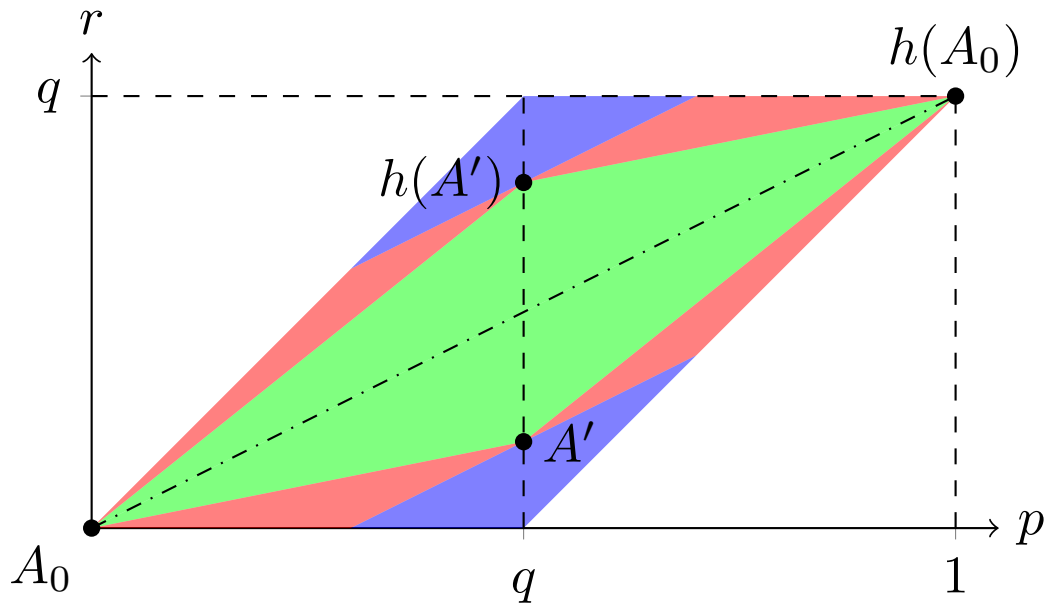
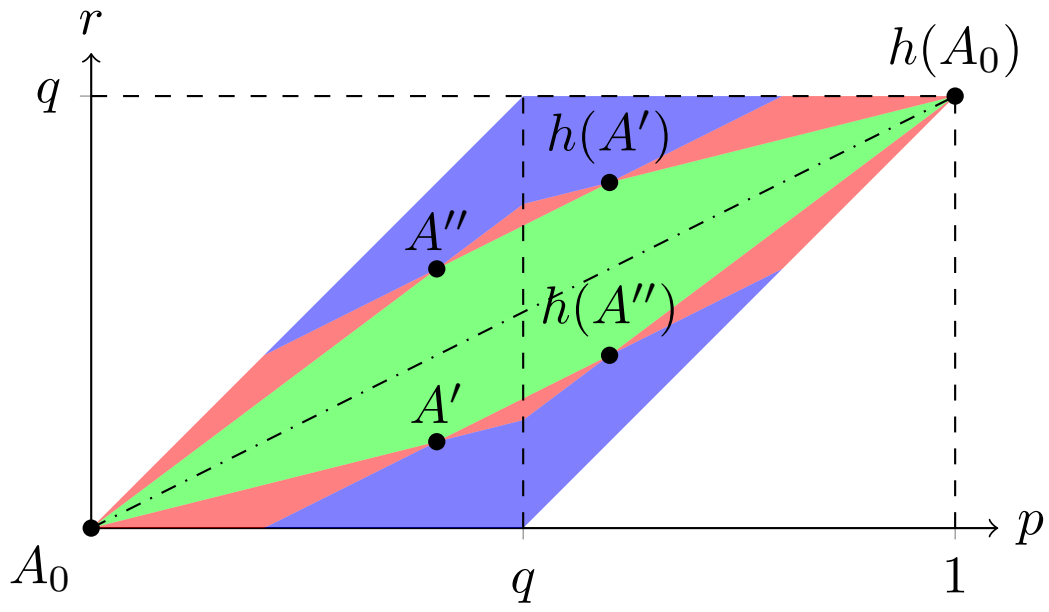

\begin{subfigure}[c]{0.49\textwidth}
\centering
\begin{RegionPlotBernoulli}{0.5}{0.1}{0.5}
  {
    \fill (0.5,0.1) circle (0.01) node[right] {$A'$};
    \fill (0.5,0.4) circle (0.01) node[left] {$h(A')$};
  }
  {
    \path[fill=\ColorCm, opacity=\OpacityPlot]
    (0,0) -- (0.5,0.1) -- (1,0.5) -- cycle;
    
    \path[fill=\ColorOut, opacity=\OpacityPlot]
    (0,0) -- (0.3,0) -- (1/2,0.1) -- cycle;

    \path[fill=\ColorOut, opacity=\OpacityPlot]
    (1/2,0.1) -- (0.7,0.2) -- (1,1/2) -- cycle;
    
    \path[fill=\ColorCp, opacity=\OpacityPlot]
    (0.3,0) -- (1/2,0) -- (0.7,0.2) -- cycle;   
  }
  
\end{RegionPlotBernoulli}

\caption{Case: $p'=1/2$ ($r' = 0.1$).}
\label{fig:q_1/2_p_prime_1/2}
\end{subfigure}
\hfill
\begin{subfigure}[c]{0.49\textwidth}
\centering
\begin{RegionPlotBernoulli}{0.4}{0.1}{0.5}
  {
    \fill (0.4,0.1) circle (0.01) node[above] {$A'$};
    \fill (0.6,0.4) circle (0.01) node[above] {$h(A')$};

    \fill (0.4,0.3) circle (0.01) node[above] {$A''$};
    \fill (0.6,0.2) circle (0.01) node[above] {$h(A'')$};
  }
  {
    \path[fill=\ColorCm, opacity=\OpacityPlot]
    (0,0) -- (0.4,0.1) -- (0.6,0.2) -- (1,0.5) -- cycle;
    
    \path[fill=\ColorOut, opacity=\OpacityPlot]
    (0,0) -- (0.2,0) -- (0.4,0.1) -- cycle;

    \path[fill=\ColorOut, opacity=\OpacityPlot]
    (0.4,0.1) -- (1/2,1/8) -- (0.6,0.2) -- cycle;

    \path[fill=\ColorOut, opacity=\OpacityPlot]
    (0.6,0.2) -- (0.8,0.3) -- (1,1/2) -- cycle;
    
    \path[fill=\ColorCp, opacity=\OpacityPlot]
    (0.2,0) -- (1/2,0) -- (0.8,0.3) -- (0.6,0.2) -- (1/2,1/8) -- (0.4,0.1) -- cycle; 
  }
\end{RegionPlotBernoulli}
\caption{Case: $p'\neq 1/2$ ($p'=0.4$, $r' = 0.1$).}
\label{fig:q1/2_p_prime_ne_1/2}
\end{subfigure}
\caption{Description of $C^+_{q, A'}$ (in green) and $C_{q,A'}^-$ (in blue) with $A'=(p', r')$ and $q=1/2$.}
\label{fig:description_cp-cas-1demi}
\end{figure}

\begin{figure}[htbp]
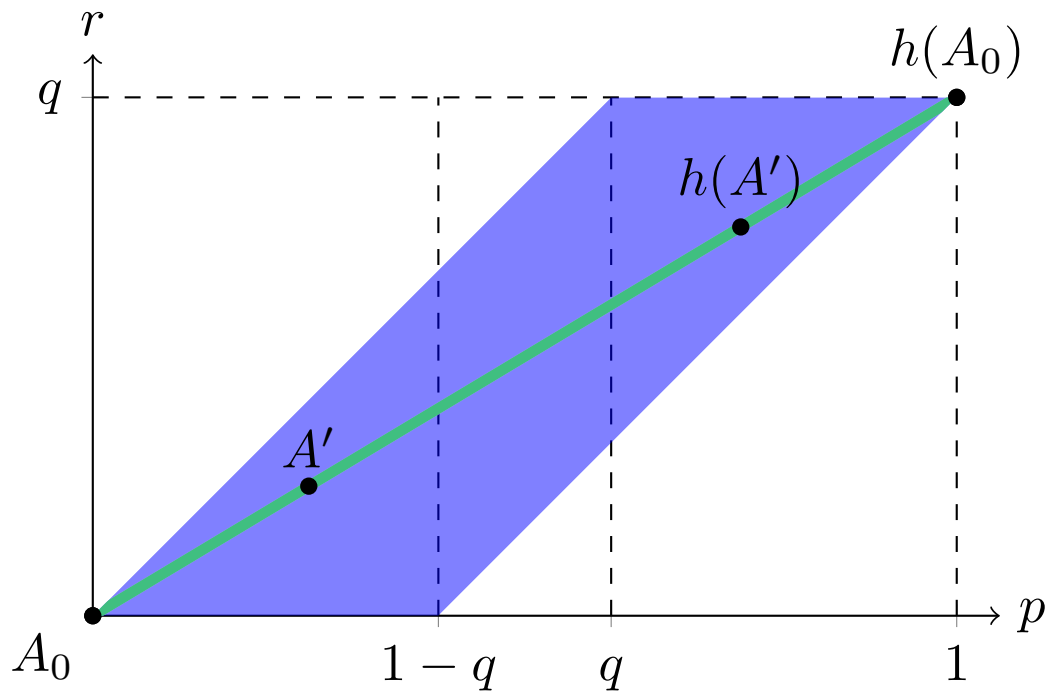
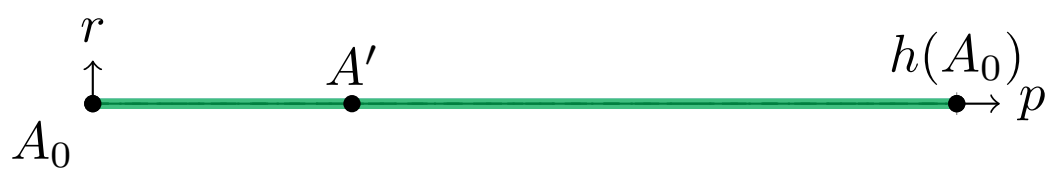
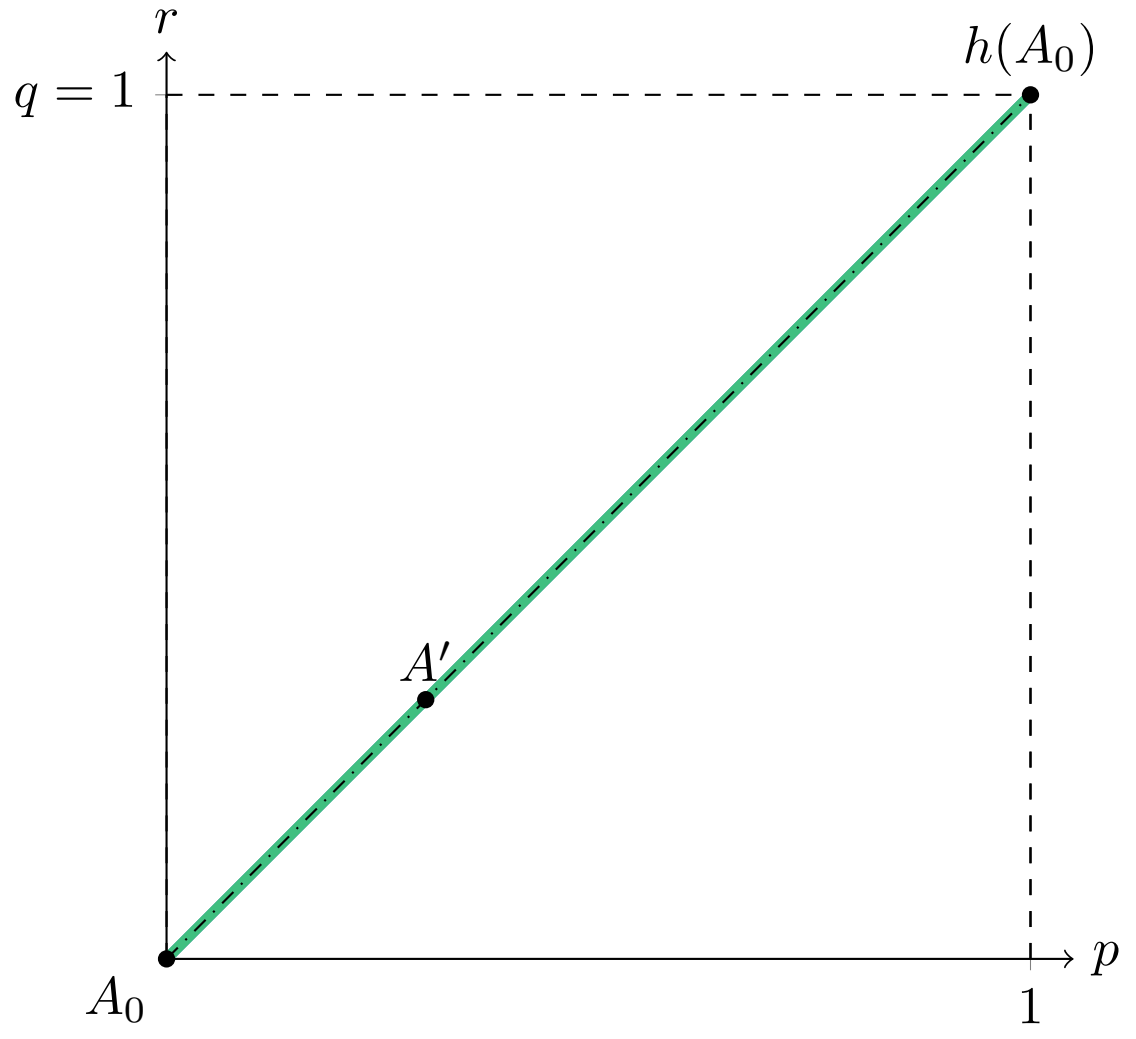

\centering

% --- Colonne gauche ---
\begin{minipage}{0.49\textwidth}
\centering

\begin{subfigure}{\textwidth}
\centering
  \begin{RegionPlotBernoulli}{0.25}{0.15}{0.6}
    {
    \begin{scope}
    \clip (0,0) -- (0.4,0) -- (1,0.6) -- (0.6,0.6) -- cycle;

    \path[fill=\ColorCp, opacity=\OpacityPlot]    (0,0) -- (0.4,0) -- (1,0.6) -- (0.6,0.6) -- cycle;

    \addplot[white, line width=2pt] coordinates {(0,0) (1,0.6)};
    \addplot[\ColorCp, opacity=\OpacityPlot, line width=2pt] coordinates {(0,0) (1,0.6)};
    \addplot[\ColorCm, opacity=\OpacityPlot, line width=2pt] coordinates {(0,0) (1,0.6)};
    \end{scope}

    \fill (1,0.6) circle (0.01) node[above] {};
    \fill (0,0) circle (0.01) node[above] {};
    
    \fill (0.25,0.15) circle (0.01) node[above] {$A' $};
    \fill (0.75,0.45) circle (0.01) node[above] {$h(A')$};

    }{}
  \end{RegionPlotBernoulli}
\caption{Case: $X'\indep Y$ ($p' = 0.25$, $q = 0.6$,  $r' = p'q = 0.15$).}
\end{subfigure}

\vspace{0.5cm}

\begin{subfigure}{\textwidth}
\centering

\begin{RegionPlotBernoulli}{0.3}{0}{0}
  {
    \pgfplotsset{xtick={1}, xticklabels={$1$}}
    \pgfplotsset{ytick=\empty, xticklabels=\empty}
    \addplot[\ColorCp, opacity=\OpacityPlot, line width=2pt] coordinates {(0,0) (1,0)};
    \addplot[\ColorCm, opacity=\OpacityPlot, line width=2pt] coordinates {(0,0) (1,0)};
    \fill (0.3,0) circle (0.01) node[above] {$A'$};

    \fill (1,0) circle (0.01) node[above] {};
    \fill (0,0) circle (0.01) node[above] {};
  }
  {}
\end{RegionPlotBernoulli}
\caption{Case $q=0$ ($p'=0.3$ and thus $r'=0$).}
\end{subfigure}

\end{minipage}
\hfill
% --- Colonne droite ---
\begin{minipage}{0.49\textwidth}
\centering

\begin{subfigure}{\textwidth}
\centering
\begin{RegionPlotBernoulli}{0.3}{0.3}{1}
  {
    \pgfplotsset{xtick={1}, xticklabels={$1$}}
    \pgfplotsset{ytick={1}, yticklabels={$q=1$}}
    \fill (0.3,0.3) circle (0.01) node[above] {$A' $};
  }
  {
    \addplot[white, line width=2pt] coordinates {(0,0) (1,1)};
    \addplot[\ColorCp, opacity=\OpacityPlot, line width=2pt] coordinates {(0,0) (1,1)};
    \addplot[\ColorCm, opacity=\OpacityPlot, line width=2pt] coordinates {(0,0) (1,1)};
  }
\end{RegionPlotBernoulli}
\caption{Case: $q=1$ ($p' =0.3$ and thus $r'=p'$).}
\end{subfigure}

\end{minipage}

\caption{Description of $C^-_{q, A'}$ (in green) and $C^+_{q, A'}$ (in blue) with $A'=(p', r')$ in the degenerated cases: $A'\in \Delta_q$ .}
\label{fig:degenerate_cases}
\end{figure}

Instances of the case~\ref{it:1-D1}, $A'\in \Delta_q$, are represented in Fig.~\ref{fig:degenerate_cases} (together with the degenerate cases $q\in \{0,1\}$); notice that $C_{q,A'}^-= \Delta_q$ if and only if $A'\in \Delta_q$. Instances of the case~\ref{it:2-D1}, $A'\not\in \Delta_q$ and $p'\in I_q$, are represented in Fig.~\ref{fig:p_prime_Iq_q_leq_1-q} and~\ref{fig:p_prime_Iq_1-q_leq_q} for $q\neq 1/2$, in Fig.~\ref{fig:q_1/2_p_prime_1/2} for $q=1/2$ (and thus $p'=1/2$), and in Fig.~\ref{fig:r_prime_eg_p_prime+q-1>0}
% and~\ref{fig:r_prime=0_p_prime_Iq}
where $A'\in \partial D_q$ and Fig.~\ref{fig:p_prime_eg_1-q_r_prime_0} where $A'\in E_q \subset \partial D_q$. Instances of the case~\ref{it:3-D1}, $A'\not\in \Delta_q$ and $p'\not\in I_q$, are represented in Fig.~\ref{fig:p_prime_<=_beta} and~\ref{fig:p_prime_>=_1-beta} for $q\neq 1/2$ in Fig.~\ref{fig:q1/2_p_prime_ne_1/2} for $q=1/2$ and in Fig.~\ref{fig:r_prime=0_p_prime<=1-q}
%\ref{fig:r_prime=p_prime+q-1>0,p_prime>=q}, 
and \ref{fig:r_prime=0,q<=1/2} for $A'\in \partial D_q$.

\begin{rem}[Maximal element]
    \label{rem:ScarsiniBernoulli}
    Following~\cite{scarsini1990ordering}, we can ask if there exists a maximal element $X'$ for the Csiszár order with a given marginal distribution, that is, given $q, p_0\in [0, 1]$ find $A'=(p_0, r')\in D_q$, such that for any $A=(p_0, r)\in D_q$, one has $X \csleq X'$ w.r.t.\ $Y$. This amounts to check that the segment $ \{p_0\}\times \R_+ \cap D_q$ belongs to $C_{q,A'}^-$; this is equivalent to $A'\in \partial D_q$ with $r'\geq p_0q$ if $(p_0-1/2)(q-1/2) \geq 0$ and $r'\leq p_0q$ if $(p_0-1/2)(q-1/2) \leq 0$. Note that if $p_0 = 1/2$ or $q=1/2$, then the maximum is no longer unique, and any $A'=(p_0, r')\in \partial D_q$ is a maximum. If $q \in \{0, 1\}$, then $\{p_0\}\times \R_+ \cap D_q$ reduces to a singleton and the maximum is unique.
    This coincides with Theorem~4.1 in~\cite{scarsini1990ordering}. 
\end{rem}

\begin{rem}[Comparison of the Csiszár order and the conditional convex order]
  \label{rem:Connection_Ccvx_Bernoulli}
According to Remark~\ref{rem:Connection_Ccvx_Gaussian}, in the Gaussian model, the ranking of the input variables $X$ and $ X'$ w.r.t.\ the one-dimensional output variable $Y$
is the same for the Csiszár order and the conditional convex order introduced by~\cite{ansari2025ordering}. This result does not hold in general: for the Bernoulli model (as the interested reader can check using Proposition 2.5 therein) taking $q=1/2$, $p=p'=1/4$, $r=0.1$, and $r'=0.15$, we get that  
 $X \csleq X'$ w.r.t.\ $Y$ whereas they are not comparable for the 
 conditional convex order. 
\end{rem}

\medskip

From the characterization of $C_{q,A'}^-$, one can derive the following description of $C_{q,A'}^+$, with $A'=(p', r')$:
\begin{equation}
    \label{eq:C-=invertC+}
C^+_{q,A'}=\left\{A= (p,r) \in D_q\, \colon\, A'\in C^-_{q,A} \right\}.
\end{equation}
From this representation it is elementary (but tedious) to derive Proposition~\ref{prop:bernoulli_cm} in Section~\ref{sec:description_Cm} on the description of $C^+_{q,A'}$ which is more complex than that of $C^-_{q,A'}$. The boundaries of $C^+_{q,A'}$ are not lines in general.
See Fig.~\ref{fig:description_ensembles_Cp} for an illustration of the boundaries, especially the curves.
For the particular case $q = 1/2$, the curves become lines and the set $C^+_{q,A'}$ becomes a set delimited by some line segments, see Section~\ref{sec:saturated_constraints} and Fig.~\ref{fig:description_cp-cas-1demi}.

In some situations, the set $C^+_{q,A'}$ is contained in the boundary of $D_q$. The following lemma provides a necessary and sufficient condition for this to occur.
The proof of the next result is in Section~\ref{sec:proof:Cm_boundaries}.

\begin{lem}[When $C_{q,A'}^+$ is a union of segments]
  \label{lem:Cm_boundaries}
  Let  $A' =(p',r') \in D_q \setminus (\Delta_q \cup E_q)$ (and thus $q\in (0,1)$). We have the following equivalence:
  \begin{equation}
    C_{q,A'}^+ \subset \partial D_q \iff A' \in \partial D_q.
  \end{equation}
\end{lem}
Notice that $C_{q,A'}^+ \cap \{A \in D_q \, \colon \, p \leq \min(q,1-q)\}$ or $C_{q,A'}^+ \cap \{A \in D_q \, \colon \, p \geq \max(q,1-q)\}$ might be a subset of $\partial D_q$ even though  $A' \notin \partial D_q$, see an instance in Fig.~\ref{fig:ber_cm_rg_segment}.

\subsection{Equality and Totality of the Csiszár Order}

In this subsection, we fix a pair of Bernoulli random variables $(X',Y)$ and characterize the set of parameters of $X$ for which $X\cseq X'$ w.r.t.\ $Y$. We also investigate the situations in which the Csiszár order is total (recall $Y$ is fixed), that is, when every pair of Bernoulli random variables $X$ and $X'$ is comparable w.r.t.\ $Y$, \emph{i.e.\ }either $X\csleq X'$ or $X\csgeq X'$.

Let us mention that $X\cseq X'$ w.r.t.\ $Y$ if and only if $A=(p,r)$ belongs to $C^+_{q, A'}\cap C_{q,A'}^-$, that is:
\[
C^+_{q, A'}\cap C_{q,A'}^-= 
\{A=(p,r)\in D_q\, \colon\, X \cseq X' \quad \text{w.r.t.} \quad Y\}.
\]
Therefore, characterizing equality in the Csiszár order amounts to characterizing the intersection $C^+_{q,A'}\cap C^-_{q,A'}$. 
The following proposition gives an explicit description of this intersection depending on the position of $A'$.
Recall that $X$ or $Y$ can be transformed by a bijection without changing the Csisz\'ar index or the Csisz\'ar order (see Proposition~\ref{prop:invariance_csiszar_index} or Axiom~\ref{axiom:transformation_invariance}). Thus, the set $C^+_{q, A'}\cap C_{q,A'}^-$ contains at least the point $A$ associated to $X$ such that $X = g(X')$ for some bijection $g$. Since $g = 1 - \id$ is a bijection, $h(A')$ belongs to $C^+_{q, A'}\cap C_{q,A'}^-$. 
The proof of the following proposition is given in Section~\ref{sec:C+-cap-C-}. 

\begin{prop}[On $C_{q,A'}^+\cap C_{q,A'}^-$]
  \label{prop:C+-cap-C-}
  Let $q\in [0, 1]$ and $A'=(p',r')\in D_q$.
  We have the following result:  
  \begin{enumerate}[(i)]
    \item \label{it:C+capC-=D} If $A'\in \Delta_q$, then $C^+_{q, A'}\cap C_{q,A'}^-=\Delta_q$.
    \item \label{it:C+capC-=2} If $A'\not \in \Delta_q$ and either $q\neq 1/2$ or $p'=q=1/2$, then $C^+_{q, A'}\cap C_{q,A'}^-$ consists in the 2 distinct points $\{A', h(A')\}$.
    \item \label{it:pneqq=1/2} If  $A'\not \in \Delta_q$ and $p'\neq q= 1/2$, then $C^+_{q, A'}\cap C_{q,A'}^-$ consists in the 4 distinct points $\{A', A'', h(A'), h(A'')\}$, where $A'' $ is defined in Proposition~\ref{prop:bernoulli_cp}.
  \end{enumerate}
\end{prop}

The first case in the previous proposition corresponds to $X'$ independent of $Y$. Thus, all the $X$ independent of $Y$ belong to $C^+_{q, A'}\cap C_{q,A'}^-$. 

The second case corresponds to the situation where $X'$ and $Y$ are not independent and either $q \neq 1/2$ or $p'=q=1/2$. In this case $C^+_{q, A'}\cap C_{q,A'}^-$ corresponds to two points. These two points in $D_q$ correspond to the random variables  $(X',Y)$, $(1-X',Y)$, which have the same Csiszár indices.

The third case corresponds to the situation where $X'$ and $Y$ are not independent and $p'\neq q =1/2$. In this case $C^+_{q, A'}\cap C_{q,A'}^-$ corresponds to four points. These four points in $D_{1/2}$ correspond to the random variables  $(X',Y)$, $(X', 1-Y)$, $(1-X', Y)$ and $(1-X', 1-Y)$ which have the same Csiszár indices. Since $q = 1/2$, we can transform $Y$ into $1-Y$ without changing the parameter $q$ since $1 - q = q$. That is why we get two new points $A''$ and $h(A'')$ compared to the second case.

\medskip

Next, to determine when the Csisz\'ar order is total, we study the set $C_{q,A'}^+ \cup C_{q,A'}^-$.
Indeed, the Csisz\'ar order is total if and only if $C_{q,A'}^+ \cup C_{q,A'}^- = D_q$, since this means every $A \in D_q$ makes $X$ either greater or smaller than $X'$ w.r.t.\ $Y$, \emph{i.e.\ }a total order.

Proposition~\ref{prop:C+-cup-C-} in Section~\ref{sec:C+-cup-C-} gives an explicit description of this union depending on the position of $A'$. In particular, we have:
\[
C_{q,A'}^+ \cup C_{q,A'}^- = D_q \iff A' \in \Delta_q \cup E_q.
\]
It gives that if $Y$ and $X'$ are Bernoulli random variables, then  there exists (on a possible enlarged probability space) a Bernoulli random variable $X$ such that $X$ and $X'$ are not comparable w.r.t.\ $Y$ if and only if  we have  neither $Y$ and $X'$ independent nor $X'$ a function of $Y$ (that is, not $Y=X'$ a.s.\ or $Y=1-X'$ a.s.).

Thus we deduce that the Csiszár order is total (that is, any  Bernoulli random variables $X$ and $X'$ are comparable w.r.t.\ $Y$)  if and only if $D_q=\Delta_q$, that is, if and only if $Y$ is constant ($q\in \{0, 1\}$). 
Hence, except in these degenerate cases, the Csiszár order is not total; furthermore $\Delta_q$ and $E_q$ are respectively the sets of minima and the set of maxima for the  Csiszár order. It is easy to check that, for $q\in (0, 1)$,   the set  $\{ P_{X|Y}\, \colon\, P_X(\{0, 1\})=1\}$ of (conditional) Bernoulli distributions with the  Csiszár order w.r.t.\ $Y$ is a lattice if and only if $q=1/2$.

\subsection{Comparison with the squared correlation order} 
For $q\in [0, 1]$ and  $A=(p,r) \in D_q$, we have for $p,q\not\in \{0, 1\}$:
\[
\corr^2(X,Y)=\frac{(r -pq)^2}{p(1-p)q(1-q)},
\]
and by convention we set $\corr(X,Y)=0$ if $p$ or $q$ belongs to $\{0, 1\}$.
One easily gets that:
\begin{equation}
    \label{eq:corr=01}
    \corr(X,Y)=0 \iff A\in \Delta_q;
    \quad\text{and}\quad
    |\corr(X,Y)|=1\iff A\in E_q \quad \text{and} \quad q \notin \{0,1\}.
\end{equation}

We compare the Csisz\'ar order between the random variables $X$ and $X'$ w.r.t.\ $Y$ and the squared correlation order (which is total), meaning that either $\corr^2(X,Y)\geq \corr^2(X',Y)$ or $\corr^2(X,Y)\leq \corr^2(X',Y)$.
The main result of this subsection is the following corollary, which follows from Proposition~\ref{prop:Cp_subset_Fp} (a technical result given in Section~\ref{sec:Cp_subset_Fp}). %It states that the Csisz\'ar order implies the same order on squared correlations.

\begin{cor}[Csiszár order implies the squared correlation order]
\label{cor:Csiszar_implies_corr}
Let $X,X',Y$ be three Bernoulli random variables. Then, we have:
\begin{equation}
  \label{eq:Csiszar_implies_corr}
  X \csleq X' \quad \text{w.r.t.} \quad Y \implies \corr^2(X,Y) \leq \corr^2(X',Y).
\end{equation}
\end{cor}

However, the converse implication is false (unlike the case where $(X,X',Y)$ is a $3$-dimensional Gaussian vector). Indeed, it is enough to consider two random variables $X$ and $X'$ that are not comparable w.r.t.\ the Csiszár order. In this case, the corresponding squared correlations with $Y$ are necessarily ordered, while neither $X \csgeq X'$ nor $X\csleq X'$ w.r.t.\ $Y$ holds.
We give more details on the comparison between the Csiszár order and the squared correlation order in Section~\ref{sec:Cp_subset_Fp}, see the technical Proposition~\ref{prop:Cp_subset_Fp}. In particular, the sets $\{A \in D_q \, \colon \, \corr(X,Y)^2 \geq \corr(X',Y)^2\}$ and $\{A \in D_q \, \colon \, \corr(X,Y)^2 \leq \corr(X',Y)^2\}$ coincide with $C_{q,A'}^+$ and $C_{q,A'}^-$, respectively, if and only if $\vert \corr(X',Y) \vert \in \{0,1\}$, that is, when $(X',Y)$ is degenerate.

\begin{figure}[htbp]
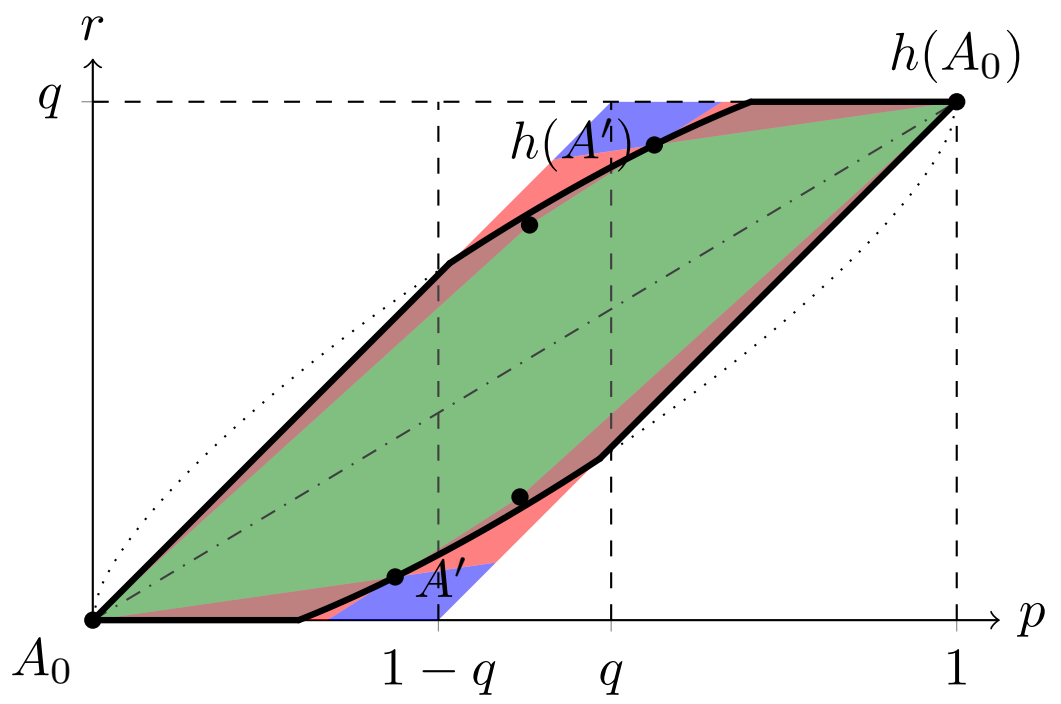
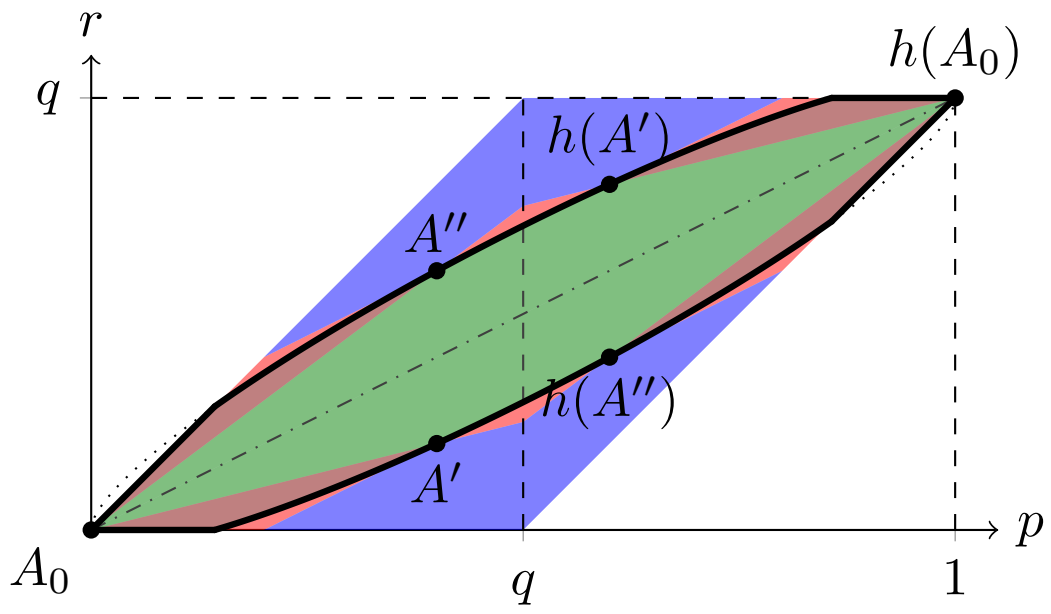

  \begin{subfigure}[c]{0.49\textwidth}
  \centering
  \begin{RegionPlotBernoulli}{0.35}{0.05}{0.6}
  {
        % \addplot[pattern={Lines[angle=45, distance=5pt, line width=0.1pt], draw=none, area legend}] 
    \path[fill=gray, opacity=0.5]
    (0,0) 
    -- plot[domain=0.412903225806452:0.761860465116279, smooth, samples = \nsamples](\x,{0.335451147751014*(\x^0.5*(1-\x)^0.5) + 0.6*\x})
    -- (1,0.6) 
    -- plot[domain=0.587096774193548:0.238139534883721, smooth, samples = \nsamples](\x,{-0.335451147751014*(\x^0.5*(1-\x)^0.5) + 0.6*\x})
    -- cycle;
    
    \addplot[black, solid, very thick, line cap=round] coordinates {(0.587096774193548, 0.18709677419354798) (1,0.6)};
    \addplot[black, solid, very thick, line cap=round] coordinates {(0, 0) (0.238139534883721,0)};
    \addplot[black, very thick, line cap=round, domain=0.238139534883721:0.587096774193548, smooth, samples = \nsamples](\x,{-0.335451147751014*(\x^0.5*(1-\x)^0.5) + 0.6*\x});
    \addplot[domain=0.587096774193548:1, smooth, samples = \nsamples, dotted](\x,{-0.335451147751014*(\x^0.5*(1-\x)^0.5) + 0.6*\x});

    \addplot[black, solid, very thick, line cap=round] coordinates {(0,0) (0.412903225806452,0.412903225806452)};
    \addplot[black, solid, very thick, line cap=round] coordinates {(1, 0.6) (0.761860465116279,0.6)};
    \addplot[black, very thick, line cap=round, domain=0.412903225806452:0.761860465116279, smooth, samples = \nsamples](\x,{0.335451147751014*(\x^0.5*(1-\x)^0.5) + 0.6*\x});
    \addplot[domain=0:0.412903225806452, smooth, samples = \nsamples, dotted](\x,{0.335451147751014*(\x^0.5*(1-\x)^0.5) + 0.6*\x});

    \fill (0.35,0.05) circle (0.01) node[right] {$A'$};
    \fill (0.4944444444444446,0.14259259259259266) circle (0.01) node[above] {};
    %\fill (0.4944444444444446,0.14259259259259266) circle (0.01) node[above] {$h(A'')$};

    \fill (0.65,0.55) circle (0.01) node[left] {$h(A')$};
    \fill (0.5055555555555554,0.4574074074074073) circle (0.01) node[below] {};
    %\fill (0.5055555555555554,0.4574074074074073) circle (0.01) node[below] {$A''$};

  }
  {
    \path[fill=\ColorCm, opacity=\OpacityPlot]
    (0,0) -- (0.35,0.05) -- (0.4944444444444446,0.14259259259259266) -- (1,0.6) -- cycle;
    
    \path[fill=\ColorOut, opacity=\OpacityPlot]
    (0,0) -- plot[domain=0.271541071319574:0.35, smooth, samples = \nsamples](\x,{0.6*\x - 0.192*(1 - \x)/(1.13 - \x)}) -- cycle;
    
    \path[fill=\ColorCp, opacity=\OpacityPlot]
    plot[domain=0.271541071319574:0.35, smooth, samples = \nsamples](\x,{0.6*\x - 0.192*(1 - \x)/(1.13 - \x)})
    -- (0.466666666666667,0.06666666666666698)
    -- (0.4,0)
    -- cycle;
    
    \path[fill=\ColorOut, opacity=\OpacityPlot]
    (0.35,0.05) -- (0.466666666666667,0.06666666666666698) -- (1,0.6) -- (0.4944444444444446, 0.14259259259259266 ) -- cycle;

  }
\end{RegionPlotBernoulli}
\caption{An instance where $C^+_{q,A'} \cap C^-_{q,A'}$ consists of two points (with $p'=0.35$, $q = 0.6$ and $r'=0.05$).}
\label{fig:boundary_F+_2pts}
\end{subfigure}
\hfill
\begin{subfigure}[c]{0.49\textwidth}
  \centering
\begin{RegionPlotBernoulli}{0.4}{0.1}{0.5}
  {

    %\addplot[pattern={Lines[angle=45, distance=5pt, line width=0.1pt]}] 
    \path[fill=gray, opacity=0.5]
    (0,0) 
    -- plot[domain=0.14285714285714302:0.857142857142857, smooth, samples = \nsamples](\x,{+0.204124145231932*(\x^0.5*(1-\x)^0.5) + 0.5*\x})
    -- (1,0.5) 
    -- plot[domain=0.857142857142857:0.142857142857143, smooth, samples = \nsamples](\x,{-0.204124145231932*(\x^0.5*(1-\x)^0.5) + 0.5*\x})
    -- cycle;

    \addplot[black, solid, very thick, line cap=round] coordinates {(1, 0.5) (0.857142857142857,0.357142857142857)};
    
    \addplot[black, solid, very thick, line cap=round, domain=0.142857142857143:0.857142857142857, smooth, samples = \nsamples](\x,{-0.204124145231932*(\x^0.5*(1-\x)^0.5) + 0.5*\x});

    \addplot[black, solid, very thick, line cap=round] coordinates {(0,0) (0.142857142857143,0)};

    \addplot[black, solid, very thick, line cap=round] coordinates {(0, 0) (0.14285714285714302,0.14285714285714302)};

    \addplot[black, solid, very thick, line cap=round, domain=0.14285714285714302:0.857142857142857, smooth, samples = \nsamples](\x,{+0.204124145231932*(\x^0.5*(1-\x)^0.5) + 0.5*\x});

    \addplot[black, solid, very thick, line cap=round] coordinates {(1, 0.5) (0.857142857142857,0.5)};

    \addplot[domain=0.857142857142857:1, smooth, samples = \nsamples, dotted](\x,{-0.204124145231932*(\x^0.5*(1-\x)^0.5) + 0.5*\x});
    \addplot[domain=0:0.14285714285714302, smooth, samples = \nsamples, dotted](\x,{+0.204124145231932*(\x^0.5*(1-\x)^0.5) + 0.5*\x});

    \fill (0.4,0.1) circle (0.01) node[below] {$A'$};
    \fill (0.6,0.4) circle (0.01) node[above] {$h(A')$};

    \fill (0.4,0.3) circle (0.01) node[above] {$A''$};
    \fill (0.6,0.2) circle (0.01) node[below] {$h(A'')$};
  }
  {
    \path[fill=\ColorCm, opacity=\OpacityPlot]
    (0,0) -- (0.4,0.1) -- (0.6,0.2) -- (1,0.5) -- cycle;
    
    \path[fill=\ColorOut, opacity=\OpacityPlot]
    (0,0) -- (0.2,0) -- (0.4,0.1) -- cycle;

    \path[fill=\ColorOut, opacity=\OpacityPlot]
    (0.4,0.1) -- (1/2,1/8) -- (0.6,0.2) -- cycle;

    \path[fill=\ColorOut, opacity=\OpacityPlot]
    (0.6,0.2) -- (0.8,0.3) -- (1,1/2) -- cycle;
    
    \path[fill=\ColorCp, opacity=\OpacityPlot]
    (0.2,0) -- (1/2,0) -- (0.8,0.3) -- (0.6,0.2) -- (1/2,1/8) -- (0.4,0.1) -- cycle; 
  }
\end{RegionPlotBernoulli}
\caption{An instance where $C^+_{q,A'} \cap C^-_{q,A'}$ consists of four points. (with $p'=0.4$, $q = 0.5$ and $r'=0.1$)}
\label{fig:boundary_F+_4pts}
\end{subfigure}
\caption{Illustrations of the fact that the correlation order is weaker than the Csiszár order. The shaded region corresponds to $\{A \in D_q \, \colon \, \corr^2(X,Y) \leq \corr(X',Y)\}$.}
\label{fig:boundary_F+}
\end{figure}

In Fig.~\ref{fig:boundary_F+}, we illustrate that the correlation order is weaker than the Csiszár order, as given by properties~\eqref{eq:inclusion_C+_F+} and~\eqref{eq:intersection_C+_F-} of Proposition~\ref{prop:Cp_subset_Fp}. The shaded region corresponds to the set $\{A \in D_q \, \colon \, \corr^2(X,Y) \leq \corr(X',Y)^2\}$. As expected, the green region, corresponding to $C_{q,A'}^-$, is contained in the shaded region. Similarly, the blue region, corresponding to $C_{q,A'}^+$, is contained in the unshaded region, which corresponds to $\{A \in D_q \, \colon \, \corr^2(X,Y) \geq \corr(X',Y)^2\}$.

\section{Csisz\'ar index for thresholded Gaussian random variables}
\label{sec:ThresholdedGaussian}

In many statistical applications, one is not interested in a variable itself but in whether it exceeds a given threshold. This is typically the case in reliability and safety analysis, where the threshold represents a physical quantity that should not be exceeded—for instance a maximal load, temperature, or stress beyond which a structure fails. Similar situations arise in environmental monitoring, hydrology, and finance. In all these settings, the relevant object is the global dependence between the exceedance indicators rather than the underlying variables; this  motivates studying the Csisz\'ar order between thresholded versions of the Gaussian vectors considered so far.

Let $(X, X',Y)\in \R^3$ be a Gaussian random vector. We recall that $X$ and $X'$ are always comparable w.r.t.\ $Y$ for the Csiszár order. 
The main result of this section states that: the thresholded random variables $\ind_{\{X \leq \alpha\}}$ and $\ind_{\{X' \leq \alpha'\}}$ are comparable w.r.t.\ $\ind_{\{Y \leq \beta\}}$ for particular values of the thresholds; and  when they are comparable then they are in the same order as $X$ and $X'$ w.r.t.\ $Y$.

\begin{theo}[Thresholded 1-dimensional Gaussian random variables] 
\label{thm:thresholded}
  For $(X,X',Y)\in \R^3$  a Gaussian random vector, we have:
  \begin{equation}
      \label{eq:thm_thresholded_eq}
      X \cseq X' \quad \text{w.r.t.}\quad Y \iff \exists \, \alpha, \alpha', \beta \in \R \quad\text{such that}\quad \ind_{\{X \leq \alpha\}} \cseq \ind_{\{X' \leq \alpha'\}} \quad \text{w.r.t.}\quad  \ind_{\{Y \leq \beta\}},
  \end{equation}
  \begin{equation}
      \label{eq:thm_thresholded_ineg}
      X \csls X' \quad \text{w.r.t.}\quad Y \iff \exists \, \alpha, \alpha', \beta \in \R \quad\text{such that}\quad \ind_{\{X \leq \alpha\}} \csls \ind_{\{X' \leq \alpha'\}} \quad \text{w.r.t.}\quad  \ind_{\{Y \leq \beta\}}.
  \end{equation}
\end{theo}

The proof of this theorem is given in Section~\ref{sec:proof_thresholded}.

\begin{rem}[The assumption that $(X,X',Y)$ is a Gaussian vector is slightly too strong]
\label{rem:thresholded_XX'Y_GV_toostrong}
% As in Section~\ref{sec:Gaussian}, we consider a Gaussian random vector $(X,X',Y)$. 
Since we are only interested in the Csisz\'ar indices between $X$ and $Y$, between $X'$ and $Y$, or their thresholded versions, the dependence between $X$ and $X'$ can be arbitrary, as explained in Remark~\ref{rem:XX'Y_GV_toostrong}. Hence the results of Theorem~\ref{thm:thresholded} hold under the slightly weaker assumption that $(X,Y) \in \R^2$ and $(X',Y) \in \R^2$ are two Gaussian random vectors.
\end{rem}

For a two dimensional Gaussian vector $(X,Y)$, we set for $\alpha, \beta \in [-\infty , \infty]$:
\[
 q = \P(Y \leq \beta), \quad p=\P(X \leq \alpha)
 \quad\text{and} \quad
 r = \P(X \leq \alpha, Y \leq \beta).
\]
These three parameters correspond to the Bernoulli pair $(\ind_{\{X \leq \alpha\}}, \ind_{\{Y \leq \beta\}})$. We define $p'$ and $r'$ in the same way.

\begin{rem}[If $(X,Y)$  or $(X',Y)$ is not a Gaussian vector, Theorem~\ref{thm:thresholded} may not hold] Suppose that $X$ and $Y$ are two standard  Gaussian random variables,
%(that is, with cdf $\Phi$)
and assume that the copula of $(X,Y)$ is given by the Clayton copula with dependence parameter equal to $1$, that is for $u,v \in [0,1]$:
\[
C(u,v) = \frac{uv}{u+v-uv}\cdot
\]
We set $X'=X$, so that obviously $X' \cseq X$ w.r.t.\ $Y$. Let $\beta= 0$ (that is, $q=1/2$) and $\alpha'$ such that $\P(X' \leq \alpha') = 3/5$ (that is, $p' = 3/5$).
Clearly we have $\ind_{\{X '\leq \alpha'\}} \cseq \ind_{\{X \leq \alpha'\}}$ w.r.t.\ $\ind_{\{Y \leq \beta\}}$. 

\medskip 

However, the ranking is not preserved for all thresholds; in particular, we claim
that $\ind_{\{X '\leq \alpha'\}} \csls \ind_{\{X \leq -\alpha'\}}$ w.r.t.\ $\ind_{\{Y \leq \beta\}}$. 
To do so, we prove that $A=(p,r)$ belongs to the interior of $C^+_{q,A'}$, where $A' = (p',r')$. By hypothesis, we have  $p = 1-p'= 2/5$ and:
\[
r = C(1-p',q) = \frac{2}{7}
\quad \text{and} \quad 
r' = C(p',q) = \frac{3}{8} \cdot
\]
Using that $q = 1/2$ and Proposition~\ref{prop:C+-cap-C-}~\ref{it:pneqq=1/2}, we have:
\[
C_{q,A'}^+ \cap C_{q,A'}^- = \{ A', A'', h(A'),  h(A'')\},
\]
with $A'' = (p'', r'')$ defined in~\eqref{eq:bernoulli_cp''} which is equal to $(1-p', r' + q(1-2p')) = (2/5, 11/40)$. Thus $A$ doesn't belong to $C_{q,A'}^+ \cap C_{q,A'}^-$, that is $\ind_{\{X' \leq \alpha'\}}$ and $\ind_{\{X \leq -\alpha'\}}$ are not equal in the Csiszár order w.r.t.\ $\ind_{\{Y \leq \beta\}}$. Moreover, since $p'' = 1-p' = p$ and $r > r'' > p'' q$ we have $A \in C_{q,A'}^+$ and 
thus $\ind_{\{X '\leq \alpha'\}} \csls \ind_{\{X \leq -\alpha'\}}$ w.r.t.\ $\ind_{\{Y \leq \beta\}}$. This proves the claim.

\medskip
Thus, assuming that the random vector $(X,Y)$ is Gaussian in Theorem~\ref{thm:thresholded} (see also Remark~\ref{rem:thresholded_XX'Y_GV_toostrong}), can not be removed even though it could certainly be weakened.
\end{rem}

\begin{rem}[Expression of $\alpha$, $\alpha'$, and $\beta$ in Theorem~\ref{thm:thresholded}]
\label{rem:expressions_aab_thresholded}
Assume that $(X,X',Y) \in \R^3$ is a  Gaussian random vector, with $X$, $X'$ and $Y$ having positive variances.
Since Csiszár indices are invariant under bijective transformations (see Proposition~\ref{prop:invariance_csiszar_index}), we may assume without loss of generality that the correlation between $X$ and $Y$ and between $X'$ and $Y$ are nonnegative by replacing $X$ and/or $X'$ with their opposites.
Recall, that the thresholded variables are in the same order as the initial variables $X$ and $X'$ or they  are not comparable.
We first give an easy to remember (partial) result on the comparison for the thresholded random variables and then a more detailed comment. 

\begin{itemize}
    \item \textbf{Simple informative choice of the thresholds.} Choosing $\alpha = \E[X]$, $\alpha' = \E[X']$ and $\beta = \E[Y]$ guarantees that the random variables $\ind_{\{X \leq \alpha\}}$ and $\ind_{\{X' \leq \alpha'\}}$ can be compared w.r.t.\ $\ind_{\{Y \leq \beta\}}$. (This corresponds to  $q=p'=p=1/2$,  and according to Proposition~\ref{prop:bernoulli_cp} and Proposition~\ref{prop:bernoulli_cm} in Section~\ref{sec:description_Cm} the vertical segment $\{(p,r) \, \colon \, r \geq 0 \} \cap D_q$  is a subset of $C^+_{q, A'} \cup C^-_{q, A'}$.)

    \item \textbf{Complement on the Equivalences~\eqref{eq:thm_thresholded_eq} and~\eqref{eq:thm_thresholded_ineg}.} We now examine (some of) the values of the parameters $\alpha$, $\alpha'$ and $\beta$ for which the thresholded variables are comparable. Let us fix $\alpha', \beta \in \R$ and define:
    \[
    \alpha_\pm^* = \pm  \,\sqrt{\frac{\Var(X)}{\Var(X')}} \left(\alpha' - \E[X']\right) + \E[X] .%\quad \text{for $\varepsilon \in \{-,+\}$}.
    \]
   In particular, $\alpha_+^* = \alpha_-^* = \E[X]$ when $\alpha'=\E[X']$. We have an explicit result for the equality case from Equivalence~\eqref{eq:thm_thresholded_eq}. Using Proposition~\ref{prop:copule_thresholded}, we get that for any $\beta \in \R$ with $\beta \neq \E[Y]$:
    \[
    X \cseq X' \quad \text{w.r.t.}\quad Y
    \iff 
    \begin{cases}
        \ind_{\{X \leq \alpha\}} \cseq \ind_{\{X' \leq \alpha'\}} \quad \text{w.r.t.}\quad  \ind_{\{Y \leq \beta\}} & \text{if $\alpha = \alpha_+^*$,} \\
        \text{$\ind_{\{X \leq \alpha\}}$ and $\ind_{\{X' \leq \alpha'\}}$ are not comparable w.r.t.\ $\ind_{\{Y \leq \beta\}}$} & \text{otherwise.}
    \end{cases}
    \]
    If  $\beta = \E[Y]$, the above equivalence still holds provided that $\alpha = \alpha_+^*$ is replaced by $\alpha \in \{\alpha_+^*, \alpha_-^*\}$.

    We now turn to Equivalence~\eqref{eq:thm_thresholded_ineg} and derive a similar result, although in a less explicit form. 
    Let $B = B(\alpha',\beta)$ be the set of $\alpha\in \R$ such that $A=(\P(X \leq \alpha), \P(X \leq \alpha, Y \leq \beta))$ belongs to $C^+_{q,A'}$ where we recall that $q = \P(Y \leq \beta)$. 
    Then Equivalence~\eqref{eq:thm_thresholded_ineg} holds if and only if $\alpha$ belongs to the set $B$. (Notice $B$ is either empty or contains $\alpha_+^*$ and is then either an interval or the union of two non empty disjoint intervals, see Fig.~\ref{fig:upsilon_rho0.7} for various instances.) 
    % \AS{Cet ensemble $B$ correspond aux courbes de la figure \ref{fig:upsilon_rho0.7} intersection la zone verte. Donc soit l'intersection est vide, soit il y a au moins un point. S'il y a un point ce sera le point de coordonnées $(\P(X \leq \alpha_+^*), \P(X \leq \alpha_+^*, Y \leq \beta))$}
    % \AS{Sortir le plot de l'annexe et détailler?}
    We also have: 
    \[
    X \csls X' \quad \text{w.r.t.}\quad Y
    \iff
    \ind_{\{X \leq \alpha\}} \csls \ind_{\{X' \leq \alpha'\}} \quad \text{w.r.t.}\quad  \ind_{\{Y \leq \beta\}}
    \quad \text{with $\alpha = \alpha_+^*$}.
    \]
    If $\beta = \E[Y]$, the above equivalence also  holds with  $\alpha_+^*$  replaced by $\alpha_-^*$.

\end{itemize}
\end{rem}

\bibliographystyle{plain}
\bibliography{references}

\newpage

\appendix

\section{A few remarkable properties on $f$-divergences and the Csiszár index}

In this section, we recall some results on $f$-divergences and the Csiszár index. 
We first extend the notion of Radon-Nikodym derivative to measures that need not be mutually absolutely continuous.
Let $\mu,\nu$ be $\sigma$-finite measures, and set $\lambda = \mu+\nu$. By the Radon-Nikodym theorem, there exist nonnegative measurable functions $p,q$ such that $\nu = p\lambda$, $\mu = q\lambda$, and $p+q=1$ $\lambda$-a.e.. Setting $h = p/q$ (with convention $0/0=\infty/\infty=0$), \textbf{even when $\mu$ and $\nu$ are not mutually absolutely continuous}, we write:
\begin{equation}
  \label{eq:Leb-dd2}
  \frac{\rd \nu}{\rd \mu} = h
  \quad\text{and}\quad
  \frac{\rd \mu}{\rd \nu} = \inv{h},
  \quad \text{with $h$ defined $(\mu+\nu)$-a.e., valued in $[0,+\infty]$.}
\end{equation}

The following corollary gives an invariance property of $f$-divergences, see~\cite[Corollary~3]{butucea2026csisz} or~\cite{liese2006divergences}.

\begin{cor}[Invariance of the $f$-divergence for random variables]
\label{cor:invariance_Df}
Let $X,Y$ be two random variables taking values in $(\Omega,\cf)$ with respective distribution $P_X$ and $P_Y$. Let $(E,\ce)$ be a measurable space and $\varphi$ be a $E$-valued measurable function defined on $\Omega$. Let $f\in \F$. Then, we have:
\[
D_f(P_{\varphi(X)} \Vert P_{\varphi(Y)}) \le D_f(P_X \Vert P_Y).
\]
If $\rd P_X/ \rd P_Y$ (in the sense of \eqref{eq:Leb-dd2}) can be written as a measurable function of $\varphi$, then we have:
\begin{equation}
 \label{eq:egalite-Df}
 D_f(P_{\varphi(X)} \Vert P_{\varphi(Y)}) = D_f(P_X \Vert P_Y).  
\end{equation}
Furthermore if the Equality~\eqref{eq:egalite-Df} holds for some $f\in \F$ strictly convex on $\R_+^*$ with $D_f(P_X \Vert P_Y)$ finite, then $\rd P_X/ \rd P_Y$ (in the sense of \eqref{eq:Leb-dd2}) is a measurable function of $\varphi$ and~\eqref{eq:egalite-Df} holds for all $f\in \F$. 
\end{cor}

By applying this corollary to Csisz\'ar indices, we recover Proposition~9 of \cite{butucea2026csisz}. 

\begin{prop}[Invariance of the Csisz\'ar index]
\label{prop:invariance_csiszar_index}
Let $f \in \F$ (and thus $f(1)=0$). For $i=1,2$, let $X_i$ be a random variable taking values in  a measurable space $(\Omega_i, \cf_i)$ and $\varphi_i$ a measurable function defined on $(\Omega_i, \cf_i)$ and taking values in a  measurable space $(E_i, \ce_i)$. Then, we have:
\[
S_f(\varphi_1(X_1), \varphi_2 (X_2)) \leq  S_f(X_1, X_2). 
\]
\begin{enumerate}[(i)]
  \item \textbf{Marginal invariance.} \label{enum:marginal_invariance}
  If $\rd P_{X_1} \otimes \, \rd P_{X_2} / \rd P _{(X_1,X_2)}$ (in the sense of~\eqref{eq:Leb-dd2}) can be written as a measurable function of $(\varphi_1, \varphi_2)$, then we have:
  \[
  S_f(\varphi_1(X_1), \varphi_2 (X_2)) = S_f(X_1, X_2). 
  \]
  \item\textbf{Invariance by injection.}
  \label{enum:bijection_invariance}
  If, for $i=1, 2$ the function  $\varphi_i$ is injective and bi-measurable, then we have:
  \[
  S_f(\varphi_1(X_1), \varphi_2 (X_2)) = S_f(X_1, X_2). 
  \]
  \item\textbf{Invariance for independent Markovian kernels} 
  \label{it:Sf-Markov} 
  Assume, that  for $i=1,2$, the random  variable $X_i$ can be  written as $X_i=(Y_i, W_i)$ and  takes values in a Polish space  (thus  the following probability kernels exist). If:
  \[
  P_{(W_1, W_2)|(Y_1,  Y_2)=(y_1, y_2)}(\rd w_1, \rd  w_2)= \prod_{i=1}^2 P_{W_i|Y_i=y_i} (\rd w_i)
  \]
  holds $P_{(Y_1, Y_2)}(\rd y_1, \rd y_2)$-a.s., then, we have:
  \[
  S_f(Y_1, Y_2) = S_f\left((Y_1, W_1),\,  (Y_2, W_2)\right). 
  \] 
\end{enumerate}
\end{prop}

The next remark corresponds to the ``forgetting'' property of independent random variables, see \cite[Remark~13~(iii)]{butucea2026csisz}.
\begin{rem}[``Forgetting'' independent random variables]
  \label{rem:forgetting_independent_variables}
  Let $X$, $Y$ and $U$ be random variables such that $U$ and $Y$ have a conditional distribution given $X$ (which is satisfied if $U$ and $Y$ take values in a Polish space). Then, for any function $f \in \F$, we have:
  \[
  S_f((X,U),Y) = S_f(X,Y)
  \quad  \text{for all} \quad  f \in \F \quad\iff\quad  \text{$U$ and $Y$ are conditionally independent given $X$}.
  \]
 \end{rem}

\section{Proofs of the results in Section \ref{sec:ConnectionConvexOrder}}
\label{sec:ProofsConnectionConvexOrder}

We first recall a technical result on the existence of the density ratio, see~\cite{butucea2026csisz}.
When there is no ambiguity on the reference $\sigma$-finite measure $\lambda$, when $P\ll \lambda$, we denote by $p$ the corresponding density. 

\begin{lem}[Existence of the density ratio]
    \label{lem:abs-cont-XY}
Let $(X,Y)$ be a random variable taking values in a product of measurable spaces.   Assume that  $P_{(X,Y)}$ is dominated by a  product of $\sigma$-finite measures   $\lambda=\lambda_X\otimes \lambda_Y$, then $P_X$ and $P_Y$ are dominated by $\lambda_X$ and $\lambda_Y$ respectively, and furthermore
$P_{(X,Y)} $ is dominated by $P_X\otimes P_Y$ with density:
\[
R_{(X,Y)}(x,y)=\frac{p_{(X,Y)}(x,y)}{p_X(x)\, p_Y(y)}\cdot
\] 
\end{lem}
Notice the ratio $R_{(X,Y)}$ is well defined on the set $A_X\times A_Y$ and can be set equal to 0 on $(A_X\times A_Y)^c$, where $A_X=\{p_X>0\}$ and $A_Y=\{p_Y>0\}$ are respectively the supports of $p_X$ and $p_Y$.

Our next result extends \cite[Theorem~3.A.1.]{shaked2007stochastic}, where a necessary and sufficient condition is established for the convex order restricted to convex functions which are finite on $\R$. Our extension to proper convex functions  is required for the study and comparison of Csisz\'ar indices. Recall that the notion of quasi-integrable random variable and convex order are defined in Section~\ref{sec:ConnectionConvexOrder}.

\begin{lem}[Characterization of the convex order]
  \label{lem:convex_order_characterization}
  Let $Z$ and $Z'$ be two real integrable random variables. Then, we have:
  \begin{equation}
    \label{eq:convex_order_characterization}
    Z \leqcv Z' 
    \quad \iff \quad
    \E[Z] = \E[Z'] \quad \text{and} \quad
    \E[(Z-t)_+] \leq \E[(Z'-t)_+] \quad \text{for all $t\in \R$}.
  \end{equation}
\end{lem}

\begin{proof}[Proof of Lemma~\ref{lem:convex_order_characterization}]
The direct implication is obvious as the functions identity and $\varphi(x) = (x-t)_+$ are a convex function for any $t\in \R$.
Assume that $Z$ and $Z'$ are integrable, with $\E[Z]=\E[Z']$, and such that $\E[(Z-t)_+] \leq \E[(Z'-t)_+] $  for all $t\in \R$. 
Theorem 3.A.1.\ in~\cite{shaked2007stochastic} states that 
$\E\left[ \varphi(Z)\right] \leq \E\left[ \varphi(Z')\right]
$ for any $\R$-valued convex function $\varphi$ (notice that $\varphi(Z)$ and $\varphi(Z')$ are quasi-integrable as $Z$ and $Z'$ are integrable). We shall prove this inequality holds also  for all proper convex function $\varphi$ such that $\varphi(Z)$ and $\varphi(Z')$ are quasi-integrable. 
To get this extension, it is enough to prove that if $\varphi$ is a proper convex function such that  $\varphi(Z)$ is quasi-integrable, then there exists a sequence $(\varphi_n)_{n\in \N}$ of real-valued convex function such that $\lim_{n\rightarrow \infty} \E[\varphi_n(Z)]=\E[\varphi(Z)]$.

Notice that $\E[(Z-t)_-] \leq \E[(Z'-t)_-] $ for all $t\in \R$. 
Let $\varphi$ be a given proper convex function. 
Set $a_0$ in the interior of $\dom(\varphi)$, $a_-=\inf \dom(\varphi)$ and $a_+=\sup \dom(\varphi)$.
Then, there exists a sequence $(\varphi_n)_{n\in \N}$ of finite positive linear combinations of maps 
$ x\mapsto c_0 x$ for some $c_0\in \R$, $x\mapsto (x-t)_+$ and $x\mapsto (x-t)_-$ for $t\in \R$ such that
$\varphi_n \leq \varphi_{n+1} \leq \varphi$ on $D_n$, with: $D_n \cap (-\infty, a_0]$ equal to $(-\infty, a_0]$ unless 
$a_- > -\infty$ and $+ \infty>\varphi(a_-)> \lim_{x \downarrow a_-} \varphi(x)$ in which case we shall assume that $\varphi_n(a_-)=\varphi(a_-)$ and $D_n
\cap (-\infty, a_0]=(-\infty, a_-]\cup[a_- + 1/n, a_0]$ (provided $n$ is small enough so that $a_-+1/n<a_0$); and similarly with 
$D_n \cap [a_0, +\infty)$ equal to $[a_0, +\infty)$ unless 
$a_+< +\infty$ and $+ \infty>\varphi(a_+)> \lim_{x \uparrow a_+} \varphi(x)$ in which case we shall assume that $\varphi_n(a_+)=\varphi(a_+)$ and $D_n
\cap [a_0, +, \infty)=[a_0, a_+ - 1/n] \cup [a_+, + \infty)$ (provided $n$ is small enough so that $a_+ - 1/n<a_0$). 
By construction, we have $\varphi_n(Z)$ is integrable and $\lim_{n\rightarrow \infty} \varphi_n(x)=\varphi(x)$ for all $x\in \R$. 
Notice that as $Z$ is integrable, then $\varphi(Z)_-$ is integrable and thus $\varphi(Z)$ is quasi-integrable.
It is then elementary to check that  $\lim_{n\rightarrow \infty} \E[\varphi_n(Z)]=\E[\varphi(Z)]$.
\end{proof}

Using this lemma, we can now prove Proposition~\ref{prop:equivalence_ConvexOrder_CsiszarOrder}.

\begin{proof}[Proof of Proposition~\ref{prop:equivalence_ConvexOrder_CsiszarOrder}]
Recall that if $f$ belongs to $ \F$, so does $f^*$. We deduce from~\eqref{eq:Sf_et_Z} that if $Z \leqcv Z'$, then $X \csleq X'$ w.r.t.\ $Y$. Suppose now that $X \csleq X'$ w.r.t.\ $Y$.
For $x,t \in \R$, we define the function $\varphi_t$ on $\R^+$ as: 
\[
\varphi_t(x) = (x-t)_+.
\]
For any $t \in \R$, the function $\tilde\varphi_t=\varphi_t-\varphi_t(1)$ belongs to $\F$. 
% Consequently, $\tilde\varphi_t^*$ also belongs to $\F$. 
% Since $f^{**}=f$ for every $f\in\F$ and $Z\geq0$ a.s., 
Equation~\eqref{eq:Sf_et_Z} gives $S_{\tilde\varphi_t}(X,Y) = \E[\tilde\varphi_t(Z)] = \E[(Z-t)_+] - (1-t)_+$, and the same identity holds for $Z'$.
Since we suppose that $X \csleq X'$ w.r.t.\ $Y$, we get that $S_{\tilde \varphi_t}(X,Y) \leq S_{\tilde \varphi_t}(X',Y)$ and so $\E[(Z-t)_+] \leq \E[(Z'-t)_+]$ for all $t \in \R$. 
As, $Z$ and $Z'$ are integrable and $\E[Z]=\E[Z']=1$, using Lemma~\ref{lem:convex_order_characterization}, we get $Z \leqcv Z'$.
This concludes the proof.

\end{proof}

\section{Proofs of the results in Section~\ref{sec:Gaussian}}
\label{sec:proof-G}

\subsection{On powers of covariance matrices}
\label{sec:MatrixDefinitions}
In this subsection, we introduce the notation used throughout this section and recall some standard facts on positive semi-definite matrices, Gaussian distributions, and correlation matrices.

With  a slight  abuse of  notation, we  denote by  $0$  any  matrix whose entries are all zero. Let $d \ge 1$ and denote by $\idmatrix_d$ the identity matrix of size $d\times d$. For a $d  \times d$ real matrix $M$, we denote its
image  and its  kernel  by $\im(M)  =  \{ Mx  : x  \in  \R^d \}$  and
$\ker(M) = \{ x \in \R^d :  Mx = 0\}$.  If furthermore the matrix $M$
is symmetric  positive semi-definite,  we denote  by $r=\dim(\im(M))$
its  rank.
If  $M=0$, that  is,  $r=0$, we  set $M^p=0$  for all  $p\in \R$.  If  $M\neq 0$,  denote  $\lambda_1\geq \ldots\geq  \lambda_d$    its
(nonnegative)   eigenvalues  ranked   in   descending   order  and   let
$\Lambda_M=  \diag(\lambda_1,  \dots,  \lambda_d)$ be  the  $d\times  d$
diagonal        matrix       of        its       eigenvalues        and
$\Lambda_M^+= \diag(\lambda_1,  \dots, \lambda_r)$  be the $r  \times r$
diagonal  matrix  of  its  positive eigenvalues.   We  then  consider  a
$d \times r$ matrix $R_M=[u_1, \ldots,  u_r]$ whose columns are given by
$r$  orthonormal  eigenvectors  associated   with  $u_i$  associated  to
$\lambda_i$.  Let $K_M$  be a $d \times (d-r)$ matrix  whose columns are
given  by  an  orthonormal  basis  of  $\ker(M)$  and  thus  the  matrix
$U_M = [R_M, K_M]$ is a  $d\times d$ orthogonal matrix. It is elementary
to check that:
\begin{equation}
   \label{eq:A-RA-LA}
   M= U_M  \Lambda_M U_M^{\transpose} = R_M \Lambda_M^+ R_M^\transpose,
   \quad
   \Lambda_M^+ = R_M ^\transpose M R_M
   \quad\text{and}\quad
   R_M^\transpose R_M= \idmatrix_r.
\end{equation}
The  matrices   $R_M  R_M^\transpose$   and  $K_M   K_M^\transpose$  are
orthogonal projectors  onto $\im(M)$ and $\ker(M)$,  respectively. We
define for $p\in \R$ the $r\times r$ diagonal matrix
$\Lambda_M^{p}  =  \diag(\lambda_1^{p}, \dots,  \lambda_r^{p})$
and the $d\times d$ matrix:
\[
  M^{p}= R_M \Lambda_M^{p} R_M^\transpose.
\]
Notice that the notation is consistent as $M^1=M$ (but $M^0=\idmatrix_d$ if and
only if $M$ is definite positive),  and that
 the matrix  $M^{p}$ is  uniquely defined  (even though  $R_M$
might  not  be  unique). The matrix   $M^{-1/2} M  M^{-1/2}= M^0 = R_M \,
R_M^\transpose $  corresponds to   the  orthogonal
projector with the same image and kernel as $M$. When $M$ is definite
positive, $M^{-1/2}$ is simply the (symmetric definite positive) square root of
$M^{-1}$.

\subsection{Proofs}

\begin{proof}[Proof of Lemma~\ref{lem:f_div_Gaussian}]
Recall that, for $i=1,2$, $P_i$ is the Gaussian distribution $\mathcal{N}(\mu_i,\Sigma_i)$, $P_1$ and $P_2$ 
have the same support (that is, $\mu_1 - \mu_2 \in \im(\Sigma_1) = \im(\Sigma_2)$), and $r = \dim(\im(\Sigma_i)) > 0$. This
gives the following density ratio: 
\begin{equation}
  \label{eq:likelihood ratio}
  \frac{\rd P_1}{\rd P_2}(x) = \sqrt{\frac{\det(\Sigma_2)}{\det(\Sigma_1)}} \exp\left( -\frac 12 \left[ (x-\mu_1)^\transpose \Sigma_1^{-1}(x-\mu_1)-(x-\mu_2)^\transpose \Sigma_2^{-1}(x-\mu_2) \right] \right), \quad (P_1 + P_2)\text{-a.s.}    
\end{equation}
According to Corollary~\ref{cor:invariance_Df}, if the ratio $\rd P_1/\rd P_2$ is $\sigma(h)-$measurable for some measurable function $h$, then:
\begin{equation}
\label{eq:invariancefdiv}
D_f(P_1 \Vert P_2) = D_f(P_{h(X_1)} \Vert P_{h(X_2)}),
\end{equation}
where $X_i \sim P_i$ and $P_{h(X_i)}$ is the distribution of $h(X_i)$ for $i=1,2$.
With $M = \Sigma_2^{-1/2} \Sigma_1 \Sigma_2^{-1/2}$ having rank $r$ we claim that:
\begin{equation}
    \label{eq:LR1}
 (x-\mu_2)^\transpose \Sigma_2^{-1} (x -\mu_2)
 =  (x-\mu_2)^\transpose \Sigma_2^{-1/2} R_M R_M^\transpose \Sigma_2^{-1/2} (x -\mu_2)
 = h(x)^\transpose h(x),
\end{equation}
where we denote $h$ the affine function  from $\R^d$ to $\R^r$ defined by
\[
h(x)= R_M^\transpose \Sigma_2^{-1/2} (x-\mu_2) .
\]
Indeed, let us recall that $ R_{\Sigma_i} R^\transpose_{\Sigma_i}$  is the
matrix of the orthogonal projector onto $\im(\Sigma_i)$. Since
$\im (\Sigma_i)=\im (M)$, we deduce that $R_{\Sigma_i}R_{\Sigma_i}^\transpose= R_{M} R_{M}^\transpose$
which implies the equality $\Sigma_2^{-1/2} R_M R_M^\transpose \Sigma_2^{-1/2} = \Sigma_2^{-1}$ used in the claim above. 
Moreover, we can write that:
\begin{align*}
\Sigma_2^{-1/2} R_M R_M^\transpose M^{-1} R_M R_M^\transpose \Sigma_2^{-1/2}
&= \Sigma_2^{-1/2} R_M R_M^\transpose \Sigma_2^{1/2} \Sigma_1^{-1} \Sigma_2^{1/2} R_M R_M^\transpose \Sigma_2^{-1/2}\\
&= \Sigma_2^{-1} \Sigma_2^{1} \Sigma_1^{-1} \Sigma_2^{1} \Sigma_2^{-1} \\
&= R_M R_M^\transpose \Sigma_1^{-1} R_M R_M^\transpose = \Sigma_1^{-1},
\end{align*}
and: 
\[
% \Sigma_1^{-1} = \Sigma_2^{-1/2} R_M R_M^\transpose M^{-1} R_M R_M^\transpose \Sigma_2^{-1/2} \quad \text{and} \quad 
R_M^\transpose M^{-1} R_M = (\Lambda_M^+)^{-1},
\]
which give:
\begin{align}    \label{eq:LR2}
    (x-\mu_1)^\transpose \Sigma_1^{-1}(x-\mu_1) &= (h(x) - h(\mu_1))^\transpose R_M^\transpose M^{-1} R_M (h(x) - h(\mu_1))\nonumber \\
&= (h(x) - h(\mu_1))^\transpose (\Lambda_M^+ )^{-1}(h(x) - h(\mu_1)).
\end{align}
By plugging \eqref{eq:LR1} and \eqref{eq:LR2} in \eqref{eq:likelihood ratio}, we conclude that the ratio $\rd P_1/\rd P_2$ is $ \sigma(h)$-measurable and use Corollary~\ref{cor:invariance_Df}.
% Corollary~3 in \cite{butucea2026csisz}. 
We now calculate the distributions of $Y_i := h(X_i)$. 

The function $h$ is affine, therefore $Y_i$ are also Gaussian vectors. 
We  have  $\E[Y_1]=  R_M^\transpose \Sigma_2^{-1/2}  (\mu_1  - \mu_2)  = \mu$  and
$\cov   (Y_1,    Y_1)   =   R_M^\transpose   M    R_M   =   \Lambda_M^+$
by~\eqref{eq:A-RA-LA}, so $Y_1 \sim \cn(\mu, \Lambda_M^+) $. We also
have $\E[Y_2]=0$ and $\cov(Y_2, Y_2)= R_M^\transpose     \Sigma_2^{-1/2}      \Sigma_2
\Sigma_2^{-1/2}                        R_M$.
According to~\eqref{eq:A-RA-LA}, we have that $ \Sigma_2^{-1/2}      \Sigma_2
\Sigma_2^{-1/2} = R_{2} R_{2}^\transpose$ which is also equal to $R_M R_M^\transpose$, and conclude that %from the last equality in~\eqref{eq:A-RA-LA} that 
$\cov(Y_2, Y_2)=\idmatrix_r$. Thus  $Y_2 \sim \cn(0, \idmatrix_r) $. Using \eqref{eq:invariancefdiv}, we obtain that:
\[
D_f(P_1\Vert P_2) %=   D_f(X_1\Vert X_2) = D_f(Y_1\VertY_2)
= D_f(P_{h(X_1)}\Vert P_{h(X_2)})
= D_f(\cn(\mu, \Lambda_M^+) \Vert \cn(0,\idmatrix_{r})).
\]
Obviously, $\mu=0$ when $\mu_1 = \mu_2$.
\end{proof}

The  next  lemma on  convex  order  is a  first  step  in the  proof  of
Proposition~\ref{prop:gaussian_gfd_properties}.                      Let
$p(\cdot,\sigma^2) $ be  the probability density function  of a centered real
Gaussian distribution with variance $\sigma^2>0$.

\begin{lem}[Convex order of Gaussian densities ratios]
\label{lem:convex_order_gaussian_density_ratio}
Let $X \sim \mathcal{N}(0,1)$, $c,c' \in \R_*^+$.  If $c \geq c' \geq 1$
or  $1 \geq c' \geq c$, we have:
\begin{equation}
    \label{eq:convex_order_gaussian_density_ratio}
\frac{p(X,c')}{p(X,1)} \, \leqcv \,  \frac{p(X,c)}{p(X,1)}\cdot
\end{equation}
\end{lem}

\begin{proof}[Proof of Lemma~\ref{lem:convex_order_gaussian_density_ratio}]
Let $X$ and $Y$ be two nonnegative random variables with finite means.
  Following \cite{shaked2007stochastic} Equation (3.A.33) and Section
  4.B.1, we say that $X$ is smaller than $Y$ in the Lorenz order, denoted by $X \, \leqlorenz \, Y$,  if:
\begin{equation}
   \label{eq:Lorentz}
\frac{X}{\E[X]} \, \leqcv \, \frac{Y}{\E[Y]}\cdot
\end{equation}

Let  $\tilde{p}(\cdot,\sigma^2)$ denote the probability density function
of the rescaled random variable $\sigma \tilde{X}$ where  $\tilde{X} = \vert X \vert$ and $X \sim
\cn(0, 1)$. Remark that $\rd P_{\sigma X} / \rd P_X (x)$ is a function of $\vert x \vert$. Thanks to Corollary~\ref{cor:invariance_Df},
we get that $D_f(P_{\sigma X} \Vert P_X ) = D_f(P_{\sigma \vert X \vert} \Vert P_{\vert X \vert})$ for any $f \in \F$, that is:
\begin{equation}
    \label{eq:equality_E_density_ratio_gaussian_absgaussian}
    \E \left[f\left( \frac{p(X,\sigma^2)}{p(X,1)}\right) \right] = \E
    \left[f\left(
        \frac{\tilde{p}(\tilde{X},\sigma^2)}{\tilde{p}(\tilde{X},1)}\right)
    \right].
\end{equation}
Let  $c,c' \in  \R_*^+$.  Assume  first  that $c \geq c' \geq 1$.  The function  $x  \mapsto  \tilde{p}(x,c)/\tilde{p}(x,c')$ is  positive  and increasing on $\R_+^*$.  Using \cite[Theorem~3.A.26]{shaked2007stochastic}, we obtain:
\begin{equation}
   \label{eq:Lorentz-N}
\frac{\tilde{p}(\tilde{X},c')}{\tilde{p}(\tilde{X},1)} \, \leqlorenz \,
\frac{\tilde{p}(\tilde{X},c)}{\tilde{p}(\tilde{X},1)} \cdot
\end{equation}
Since both sides of~\eqref{eq:Lorentz-N}  have the same mean  equal to
$1$, we deduce from~\eqref{eq:Lorentz} that~\eqref{eq:Lorentz-N}  holds
also for the  convex order. Then,
use~\eqref{eq:equality_E_density_ratio_gaussian_absgaussian} to
get that:
\begin{equation}
    \label{eq:convex_order_gaussian_density_ratio_case1}
\frac{p(X,c')}{p(X,1)} \, \leqcv \,  \frac{p(X,c)}{p(X,1)}\cdot
\end{equation}
We now assume that $1 \geq c' \geq c > 0$, and we set $\theta = 1/c$ and $\theta' = 1/c'$, so that $\theta \geq \theta' \geq 1$. For any $a \in \R_+^*$ and any $f \in \F$, we get (with the change of variable $y=a x$):
\begin{equation}
  \label{eq:equality_E_density_ratio_gaussian}
 \E\left[f^*\left(\frac{p(X,1/a)}{p(X,1)}\right)\right]  
%  = \E\left[f^*\left(\frac{p(\sqrt{a}\, X,1)}{p(\sqrt{a}\, X,a)}\right)\right]
  = \int_\R p\left(x,1/a\right) \, f \left(\frac{p(x,1)}{p(x,1/a)} \right)\,  \rd x 
  =\E\left[f\left(\frac{p(X,a)}{p(X,1)}\right)\right].
\end{equation}
Applying~\eqref{eq:equality_E_density_ratio_gaussian} with $a = c'$ and $a = c$, and then~\eqref{eq:convex_order_gaussian_density_ratio_case1}   to the convex function~$f^*$, we obtain: 
\[
  \E\left[f\left(\frac{p(X,c')}{p(X,1)}\right)\right]
  = \E\left[f^*\left(\frac{p(X,\theta')}{p(X,1)}\right)\right]
  \leq \E\left[f^*\left(\frac{p(X,\theta)}{p(X,1)}\right)\right]
  = \E\left[f\left(\frac{p(X,c)}{p(X,1)}\right)\right].
\]
This complete the proof of the Lemma. 
\end{proof}

\begin{proof}[Proof of Proposition~\ref{prop:gaussian_gfd_properties}]
  Let $x=(x_1, \ldots, x_d) \in [0, 1]^d$. 
  Since for any probability distributions $P,Q,P',Q'$ we have $D_f(P \otimes P' \Vert Q \otimes Q') = D_f(P' \otimes P \Vert Q' \otimes Q)$, the proof of~\ref{enum:gaussian_gfd_properties_permutation} is immediate.

  We now prove~\ref{enum:gaussian_gfd_properties_min}: if $x=0$ is the null vector,  then
  $\Gamma_x= \idmatrix_{2d}$ and by the definition of $\kappa_f^{(d)}$, we get
  $\kappa_f^{(d)}(x)=0$. 
  
  We now prove~\ref{enum:gaussian_gfd_properties_max}. If the vector $x$ has at least one entry equal to the maximum possible value 1, 
  then $\Gamma_x$ cannot be full rank, \emph{i.e.\ }$\rank (\Gamma_x)\leq  2d-1 .$ Therefore, the distributions $\cn(0,
  \Gamma_x)$ and $\cn(0, \idmatrix_{2d})$ are mutually singular and thus
  $\kappa_f^{(d)}(x)=f(0) + f^*(0)$, which is the largest attainable value for the $f$-divergence.
  
  To prove~\ref{enum:gaussian_gfd_properties_dim}, as
  $\kappa_f^{(d)}$ is  symmetric, we  can assume  without loss  of generality
  that it is the last coordinate which is null: $x_d=0$. Then we get that
  $\cn(0, \Gamma_x)=\cn (0, \Gamma_{x'}) \otimes \cn(0, \idmatrix_{2})$,
  with $x'=(x_1, \ldots, x_{d-1})$.  
  It is elementary to check that $D_{f}(\cn(0, \Gamma_{x'}) \otimes \cn(0, \idmatrix_{2}) \, \Vert \, \cn(0, \idmatrix_{2(d-1)}) \otimes \cn(0, \idmatrix_{2})) = D_{f}(\cn(0, \Gamma_{x'}) \, \Vert \, \cn(0, \idmatrix_{2(d-1)}))$.
  Then we get \ref{enum:gaussian_gfd_properties_dim}.

\bigskip

  Lastly,  we  prove~\ref{enum:gaussian_gfd_properties_monoticity}.  Let
  $y=(y_1, \ldots,  y_d)\in [0,  1]^d$ such that  $y_i\leq x_i$  for all
  $i    \in    \{1,    \dots,     d\}$.     If  $x$ has at least one entry equal to 1, we have seen that $\kappa_f^{(d)}(x)$ already attains the maximal value of the divergence, which equals $f(0)+f^*(0)$. Thus,   we        get
  $\kappa_f^{(d)}(y)  \leq  \kappa_f^{(d)}(x)$, which  proves  the  claim in  that
  case.  
  
  From now  on, we  assume that  $\max y  \leq \max  x <1$.   Let
  $X_1  , \dots  , X_{2d}$ be independent centered  Gaussian random
  variable      with      variance      equal to  1.   Then, the likelihood of this vector under the distribution $\cn(0,\Gamma_x)$ can be written as:
  \[
  \prod_{i=1}^d p(X_i, 1-\sqrt{x_i}) \cdot p(X_{d+i}, 1+\sqrt{x_i}),
  \]
   where $p(\cdot,\sigma^2)$ denotes the probability density  of $\cn(0, \sigma^2)$. 
  Let us write $\kappa_f^{(d)}(x) = \E [f (V_x)]$, where $V_x$ is the likelihood ratio: 
  \[
  V_x=  \prod_{i = 1}^d \frac{p(X_i, 1-\sqrt{x_i})}{p(X_i,1)} \cdot \frac{p(X_{d+i}, 1+\sqrt{x_i})}{p(X_{d+i},1)} \cdot % \quad \text{and} \quad X \sim \cn(0, \idmatrix_{2d}).
  \]
  Analogously, $\kappa_f^{(d)}(y) = \E[f(V_y)]$. 
  From   Lemma~\ref{lem:convex_order_gaussian_density_ratio},  we   can  deduce   that for all $i \in \{1, \dots, d\}$:
    \[
    \frac{p(X_i, 1-\sqrt{y_i})}{p(X_i,1)} \, \leqcv \, \frac{p(X_i,
      1-\sqrt{x_i})}{p(X_i,1)}
    \quad\text{and}\quad
        \frac{p(X_{d+i}, 1+\sqrt{y_i})}{p(X_{d+i},1)} \, \leqcv \, \frac{p(X_{d+i}, 1+\sqrt{x_i})}{p(X_{d+i},1)}\cdot
    \]
Since the random variables $X_1, \ldots, X_{2d}$ are independent, we
deduce from  \cite[Corollary~3.A.22]{shaked2007stochastic} that 
\[
V_y \, \leqcv \, V_x.
\]
    Then, for all $f \in \F$, we get that
    $    \kappa_f^{(d)}(y )= \E[f(V_y)] \leq  \E[f(V_x)]=  \kappa_f^{(d)}(x)$, which proves~\ref{enum:gaussian_gfd_properties_monoticity}.
\end{proof}

\begin{proof}[Proof   of  Theorem~\ref{thm:csiszar_index_Gaussian}]   To
  simplify, we set $X_1=X$ and  $X_2=Y$. We consider the dimension $d_i$
  of $X_i$,  the covariance  matrix $\Sigma_i$ of  $X_i$, $r_i$ its
  rank, and if  $r_i\geq 1$ consider the projector  $Q_i=R_i R_i^\transpose$ with
  the shorter notation $R_i=R_{\Sigma_i}$.  If  $\min(r_1,r_2)=0$, then $S_f(X_1,  X_2)=0$ as
  $X_1$  and  $X_2$  are  independent.   In  this  case,  by  convention
  $C=\corr(X_1, X_2)=0$  and thus  $k=\rank(C)=0$. As  $\kappa^{(0)}_f=0$ by convention, we
  get that the equality~\eqref{eq:Sf=gf} holds.

\medskip

We    now    assume    that   $\min(r_1,r_2)\geq    1$.    
By considering injective and bi-measurable functions of the random variables 
$X_1,\, X_2$, the Csisz\'ar index is unchanged (see Proposition~\ref{prop:invariance_csiszar_index}). In particular, apply translations to the Gaussian variables to see that the value of $S_f(X_1, X_2)$
does not depend on the mean of $X_i$; so,  without loss of generality,
we assume that $\E[X_i]=0$ for both $i=1$ and $i=2$.  The Gaussian distributions $P_{(X_1, X_2)}$
and  $P_{X_1}\otimes P_{X_2}$  are centered  with respective  covariance
matrix:
\[
  \Sigma=  \begin{pmatrix}
    \Sigma_1 & \cov(X_1,X_2) \\
    \cov(X_2,X_1)  & \Sigma_2
  \end{pmatrix}
  \quad\text{and}\quad
  \tilde \Sigma= \begin{pmatrix} \Sigma_1 & 0 \\ 0 & \Sigma_2 \end{pmatrix}.
\]
We consider the correlation
matrix  between $X_1$ and $X_2$:
\[
  C=\corr(X_1, X_2)= \Sigma_{1}^{-1/2} \cov (X_1, X_2)
  \Sigma_{2}^{-1/2},
\]
and the matrix $ M=\tilde \Sigma^{-1/2} \Sigma \tilde \Sigma^{-1/2} $. Using  $
  \tilde \Sigma^p =
  \begin{pmatrix} \Sigma_1^p & 0 \\ 0 & \Sigma_2^p \end{pmatrix}$
for any real number $p$ and taking into account that $\Sigma_1$ and $\Sigma_2$ are non-degenerate because we assumed that  $\min(r_1,r_2)\geq 1$,   we get:
\[
M  =  \begin{pmatrix} Q_1 & C \\ C^\transpose & Q_2 \end{pmatrix}.
\]

\medskip

Let $x=(x_1,  x_2)\in \ker  (\tilde \Sigma)$, that is $\tilde{\Sigma} x = 0$. It implies that  $x_i\in \ker(\Sigma_i)$
and thus  a.s.\ $X_i^\transpose x_i=0$  for $i=1, 2$. This  implies that
$\Sigma  x=0$,   that  is  $\ker(\tilde \Sigma)  \subset   \ker(\Sigma)$. % and thus $\im(\Sigma)\subset   \im   (\tilde \Sigma)$.  

We claim that $\im(\tilde \Sigma) \cap \ker(\Sigma) \neq \{0\}$ if and only if $1$ is an eigenvalue of $C C^\transpose$.
Suppose that $\im(\tilde \Sigma) \cap \ker(\Sigma) \neq \{0\}$.
Let $x=(x_1, x_2) \in  \im(\tilde \Sigma) \cap \ker(\Sigma) \setminus \{0\}$. There exists $y=(y_1, y_2)$ non-zero such that $y_i$ belongs to $ \im(\Sigma_i)$ which is the same space as $\im (\Sigma_i^{1/2})$ and $x_i= \Sigma_i^{-1/2} y_i$, for $i=1,2$. This gives $x= \tilde \Sigma^{-1/2} y$ and thus
$My= \tilde \Sigma^{-1/2} \Sigma x = 0$. As $Q_i y_i=y_i$, we deduce that $C C^\transpose y_i=y_i$ and $1$ is an eigenvalue of $C C^\transpose$ as $y_1$ or $y_2$ is non-zero.
We now prove the converse implication. Assume that $1$ is an eigenvalue of $C C^\transpose$. Let $y_1$ be an eigenvector of $C C^\transpose$ associated to the eigenvalue $1$. Let $y_2 = C^\transpose y_1$, $x_1 = \Sigma_1^{-1/2} y_1$ and $x_2 = -\Sigma_2^{-1/2} y_2$. 
If $x_1 = 0$, then $C^\transpose y = 0$ and then $CC^\transpose y_1=0$, which is impossible as $y_1$ is an eigenvector; thus we get that $x_1 \neq 0$ and so $x \neq 0$. 
Moreover, we have:
\[
\Sigma x = 
\tilde \Sigma^{1/2} M 
\begin{pmatrix}
  y_1 \\ -y_2
\end{pmatrix}
=
\tilde \Sigma^{1/2}
\begin{pmatrix}
  y_1 - C y_2 \\ C^\transpose y_1 - y_2
\end{pmatrix}
= 
\tilde \Sigma^{1/2}
\begin{pmatrix}
  y_1 - C C^\transpose y_1 \\ C^\transpose y_1 - C^\transpose y_1
\end{pmatrix}
= 
\begin{pmatrix}
  0 \\ 0
\end{pmatrix},
\]
that is $x \in \ker(\Sigma)$. We also have by construction that $x \in \im (\tilde \Sigma^{-1/2})= \im(\tilde \Sigma)$. 
Thus $\im(\tilde \Sigma) \cap \ker(\Sigma) \neq \{0\}$. This proves the claim.

In conclusion, since $\ker(\tilde \Sigma) \subset \ker(\Sigma)$, we get that $\ker(\Sigma)=\ker(\tilde \Sigma)$ (that is $\im(\Sigma)=\im(\tilde \Sigma)$) if and only if $1$ is not an eigenvalue of $C C^\transpose$ (or equivalently of $C^\transpose C$).

\medskip If $1$  is an eigenvalue of $C C^\transpose$,  then the supports of
the probability measures $P_{(X_1, X_2)}$ and $ P_{X_1}\otimes P_{X_2}$,
that  is,   $\im(\Sigma)$  and  $\im(\tilde \Sigma)$ are   distinct,  so  the
probability     measures    are     mutually    singular     and    thus
$S_f(X_1,
X_2)=f(0)+f^*(0)$. Using~Proposition~$\ref{prop:gaussian_gfd_properties}~\ref{enum:gaussian_gfd_properties_max}$
as  $1$  is  an  eigenvalue  of  $C  C^\transpose$,  we  also  get  that $\kappa^{(d)}_f(\lambda)=f(0)+ f^*(0)$ and thus~\eqref{eq:Sf=gf} holds.

\medskip From now on we assume that $1$ is not an eigenvalue of $C
C^\transpose$, and thus that $\im(\Sigma)=\im(\tilde \Sigma)$ and
$\ker(\Sigma)=\ker(\tilde \Sigma)$.
The rank of $\Sigma$ is equal to  $r=r_1+r_2$ and the matrix $M$ is
semi-definite positive with the same rank.
   Using Lemma~\ref{lem:f_div_Gaussian}, we get:
\[
S_f(X,Y) =D_{f}\left(P_{(X_1,X_2)} \Vert P_{X_1} \otimes
  P_{X_2}\right) =  D_{f}\left( \cn(0, \Lambda^+) \Vert
\cn(0, \idmatrix_{r})\right),
\]
where $\Lambda^+$ is a diagonal matrix whose diagonal is given by
the $r$ positive eigenvalues of $M$. The eigenvalues of
$M$ are the roots of $\det (t \idmatrix_{d} - M)$ with $d=d_1+d_2$.

Recall that the determinant of the block matrix 
$\begin{pmatrix}
A & B \\ C & D
\end{pmatrix}$ 
is equal to $\det(D) \cdot \det (A- B D^{-1} C)$ provided
$\det(D)\neq 0$. Thus, provided $\det(t \idmatrix_{d_2} - Q_2) \neq 0$, we get:
\begin{align*}
    \det(t \idmatrix_{d} - M) &= \det \begin{pmatrix}
        t \idmatrix_{d_1} - Q_1 & C\\
        C^\transpose & t \idmatrix_{d_2} - Q_2 
    \end{pmatrix}\\
    & = \det(t \idmatrix_{d_2} - Q_2) \cdot \det ( t \idmatrix_{d_1} -  Q_1 -  C  (t \idmatrix_{d_2} - Q_2)^{-1} C^\transpose ).
\end{align*}
Notice that if the matrix $Q'$ is a projector in $\R^{d'}$, then $\det(t\idmatrix_{d'} -Q')= t^{d'-r'} (t-1)^{r'}$ is non-zero for $t\not \in \{0,1\}$, with $r'=\rank (Q')$, and the inverse of the matrix $t\idmatrix_{d'} - Q'$ exists and is given by
$t^{-1} \idmatrix_{d'} + (t(t-1))^{-1} Q'$.
We get for $t\not \in \{0,
1\}$:
\begin{align*}
  \det (t \idmatrix_{d} - M)
  & =t^{d_2 -r_2} (t-1)^{r_2} \det \left(t \idmatrix_{d_1} -  Q_1 - \inv{t} C \left(\idmatrix_{d_2} + \inv{t-1} Q_2\right) C^\transpose\right)\\
  & =t^{d_2 -r_2} (t-1)^{r_2} \det \left(t \idmatrix_{d_1} -  \left(Q_1 + \inv{t-1} C C^\transpose\right)\right).
\end{align*}
Indeed, $\Sigma_2^{-1/2} Q_2 \Sigma_2^{-1/2}= \Sigma^{-1}_2$ thanks to~\eqref{eq:A-RA-LA} and thus $C Q_2 C^\transpose =C C^\transpose$. 
Now, let us use that $Q_1$ and $ \idmatrix_{d_1} -  Q_1 $ are orthogonal projectors on $\im(\Sigma_1)$ and $\ker(\Sigma_1)$ respectively and that $\im (C C^\transpose) \subset \im (\Sigma_1)=\im (Q_1)$ to get:
\begin{align*}
  \det (t \idmatrix_{d} - M)
  &= t^{d_2 -r_2} (t-1)^{r_2 - d_1} \det \left( (t-1) 
  \begin{pmatrix}  
    (t-1) \idmatrix_{r_1} & 0 \\ 
    0 & t \idmatrix_{d_1-r_1}  
  \end{pmatrix} - 
  \begin{pmatrix}
    \diag(\lambda_1, \dots, \lambda_{r_1}) & 0
     \\ 0 & 0
  \end{pmatrix} \right)  \\
  & =t^{d -r} (t-1)^{r_2 - r_1} \det \left((t-1)^2 \idmatrix_{r_1} -  \diag(\lambda_1, \dots, \lambda_{r_1}) \right)\\
  &= t^{d -r} (t-1)^{r_2 - r_1}  \prod_{j=1}^{r_1} \left((t-1)^2 - \lambda_j\right)\\
  &= t^{d -r} (t-1)^{r-2k}  \prod_{j=1}^{k} \left(t-(1+\sqrt{ \lambda_j})\right)\,  \cdot  \left(t-(1- \sqrt{\lambda_j})\right),
\end{align*} 
with the notation $\lambda_1\geq  \ldots \geq  \lambda_{r_1}$ for  the (nonnegative) eigenvalues of $C C^\transpose$ ranked in descending order, such that only
the $k=\rank(C C^\transpose)$ first are positive.
By continuity in $t$, we deduce that the expression of $\det (t \idmatrix_{d} - M)$ holds for all $t\in \R$,
and thus the eigenvalues of $M$ are $0$ with multiplicity $d-r$, $1$
with multiplicity $r-2k$ and $1+\sqrt{ \lambda_j}$ as well as $1-\sqrt{
  \lambda_j}$ for $j\in \{1, \ldots k\}$. We consider
$\Gamma_\lambda$ defined in~\eqref{eq:def-Gamma-x} with
$x=\lambda=(\lambda_1, \ldots, \lambda_k)$. Since  one can choose
$\Lambda^+$ equal to $\begin{pmatrix}
   \idmatrix_{r-2k} & 0\\ 0 & \Gamma_\lambda
\end{pmatrix}$, we get:
\[
   D_{f}\left( \cn(0, \Lambda^+) \Vert
     \cn(0, \idmatrix_{r})\right)
=  D_{f}\left( \cn(0, \Gamma_\lambda) \Vert
     \cn(0, \idmatrix_{2k})\right)=
   \kappa_f^{(k)}(\lambda),
\]
where we used the definition of $\kappa_f^{(k)}$ for the last equality.  This ends the proof.
 \end{proof}

\begin{proof}[Proof of Corollary~\ref{cor:gaussian_case_dim_1}]
  Proof   of~$\ref{enum:gaussian_dim1_case1}$.   Suppose   that
  $d_Y = d_X = d_{X'} = 1$.  From the definition of the matrix $C$ given
  in Theorem~\ref{thm:csiszar_index_Gaussian}  for the pair  $(X,Y)$, we
  get    that    $C$     is    a    $1    \times     1$    matrix    and
  $C^{\transpose}  C=\corr(X,Y)^2$.  For  any  $f \in  \F$,  this  gives
  $S_f(X,Y)  = \kappa_f^{(1)}\left(\corr(X,Y)^2\right)$.  Since for example the function
  $\kappa_\mathrm{KL}^{(1)}$ is increasing (see Table~\ref{tab:expression_gf1}), the result follows.

Proof of~$\ref{enum:gaussian_dim1_case2}$.  If $d_Y = 1$
(resp. if $d_X = d_{X'} = 1$), then $C^{\transpose} C$ (resp. $C
C^{\transpose}$) is a $1 \times 1$ matrix which is   equal to $h(X,Y)$ (with $h$ given by \eqref{eq:def-h}),
and consequently $S_f(X,Y) = \kappa_f^{(1)}\left(h(X,Y)\right)$.
We conclude again using that the function $\kappa_\mathrm{KL}^{(1)}$ is increasing.
\end{proof}

\begin{proof}[Proof of Lemma~\ref{lem:cvx-kappaf} \ref{enum:kappa_KL} and \ref{enum:kappa_alpha}]

We first prove~\ref{enum:kappa_KL}. From the expressions of $\kappa_{\mathrm{KL}}^{(d)}$ and $\kappa_{\mathrm{KL}^*}^{(d)}$
given in Table~\ref{tab:expression_gf1}, we immediately get they are convex as sums of convex functions.

We now prove~\ref{enum:kappa_alpha}. Recall that we have the convention $1/0 = +\infty$. 
Let $d \in \N^*$, $\alpha \in \R \setminus \{0,1\}$ and define the function $u_\alpha$ on $I_\alpha = [0,\min(1/\alpha^2,1))$ by:
\begin{equation}
  \label{eq:ualpha_lemmeconvexite}
  u_\alpha(t) = \frac{(1-t)^\alpha}{1-\alpha^2 t} \cdot
\end{equation}
Then, using expressions from Table~\ref{tab:expression_gf1}, for $x \in I_\alpha^d$ we have:
\[
\kappa_{1 - \alpha}^{(d)}(x) = c_\alpha \left( 1 - \prod_{i=1}^d \sqrt{u_\alpha(x_i)} \right),
\]
where $c_\alpha = (\alpha(1-\alpha))^{-1}$.
If $c_\alpha < 0$ (that is, $\alpha \notin [0,1]$), then $\kappa_{1 - \alpha}^{(d)}$ is convex if and only if $x \mapsto \prod_{i=1}^d \sqrt{u_\alpha(x_i)}$ is convex on $I_\alpha^d$.
If $c_\alpha > 0$ (that is, $\alpha \in (0,1)$), then $\kappa_{1 - \alpha}^{(d)}$ is convex if and only if
$x \mapsto \prod_{i=1}^d \sqrt{u_\alpha(x_i)}$ is concave on $I_\alpha^d$.

\medskip
\noindent\textbf{Case 1: $\alpha \notin [0,1]$ (\emph{i.e.\ }$c_\alpha < 0$).}
In this subcase, we prove that $\kappa_{1 - \alpha}^{(d)}$ is convex for every $d \in \N^*$.
Let $v_\alpha = \log(\sqrt{u_\alpha})$. This function is twice differentiable on the interval $I_\alpha$, with:
\[
v_\alpha''(t) = \frac{1}{2} \left( \frac{\alpha^4}{(1-\alpha^2 t)^2} - \frac{\alpha}{(1-t)^2}\right).
\]
For $\alpha < 0$, it is clear that $v_\alpha''(t) > 0$, so $v_\alpha$ is convex on $I_\alpha$. For $\alpha > 1$ we have:
\[
v_\alpha''(t) > 0 \iff \frac{\alpha^3}{(1-\alpha^2 t)^2} - \frac{1}{(1-t)^2} > 0
\iff (\alpha^{3/2} - 1) + \alpha^{3/2}\, t\, (\sqrt{\alpha} - 1) > 0,
\]
and the right-hand side is positive since $\alpha > 1$; hence $v_\alpha''(t) > 0$ and $v_\alpha$ is convex on $I_\alpha$. 
% Thus, for every $\alpha \notin [0,1]$, $v_\alpha$ is convex, so $\sqrt{u_\alpha}$ is log-convex. 
For $x \in I_\alpha^d$ we have $\prod_{i=1}^d \sqrt{u_\alpha(x_i)} = \exp ( \sum_{i=1}^d v_\alpha(x_i) )$ which is convex, and therefore $\kappa_{1 - \alpha}^{(d)}$ is convex.

\medskip
\noindent\textbf{Case 2: $\alpha \in (0,1)$ (\emph{i.e.\ }$c_\alpha > 0$).}
In this subcase, we suppose that $d \in \{1,2\}$ and we prove that $\kappa_{1 - \alpha}^{(d)}$ is convex.
Notice that $I_\alpha = [0,1)$. 
The function $u_\alpha$ is continuous and twice differentiable on $(0,1)$, with:
\[
u_\alpha''(t) = -\frac{\alpha \left(1 - t\right)^{\alpha} \left(1-\alpha \right) \left(t \alpha^2 (\alpha - 2)(t\alpha + 2) + 2 \alpha^2 + 1\right)}{\left(t - 1\right)^{2} \left(1-\alpha^{2} t\right)^{3}} \cdot
\]
The polynomial $t \alpha^2 (\alpha - 2)(t\alpha + 2) + 2 \alpha^2 + 1$ is a quadratic in $t$ with negative leading coefficient $\alpha^{3}(\alpha-2)$, hence concave; since its value at $t=0$ is $2 \alpha^2 + 1 > 0$ and its value at $t=1$ is $(\alpha^2-1)^2 > 0$, it is positive on all of $[0,1]$ for every $\alpha \in (0,1)$. The factors $\alpha$, $(1-t)^{\alpha}$, $(1-\alpha)$, $(t-1)^2$ and $(1-\alpha^2 t)^3$ are also positive for all $t \in [0,1]$ and $\alpha \in (0,1)$. Hence $u_\alpha''(t) \leq 0$ for all $t \in (0,1)$, so $u_\alpha$ is concave on $[0,1]$. Let $v(x) = \prod_{i=1}^d \sqrt{x_i}$ for $x \in [0,1]^d$. The function $v$ is concave and coordinate-wise nondecreasing (recall that $d \in \{1,2\}$). 
% By page~86 of~\cite{boyd2004convex}, 
Since $v$ is concave and coordinate-wise nondecreasing and $u_\alpha$ is concave, the function $w$ defined on $[0,1]^d$ by:
\[
w(x) = v(u_\alpha(x_1), \dots, u_\alpha(x_d)) = \sqrt{\prod_{i=1}^d \frac{(1-x_i)^{\alpha}}{1-\alpha^2 x_i}},
\]
is concave
and therefore $\kappa_{1 - \alpha}^{(d)}$ is convex for $d \in \{1,2\}$.

Thus, we have shown that, in the case $\alpha \notin [0,1]$, the function $\kappa_{1 - \alpha}^{(d)}$ is convex for every dimension $d \in \N^*$, and in the case $\alpha \in (0,1)$, the function $\kappa_{1 - \alpha}^{(d)}$ is convex for $d \in \{1,2\}$.
Since the sets $\R \setminus [0,1]$ and $(0,1)$ are stable under the transformation $\alpha \mapsto 1-\alpha$, then~\ref{enum:kappa_alpha} follows.

\end{proof}

\begin{rem}[The function $\kappa_\alpha^{(d)}$ is not convex for $\alpha d >2$ and $\alpha\in (0, 1)$]
Let $\alpha\in(0,1)$ and $d \in \N^*$. For $z \in [0,1]$, we have:
\[
\kappa_{\alpha}^{(d)}(z\mathds{1}) = c_\alpha (1 - v_\alpha(z)),
\]
with $v_\alpha = u_\alpha^{d/2}$ and $u_\alpha$ defined in~\eqref{eq:ualpha_lemmeconvexite}.
We show that $c_\alpha (1 - v_\alpha)$ is convex if and only if $\alpha d \leq 2$. 
In particular, if $\alpha d > 2$ (which can happen as soon
as $d\geq 3$), $\kappa_{\alpha}^{(d)}$ is not convex on the diagonal, and thus not convex on $[0,1]^d$.
 The function $v_\alpha$ is continuous on $[0,1]$ and twice differentiable on $(0,1)$. We have $v''_\alpha = w_1 w_2$ with:
\begin{align*}
w_1(t) &= \frac{d}{2} \alpha (1-\alpha) (1-t)^{\frac{d}{2} \alpha - 2} (1-\alpha^2 t)^{-\frac{d}{2} - 2}  \\
w_2(t) &= \frac{d}{2} \alpha (1-\alpha) (1+\alpha t)^2 + \alpha^3 t^2 + 2 \alpha^2 t - (1 + \alpha + \alpha^2).
\end{align*}
The function $w_1$ is positive for $t \in (0,1)$ (recall that $\alpha \in (0,1)$).
Note that $w_2$ is a quadratic polynomial in $t$, with $w_2'(t) = \alpha^2 (1+t \alpha) ((1-\alpha) d + 2)$, which is positive for $t \in [0,1]$; thus the function $w_2$ is increasing.
We have $w_2(1) = (1-\alpha)(1+\alpha)^2(\alpha d - 2)/2$.
If $\alpha d \leq 2$, then $w_2(1)\le 0$ and 
 $v_\alpha''(t) = w_1(t) w_2(t) < 0$ for all $t \in (0,1)$, so $v_\alpha$ is strictly concave.
If $\alpha d > 2$, then $w_2(1)>0$, and by continuity $w_2>0$ on some interval $(1-\varepsilon,1)$; thus $v_\alpha''>0$ on $(1-\varepsilon,1)$, so $v_\alpha$ is not concave.
 In conclusion, the function  $c_\alpha (1 - v_\alpha)$ is not convex if and only if $\alpha d > 2$.
So,  for $\alpha d>2$, the function $z \mapsto \kappa_{\alpha}^{(d)}(z\mathds{1})$
is not convex, and thus the function $\kappa_{\alpha}^{(d)}$ is not convex on $[0, 1]^d$. 
\end{rem}

The proof of Lemma~\ref{lem:cvx-kappaf}~\ref{enum:kappa_TV} relies on the following technical lemma. Its proof is given at the end of this section.

\begin{lem}[Taylor expansion of $\kappa_{\mathrm{TV}}^{(d)}$ at $0$]
\label{lem:dev_lim_kappatv}
Let $d \geq 1$, $x =(x_1, \dots , x_d)\in [0,1]^d$ and $Z_1, \dots, Z_d, Z_1', \dots, Z_d'$ being independent standard Gaussian random variables. As $\max(x)$ goes down to zero, we have:
\begin{equation}
    \label{eq:dev_lim_kappatv}
      \kappa_{\mathrm{TV}}^{(d)}(x) 
      = \frac{1}{2} \E \left[ \left\vert \sum_{i=1}^d \sqrt{x_i} (Z_i^2 - Z_i'^2) \right\vert \right] + O \left( \sum_{i=1}^d x_i \right).
      %= \E \left[ \left\vert \sum_{i=1}^d \sqrt{x_i} Z_i Z_i' \right\vert \right] + O \left( \sum_{i=1}^d x_i \right),
\end{equation}
and for $x = z\mathds{1} $ with $\mathds{1}\in \R^d$ being the vector of ones, we have  for $z \in [0,1]$ going to zero:
\begin{equation}
    \label{eq:dev_lim_kappatv_vec1}
    \kappa_{\mathrm{TV}}^{(d)}(z\mathds{1}) = \sqrt{z} \, \frac{2}{\sqrt{\pi}} \, \frac{\Gamma((d+1)/2)}{\Gamma(d/2)} + O(z).
\end{equation}
\end{lem}

\begin{proof}[Proof of Lemma~\ref{lem:cvx-kappaf}~\ref{enum:kappa_TV}]
Using Lemma~\ref{lem:dev_lim_kappatv}, we get in dimension $1$ that:
\begin{equation}
    \label{eq:dev_lim_kappatv_1}
    \kappa_{\mathrm{TV}}^{(1)}(z) = \frac{2}{\pi} \sqrt{z} + O\left(z \right).
\end{equation}
Thus, we get that $\kappa^{(1)}_{\mathrm{TV}}$ is not convex. 
In dimension $d$, we get:
\[
\kappa_{\mathrm{TV}}^{(d)}(z \mathds{1}) = \sqrt{z} \frac{2}{\sqrt{\pi}} R_d + O(z ), \quad \text{with} \quad R_d = \frac{\Gamma((d+1)/2)}{\Gamma(d/2)} \cdot
\]
We claim that $R_d > \sqrt{d/\pi}$ for $d \geq 2$. To prove it, we bound $R_d$ from both sides. Using Cauchy-Schwarz inequality, we have $\Gamma(x + 1/2)^2 \leq \Gamma(x) \Gamma(x+1) = x \Gamma(x)^2$ for $x > 0$. Applying this with $x = d/2$, we get $R_d \leq \sqrt{d/2}$. Using the previous inequality with $d+1$ instead of $d$, we have:
\[
R_d = \frac{d/2}{R_{d+1}} \geq \frac{d/2}{\sqrt{(d+1)/2}} > \sqrt{\frac{d}{\pi}} \cdot
\]
Thus we get:
\[
\kappa_{\mathrm{TV}}^{(d)}(z \mathds{1}) - \kappa_{\mathrm{TV}}^{(1)}(z d) = \sqrt{z} \frac{2}{\sqrt{\pi}} \left( R_d - \sqrt{\frac{d}{\pi}} \right) + O(z).
\]
We obtain that $\kappa_{\mathrm{TV}}^{(d)}(z \mathds{1}) > \kappa_{\mathrm{TV}}^{(1)}(z d)$ for $d \geq 2$ and $z$ small enough. 
Using the fact that $\kappa_{\mathrm{TV}}^{(d)}$ is symmetric and 
$\kappa_{\mathrm{TV}}^{(d)}(z \mathds{1})=\kappa_{\mathrm{TV}}^{(d)}(z d, 0, \ldots, 0)$ from 
the dimension reduction property (Proposition~\ref{prop:gaussian_gfd_properties}~\ref{enum:gaussian_gfd_properties_dim}), we get that $\kappa_{\mathrm{TV}}^{(d)}$ is not quasi-convex for $d \geq 2$.

To see it more clearly, let us consider the case $d=2$. 
Let $h$ be a non-negative quasi-convex symmetric function on $[0,1]^2$ such that $h(0,0) = 0$. Let  $z \in (0,1)$. By definition of quasi-convexity, we have $h(z, z) \leq \max(h(0, 2z), h(2z, 0)) = h(2z, 0)$. Thus, if we plot the level lines of $h$, in particular the one passing through $(z,z)$, we obtain Fig.~\ref{fig:non_convex_level_set} where the point $A$, defined as the intersection of the level line of $h$ through $(z,z)$ with the $x_1$-axe, 
is on  the left of  $(2z,0)$. 

\begin{figure}[htbp]
\centering
\begin{tikzpicture}[scale=3.5]
  \draw[->] (-0.1,0) -- (1.2,0) node[right] {$x_1$};
  \draw[->] (0,-0.1) -- (0,1.2) node[above] {$x_2$};
  \draw[dashdotted, gray] (0,0) -- (1,1);

  \fill (0.4,0.4) circle (0.01) node[right] {$(z,z)$};

  \draw (0.8,0.02) -- (0.8,-0.02);
  \fill (0.8,0) circle (0) node[below] {$2z$};

  \draw (0.02,0.8) -- (-0.02,0.8);
  \fill (0,0.8) circle (0) node[left] {$2z$};

  \draw[dashed, gray] (0.8,0) -- (0,0.8);

  \draw (0.4,0.4) arc (45:45-27.12:1);
  \draw (0.4,0.4) arc (45:45+27.12:1);

  \fill (0.6446,0) circle (0.01) node[below] {$A$};
  % \fill (0,0.6446) circle (0.01) node[left] {$B$};

  \fill (0,0) circle (0.01) node[below left] {$(0,0)$};

\end{tikzpicture}
\caption{An example of a level set of a non-negative quasi-convex symmetric function $h$ defined on $[0,1]^2$ and such that $h(0,0) = 0$.}
\label{fig:non_convex_level_set}
\end{figure}
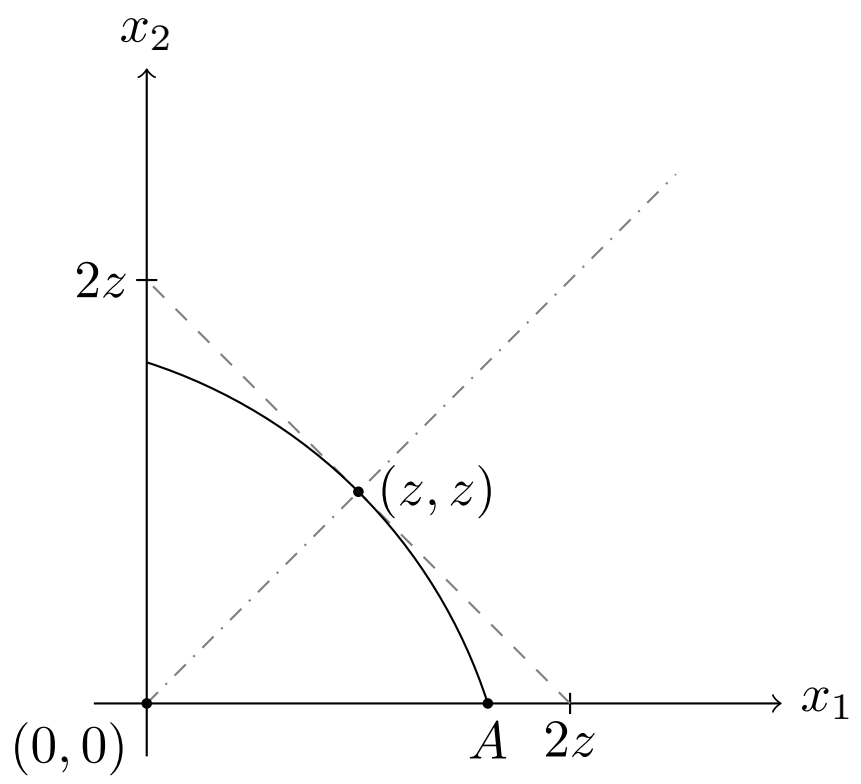

\end{proof}

We now give a proof of Lemma~\ref{lem:dev_lim_kappatv}. 

\begin{proof}[Proof of Lemma~\ref{lem:dev_lim_kappatv}]
Let $d\geq 1$ and $x \in (0,1)^d$. Let $Z_1, \dots, Z_d, Z_1', \dots, Z_d'$ be independent standard Gaussian random variables. Let $\eta_1$ and $\eta_2$ be functions defined on $[0,1)$ by:
\[
\eta_1(t) = \frac{\sqrt{t}}{1-\sqrt{t}} \quad\text{and} \quad \eta_2(t) = \frac{\sqrt{t}}{1+\sqrt{t}} \cdot
\]
The functions $\eta_1$ and $\eta_2$ are increasing and $\eta_2 \leq \eta_1$ on $[0,1)$. 
Recall that the matrix $\Gamma_x$ is defined in~\eqref{eq:def-Gamma-x} and we denote by $p(.,\sigma^2)$ the density of $X$ where $X \sim \cn(0,\sigma^2)$. 
For $x \in [0,1)^d$, the function $\kappa_{\mathrm{TV}}^{(d)}$ is given by:
\begin{align*}
\kappa_{\mathrm{TV}}^{(d)}(x) 
&= D_\mathrm{TV} \left( \cn(0,\Gamma_x) \Vert \cn(0,\idmatrix_{2d})\right) \\
&= \E \left[ \left\vert \prod_{i=1}^d \frac{p(Z_i, 1-\sqrt{x_i}) p(Z_i', 1+\sqrt{x_i})}{p(Z_i, 1) p(Z_i', 1)} - 1 \right\vert \right] \\
&= \E \left[ \left\vert \prod_{i=1}^d \frac{1}{\sqrt{(1-\sqrt{x_i})(1+\sqrt{x_i})}} \exp\left( - \frac{Z_i^2}{2 (1 - \sqrt{x_i})} + \frac{Z_i^2}{2} - \frac{Z_i'^2}{2 (1 + \sqrt{x_i})} + \frac{Z_i'^2}{2} \right) - 1 \right\vert \right] \\
% &= \E \left[ \left\vert \prod_{i=1}^d \inv{\sqrt{1-x_i}} \exp\left( \frac{1}{2} \sum_{i=1}^d \eta_2(x_i) Z_i'^2 - \eta_1(x_i) Z_i^2 \right) -1\right\vert \right].
&= \E \left[ \left\vert g(x) \exp\left( S_x \right) -1\right\vert \right],
\end{align*}
with $g(x) = \prod_{i=1}^d (1-x_i)^{-1/2}$ and $S_x = \frac{1}{2} \sum_{i=1}^d \left(\eta_2(x_i) Z_i'^2 - \eta_1(x_i) Z_i^2\right)$.
We consider the error term:
\[
R = g(x) \exp\left(S_x \right) -1 - \sum_{i=1}^d \frac{\sqrt{x_i}}{2} (Z_i'^2 - Z_i^2)
\]
and claim it is small in $L^1$:
\begin{equation}
    \label{eq:error_term_gaussian}
    \E \left[ \vert R \vert \right] = O\left( \sum_{i=1}^d x_i\right).
\end{equation}
From this claim, we deduce that~\eqref{eq:dev_lim_kappatv} holds.
We also prove~\eqref{eq:dev_lim_kappatv_vec1}. Let $z \in [0,1]$. We have:
\[
\kappa_\mathrm{TV}^{(d)}(z\mathds{1}) = \frac{\sqrt{z}}{2} \E \left[ \left\vert \sum_{i=1}^d \left(Z_i^2 - Z_i'^2 \right) \right\vert \right] + O(z) = \frac{\sqrt{z}}{2} \E \left[ \left\vert U - V \right\vert \right] + O(z),
\]
where $U = \sum_{i=1}^d Z_i^2$ and $V = \sum_{i=1}^d Z_i'^2$ follow a $\chi^2(d)$ distributions.
Let $U' = U + V$ and $V' = U/U'$. It is known that $U'$ and $V'$ are independent, $U'$ follows a $\chi^2(2d)$ distribution, and $V'$ follows a Beta distribution with parameters $(d/2,d/2)$. We have $U - V = U' (2V' - 1)$ and thus:
\[
\E \left[ \left\vert U - V \right\vert \right] = \E \left[ \left\vert U' \right\vert \right] \E \left[ \left\vert 2V' - 1 \right\vert \right] = 2d \E \left[ \left\vert 2V' - 1 \right\vert \right].
\]
We denote by $p_{V'}$ the density of $V'$. By making a change of variables and then using the symmetry of $p_{V'}$, we obtain:
\[
\E \left[ \left\vert 2V' - 1 \right\vert \right] 
= \int_{0}^{1} \vert 2v - 1 \vert p_{V'}(v) \rd v 
= \frac{1}{2} \int_{-1}^{1} \vert w \vert \, p_{V'} \left( \frac{w+1}{2} \right) \rd w 
= \int_{0}^{1} w \, p_{V'} \left( \frac{w+1}{2} \right) \rd w.
\]
Using the expression of $p_{V'}$, we have:
\[
\E \left[ \left\vert 2V' - 1 \right\vert \right] 
= \frac{\Gamma(d)}{\Gamma(d/2)^2} \int_{0}^{1} w \left( \frac{1-w^2}{4}\right)^{d/2-1} \rd w 
= \frac{\Gamma(d)}{\Gamma(d/2)^2}\, \inv{d 2^{d - 2} } \cdot
\]
Using that $\Gamma(d) = \Gamma(d/2) \Gamma((d+1)/2)\,  2^{d-1}/\sqrt{\pi}$, we get
$
\E \left[ \left\vert U - V \right\vert \right] = \frac{4}{\sqrt{\pi}} \frac{\Gamma((d+1)/2)}{\Gamma(d/2)},
$
and thus~\eqref{eq:dev_lim_kappatv_vec1} holds. 

\medskip

We now check the claim, that is,  the error term satisfies~\eqref{eq:error_term_gaussian}. Firstly, we have $\E \left[ \vert R \vert \right] \leq A + B$ with:
\[
A = \E\left[ \left\lvert g(x) \left(  1 + S_x \right) - 1 - \sum_{i=1}^d \frac{\sqrt{x_i}}{2} (Z_i'^2 - Z_i^2) \right\rvert \right] 
\quad \text{and} \quad
B =  g(x) \E\left[ r\left(S_x \right) \right],
\]
where $r(t) = \expp{t} - 1 - t$ for $t \in \R$.
We first focus on $A$. We have:
\begin{align*}
A &= \E\left[ \left\vert g(x) - 1  +  \inv{2} g(x) \left( \sum_{i = 1}^d \left(- \eta_1(x_i) + \sqrt{x_i}\right) Z_i^2 + \left(\eta_2(x_i)-\sqrt{x_i} \right) Z_i'^2 \right) \right\vert  \right]\\
&\leq \left\vert g(x) - 1 \right\vert + \inv{2}g(x) \left( \sum_{i = 1}^d \left\vert - \eta_1(x_i) + \sqrt{x_i} \right\vert  + \left\vert \eta_2(x_i)-\sqrt{x_i}   \right\vert   \right) \\
&\leq \left\vert g(x) - 1 \right\vert + g(x) \left( \sum_{i = 1}^d \max_{j \in \{1,2\}} \left\vert \eta_j(x_i) - \sqrt{x_i} \right\vert \right).
\end{align*}
Since $g(x) - 1 = O\left(\sum_{i=1}^d x_i\right)$ and $\left\vert \eta_j(x_i) - \sqrt{x_i} \right\vert = O(x_i)$ for all $(i,j) \in \{1, \dots,d \} \times \{1,2 \}$,
we deduce that $A = O\left(\sum_{i=1}^d x_i \right)$.
We now focus on the second term $B$. 
We have $B = g(x) \E[ r(S_x) ]$. Since $g(x)$ goes to 1 as $\max(x) \to 0$, it is enough to show that $\E[r(S_x)] = O\left(\sum_{i=1}^d x_i\right)$. We have $\E[r(S_x)] = \E \left[ \expp{S_x} \right] - 1 - \E[S_x]$. We first compute $\E \left[ \expp{S_x} \right]$ using the moment generating function of a $\chi^2(1)$ distribution. We have:
\[
\E \left[\expp{S_x}\right] 
= \prod_{i=1}^d \E \left[\expp{\frac{1}{2} \eta_2(x_i) Z_i'^2} \right] \E \left[\expp{-\frac{1}{2} \eta_1(x_i) Z_i^2 } \right] = \prod_{i=1}^d \frac{1}{\sqrt{1-\eta_2(x_i)}} \frac{1}{\sqrt{1+\eta_1(x_i)}} = \prod_{i=1}^d \sqrt{1-x_i}.
\]
Next, we compute $\E[S_x]$:
\[
\E[S_x] = \sum_{i=1}^d \frac{1}{2} \left( \eta_2(x_i) - \eta_1(x_i) \right) = \sum_{i=1}^{d} \frac{x_i}{x_i-1} \cdot
\]
Thus, we have:
\[
  \E\left[ r(S_x) \right] = \prod_{i=1}^d \sqrt{1-x_i} - 1 + \sum_{i=1}^{d} \frac{x_i}{1-x_i} = O\left(\sum_{i=1}^d x_i\right).
\]
This proves the claim~\eqref{eq:error_term_gaussian}. 
\end{proof}

\section{Proofs of the results in Section~\ref{sec:Bernoulli}}

When there is no ambiguity on the parameters $q$ and $A'$, we shall write $C^\bullet$ for $C^\bullet_{q,A'}$. 
We recall that for $p,q \in (0,1)$, the density ratio $Z$ (defined in~\eqref{eq:Sf_et_Z}) takes values in the set
$\{a_i : i=1,\ldots,4\}$, with:
\begin{equation}
  \label{eq:def-a}
  a_1=\frac{r}{\rho_1}, \quad a_2=\frac{1-p-q+r}{\rho_2}, \quad
  a_3=\frac{q-r}{\rho_3} \quad\text{and}\quad a_4=\frac{p-r}{\rho_4},
\end{equation}
where the $\rho_i$ are defined as:
\begin{equation}
    \label{eq:def-rho}
\rho_1=pq, \quad \rho_2=(1-p)(1-q), \quad \rho_3=(1-p)q \quad\text{and}\quad \rho_4=p(1-q).
\end{equation}
We also define two sets:
\[
D_q^- = \left\{ A = (p,r) \in D_q \, \colon \, r \leq pq \right\} 
\quad \text{and} \quad
D_q^+ = \left\{ A = (p,r) \in D_q \, \colon \, r \geq pq \right\}.
\]
Note that we have $D_q^- \cup D_q^+ = D_q$ and $D_q^- \cap D_q^+ = \Delta_q$.
We also recall an equivalence we will use several times. Since we have $S_f(X, Y)=\E[f(Z)]$, we have that $X$ and $Y$ are independent if and only if $Z=1$ a.s. with respect to the product measure:
\begin{equation}
\label{eq:Sf=0iff}
S_f(X, Y)=0 \quad\text{for all $f\in \F$} \quad \iff \quad Z =1 \quad \text{ $P_X\otimes P_Y$-a.s.}  %\quad  \iff A\in \Delta_q. 
\end{equation}

\subsection{Proof of Proposition~\ref{prop:bernoulli_cp}}
\label{sec:proof-bernoulli_cp}

The technical proof is divided as follows: in Section~\ref{sec:deg-Dq} we consider the degenerate case $A'\in \Delta_q$; in Section~\ref{sec:prelim-proof} we give an elementary characterization of $C^-$; we then consider the cases $q\in [1/2, 1)$, $(p', r')\in D_q$ with
$r'<p'q$, and $p'\in (0, 1-q)$
in Section~\ref{sec:proof_case1}, $p'\in [1-q, q]$ in  Section~\ref{sec:proof_case2}, and 
$p'\in(q, 1)$ in  Section~\ref{sec:proof_case3}; eventually using symmetries we end the proof in Section~\ref{sec:proof_bernoulli_sym_expression_finale}.

\subsubsection{A preliminary result} 
\label{sec:prelim-proof}
Let $q\in [0, 1]$, $(p,r), (p', r')\in D_q$. 
Recall the density ratio $Z=(\rd P_{(X,Y)}/\rd P_X\otimes \rd P_Y)(X,Y)$ (under $P_X\otimes P_Y$) and $Z'$ is defined similarly with $X$ replaced by $X'$.  
We define the function $\Gamma$ on $\R$ by:
\[
\Gamma(t)=\E_{P_{X'}\otimes P_Y}[(t-Z')_+] - \E_{P_X\otimes P_Y}[(t-Z)_+].
\]
As $\E_{P_X\otimes P_Y}[Z]=\E_{P_{X'}\otimes P_Y}[Z']=1$, we also have:
\begin{equation}
    \label{eq:Gamm=E-}
\Gamma(t) = \E_{P_{X'}\otimes P_Y}[(Z'-t)_+] - \E_{P_X\otimes P_Y}[(Z-t)_+].
\end{equation}
According to~\cite[Theorem 3.A.1.]{shaked2007stochastic} on the characterization of the convex order (which is easily extended from finite convex function to proper convex functions),  since $Z$ and $Z'$ are integrable with mean 1, we have
$  Z'\, \geqcv \,  Z$ (that is, 
$(p,r)\in C^+$) if and only if $\Gamma$ is nonnegative. In the particular case of this section, we can simplify this criterion.

\begin{lem}[An elementary criterion for the convex order]
\label{lem:equivalence_convex_order_int_cdf}
 Let $q\in (0, 1)$ and $(p,r), (p', r')\in D_q$ with  $p',p\in (0, 1)$.   We have $(p,r)\in C^-$ (resp. $(p,r) \in C^+$) if and only if $\Gamma(a_i')\geq 0$  (resp. $\Gamma(a_i)\leq 0$) for all $i=1, \ldots, 4$.
\end{lem}

We shall use the following notations in the proof of this lemma and later on. 
Let  $a_{(1)}\leq a_{(2)}\leq a_{(3)}\leq a_{(4)}$ be the ordered sequence of $a_1, a_2, a_3, a_4$, defined in~\eqref{eq:def-a}, and, following the definition of  $\rho_i$  in~\eqref{eq:def-rho}, we set $\rho_{(i)} = \P \left(X = a_{(i)}\right)$ for  $i = 1,\dots,4$. We define similarly $a'_{(i)}$ and $ \rho'_{(i)}$ with $X$ replaced by  $X'$.

\begin{proof}[Proof of Lemma~\ref{lem:equivalence_convex_order_int_cdf}]
We consider the function $\psi(t)=
\E[(t-Z)_+]$ defined on $\R$ and similarly for $\psi'$ with $Z$ replaced by $Z'$. The functions $\psi$ and $\psi'$ are convex and  $\Gamma = \psi' - \psi$.
Notice that  $\E[Z] = \E[Z']=1$.
According to~\cite[Theorem 3.A.1.]{shaked2007stochastic}, $  Z'\, \geqcv \,  Z$ (that is, $(p,r)$ belongs to $ C^-$) if and only if the function $\Gamma$ is nonnegative. Notice that $\Gamma=0$ outside the interval $(a_{(1)} \wedge a'_{(1)}, a_{(4)} \vee a'_{(4)})$. Since the function $\Gamma$ is piecewise affine, with possible changes of slope occurring at the points $a_i$ and $a_i'$, we get that $\Gamma$ is nonnegative if and only if  $\Gamma$ is nonnegative at $a_i$ and $a'_i$ for all $i=1, \ldots, 4$. 
In particular, if  $Z'\, \geqcv \, Z$, then one has $\Gamma(a_i') \geq 0$ for all $i=1, \ldots, 4$. 

We shall now prove the converse. First, notice that:
\begin{equation}
  \label{eq:gamma1_gamma4}
\begin{aligned}
    &\Gamma(a_{(1)}') = 0 \iff  \Gamma(a_{(1)}') \geq 0 \iff a_{(1)}' \leq a_{(1)}, \\
     &\Gamma(a_{(4)}') = 0 \iff   \Gamma(a_{(4)}') \geq 0 \iff a_{(4)} \leq a'_{(4)},
\end{aligned}  
\end{equation} 
where for the second line (the proof of the first line is similar), we used that 
 $\E[(a_{(4)}'-Z')_-] = 0$, and  if $\Gamma(a_{(4)}') \geq 0$, then  we get $-\E[(a_{(4)}'-Z)_-] \geq 0$ and thus  $\E[(a_{(4)}'-Z)_-] = 0$, which implies that  $a_{(4)}' \geq a_{(4)}$.

We now assume that $\Gamma(a_i') \geq 0$ for all $i=1, \ldots, 4$.  To get $  Z'\, \geqcv \,  Z$, it is enough to prove that $\Gamma(a_i)\geq 0$ for all $i=1, \ldots, 4$. Fix $i \in \{1,\dots,4\}$.
We deduce from~\eqref{eq:gamma1_gamma4}
%Since  $\Gamma(a_{(1)}') \geq 0$ and $\Gamma(a_{(4)}') \geq 0$, we deduce 
that  $a_{(1)}' \leq a_i \leq a'_{(4)}$. 
Thus there exists $j \in \{1,2,3\}$ such that $a_i$ belongs to $I_j= [a_{(j)}', a_{(j+1)}']$, that is, there exists $\alpha \in [0,1]$ such that $a_i = \alpha a_{(j)}' + (1-\alpha) a_{(j+1)}'$. Since the function $\psi'$ is affine on $I_j$, we have $\psi'(a_i) = \alpha \, \psi'(a_{(j)}') + (1-\alpha) \, \psi'(a_{(j+1)}')$. Using the convexity of $\psi$ and the fact that $\Gamma$ is nonnegative at $a_{(j)}'$ and $a_{(j+1)}'$, we get:
\[
\psi(a_i) \leq \alpha \psi(a_{(j)}') + (1-\alpha) \psi(a_{(j+1)}') \leq \alpha \psi'(a_{(j)}') + (1-\alpha) \psi'(a_{(j+1)}') = \psi'(a_i).
\]
We thus obtain that $\Gamma(a_i) \geq 0$. Since  $i \in \{1,\dots,4\}$ is arbitrary, we deduce the converse holds. 

The proof for $C^+$ is similar. 
\end{proof}

We end this section with a very usefull remark on the order statistics of $\{a_1, a_2, a_3, a_4\}$. 

\begin{rem}[On the $a_{(i)}$]
    \label{rem:stat-ordre}
Assume $q\in [1/2, 1)$, $p\in (0,1)$. Recall the $a_i$ defined in~\eqref{eq:def-a}. For any $i\in \{1, \ldots, 4\}$, we have $a_i=1$ if and only if $r=pq$ (and thus all the $a_j$ are equal to 1). For $r<qp$, we  obtain that:
\[
\{a_{(1)},a_{(2)}\}= \{a_1, a_2 \}, \quad
\{a_{(3)}, a_{(4)}\}= \{a_3, a_4 \} \quad\text{and}\quad 
a_{(2)} <1 < a_{(3)},
\]
as well as the following:
\begin{enumerate}[(i)]
    \item  $a_1=a_2$ if  and only if $p=1-q$;
    \item $a_3=a_4$ if and only if $p=q$;
   \item if $p\in (0, 1-q]$, then  $a_{(1)}=a_1$ and $a_{(4)}=a_4$;
    \item if $p\in [1-q, q]$, then   $a_{(1)}=a_2$ and $a_{(4)}=a_4$;
    \item if $p\in [q, 1)$, then   $a_{(1)}=a_2$ and $a_{(4)}=a_3$.
     \end{enumerate}
\end{rem}

\subsubsection{The case $(p',r')$ or  $(p,r)$ in $\Delta_q$}
\label{sec:deg-Dq}
Notice that  $(p, r)\in \Delta_q$ if and only if  $X$ and $Y$ are independent, that is, if and only if $S_f(X,Y)=0$ for one $f\in \F$ or equivalently for all $f\in \F$. 
We deduce that $\Delta_q\subset C^-_{q, A'}$ for all $A'\in D_q$.

Let $A'=(p',r')\in \Delta_q$. In particular, we get from the above comment  that $C^-_{q, A'}=\Delta_q$. Since $S_f(X,Y) \geq S_f(X', Y)=0$ for all  Bernoulli random variable $X$ and for all $f\in \F$; we also deduce that  $C^+_{q, A'}=D_q$. 

\medskip 

Notice also that $q\in \{0, 1\}$ if and only if $D_q=\Delta_q$. This degenerated case corresponds to  $Y$ constant and thus $S_f(X',Y)=S_f(X,Y)=0$ for all Bernoulli random variables $X$ and $X'$ and for all $f\in \F$; and then $C^+_{q, A'}=C^-_{q, A'}=\Delta_q=D_q$.

\medskip

Notice that if $A=(p,r) \in D_q$ with  $p\in \{0, 1\}$ implies that  $A\in \Delta_q \subset C^-_{q, A'}$ for any $q\in [0, 1]$ and $A'\in D_q$. 

\medskip

In conclusion, we shall now consider the case  $q\in (0,1)$, that is $\Delta_q \neq D_q$, and  $A'=(p', r')$ and $A=(p,r ) $ in $D_q \setminus  \Delta_q$, that is $r'\neq qp'$ and $r\neq qp$, which implies also that $p'$ and $p$ belong to $(0, 1)$.

\subsubsection{The case $q\in [1/2, 1)$, $(p', r')\in D_q$, $p'\in (0, 1-q)$ and $r'<p'q$}
\label{sec:proof_case1}
This is  part of case~\ref{it:3-D1} in Proposition~\ref{prop:bernoulli_cp} as the latter condition is equivalent to $A'\not \in \Delta_q$. The proof  is divided in five parts: after some preliminary remarks, we consider the cases $p$ belonging to $(0, p']$, $[p', \tilde p '']$ and $[\tilde p'', 1)$ for which different binding constraints apply, and then give the conclusion.

By symmetry, it is enough to consider the case $r< pq$ (recall the equality is studied in Section~\ref{sec:deg-Dq}).

\noindent\textbf{a) Preliminaries.}
Notice that $2(1-q)-p'$ is positive. 
As the point $A''$ lies above $\Delta_q$, we shall consider $\tilde A'' = h_q(A'')$, with coordinates $(\tilde p'', \tilde r'')$ given by:
\begin{equation}
    \label{eq:def-tp''}
\tilde p'' = 1-\frac{p'(1-p')}{2(1-q)-p'} \quad \text{and} \quad
\tilde r'' = \tilde p'' q - (p'q- r') \frac{(1-p')(1-q)}{q(2(1-q)-p')}\cdot
\end{equation}
Since $0<p'< 1-q$ and $q\geq 1/2$, we easily deduce that  $\tilde p'' \in (1-q,1)$.
Since  $r' < p'q$, we also have $\tilde r'' < \tilde p'' q$. It is also easy to check that $0<\tilde p'' + q-1 \leq \tilde r''$. We deduce that $\tilde A''\in D_q$ (and 
also that $A''=h_q(\tilde A'')\in D_q$). 

Let $\Delta_1$ be the line passing through $A_0$ and $A'$, $\Delta_2$ the line passing through $A'$ and $\tilde A''$, and $\Delta_3$ the line passing through $\tilde A''$ and $h_q(A_0)=(1,q)$. Let $s_i$ be the slope of $\Delta_i$ for $i=1,2,3$. We have:
\[
s_1= \frac{r'}{p'}, \quad s_2=\frac{\tilde r''-r'}{\tilde p'' - p'}=\frac{r'}{p'}+ (p'q- r')  \frac{2q -p'}{2q p'(1-p')}  \quad \text{and} \quad s_3= \frac{q-\tilde r''}{1-\tilde p''}=\frac{r'}{p'}+ (p'q- r')\frac{1}{p'q}\cdot
\]
Since  $p' < 1 - q$ and $r'< p'q$, we get that  $s_1 < s_2 < s_3$. 
Thus the points $A_0$, $A'$, $\tilde A''$, and $h_q(A_0)$ are in convex position. Using the symmetry, we deduce that the set $\ca=\ca_0\cup h_q(\ca_0)$, with $\ca_0=\{A_0, A', A''\}\subset D_q$, consists in 6 distinct points and that they are in convex position. 

\medskip

Let $r=\Delta_i(p)$ be the equation of the line $\Delta_i$, that is, $\Delta_i(p)=p_i+ s_i p$, with $p_1=0$, $p_2=r' - s_2 p'$ and $p_3=q-s_3$. 
As $s_1<s_2<s_3$ and since $\Delta_1$ and $\Delta_2$ intersect at $(p',r')$ and $\Delta_2$ and $\Delta_3$ intersect at $(\tilde p'',\tilde r'')$, we deduce that:
\begin{equation}
  \label{eq:gamma_inequalities_case1}
\begin{aligned}
     \Delta_1(p) > \Delta_2(p) \iff p < p'; \quad \text{and} \quad \Delta_2(p) > \Delta_3(p) \iff p < \tilde p'',\\
      \Delta_1(p) =\Delta_2(p) \iff p = p'; \quad \text{and} \quad \Delta_2(p) = \Delta_3(p) \iff p = \tilde p''.
\end{aligned} 
\end{equation}

\medskip
We now derive some properties related to the atoms $(a'_i)_{i\in \{1, \ldots, 4\}}$ of $Z'$. Remark~\ref{rem:stat-ordre} implies that:
\begin{equation}
    \label{eq:case1-order-a'}
    0\leq a'_1= a'_{(1)}<a'_2 = a'_{(2)} < 1 < a'_3= a'_{(3)} < a'_4= a'_{(4)} \leq \frac{1}{1-q}\cdot
\end{equation}
We shall use many time that:
\begin{align}
   \label{eq:s3}
 &s_3= 1 - \frac{1-q}{q}\frac{r'}{p'} ,\\
    \label{eq:G(a'3)}
&(a'_4 -a'_3) \rho'_4 - (a_4 -a'_3) \rho_4 - (a_3 -a'_3) \rho_3
= 2(r-\Delta_2(p)). 
 \end{align}
Using this representation of $s_3$ for~\eqref{eq:case1_cns_delta3}, we get:
\begin{align}
  \label{eq:case1_cns_delta1} r \geq \Delta_1(p) &\iff a_1 \geq a_1' \iff a_4 \leq a_4', \\
  \label{eq:case1_cns_delta3} r \geq \Delta_3(p) &\iff a_2 \geq a_1'.
\end{align}

\medskip

\noindent 
\textbf{b) The case  $p \in (0,p']$.} 
 As $0<p\leq p'< 1-q\leq 1/2$ and $r'<qp'$, we get from Remark~\ref{rem:stat-ordre} that $a_1=a_{(1)}< a_2<1$ and $1<a_3<a_4=a_{(4)}$. Recall and set:
 \[
 \Delta_1(p)=\frac{r'}{p'}\, p
 \quad\text{and}\quad
 \cstb= \frac{(p'-p)(qp'-r')}{p'(1-p')}\geq 0.
 \]
 We first consider the binding constraint. We deduce from~\eqref{eq:gamma1_gamma4},~\eqref{eq:case1-order-a'} and~\eqref{eq:case1_cns_delta1} that:
\begin{equation}
  \label{eq:gamma1_subcase1i}
  \Gamma (a_{(1)}') \geq 0 \iff a_{(1)}' \leq a_{(1)} \iff a_1' \leq a_1 \iff r \geq \Delta_1(p)
  \iff a_4 \leq a'_4 \iff \Gamma (a_{(4)}') \geq 0.
\end{equation}
We shall now assume that $r\geq \Delta_1(p)$ and prove that $\Gamma(a'_{i})\geq 0$ for $i=2,3$. First notice that:
\begin{equation}
    \label{eq:subcase1i_eqa2a3}
    \rho_2 (a_2 -a'_2) = \rho_3 (a_3' - a_3) = r - \Delta_1(p) + \cstb \geq 0.
\end{equation}
For $i=2$, with~\eqref{eq:subcase1i_eqa2a3} notice that we have $a_2 \geq a_2'$.
So we have either  $a'_{2}\in (a'_1, a_1) $ or   $a'_2\in [a_1, a_2]$. 
We have: for $a'_{2}\in [a'_1, a_1] $ that $\Gamma(a'_2)= \E[(a'_2 - Z')_+]\geq 0$; and 
for $a'_2\in [a_1, a_2]$ that:
\[
\Gamma(a'_2)= (a'_2 -a'_1) \rho'_1 - (a'_2 - a_1) \rho_1
= r - \Delta_1(p)  +\frac{1-q -p'}{1-q} \, \cstb \geq 0. 
\]
For $i=3$, with~\eqref{eq:subcase1i_eqa2a3} notice that we have $a_3 \leq a_3'$. So we have either $a'_3\in [a_3, a_4]$ or $a'_3\in (a_4, a'_4)$. We shall use the representation of  $\Gamma$ given in~\eqref{eq:Gamm=E-}. We have: for $a'_3\in (a_4, a'_4)$ that $\Gamma(a'_3)= \E[(Z'-a'_3)_+]\geq 0$; 
for $a'_3\in [a_3, a_4]$ that:
\begin{equation}
    \label{eq:cas-b-a3}
    \Gamma(a'_3)= (a'_4 -a'_3) \rho'_4 - (a_4 -a'_3) \rho_4= r-\Delta_1(p) + \frac{q-p'}{q}\, \cstb\geq 0.
\end{equation}
So, as $\cstb\geq 0$,  we have obtained that: 
\begin{equation}
    \label{eq:case1i-conclusion}
    \Gamma(a'_1)\geq 0 \iff \Gamma (a'_4) \geq 0 \iff r\geq \Delta_1(p) \implies \Gamma(a'_i)\geq 0 \quad \text{for $i=2,3$}.
\end{equation}

\medskip

\noindent
\textbf{c) The case  $p \in  [p',\tilde p'']$.} 
Recall $\tilde p''$ in~\eqref{eq:def-tp''} belongs to $(1 - q, 1)$. 
As $0< p'< 1-q$, we get that $p'<1-q<\tilde p''$. 
 As  $r<qp$,  we get from Remark~\ref{rem:stat-ordre} that $\{a_3, a_4\}=\{a_{(3)}, a_{(4)}\}$. Recall:
\[
\Delta_2(p)= r'+ \left( \frac{r'}{p'} +\frac{2q-p'}{p'}\, \cstc\right) \, (p-p') 
\quad\text{with}\quad
\cstc=\frac{qp'-r'}{2q(1-p')}\cdot
\]
Elementary computations give:
\begin{align}
\label{eq:case1i-3}
    \rho_3(a'_3-a_3)&=r - \Delta_2(p)  - \cstc' 
    \quad\text{with}\quad \cstc'=\cstc\, (p-p')\geq 0, \\
 \label{eq:case1i-4}
   \rho_4 (a'_3 - a_4) &=r - \Delta_2(p)
 - \cstc''
 \quad\text{with}\quad \cstc''=\cstc\, \left( (\tilde p'' -p) + 2\frac{1- \tilde p''}{p'(1-p')} (1-q)(2q-1)\right)
 \geq 0.
\end{align}

We  first check the binding constraint:
\begin{equation}
  \label{eq:gamma3_subcase1ii}
  \Gamma(a_{3}') \geq 0 \iff  r \geq \Delta_2(p).
\end{equation}
To do so we distinguish three  cases $a'_3\in (1, a_{(3)})$, $a'_3\in [a_{(3)}, a_{(4)}]$ and $a'_3\in (a_{(4)}, a'_4)$.
In the latter case, we have $\Gamma(a'_3)=\E[(Z'-a'_3)_+] \geq 0$ and 
we deduce from~\eqref{eq:case1i-3} and~\eqref{eq:case1i-4} (depending on  $a_{(4)}$ being equal to $a_3$ or $a_4$) that $r\geq \Delta_2(p)$; and thus~\eqref{eq:gamma3_subcase1ii} holds.
For the first case $a'_3\in (1, a_{(3)})$, we have using~\eqref{eq:G(a'3)}:
\[
\Gamma(a_{3}')=(a_{4}' -a_{3}') \rho_{4}' - (a_{4} - a_{3}')\rho_{4} - (a_{3} - a_{3}') \rho_{3}=2(r -\Delta_2(p));
\]
and thus~\eqref{eq:gamma3_subcase1ii} holds also. 
For the case $a'_3\in [a_{(3)}, a_{(4)}]$, we have: 
for $p\leq q$ that $a_{(4)}=a_4$ by Remark~\ref{rem:stat-ordre} and thus $a'_3\in [a_3, a_4]$ which thanks to~\eqref{eq:case1i-3} implies $r\geq \Delta_2(p)$ as well as (subtract~\eqref{eq:case1i-3} to~\eqref{eq:G(a'3)}):
\[
\Gamma(a'_3)=(a_4' - a_3') \rho_4' - (a_4 - a_3') \rho_4 = r -\Delta_2(p) +  \cstc'  \geq 0;
\]
and for  $p\geq q$
that $a_{(4)}=a_3$  and thus $a'_3\in [a_4, a_3]$ which thanks to~\eqref{eq:case1i-4} implies also $r\geq \Delta_2(p)$ as well as (subtract~\eqref{eq:case1i-4} to~\eqref{eq:G(a'3)}):
\[
\Gamma(a'_3)=(a_4' - a_3') \rho_4' - (a_3 - a_3') \rho_3 = r -\Delta_2(p) + \cstc'' \geq 0;
\]
and thus~\eqref{eq:gamma3_subcase1ii} holds also. 
In conclusion Equation~\eqref{eq:gamma3_subcase1ii}  holds. 

\medskip

We  assume that $r\geq \Delta_2(p)$ and shall now prove that 
 $\Gamma(a'_{i})\geq 0$ for $i=1,2,4$. 
 To prove that  $\Gamma(a'_1) \geq 0$, it is enough by~\eqref{eq:gamma1_gamma4}, to prove that  $a_1' \leq a_{(1)}$ (as $a'_1=a'_{(1)}$).
We have:
\[
r\geq \Delta_2(p) \implies 
\begin{cases}
   r \geq \Delta_1(p)
\iff a_{1} \geq a'_{1}
\iff  a_{(1)} \geq a'_1 & \quad \text{for} \quad p\in [p', 1-q],\\
  r \geq \Delta_3(p)
\iff a_{2} \geq a'_{1}
\iff  a_{(1)} \geq a'_1 & \quad \text{for} \quad p\in [1-q, \tilde p''],
\end{cases}
\]
where we used~\eqref{eq:gamma_inequalities_case1} for the implication, 
\eqref{eq:case1_cns_delta1} and~\eqref{eq:case1_cns_delta3} for the first implications and the value of $a_{(1)}$ given in 
Remark~\ref{rem:stat-ordre} for the second implications.
In conclusion we get $\Gamma(a'_1)\geq 0$.

Similarly, to prove that  $\Gamma(a'_4) \geq 0$, it is enough by~\eqref{eq:gamma1_gamma4}, to prove that  $a_{(4)} \leq a'_4$ (as $a'_4=a'_{(4)}$). First notice that:
\begin{equation}
    \label{eq:a'4-a3}
    \rho_3(  a_4' - a_3) = r- \Delta_3(p)+ (1-p)\frac{(qp'-r')(2q-1)}{p'q(1-q)}\cdot 
\end{equation}
We have:
\[
r\geq \Delta_2(p) \implies 
\begin{cases}
   r \geq \Delta_1(p)
\iff a_{4} \leq a'_{4}
\iff  a_{(4)} \leq a'_4 & \quad \text{for} \quad p\in [p', q],\\
  r \geq \Delta_3(p)
\implies a_{3} \leq a'_{4}
\iff  a_{(4)} \leq a'_4 & \quad \text{for} \quad p\in [q, \tilde p''],
\end{cases}
\]
where we used~\eqref{eq:gamma_inequalities_case1} for the first implication, then
\eqref{eq:case1_cns_delta1} and~\eqref{eq:a'4-a3} and last the value of $a_{(4)}$ given in 
Remark~\ref{rem:stat-ordre}. 
In conclusion we get $\Gamma(a'_4)\geq 0$.

For $i = 2$, notice that etiher $a_2' \in (a_1',a_{(1)})$ or $a_2' \in [a_{(1)},a_{(2)}]$ or $a_2' \in (a_{(2)},1)$. We have: for $a_2' \in (a_1',a_{(1)})$ that $\Gamma(a_2')= \E[(a_2' - Z')_+]\geq 0$; for $a_2' \in [a_{(1)},a_{(2)}]$ that:
\[
\Gamma(a_2') = (a_2' -a_1')\rho_1' -   (a_2' -a_{(1)})\rho_{(1)}
= (a_2' -a_1')\rho_1' -
\begin{cases}
  (a_2' -a_1)\rho_1 & \quad\text{if $p \leq 1-q$,}\\
  (a_2' -a_2)\rho_2 & \quad\text{if $p \geq 1-q$,} 
\end{cases} 
\]
that is:
\[
\Gamma(a_2') =  r - \Delta_2(p) + 
\frac{\cstc}{1-q} *\begin{cases}
  (3q-1)(p-p')  & \quad\text{if $p \leq 1-q$,}\\
   (1-q)(\tilde p'' - p) + 2\frac{1- \tilde p''}{1-p'}\, 
  \frac{(2q-1)(1-q-p')^2}{p'} & \quad \text{if $p \geq 1-q$}; 
\end{cases} 
\]
and for $a_2' \in (a_{(2)},1)$ that:
\[
\Gamma(a_2') = (a_2' -a_1')\rho_1' - (a_2' -a_1)\rho_1 -  (a_2' -a_2)\rho_2 = 2(r - \Delta_2(p)) +  2  \frac{2q-1}{1-q} \, \cstc'.
\]
In conclusion, as $r\geq \Delta_2(p)$, we get $\Gamma(a'_2)\geq 0$.
We have obtained that: 
\begin{equation}
    \label{eq:case1ii-conclusion}
 \Gamma(a'_3)\geq 0 \iff r\geq \Delta_2(p) \implies \Gamma(a'_i)\geq 0
 \quad \text{for $i=1,2,4$.}
\end{equation}

\medskip
\noindent
\textbf{d) The case  $p \in [\tilde p'',1)$.}
 Recall $\tilde p''$ in~\eqref{eq:def-tp''} belongs to $(1-q, 1)$. 
 Since $p\geq \tilde p''>1-q$ (and also $r<qp$),
 we get from Remark~\ref{rem:stat-ordre} that $a_2=a_{(1)}< a_1<1$. Recall (using $s_3$ from~\eqref{eq:s3}) and set:
 \[
 \Delta_3(p)=\left( 1- \frac{1-q}{q} \frac{r'}{p'} \right) \, (p -1) +q
 \quad\text{and}\quad
 \cstd^0= \frac{qp'-r'}{qp'(1-p')}\geq 0.
 \]
 We first consider the binding constraint. Using~\eqref{eq:s3}, we deduce from~\eqref{eq:gamma1_gamma4} and~\eqref{eq:case1_cns_delta3} that:
\[
\Gamma(a_{(1)}') \geq 0 \iff 
%\Gamma(a_{1}') \geq 0 \iff 
a'_{(1)} \leq a_{(1)}
\iff a'_1 \leq a_2
\iff
r \geq \Delta_3(p).
\]
We now assume that 
 $r\geq \Delta_3(p)$ and prove that  $\Gamma(a'_i) \geq 0$ for $i=2,3,4$. 
For $i = 4$, since $a'_{(4)}=a'_4$, we  get from~\eqref{eq:gamma1_gamma4}   that:
\[
\Gamma(a_4') \geq 0 \iff 
\Gamma(a_{(4)}') \geq 0 \iff 
a_{(4)}'\geq a_{(4)} \iff
a_4'- a_{(4)}\geq 0. 
\]
We have 
either $p \in [p',q]$ or $p \in [q,1)$. Thanks to 
Remark~\ref{rem:stat-ordre}, this gives 
$a_{(4)}=a_4$ in the former case and 
$a_{(4)}=a_3$ in the latter case.
Using~\eqref{eq:gamma_inequalities_case1} (with $p\geq \tilde p''$), elementary computation  gives:
\[
\rho_4(a'_4-a_4)= r - \Delta_1(p)  \geq r - \Delta_3(p)\geq 0
\quad\text{and}\quad 
\rho_3(  a_4' - a_3) = r- \Delta_3(p)+ (1-p)\frac{(qp'-r')(2q-1)}{p'q(1-q)}\geq 0,
\]
as $q\geq 1/2$. Thus, we have obtained that $r \geq \Delta_3(p) $ implies $ \Gamma(a'_4)\geq 0$.
For $i = 2$, notice that $\rho_1 (a_1 - a_2') = r - \Delta_3(p)+ \cstd$ 
with:
\[
    \cstd= \frac{\cstd^0}{1-q}
\left( (p- \tilde p'') \left( (1-q)^2 + (2q -1) p'\right) + \frac{1-\tilde p''}{1-p'} (2q-1) \left((1-p')^2 - q^2\right)\right)\geq 0. 
\]
So we have either $a_2' \in (a_1',a_2)$ or $a_2' \in [a_2,a_1]$. We have: for $a_2' \in (a_1',a_2)$ that $\Gamma(a_2')= \E[(a_2' - Z')_+]\geq 0$; and  for $a_2' \in [a_2,a_1]$ that:
\[
\Gamma(a_2') = (a_2' -a_1')\rho_1' - (a_2' -a_2)\rho_2 = r - \Delta_3(p) + \cstd '\geq 0,
\]
where:
\[
    \cstd' = (1-p'-q)\,\frac{ \cstd^0}{1-q} \left( (p- \tilde p'') (1-q) + \frac{1- \tilde p''}{1-p'} (2q-1) (1-q-p') \right)>0.
\]

For $i = 3$, notice that either $a'_3 \in (a_{(4)},a_4')$ or $a'_3 \in [a_{(3)},a_{(4)}]$ or $a'_3 \in (1,a_{(3)})$.  Using~\eqref{eq:Gamm=E-}, we have: for $a'_3 \in(a_{(4)},a_4')$ that $\Gamma(a'_3)= \E[(Z'-a'_3)_+]\geq 0$; for $a'_3 \in [a_{(3)},a_{(4)}]$ that:
\[
\Gamma(a'_3) =\begin{cases}
   (a_4' - a_3') \rho_4' - (a_4 - a_3') \rho_4 = r - \Delta_3(p) +\cstd'' 
   & \text{if $p \leq q$,} \\
  (a_4' - a_3') \rho_4' -  (a_3 - a_3') \rho_3
  = r - \Delta_3(p) +\cstd'''
  & \text{if $p \geq q$,}
\end{cases} 
\]
where:
\[
\begin{aligned}
\cstd''&=\cstd^0\, \left( (p-\tilde p'') (1-q) + (1-\tilde p'')(1-q-p') \right) \geq 0, 
\\
  \cstd'''&= \cstd^0\,\left( (p-\tilde p'') (1-q-p') + \frac{1-\tilde p''}{1-p'} (2q-1)(1-q) \right) \geq 0; 
\end{aligned}
\]
and for $a'_3 \in (1,a_{(3)})$  that
$\Gamma(a'_3) = (a_4' - a_3') \rho_4' - (a_4 - a_3') \rho_4 - (a_3 - a_3') \rho_3 $, that is,  with~\eqref{eq:G(a'3)} and~\eqref{eq:gamma_inequalities_case1} (with $p\geq \tilde p''$):
\[
\Gamma(a'_3)= 2( r- \Delta_2(p))\geq 2 (r - \Delta_3(p)).
\]
So we have obtained:
\begin{equation}
    \label{eq:case1iii-conclusion}
 \Gamma(a'_1)\geq 0 \iff r\geq \Delta_3(p) \implies \Gamma(a'_i)\geq 0
 \quad \text{for $i=2,3,4$.}
\end{equation}

\medskip
\noindent
\textbf{e) Conclusion.} 
Recall that  $q\in [1/2, 1)$, $p'\in (0, 1-q)$, $r'<p'q$. We have shown that the set $\ca=\ca_0\cup h_q(\ca_0)$, with $\ca_0=\{A_0, A', A''\}\subset D_q$, consists in $6$ distinct points and that they are in convex position. Moreover, for  $(p,q)\in D_q$  and $r<pq$, we deduce from Lemma~\ref{lem:equivalence_convex_order_int_cdf} that $(p,r)$ belongs to $ C^-=C^-_{q, A'}$ if and only if:
\begin{enumerate}[(i)]
    \item $r \geq \Delta_1(p)$
    if $p\in (0, p']$, see~\eqref{eq:case1i-conclusion},
    \item  $r \geq \Delta_2(p)$
    if $p\in [p', \tilde p'']$, see~\eqref{eq:case1ii-conclusion},
    \item  $r \geq \Delta_3(p)$
    if $p\in [ \tilde p'',1)$, see~\eqref{eq:case1iii-conclusion}. 
\end{enumerate}
Since $\Delta_q$ belongs also to $C^-$, see Section~\ref{sec:deg-Dq}, 
we deduce that $C^- \cap D_q^-$ is the convex hull of $\ca_0$, $\Cv(\ca_0)$. 
Since $C^-$ is invariant by the symmetry $h_q$, see~\eqref{eq:sym-C}, and since $D_q=D_q^- \cup h_q(D_q^-)$, we deduce that $C^-=\Cv(\ca_0) \cup \Cv(h_q(\ca_0))$ and thus that:
\begin{equation}
    \label{eq:C+-p'<beta}
    q\in [1/2, 1), \quad p'\in (0, 1-q) \quad\text{and}\quad r'<p'q
    \, \implies\,  C^-=\Cv(\ca).
\end{equation}

\subsubsection{The case $q\in [1/2, 1)$, $p'\in [1-q, q]$ and $r'<qp'$}
\label{sec:proof_case2}

This is  part of case~\ref{it:2-D1} in Proposition~\ref{prop:bernoulli_cp} as the latter condition is equivalent to $A'\not \in \Delta_q$. The proof  is divided in four parts: after some preliminary remarks, we consider the cases $p$ belonging to $(0, p']$,  and $[p', 1)$ for which different binding constraints apply, and then give the conclusion.

By symmetry, it is enough to consider the case $r< pq$ (recall the equality is studied in Section~\ref{sec:deg-Dq}).

\medskip
\noindent
\textbf{a) Preliminaries.} Let $\Delta_1$ be the line passing through $A_0$ and $A'$, $\Delta_2$ the line passing through $A'$ and $h_q(A_0)$. Let $s_i$ be the slope of $\Delta_i$ for $i=1,2$. We have:
\[
s_1= \frac{r'}{p'}, \quad \text{and} \quad s_2=\frac{q-r'}{1-p'} \cdot
\]
Since we assume that $r'<p'q$, the point $A'$ is strictly under $\Delta_q$ and we have $s_1 < s_2$. 
We thus obtain that the points $A_0$, $A'$, and $h_q(A_0)$ are in convex position. Using the symmetry, we deduce that the set $\ca=\ca_0\cup h_q(\ca_0)$, with $\ca_0=\{A_0, A'\}\subset D_q$, consists in $4$ distinct points and that they are in convex position.

As $\Delta_1$ and $\Delta_2$ intersect at $A'$, we deduce that:
\begin{equation}
  \label{eq:gamma_inequalities_case2}
  \Delta_1(p) > \Delta_2(p) \iff p > p'; \quad \text{and} \quad \Delta_1(p) = \Delta_2(p) \iff p = p'.
\end{equation}
It is elementary to check that:
\begin{align}
  \label{eq:case2_cns_delta1} r \geq \Delta_1(p) &\iff a_1 \geq a_1' \iff a_4 \leq a_4', \\
  \label{eq:case2_cns_delta2} r \geq \Delta_2(p) &\iff  a_2 \geq a_2' \iff a_3 \leq a_3'.
\end{align}
Furthermore, we deduce from Remark~\ref{rem:stat-ordre} that:
\begin{equation}
    \label{eq:case2-order-a'}
    0\leq a'_2= a'_{(1)}<a'_1 = a'_{(2)}<1 < a'_3= a'_{(3)} < a'_4= a'_{(4)} \leq \frac{1}{1-q}\cdot
\end{equation}

\medskip
\noindent
\textbf{b) The case $p\in (0,p']$.} As $p \leq p' \leq q$, we get from Remark~\ref{rem:stat-ordre} that $a_{(4)} = a_4$.  
Recall:
\[
\Delta_1(p)=\frac{r'}{p'}\, p.
\]
We first consider the binding constraint. We deduce from~\eqref{eq:gamma1_gamma4},~\eqref{eq:case2-order-a'} and~\eqref{eq:case2_cns_delta1} that:
\begin{equation}
  \label{eq:gamma4_subcase2i}
  \Gamma (a_{(4)}') \geq 0 \iff a_{(4)} \leq a_{(4)}' \iff a_4 \leq a_4' \iff r \geq \Delta_1(p).
\end{equation}
We shall now assume that $r\geq \Delta_1(p)$ and prove that $\Gamma(a'_{i})\geq 0$ for $i=1,2,3$. For $i=1$, we have:
\[
\Gamma(a_{(1)}') \geq 0 \iff a_{(1)} - a_2' \geq 0 \quad \text{with} \quad \rho_{(1)} (a_{(1)} - a_2') =  r - \Delta_1(p) + \begin{cases}
  \frac{p(qp'-r')(p'-(1-q))}{p'(1-q)(1-p')} & \text{if $p \in (0,1-q]$,} \\
  \Delta_1(p) - \Delta_2(p) & \text{if $p \in [1-q,p']$.} 
\end{cases}
\]
For $i=2$, notice that $
\rho_1 (a_1 - a_1') = r - \Delta_1(p) \geq 0$.
So we have either $a_1'\in (a_2', a_{(1)})$ or $a_1' \in [a_{(1)}, a_{(2)}]$. We have: for $a_1'\in (a_2', a_{(1)})$ that $\Gamma(a'_{(2)})= \E[(a'_2 - Z')_+]\geq 0$; and for $a_1' \in [a_{(1)}, a_{(2)}]$ that:
\[
\Gamma(a_{(2)}') = (a_1' - a_2') \rho_2' - (a_1' - a_2) \rho_2 = r - \Delta_1(p) + \frac{(p'-p)(qp'-r')}{qp'} \cdot
\]
For $i = 3$, the proof proceeds exactly as in the case b) $q \in [1/2,1)$, $p' \in (0,1-q]$, and $r' < q p'$, when $p \in (0,p']$, as the argument remains unchanged, see~\eqref{eq:cas-b-a3}
and the equation before. We obtain:
\begin{equation}
  \label{eq:case2i-conclusion}
\Gamma(a_{(4)}') \geq 0 \iff r \geq \Delta_1(p) \implies \Gamma(a_{(i)}') \geq 0 \quad \text{for $i = 1,2,3$}.
\end{equation}

\medskip
\noindent
\textbf{c) The case $p\in [p',1)$.} As $p \geq p' \geq 1-q$, we get from Remark~\ref{rem:stat-ordre} that $a_{(1)} = a_2$. Recall:
\[
\Delta_2(p) = \frac{(p-p')(q-r')}{1-p'} + r'.
\]
We first consider the binding constraint. We deduce from~\eqref{eq:gamma1_gamma4},~\eqref{eq:case2-order-a'} and~\eqref{eq:case2_cns_delta2} that:
\[
\Gamma(a_{2}') \geq 0 \iff 
\Gamma(a_{(1)}') \geq 0 \iff a_{(1)}' \leq a_{(1)} \iff a_2'\leq a_2 \iff r \geq \Delta_2(p).
\]
We shall now assume that $r\geq \Delta_2(p)$ and prove that $\Gamma(a'_{(i)})\geq 0$ for $i=2,3,4$. For $i=4$, we have:
\[
\Gamma(a_{(4)}') \geq 0 \iff a_4' - a_{(4)} \quad \text{with} \quad \rho_{(4)} (a_4' - a_{(4)}) = r - \Delta_2(p) + \begin{cases}
  \Delta_2(p) - \Delta_1(p) & \text{if $p \leq q$,} \\
  \frac{(p-p')(qp'-r')}{p'(1-p')} & \text{if $p \geq q$.}
\end{cases}
\]
For $i=2$, notice that $\rho_1 (a_1 - a_1') = r - \Delta_2(p) + \Delta_2(p) - \Delta_1(p) \geq 0$. So we have either  $a_1'\in (a_2', a_2) $ or $a_1'\in [a_2, a_1]$. We have: for $a_1'\in (a_2', a_2)$ that $\Gamma(a'_{(2)})= \E[(a'_2 - Z')_+]\geq 0$; for $a_1'\in [a_2, a_1]$ that:
\[
\Gamma(a'_{(2)})= (a_1' - a_2') \rho_2' - (a_1' - a_2) \rho_2  = r - \Delta_2(p) + \frac{(p-p')(qp'-r')(p' -(1-q))}{qp'(1-p')}\cdot
\]
For $i=3$, notice that $
\rho_3 (a_3' - a_3) = r - \Delta_2(p)$.
So we have either $a_3' \in [a_{(3)}, a_{(4)}]$ or $a_3' \in (a_{(4)}, a_4')$. We have: for $a_3' \in [a_{(3)}, a_{(4)}]$ that:
\[
\Gamma(a_{(3)}') = (a_4' - a_3') \rho_4' - (a_4 - a_3') \rho_4 = r - \Delta_2(p) + \frac{(p-p')(qp'-r')}{q(1-p')};
\]
and for $a_3' \in (a_{(4)}, a_4')$ that $\Gamma(a_{(3)}') = \E[(a_3' - Z')_+]\geq 0$. We obtain:
\begin{equation}
  \label{eq:case2ii-conclusion}
  \Gamma(a_2') \geq 0 \iff r \geq \Delta_2(p) \implies \Gamma(a_{(i)}') \geq 0 \quad \text{for $i = 2,3,4$}.
  \end{equation}

\medskip
\noindent
\textbf{d) Conclusion.} 
Recall that $q\in [1/2, 1)$, $p'\in [1-q,q]$, $r'<p'q$. 
We have shown  that the set $\ca=\ca_0\cup h_q(\ca_0)$, with $\ca_0=\{A_0, A'\}\subset D_q$, consists in $4$ distinct points and that they are in convex position. Moreover, for $(p,q)\in D_q$  and $r<pq$, we deduce from Lemma~\ref{lem:equivalence_convex_order_int_cdf} that $(p,r)$ belongs to $ C^-=C^-_{q, A'}$ if and only if:
\begin{enumerate}[(i)]
    \item $r \geq \Delta_1(p)$
    if $p\in (0, p']$, see~\eqref{eq:case2i-conclusion},
    \item  $r \geq \Delta_2(p)$
    if $p\in [p', 1)$, see~\eqref{eq:case2ii-conclusion}.
\end{enumerate}
Since $\Delta_q$ belongs also to $C^-$, see Section~\ref{sec:deg-Dq}, 
we deduce that $C^- \cap D_q^-$ is the convex hull of $\ca_0$, $\Cv(\ca_0)$. 
Since $C^-$ is invariant by the symmetry $h_q$, see~\eqref{eq:sym-C}, and since $D_q=D_q^- \cup h_q(D_q^-)$, we deduce that $C^-=\Cv(\ca_0) \cup \Cv(h_q(\ca_0))$ and thus that: 
\begin{equation}
    \label{eq:C+-p'>beta}
    q\in [1/2, 1), \quad p'\in [ 1-q,q] \quad\text{and}\quad r'<p'q
    \, \implies\,  C^-=\Cv(\ca).
\end{equation}

\subsubsection{The case $q\in [1/2, 1)$, $p'\in (q, 1)$ and $r'<p'q$}
\label{sec:proof_case3}

This is  part of case~\ref{it:3-D1} in Proposition~\ref{prop:bernoulli_cp} as the latter condition is equivalent to $A'\not \in \Delta_q$. The proof is similar 
to the one of Section~\ref{sec:proof_case1}; it 
is divided in five parts: after some preliminary remarks, we consider the cases $p$ belonging to $(0, p'']$, $[p'', p']$ and $[p', 1)$ for which different binding constraints apply, and then give the conclusion.

By symmetry, it is enough to consider the case $r< pq$ (recall the equality is studied in Section~\ref{sec:deg-Dq}).

\medskip \noindent \textbf{a) Preliminaries.} Recall: 
\begin{equation}
  \label{eq:def-p''}
  p'' = \frac{p'(1-p')}{2(1-q) - (1-p')} \quad \text{and} \quad r'' = q p'' - (qp' - r') \frac{p'(1-q)}{q(2(1-q) - (1-p'))} \cdot
\end{equation}
Notice that $2(1-q) - (1-p')$ is positive and $r' < p'q$, then we have $r'' < qp''$. Since $q < p' < 1$ and $q \geq 1/2$, we easily deduce that $p'' \in (0,q)$. 
It is also easy to check that  $r''\in [0, q)$ and $r'' > p''+q-1$. 
In particular, we have $A''\in D_q^-$. 
We deduce that $A''\in D_q$ (and also that $\tilde A''=h_q(A'')\in D_q$). 
Let $\Delta_1$ be the line passing through $A_0$ and $A''$, $\Delta_2$ the line passing through $A''$ and $A'$, and $\Delta_3$ the line passing through $A'$ and $h_q(A_0)=(1,q)$. Let $s_i$ be the slope of $\Delta_i$ for $i=1,2,3$. We have:
\[
s_1= \frac{r''}{p''} = \frac{r'}{p'} + \frac{(qp'-r')(q-p')}{qp'(1-p')} , \quad 
s_2= \frac{r' - r''}{p' - p''} = \frac{r'}{p'} + \frac{qp'-r'}{2qp'} \quad \text{and} \quad 
s_3= \frac{q-r'}{1-p'} = \frac{r'}{p'} + \frac{qp'-r'}{p'(1-p')} \cdot
\]
Since we assume that $p' > q \geq 1/2$ and $r' < p'q$, we have $s_1 < s_2 < s_3$. 
Thus the points $A_0$, $A'$, $\tilde A''$, and $h_q(A_0)$ are in convex position. Using the symmetry, we deduce that the set $\ca=\ca_0\cup h_q(\ca_0)$, with $\ca_0=\{A_0, A', A''\}\subset D_q$, consists in 6 distinct points and that they are in convex position. 

\medskip

Let $r=\Delta_i(p)$ be the equation of the line $\Delta_i$, that is, $\Delta_i(p)=p_i+ s_i p$, with $p_1=0$, $p_2=r'' - s_2 p''$ and $p_3=q-s_3$. 
As $s_1<s_2<s_3$ and since $\Delta_1$ and $\Delta_2$ intersect at $A''$ and $\Delta_2$ and $\Delta_3$ intersect at $A'$, we deduce that:
\begin{equation}
  \label{eq:gamma_inequalities_case3}
\begin{aligned}
     \Delta_1(p) > \Delta_2(p) \iff p < p''; \quad \text{and} \quad \Delta_2(p) > \Delta_3(p) \iff p < \tilde p',\\
      \Delta_1(p) =\Delta_2(p) \iff p = p''; \quad \text{and} \quad \Delta_2(p) = \Delta_3(p) \iff p = \tilde p'.
\end{aligned} 
\end{equation}
We now derive some properties related to the atoms $(a'_i)_{i\in \{1, \ldots, 4\}}$ of $Z'$. Remark~\ref{rem:stat-ordre} implies that:
\begin{equation}
    \label{eq:case3-order-a'}
    0\leq a'_1= a'_{(2)}<a'_2 = a'_{(1)}< 1 <a'_3= a'_{(4)} < a'_4= a'_{(3)} \leq \frac{1}{1-p'}\cdot
\end{equation}
We shall use many time that: 
\begin{equation}
  \label{eq:case3_delta2}
  (a_1' - a_2') \rho_2'  -  (a_1' - a_1) \rho_1 - (a_1' - a_2) \rho_2 = 2(r-\Delta_2(p)). 
\end{equation}
Moreover we have:
\begin{align}
  \label{eq:case3_cns_delta1} r \geq \Delta_1(p) &\iff a_3' \geq a_4, \\
  \label{eq:case3_cns_delta3} r \geq \Delta_3(p) &\iff a_3' \geq a_3 \iff a_2' \leq a_2. 
\end{align}

\medskip \noindent \textbf{b) The case  $p \in (0,p'']$.} 
Recall $p''$ in~\eqref{eq:def-p''} belongs to $(0, q)$. 
As $0<p\leq p''< q$ and $r'<qp'$, we get from Remark~\ref{rem:stat-ordre} that $a_1=a_{(2)}< a_1<1$ and $1<a_3<a_4=a_{(4)}$. Recall and set:
\[
\Delta_1(p)=\frac{r''}{p''}\, p
\quad\text{and}\quad
\cstb= \frac{qp'-r'}{qp'(1-p')}\geq 0.
\]
We first consider the binding constraint. We deduce from~\eqref{eq:gamma1_gamma4} and~\eqref{eq:case3_cns_delta3} that:
\[
\Gamma(a_{(4)}') \geq 0 \iff \Gamma(a_{3}') \geq 0 \iff a_{(4)}' \geq a_{(4)} \iff a_3' \geq a_4 \iff r \geq \Delta_1(p).
\]
We now assume that $r\geq \Delta_1(p)$ and prove that  $\Gamma(a'_{(i)}) \geq 0$ for $i=1,2,3$. 
For $i = 1$, since $a'_{(1)}=a'_2$, we  get from~\eqref{eq:gamma1_gamma4}  that:
\[
\Gamma(a_{(1)}') \geq 0 \iff 
\Gamma(a_{2}') \geq 0 \iff 
a_{(1)}'\leq a_{(1)} \iff
a_{(1)} - a_2'\geq 0. 
\]
We have either $p \in (0,1-q]$ or $p \in [1-q,p'']$.  Thanks to Remark~\ref{rem:stat-ordre}, this gives 
$a_{(1)}=a_1$ in the former case and $a_{(1)}=a_2$ in the latter case. Using~\eqref{eq:gamma_inequalities_case3}, elementary computation  gives:
\[
\rho_1(a_1 - a_2')= r - \Delta_1(p) + \frac{p(2q-1)(qp'-r')}{q(1-p')(1-q)} \geq 0
\quad\text{and}\quad 
\rho_1(a_1 - a_2') = r - \Delta_3(p) \geq r - \Delta_1(p) \geq 0,
\]
as $q \geq 1/2$.  Thus, we have obtained that $r \geq \Delta_1(p) $ implies $ \Gamma(a'_{(1)}) = \Gamma(a'_{2}) \geq 0$.

For $i = 2$, notice that either $a_1' \in (a_2' , a_{(1)})$ or $a_1' \in [a_{(1)} , a_{(2)}]$ or $a_1' \in (a_{(2)} ,1)$. We have: for  $a_1' \in (a_2' , a_{(1)})$ that $\Gamma(a'_1)= \E[(Z'-a'_1)_+]\geq 0$; for $a_1' \in [a_{(1)} , a_{(2)}]$ that:
\[
\Gamma(a'_1) =\begin{cases}
  (a_1' - a_2') \rho_2' - (a_1' - a_1) \rho_1= r - \Delta_1(p) + \cstb \left( (p'-q)(p''-p) + (1-q) \cstb' \right)
  & \text{if $p \leq 1-q$,} \\
  (a_1' - a_2') \rho_2' - (a_1' - a_2) \rho_2 = r - \Delta_1(p) + \cstb \left( (1-q)(p''-p) + \cstb'' \right)
  & \text{if $p \geq 1-q$,}
\end{cases} 
\]
with:
\[
\cstb' = \frac{(1-p')(2q-1)}{2(1-q) - (1-p')} \geq 0, \quad \text{and} \quad \cstb'' = \frac{p'(1-p')(p'-q)}{2(1-q)-(1-p')} \geq 0,
\]
and for $a_1' \in (a_{(2)} ,1)$ that $\Gamma(a'_1) = (a_1' - a_2') \rho_2' - (a_1' - a_1) \rho_1 -  (a_1' - a_2) \rho_2$, that is, with \eqref{eq:case3_delta2} and \eqref{eq:gamma_inequalities_case3}:
\[
\Gamma(a'_1) = 2(r-\Delta_2(p)) \geq 2(r-\Delta_1(p)) \geq 0.
\]
For $i = 3$, notice that: 
\[
\rho_3 (a_4' - a_3) = r - \Delta_1(p) + \frac{\cstb}{(1-q)} \left( ((2q-1)(1-p') + (1-q)^2)(p'' -p) + \cstb' (p'+q)(p'-q)  \right) \geq 0.
\]
So we have either $a_4' \in [a_3,a_4]$ or $a_4' \in (a_4,a_3')$. We have: for $a_4' \in [a_3,a_4]$ that:
\[
\Gamma(a_4') = (a_3' - a_4') \rho_3' - (a_4 - a_4') \rho_4 = r - \Delta_1(p) + \cstb \frac{p'-q}{1-q} \left( (1-q)(p''-p) + \cstb' (p'-q)\right),
\]
and for $a_4' \in (a_4,a_3')$ that $\Gamma(a_4') = \E[(Z'-a_4')_+]\geq 0$. So we have obtained:
\begin{equation}
    \label{eq:case3i-conclusion}
 \Gamma(a'_3)\geq 0 \iff r\geq \Delta_1(p) \implies \Gamma(a'_i)\geq 0 \quad \text{for $i=1,2,4$.}
\end{equation}

\medskip \noindent \textbf{c) The case  $p \in  [p'',p']$.} Recall $p''$ in~\eqref{eq:def-p''} belongs to $(0, q)$. 
As $q< p'<1$, we get that $p''<q<p'$.  As  $r<qp$,  we get from Remark~\ref{rem:stat-ordre} that $\{a_1, a_2\}=\{a_{(1)}, a_{(2)}\}$. Recall:
\[
\Delta_2(p)= (p-p'') \left(\frac{r'}{p'} + \cstc \right) + r'' 
\quad\text{with}\quad
\cstc= \frac{qp'-r'}{2qp'} \geq 0.
\]
Elementary computations give:
\begin{align}
\label{eq:case3i-1}
    \rho_1(a_1-a_1')&=r - \Delta_2(p) - \cstc' \quad \text{with} \quad \cstc' = \cstc (p-p') \geq 0,\\
 \label{eq:case3i-2}
    \rho_2(a_2-a_1')&=r - \Delta_2(p) - \cstc'' \quad \text{with} \quad \cstc'' = \cstc \left( p-p'' + \frac{2(1-q)(2q-1)}{2(1-q)-(1-p')} \right) \geq 0. 
\end{align}
We  first check the binding constraint:
\begin{equation}
  \label{eq:gamma1_subcase3ii}
  \Gamma(a_1') \geq 0 \iff  r \geq \Delta_2(p).
\end{equation}
To do so we distinguish three cases $a_1' \in (a_2', a_{(1)})$, $a_1' \in [a_{(1)},a_{(2)}]$ and $a_1' \in (a_{(2)},1)$. For the first case $a_1' \in (a_2', a_{(1)})$, we have $\Gamma(a_1') = \E[(a_1'-Z')_+]\geq 0$ and we deduce from~\eqref{eq:case3i-1} and~\eqref{eq:case3i-2} (depending on  $a_{(1)}$ being equal to $a_1$ or $a_2$) that $r\geq \Delta_2(p)$; and thus~\eqref{eq:gamma1_subcase3ii} holds. For the third case $a_1' \in (a_{(2)},1)$, we have using~\eqref{eq:case3_delta2}:
\[
\Gamma(a_1') = (a_1' - a_2') \rho_2'  -  (a_1' - a_1) \rho_1 - (a_1' - a_2) \rho_2 = 2(r-\Delta_2(p));
\]
and thus~\eqref{eq:gamma1_subcase3ii} holds. For the second case $a_1' \in [a_{(1)},a_{(2)}]$, we have for $p \leq 1-q$ that $a_{(1)}=a_1$ by Remark~\ref{rem:stat-ordre} and thus $a'_1\in [a_1, a_2]$ which thanks to~\eqref{eq:case3i-2} implies that $r \geq \Delta_2(p)$ as well as (subtract~\eqref{eq:case3i-1} to~\eqref{eq:case3_delta2}):
\[
\Gamma(a'_1)=(a_1' - a_2') \rho_2' - (a_1' - a_1) \rho_1 = r -\Delta_2(p) +  \cstc'  \geq 0;
\]
and for  $p\geq 1-q$
that $a_{(1)}=a_2$  and thus $a'_1\in [a_2, a_1]$ which thanks to~\eqref{eq:case3i-1} implies also $r\geq \Delta_2(p)$ as well as (subtract~\eqref{eq:case3i-2} to~\eqref{eq:case3_delta2}):
\[
\Gamma(a'_1)=(a_1' - a_2') \rho_2' - (a_1' - a_2) \rho_2 = r -\Delta_2(p) +  \cstc''  \geq 0;
\]
and thus~\eqref{eq:gamma1_subcase3ii} holds also. In conclusion Equation~\eqref{eq:gamma1_subcase3ii}  holds. 

\medskip

We  assume that $r\geq \Delta_2(p)$ and shall now prove that $\Gamma(a'_{i})\geq 0$ for $i=2,3,4$. 
To prove that  $\Gamma(a'_3) \geq 0$, it is enough by~\eqref{eq:gamma1_gamma4}, to prove that  $a_3' \geq a_{(4)}$ (as $a'_3=a'_{(4)}$).
We have:
\[
r\geq \Delta_2(p) \implies 
\begin{cases}
  r \geq \Delta_1(p) \iff a_3' \geq a_4 \iff  a_3' \geq a_{(4)} & \quad \text{for} \quad p\in [p'', q],\\
  r \geq \Delta_3(p) \iff a_3' \geq a_3 \iff  a_3' \geq a_{(4)} & \quad \text{for} \quad p\in [q, p'],
\end{cases}
\]
where we used~\eqref{eq:gamma_inequalities_case3} for the implication, \eqref{eq:case3_cns_delta1} and~\eqref{eq:case3_cns_delta3} for the first implications and the value of $a_{(4)}$ given in Remark~\ref{rem:stat-ordre} for the second implications. In conclusion we get $\Gamma(a'_3)\geq 0$.

Similarly, to prove that  $\Gamma(a'_2) \geq 0$, it is enough by~\eqref{eq:gamma1_gamma4}, to prove that $a_2' \leq a_{(1)}$ (as $a'_2=a'_{(1)}$). First notice that:
\begin{equation}
  \label{eq:a'2-a1}
    \rho_1(  a_1- a_2' ) = r - \Delta_1(p)+ \frac{p(2q-1)(qp'-r')}{q(1-p')(1-q)} \cdot
\end{equation}
We have:
\[
r\geq \Delta_2(p) \implies 
\begin{cases}
  r \geq \Delta_1(p) \implies a_2' \leq a_1 \iff  a_2' \leq a_{(1)} & \quad \text{for} \quad p\in [p'', 1-q],\\
  r \geq \Delta_3(p) \iff a_2' \leq a_2 \iff  a_2' \leq a_{(1)} & \quad \text{for} \quad p\in [1-q, p'],
\end{cases}
\]
where we used~\eqref{eq:gamma_inequalities_case3} for the first implication, then \eqref{eq:a'2-a1} and~\eqref{eq:case3_cns_delta3} and last the value of $a_{(1)}$ given in Remark~\ref{rem:stat-ordre}. In conclusion we get $\Gamma(a'_2)\geq 0$.

For $i=4$, notice that either  $a'_{4}\in (1, a_{(3)}) $ or $a'_4\in [a_{(3)}, a_{(4)}]$ or $a'_4\in (a_{(4)}, a'_4)$. We have: for $a'_{4}\in (1, a_{(3)}) $ that:
\[
\Gamma(a_4') = (a_3' - a_4') \rho_3' - (a_4 - a_4') \rho_4 - (a_3 - a_4') \rho_3 = 2(r-\Delta_2(p)) + 2 \cstc' \frac{2q-1}{1-q};
\]
for $a'_4\in [a_{(3)}, a_{(4)}]$ that:
\[
\Gamma(a_4') = (a_3' - a_4') \rho_3' - (a_{(4)}- a_4') \rho_{(4)} = (a_3' - a_4') \rho_3' - 
\begin{cases}
(a_{4}- a_4') \rho_{4} \quad\text{if $p \leq q$,} \\
(a_{3}- a_4') \rho_{3} \quad\text{if $p \geq q$,}
\end{cases}
\]
that is:
\[
\Gamma(a_4') = r - \Delta_2(p) + \frac{\cstc}{1-q}
\begin{cases}
(1-q)(p-p'') + \cstc'''\quad\text{if $p \leq q$,} \\
(p'-p)(3q-1) \quad\text{if $p \geq q$,}
\end{cases}
\]
with:
\[
\cstc''' = \frac{2(p'-q)^2(2q-1)}{2(1-q)-(1-p')} \geq 0,
\]
and for $a'_4\in (a_{(4)}, a'_4)$ that $\Gamma(a_4')= \E[(Z'-a_4')_+]\geq 0$. In conclusion, as $r\geq \Delta_2(p)$, we get $\Gamma(a'_4)\geq 0$.
We have obtained that: 
\begin{equation}
    \label{eq:case3ii-conclusion}
 \Gamma(a'_1)\geq 0 \iff r\geq \Delta_2(p) \implies \Gamma(a'_i)\geq 0 \quad \text{for $i=2,3,4$.}
\end{equation}

\medskip \noindent \textbf{d) The case  $p \in [p',1)$.} As $1/2 \leq q < p' \leq p <1$, $r' < qp'$, and $r < qp$, we get from Remark~\ref{rem:stat-ordre} that $a_{(1)}' = a_2'$, $a_{(4)}' = a_3'$,  $a_{(1)} = a_2$ and $a_{(4)} = a_3$. Recall and set:
\[
\Delta_3(p) = (p-p') \left(\frac{r'}{p'} + \cstd \right)+ r'
\quad\text{with}\quad
\cstd = \frac{qp'-r'}{p'(1-p')} \geq 0.
\]
We first consider the binding constraint. We deduce from~\eqref{eq:gamma1_gamma4},~\eqref{eq:case3-order-a'} and~\eqref{eq:case3_cns_delta3} that:
\begin{equation}
  \label{eq:gamma1_subcase3iii}
  \Gamma (a_{(1)}') \geq 0 \iff a_{(1)}' \leq a_{(1)} \iff a_2' \leq a_2 \iff r \geq \Delta_3(p) 
  \iff a_3 \leq a'_3 \iff \Gamma (a_{(4)}') \geq 0.
\end{equation}
We shall now assume that $r\geq \Delta_3(p)$ and prove that $\Gamma(a'_{(i)})\geq 0$ for $i=2,3$. 
First notice that:
\begin{equation}
    \label{eq:subcase3iii_eqa1a4}
    \rho_1 (a_1 -a'_1) = \rho_4 (a_4' - a_4) = r - \Delta_3(p) + \cstd (p-p') \geq 0.
\end{equation}
For $i=2$, with~\eqref{eq:subcase3iii_eqa1a4} notice that we have $a_1 \geq a_1'$.
So we have either $a_1' \in (a_2',a_2)$ or $a_1' \in [a_2,a_1]$. We have: for $a_1' \in (a_2',a_2)$ that $\Gamma(a'_1)= \E[(a'_1 - Z')_+]\geq 0$; for $a_1' \in [a_2,a_1]$ that:
\[
\Gamma(a_1') = (a_1' - a_2') \rho_2' - (a_1' - a_2) \rho_2 = r - \Delta_3(p) +\cstd  \,  \frac{(p-p')(p' - (1-q))}{q}  \geq 0.
\]
For $i=3$, with~\eqref{eq:subcase3iii_eqa1a4} notice that we have $a_4 \leq a_4'$.
So we have either $a_4' \in [a_4,a_3]$ or $a_4' \in (a_3, a_3')$. We have: for $a_4' \in [a_4,a_3]$ that:
\[
\Gamma(a_4') = (a_3'-a_4') \rho_3' - (a_3-a_4') \rho_3 = r - \Delta_3(p) + \cstd \, \frac{(p-p')(p'-q)}{(1-q)} \geq 0;
\]
and for $a_4' \in (a_3, a_3')$ that $\Gamma(a_4') = \E[(Z'-a_4')_+]\geq 0$. So we have obtained that: 
\begin{equation}
    \label{eq:case3iii-conclusion}
    \Gamma(a'_{(1)})\geq 0 \iff \Gamma (a'_{(4)}) \geq 0 \iff r\geq \Delta_3(p) \implies \Gamma(a'_{(i)})\geq 0 \quad \text{for $i=2,3$}.
\end{equation}

\medskip \noindent \textbf{e) Conclusion.} 
Recall that $q\in [1/2, 1)$, $p'\in (q,1)$, $r'<p'q$. 
We have shown  that the set $\ca=\ca_0\cup h_q(\ca_0)$, with $\ca_0=\{A_0, A', A''\}\subset D_q$, consists in $6$ distinct points and that they are in convex position. Moreover,  for  $(p,q)\in D_q$  and $r<pq$, we deduce from Lemma~\ref{lem:equivalence_convex_order_int_cdf} that $(p,r)$ belongs to $ C^-=C^-_{q, A'}$ if and only if:
\begin{enumerate}[(i)]
    \item $r \geq \Delta_1(p)$
    if $p\in (0, p'']$, see~\eqref{eq:case3i-conclusion},
    \item  $r \geq \Delta_2(p)$
    if $p\in [p'', p']$, see~\eqref{eq:case3ii-conclusion},
    \item  $r \geq \Delta_3(p)$
    if $p\in [p',1)$, see~\eqref{eq:case3iii-conclusion}. 
\end{enumerate}
Since $\Delta_q$ belongs also to $C^+$, see Section~\ref{sec:deg-Dq}, 
we deduce that $C^- \cap D_q^-$ is the convex hull of $\ca_0$, $\Cv(\ca_0)$. 
Since $C^-$ is invariant by the symmetry $h_q$, see~\eqref{eq:sym-C}, and since $D_q=D_q^- \cup h_q(D_q^-)$, we deduce that $C^-=\Cv(\ca_0) \cup \Cv(h_q(\ca_0))$ and thus that:
\begin{equation}
    \label{eq:C+-p'>1-beta} 
    q\in [1/2, 1), \quad p'\in (q,1) \quad\text{and}\quad r'<p'q
    \, \implies\,  C^-=\Cv(\ca).
\end{equation}

\subsubsection{Using the symmetries to get the expression of $C^-$} 
\label{sec:proof_bernoulli_sym_expression_finale} 
Recall the symmetry $h_q$ defined in Section~\ref{sec:Bernoulli-intro} leaves $C^-_{q, A'}$,  $D_q$ and $\Delta_q$ invariant.
Using~\eqref{eq:C+-p'<beta}, \eqref{eq:C+-p'>beta}
and~\eqref{eq:C+-p'>1-beta}, and Section~\ref{sec:deg-Dq} for the cases $q=1$ or $p'\in \{0, 1\}$ or $r'=qp'$, we deduce that Proposition~\ref{prop:bernoulli_cp} holds for $q\in [1/2, 1]$ and $A'$ in $D_q^-$.

Assume $q\in [1/2, 1]$ and $A'=(p', r')$ in $\{(p,r)\in D_q\, \colon\, r> pq\}$. Notice that $\tilde A'=h_q(A')=(\tilde p'=1-p', \tilde r'= q-r')$ belongs to $D_q^-$. Since 
 $D_q$, $\Delta_q$ and $C_{q,A'}^-$ 
  are invariant by $h_q$ and since $h_q$ is an involution, 
  we get that Properties~\ref{it:1-D1}-\ref{it:2-D1} 
  % and case $p'\in I_q$ of Property~\ref{it:3-D1} 
  from
Proposition~\ref{prop:bernoulli_cp} are obvious. To prove Property~\ref{it:3-D1}, that is,  $C^-_{q, A'}=\Cv( \{A_0, A', A'', h_q(A_0), h_q(A'), h_q(A'')\}$, 
 it is thus enough to check that $ \tilde  A'' =A''$, 
 %\in \{A'', h_q(A'')\})$, 
  where $\tilde A''= (\tilde p'', \tilde r'')$
 is defined by~\eqref{eq:bernoulli_cp''} with $A'=(p', r')$ replaced by $\tilde A'=(\tilde p', \tilde  r')$. Using that 
  $\tilde \alpha'= \tilde p' \wedge (1- \tilde p')=\alpha'$ and  
$\tilde \varepsilon= \ind_{\{\tilde p' <  \beta\}} - \ind_{\{\tilde  p' > 1 - \beta\}}=-\varepsilon$, we get:
\[
\begin{aligned}
 \tilde p''&= \frac{\tilde \alpha'(1-\tilde \alpha')}{2 \beta- \tilde \alpha'}
= p'',\\
\tilde r''&= q\tilde p'' + \tilde \varepsilon \tilde p'' ( q\tilde p' - \tilde r') \frac{\beta}{\alpha' (1- \beta)}
= q p'' + \varepsilon p'' (qp' -r') \frac{\beta}{\alpha' (1- \beta)}= r''. 
\end{aligned}
\]
This gives $\tilde A''=A''$, which finishes the proof for $q\in [1/2, 1]$ and $A'\in 
\{(p,r)\in D_q\, \colon\, r> pq\}$.

\medskip

To handle the case $q\in [0, 1/2)$, we first 
consider $\hat Y=1-Y$ and notice that 
$S_f(X,Y)=S_f(X,\hat Y)$ and $S_f(X',Y)=S_f(X',\hat Y)$ thanks to the invariance properties of $S_f$ given in Proposition~\ref{prop:invariance_csiszar_index}.
Furthermore, if we denote by  $(q,p',r', p,r)$ the parameters associated to the Bernoulli model $(Y,X',X)$, we get that the parameters associated to the Bernoulli model $(\hat Y,X',X)$, say $(\hat q, p', \hat r', p, \hat r)$, are given by 
$\hat q = \P( \hat Y=1) =1-q$, $\hat r' = \P (X'=\hat Y=1) = p'-r'$ and similarly for $\hat r$. We define
the involution  $\hk: (x,y) \mapsto (x, x-y)$. 
Recall $A_0=(0,0)$ and $h_q(A_0)=(1, q)$. 
Then as $\hat q= 1- q$, notice that 
$\hk( A_0)=A_0$, $\hk(h_q(A_0))=h_{\hat q} (A_0)$, and that 
$A$ belongs to $D_q$ (resp. $\Delta_q)$ if and only if $\hk( A)$
belongs to $D_{\hat q}$ (resp. $\Delta_{\hat q}$). 
In conclusion, we deduce that, with the obvious notation $\hk( A')=(\hat p', \hat r')=(p', p'-r')$:
\[
A \in C_{q,A'}^- \iff \hk( A) \in C_{\hat q,\hk( A')}^-.
\]
Notice also that $\hat I_q=[\hat \beta, 1- \hat \beta]$ is equal to $I_q=[ \beta, 1-  \beta]$, as $\hat \beta= \hat q \wedge (1-\hat q)= (1-q) \wedge q=\beta$. 
Since $\hk$ is an involution, we get that Properties~\ref{it:1-D1}-\ref{it:2-D1} of Property~\ref{it:3-D1}
  from
Proposition~\ref{prop:bernoulli_cp} are obvious. To prove Property~\ref{it:3-D1}, it is thus enough to check that $\hk( \hat A'' )=A''$, where 
$\hat A''=(\hat p'',  \hat r'') $ is defined as in~\eqref{eq:bernoulli_cp''} with $(q, p', r')$ replaced by $(\hat q, \hat p', \hat r')=(1-q, p', p'-r')$. Since $\hat \beta=
\beta$,  $\hat \alpha'=\hat p' \wedge (1- \hat p')=\alpha'$ and $\hat \varepsilon = \ind_{\{\hat p' < \hat \beta\}} - \ind_{\{\hat p' > 1 - \beta\}}=\varepsilon$. We have: 
\[
\begin{aligned}
    \hat p''& = \frac{\hat \alpha'(1-\hat \alpha')}{2 \hat \beta- \hat  \alpha'}
= p'',\\
\hat r''&= \hat q\hat  p'' + \hat \varepsilon \hat p'' ( \hat q\hat p' - \hat r') \frac{\hat \beta}{\hat \alpha' (1- \hat\beta)}
= (1-q) p'' + \varepsilon p'' (r'-q  p') \frac{\beta}{\alpha' (1- \beta)}= q-r''. 
\end{aligned}
\]
This gives $\hk(\hat A'')=A''$, which finishes the proof for $q\in [0, 1/2)$ and $A'\in 
 D_q$. 
Thus, the proof of Proposition~\ref{prop:bernoulli_cp} is complete.

\subsection{Description of $C^+_{q,A'}$}
\label{sec:description_Cm}
In this section, we give a description of $C^+_{q,A'}$.

There are straightforward cases such as: $A'\in\Delta_q$ and $A'\in E_q$. One has $C^+_{q,A'}=D_q$ in the former case and $C^+_{q,A'}=E_q$ in the latter. Further details are provided in~\eqref{eq:equiv-Dq} and~\eqref{eq:equiv-Eq}. 

Next, we assume that  $A'\not \in \Delta_q \cup E_q$ (which implies that $q\in (0, 1)$) and introduce more notation before the result. 
Set  $\beta = \min(q ,1-q)\in (0, 1)$ and recall $I_q=[\beta, 1-\beta]$.  We  define the function $s$ on $ D_q \setminus \Delta_q$  by:
\begin{equation}
    \label{eq:def-s}
s: A=(p,r) \mapsto q - \frac{1- \beta}{\beta} \, \frac{qp-r}{1-p}\cdot
\end{equation}
We also define the following convex subsets of $D_q$:
\[
R_\m=\{A\in D_q\, \colon\, p\in I_q\},
\quad
R_\rg=\{A\in D_q\, \colon\, p\leq \beta\}
\quad\text{and}\quad
R_\rd=\{A\in D_q\, \colon\, p\geq 1-\beta\}
\quad\text{where $A=(p,r)$}.
\]

These three sets partition $D_q$ according to the position of $p$ relative to $I_q$, and will play the role of the "middle", "left" and "right" regions in the description of $C_{q,A'}^+$ below.

Let $\ch_\rg$ (resp.\ $\ch'_\rg$) be the closed half-plane
delimited by the line passing through $A_0$ and $A'$ (resp. through $A_0$ with slopes $s(A')$) and  which does not contain the point   $A_1$ if $A' \in D_q^-$, and does not  contain the point
$h(A_1)=(1-q,0)$ if $A'\in D_q^+$. 
Define the closed half-planes $\ch_\rd$ and $\ch'_\rd$ similarly as $\ch_\rg$ and $\ch'_\rg$ but with $A_0$ replaced by $h(A_0)=(1,q)$ and with the slope $s(A')$ replaced by $s(h(A'))$.

\begin{figure}[htbp]
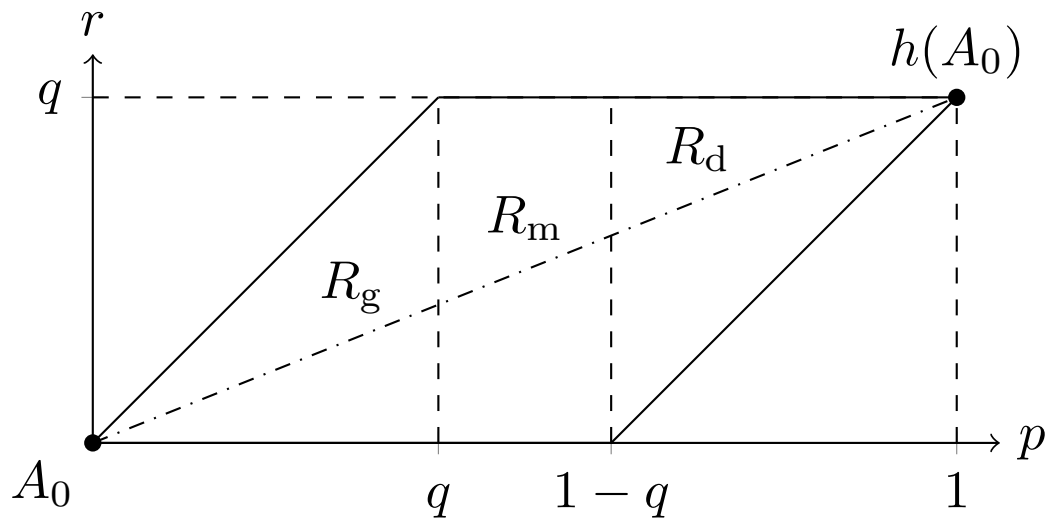
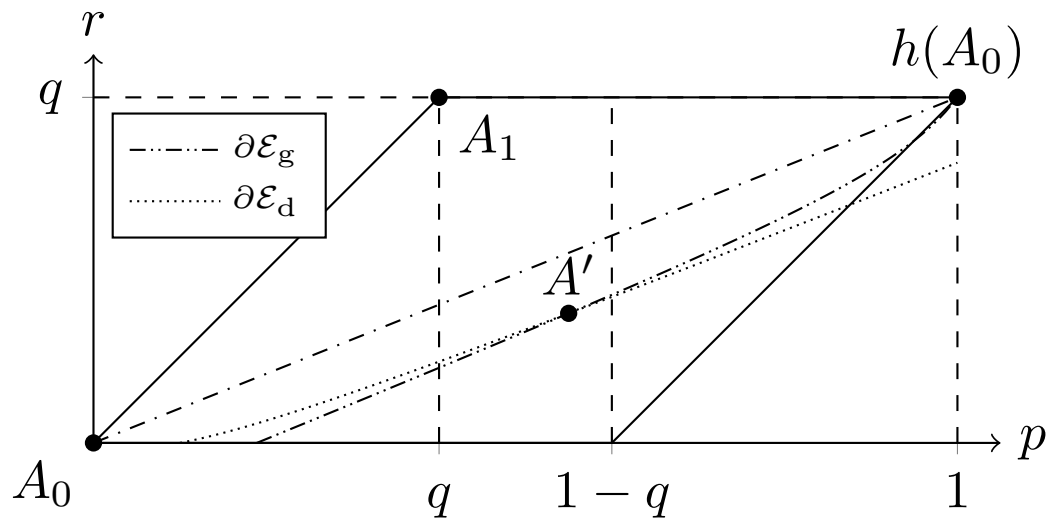
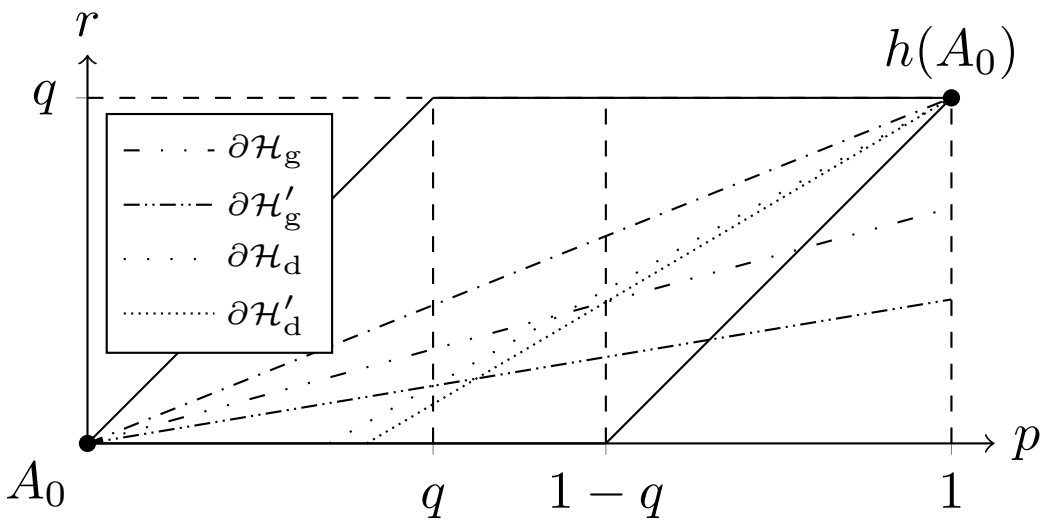

\centering
\begin{subfigure}[c]{0.49\textwidth}
\centering
\begin{RegionPlotBernoulli}{0.55}{0.15}{0.4}
  {
    \addplot[black] coordinates {(0,0) (0.4,0.4)};
    \addplot[black] coordinates {(0.6,0) (1,0.4)};
    \addplot[black] coordinates {(1,0.4) (0.4,0.4)};
    \addplot[black] coordinates {(0,0) (0.6,0)};

    \fill (0.3,0.18) circle (0) node {$R_\rg$};
    \fill (0.5,0.26) circle (0) node {$R_\m$};
    \fill (0.7,0.34) circle (0) node {$R_\rd$};
  }
  {

  }
\end{RegionPlotBernoulli}
\caption{The regions $R_\rg$, $R_\m$, and $R_\rd$ (in white).}
%\label{fig:}
\end{subfigure}
\hfill
\begin{subfigure}[c]{0.49\textwidth}
\centering
\begin{RegionPlotBernoulli}{0.55}{0.15}{0.4}
  { 
    \pgfplotsset{legend style={font=\tiny, at={(0.02,0.85)}, anchor=north west}}

    \addplot[black, forget plot] coordinates {(0,0) (0.4,0.4)};
    \addplot[black, forget plot] coordinates {(0.6,0) (1,0.4)};
    \addplot[black, forget plot] coordinates {(1,0.4) (0.4,0.4)};
    \addplot[black, forget plot] coordinates {(0,0) (0.6,0)};

    \addplot[domain=0.189022777135356:1, smooth, samples = \nsamples, densely dashdotdotted](\x,{0.4*\x - 0.084*(1 - \x)/(1.09 - \x)});
    \addlegendentry{$\partial \ce_\rg$}

    \addplot[domain=0.1:1, smooth, samples = \nsamples, densely dotted](\x,{0.4*\x - 0.084*\x/(\x + 0.11)});
    \addlegendentry{$\partial \ce_\rd$}

    \fill (0.55,0.15) circle (0.01) node[above] {$A'$};
    \fill (0.4,0.4) circle (0.01) node[below right] {$A_1$};
  }
  {

  }
\end{RegionPlotBernoulli}
\caption{Boundaries of $\ce_\rg$ and $\ce_\rd$.}
\end{subfigure}
\vspace{0.5cm}

\begin{subfigure}[c]{0.49\textwidth}
\centering
\begin{RegionPlotBernoulli}{0.55}{0.15}{0.4}
  {
    \pgfplotsset{legend style={font=\tiny, at={(0.02,0.85)}, anchor=north west}}

    \addplot[black, forget plot] coordinates {(0,0) (0.4,0.4)};
    \addplot[black, forget plot] coordinates {(0.6,0) (1,0.4)};
    \addplot[black, forget plot] coordinates {(1,0.4) (0.4,0.4)};
    \addplot[black, forget plot] coordinates {(0,0) (0.6,0)};

    \addplot[domain=0:1, smooth, samples = \nsamples, loosely dashdotdotted](\x,{0.272727272727273*\x});
    \addlegendentry{$\partial \ch_\rg$}
    \addplot[domain=0:1, smooth, samples = \nsamples, densely dashdotdotted](\x,{1/6 * \x});
    \addlegendentry{$\partial \ch_\rg'$}

    \addplot[domain=0.28:1, smooth, samples = \nsamples, loosely dotted](\x,{0.555555555555556*\x - 0.155555555555556});
    \addlegendentry{$\partial \ch_\rd$}
    \addplot[domain=0.32307692307692315:1, smooth, samples = \nsamples, densely dotted](\x,{0.590909090909091*\x - 0.190909090909091});
    \addlegendentry{$\partial \ch_\rd'$}
  }
  {

  }
\end{RegionPlotBernoulli}
\caption{Boundaries of the half planes $\ch_\rg$, $\ch_\rg'$, $\ch_\rd$ and $\ch_\rd'$ (not containing $A_1$).}
\end{subfigure}
\caption{An instance  of the boundaries of the sets describing $C_{q,A'}^+\cap D_q^-$ where $A'\in D_q^-$ and $A'\not  \in \Delta_q \cup E_q$ (with  $q=0.4$, $p' = 0.55\in I_q$, and $r'=0.15$).}
\label{fig:description_ensembles_Cp}
\end{figure}

Notice that $(1-p) + (1-p')(1-2\beta)\geq 0$ (resp.\ $p + p'(1-2\beta)\geq 0$) is null if and only if $\beta=1/2$ and $p=1$ (resp.\ $p=0$). We define the continuous functions
 $\psi_\rg$ and $\psi_\rd$  on $[0, 1]$ by:
\begin{equation}
\label{eq:bernoulli_cm_psi}
    \psi_\rg(p) = qp - \frac{c(1-p)}{(1-p) + (1-p')(1-2\beta)}
    \quad \text{and} \quad
    \psi_\rd(p) = qp - \frac{cp}{p + p'(1-2\beta)},
\end{equation}
with $c = 2(1-\beta)(qp'-r')$, and the conventions $\psi_\rg(1)=r'+(1-p')/2$
and $\psi_\rd(0)=r'-p'/2$ when $\beta =1/2$. For the particular case $q = 1/2$, those two functions reduce to: 
\begin{equation}
    \label{eq:bernoulli_psi_1/2}
    \psi_\rg(p) = \psi_\rd(p) = \frac{p-p'}{2}+r'. 
\end{equation}
For $i=\rg, \rd$, let $\ce_i$  be the set defined by:
\[
\ce_i = \left\{A= (p,r) \in D_q \, \colon \, \varepsilon (r- \psi_i(p)) \geq 0\right\}
\quad\text{with}\quad
\varepsilon= \begin{cases}
    + 1 & \text{ if $A'\in D_q^+$,}\\
    - 1& \text{ if $A'\in D_q^-$.}
\end{cases}
\]
We mention that the set $\ce_i$ is closed but not convex in general. 
In Fig.~\ref{fig:description_ensembles_Cp}, the various sets $R_\rg$, $R_\m$, $R_\rd$, together with the boundaries of $\ce_\rg$, $\ce_\rd$, $\ch_\rg$, $\ch_\rd$, $\ch'_\rg$ and $\ch'_\rd$, are shown for an example of $A'$.

The proof of the next result is in Section~\ref{sec:proof-bernoulli_cm}. Its proof relies simply on Proposition~\ref{prop:bernoulli_cp} and relation~\eqref{eq:C-=invertC+}.

\begin{prop}[Description of $C_{q,A'}^+$]\label{prop:bernoulli_cm}
Let $q \in [0,1]$ and $A' = (p',r') \in D_q$.  Set $\beta=\min(q,1-q)$. We have:
\begin{enumerate}[(i)]
    \item If $A' \in \Delta_q$, then $C_{q,A'}^+ = D_q$.
    \item If $A'\in E_q$, then $C_{q,A'}^+ =E_q$.
    \item \label{it:C-3}
    If $A' \notin \Delta_q\cup E_q$ (and thus $q\in (0, 1)$), we consider the closed sets (possibly empty): 
    \[
    \cc_i =  R_i \cap \ch_i \cap \ch_i' \cap \ce_i \quad \text{for $i=\rg, \rd$, and} \quad
    \cc_q =  R_\m \cap \ch_\rg \cap \ch_\rd.
    \]
    Then we have $C_{q,A'}^+ = \cc \cup h(\cc)$, with $\cc = \cc_\rg \cup \cc_q \cup \cc_\rd$.
\end{enumerate}
\end{prop}

\begin{figure}
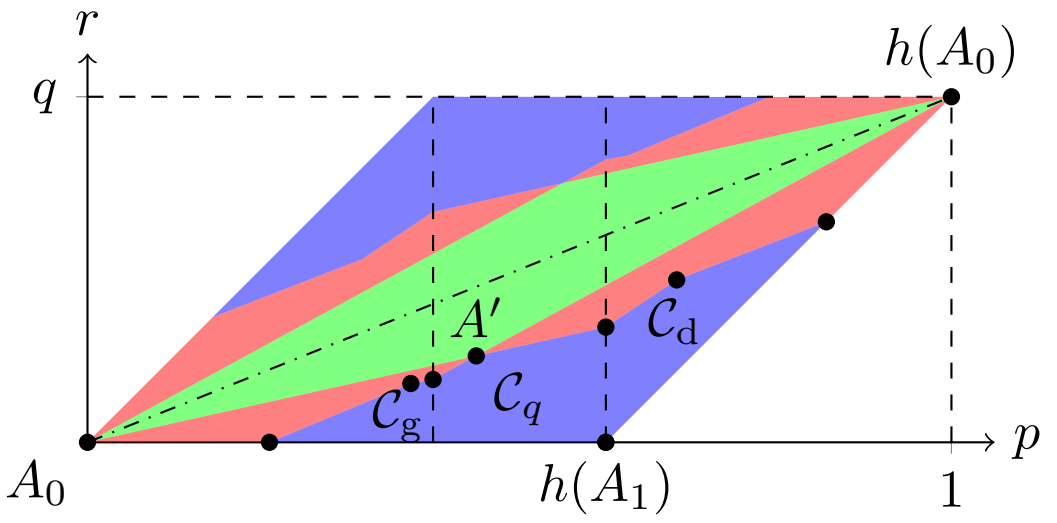
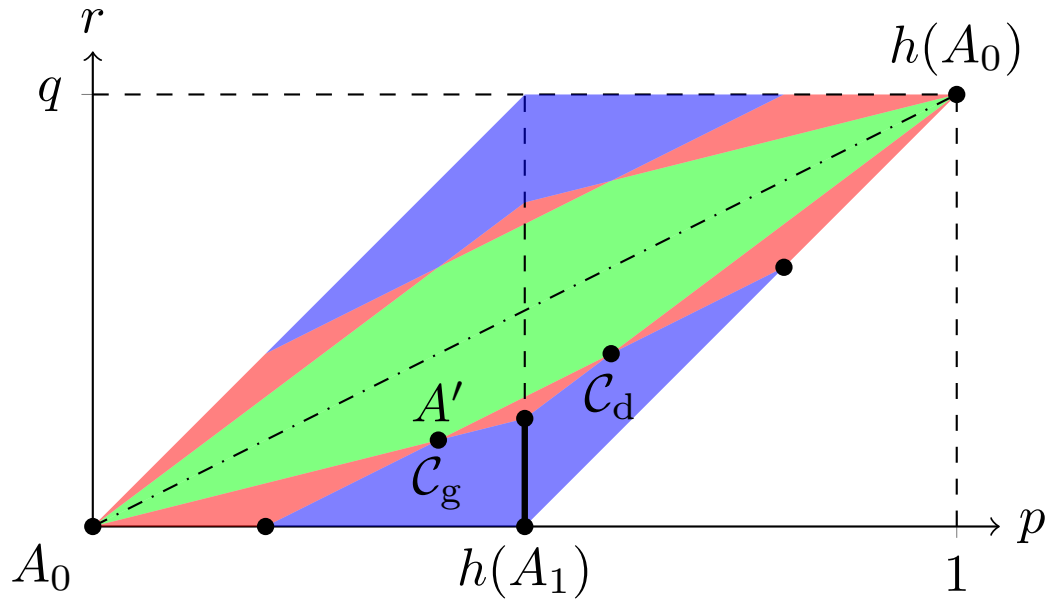
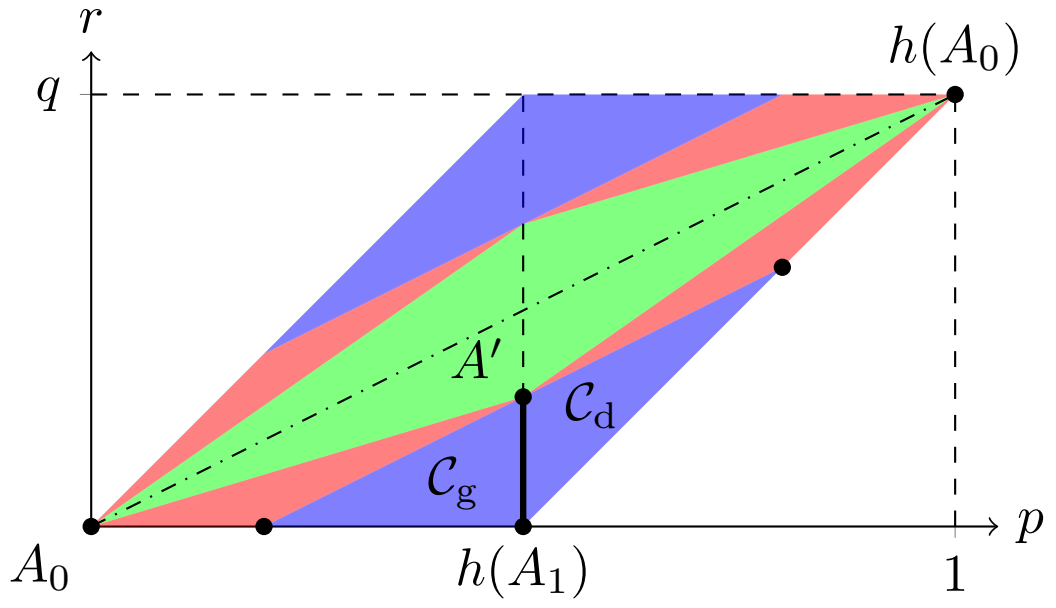
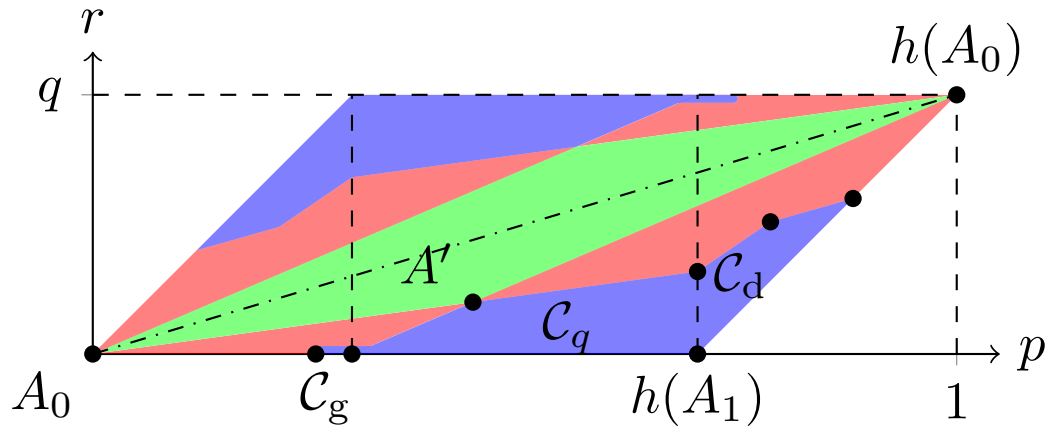
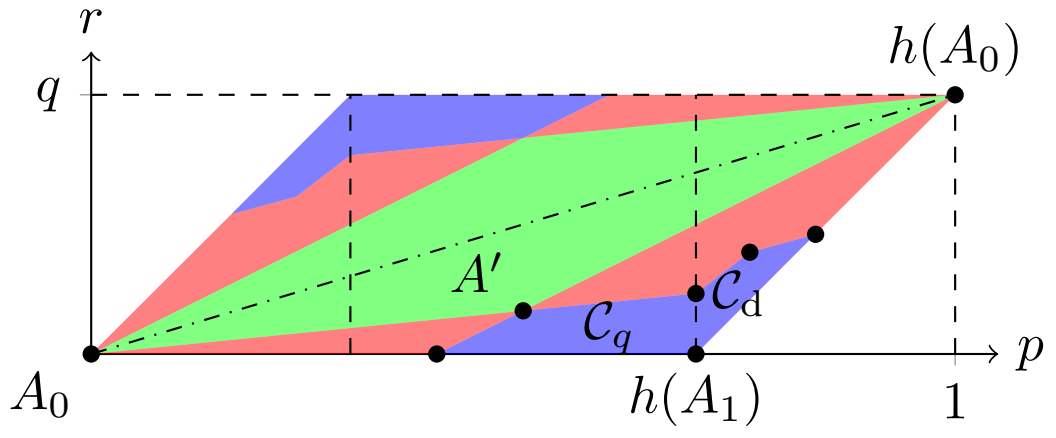

    \centering

\begin{subfigure}[c]{0.49\textwidth}
  \centering
  \begin{RegionPlotBernoulli}{0.3}{0.05}{0.45}
  {
    \pgfplotsset{xtick={1}, xticklabels={$1$}}

    \fill (0.55,0) circle (0.01) node[below] {$h(A_1)$};
    \fill (0.191226924268372,0) circle (0.01);
    \fill (0.3,0.05) circle (0.01) node[above] {$A'$};
    \fill (0.55,0.09166666666666669) circle (0.01);
    \fill (0.7405124837953327,0.2433710519110983) circle (0.01);
    \fill (0.835889894354067,0.285889894354067) circle (0.01);

    \fill (0.37,0.03) circle (0) node {$\cc_\rg$};
    \fill (0.5,0.04) circle (0) node {$\cc_q$};
    \fill (0.62,0.11) circle (0) node {$\cc_\rd$};
  }
  {
    \path[fill=\ColorCm, opacity=\OpacityPlot]
    (0,0) -- (0.3,0.05) -- (0.65,0.21136363636363642) -- (1,0.45) -- cycle;
    
    \path[fill=\ColorOut, opacity=\OpacityPlot]
    (0,0) -- plot[domain=0.191226924268372:0.3, smooth, samples = \nsamples](\x,{0.45*\x - 0.0935*(1 - \x)/(1.07 - \x)}) -- cycle;
    
    \path[fill=\ColorCp, opacity=\OpacityPlot]
    plot[domain=0.191226924268372:0.3, smooth, samples = \nsamples](\x,{0.45*\x - 0.0935*(1 - \x)/(1.07 - \x)})
    -- (0.55,0.09166666666666669)
    -- plot[domain=0.7405124837953327:0.835889894354067, smooth,samples = \nsamples](\x,{0.45*\x - 0.0935*\x/(\x + 0.03)})
    -- (0.55,0)
    -- cycle;
    
    \path[fill=\ColorOut, opacity=\OpacityPlot]
    (0.3,0.05) -- (0.55,0.09166666666666669) -- plot[domain=0.7405124837953327:0.835889894354067, smooth,samples = \nsamples](\x,{0.45*\x - 0.0935*\x/(\x + 0.03)}) -- (1,0.45) -- (0.65,0.21136363636363642) -- cycle;
    }
\end{RegionPlotBernoulli}
  \caption{Case: $p' < \beta$ ($p' = 0.3$, $r' = 0.05$, $q = 0.45$).}
  \label{fig:ber_cm_p_in_Rg}
\end{subfigure}
\hfill
\begin{subfigure}[c]{0.49\textwidth}
\centering
  \begin{RegionPlotBernoulli}{0.7}{0.3}{0.55}
  {
    \pgfplotsset{xtick={1}, xticklabels={$1$}}

    \fill (0.45,0) circle (0.01) node[below] {$h(A_1)$};
    \fill (0.164110105645933,0) circle (0.01);
    \fill (0.25948751620466726, 0.0528585681157655) circle (0.01);
    \fill (0.45,0.09166666666666667) circle (0.01);
    \fill (0.7,0.3) circle (0.01) node[left] {$A'$};
    \fill (0.808773075731628,0.35877307573162803) circle (0.01);

    \fill (0.35,0.035) circle (0) node {$\cc_\rg$};
    \fill (0.5,0.09) circle (0) node {$\cc_q$};
    \fill (0.6,0.18) circle (0) node {$\cc_\rd$};
  }
  {
    \path[fill=\ColorCm, opacity=\OpacityPlot]
    (0,0) -- (0.35,0.1113636363636364) -- (0.7,0.3) -- (1,0.55) -- cycle;
    
    \path[fill=\ColorOut, opacity=\OpacityPlot]
    (0,0) 
    -- plot[domain=0.164110105645933:0.25948751620466726, smooth,samples = \nsamples](\x,{0.55*\x - 0.0935*(1 - \x)/(1.03 - \x)})
    -- (0.45,0.09166666666666667)
    -- (0.7,0.3)
    -- (0.35,0.1113636363636364)
    -- cycle;

    \path[fill=\ColorOut, opacity=\OpacityPlot]
    plot[domain=0.7:0.808773075731628, smooth,samples = \nsamples](\x,{0.55*\x - 0.0935*\x/(\x + 0.07)})
    -- (1,0.55) 
    -- cycle;
    
    \path[fill=\ColorCp, opacity=\OpacityPlot] 
    plot[domain=0.164110105645933:0.25948751620466726, smooth,samples = \nsamples](\x,{0.55*\x - 0.0935*(1 - \x)/(1.03 - \x)})
    -- (0.45,0.09166666666666667)
    -- (0.7,0.3)
    -- plot[domain=0.7:0.808773075731628, smooth,samples = \nsamples](\x,{0.55*\x - 0.0935*\x/(\x + 0.07)})
    -- (0.45,0) 
    -- cycle;
  }
  \end{RegionPlotBernoulli}
\caption{Case: $p' > 1-\beta$ ($p' = 0.7$, $r' = 0.3$, $q = 0.55$).}
\label{fig:ber_cm_p_in_Rd}
\end{subfigure}

\vspace{0.5cm}

\begin{subfigure}[c]{0.49\textwidth}
\centering
  \begin{RegionPlotBernoulli}{0.45}{0.1}{0.4}
    {
    \pgfplotsset{xtick={1}, xticklabels={$1$}}

    \fill (0.6, 0) circle (0.01);
    \fill (1-0.789354390525168, 0) circle (0.01);
    \fill (0.37421951145296506, 0.06803991117326635) circle (0.01);
    \fill (0.4, 0.07272727272727275) circle (0.01);
    \fill (0.6, 0.13333333333333333) circle (0.01); 
    \fill (0.6819705149024926, 0.1879803432683284) circle (0.01);
    \fill (1-0.144765682192536, 0.4-0.144765682192536) circle (0.01);

    \fill (0.45, 0.1) circle (0.01) node[above] {$A'$};
    \fill (0.6, 0) circle (0.01) node[below] {$h(A_1)$};

    \fill (0.36,0.03) circle (0) node {$\cc_\rg$};
    \fill (0.5,0.05) circle (0) node {$\cc_q$};
    \fill (0.68,0.14) circle (0) node {$\cc_\rd$};
    }
    {
    \path[fill=\ColorCm, opacity=\OpacityPlot]
    (0,0) -- (0.55,0.3) -- (1,0.4) -- cycle;
    
    \path[fill=\ColorOut, opacity=\OpacityPlot]
    (0,0) -- plot[domain=0.144765682192536:0.3180294850975073, smooth, samples = \nsamples](\x,{0.4*\x + 0.0959999999999999*(1 - \x)/(1.09 - \x)}) -- (0.4, 0.2666666666666666) -- (0.55,0.3)  -- cycle;

    \path[fill=\ColorOut, opacity=\OpacityPlot]
    plot[domain=0.6257804885470349:0.789354390525168, smooth, samples = \nsamples](\x,{0.4*\x + 0.0959999999999999*\x/(\x + 0.11)})
    -- (1,0.4) 
    -- (0.55,0.3)
    -- (0.6, 0.3272727272727272)
    -- cycle;
    
    \path[fill=\ColorCp, opacity=\OpacityPlot]
    plot[domain=0.144765682192536:0.3180294850975073, smooth, samples = \nsamples](\x,{0.4*\x + 0.0959999999999999*(1 - \x)/(1.09 - \x)}) -- (0.4, 0.2666666666666666) 
    -- (0.4, 0.2666666666666666)
    -- (0.55,0.3)
    -- (0.6, 0.3272727272727272)
    -- plot[domain=0.6257804885470349:0.789354390525168, smooth, samples = \nsamples](\x,{0.4*\x + 0.0959999999999999*\x/(\x + 0.11)})
    -- (0.4,0.4)
    -- cycle;   
    }
  \end{RegionPlotBernoulli}
\caption{Case: $p' \in I_q$ ($p' = 0.45$, $r' = 0.1$, $q=0.4$).}
\label{fig:ber_cm_p_in_Iq}
\end{subfigure}

\vspace{0.5cm}

\begin{subfigure}[c]{0.49\textwidth}
\centering
\begin{RegionPlotBernoulli}{0.4}{0.1}{0.5}
  {
    \pgfplotsset{xtick={1}, xticklabels={$1$}}

    \addplot[white, line width=1.5pt] coordinates {(0.5,0) (0.5,0.125)};

    \path[fill=\ColorCp, opacity=\OpacityPlot]
    (0.2,0) -- (1/2,0) -- (0.8,0.3) -- (0.6,0.2) -- (1/2,0.125) -- (0.4,0.1) -- cycle;
    \path[fill=\ColorCp, opacity=\OpacityPlot]
    (1-0.2,0.5) -- (1-1/2,0.5) -- (1-0.8,0.2) -- (1-0.6,0.3) -- (1-1/2,0.5-0.125) -- (1-0.4,0.4) -- cycle;

    \fill (0.4,0.1) circle (0.01) node[above]{$A'$};
    \fill (0.6,0.2) circle (0.01);

    \fill (0.5,0) circle (0.01) node[below] {$h(A_1)$};
    \fill (0.2,0) circle (0.01);
    \fill (1/2,0.125) circle (0.01);
    \fill (0.8,0.3) circle (0.01);

    \addplot[black, very thick] coordinates {(0.5,0) (0.5,0.125)};

    \fill (0.4,0.05) circle (0) node {$\cc_\rg$};
    \fill (0.6,0.15) circle (0) node {$\cc_\rd$};
  }
  {
    \path[fill=\ColorCm, opacity=\OpacityPlot]
    (0,0) -- (0.4,0.1) -- (0.6,0.2) -- (1,0.5) -- cycle;

    \path[fill=\ColorOut, opacity=\OpacityPlot]
    (0,0) -- (0.2,0) -- (0.4,0.1)  -- cycle;

    \path[fill=\ColorOut, opacity=\OpacityPlot]
    (0.4,0.1) -- (1/2,0.125) -- (0.6,0.2)  -- cycle;

    \path[fill=\ColorOut, opacity=\OpacityPlot]
    (0.6,0.2) -- (0.8,0.3) -- (1,1/2)  -- cycle;
  }
\end{RegionPlotBernoulli}
\caption{Case: $q = 1/2$ ($p' = 0.4$, $r' = 0.1$). The vertical bold line starting from $h(A_1)$ corresponds to $\cc_q$.}
\label{fig:ber_cm_q=1/2}
\end{subfigure}
\hfill
\begin{subfigure}[c]{0.49\textwidth}
\centering
\begin{RegionPlotBernoulli}{0.5}{0.15}{0.5}
  {
    \pgfplotsset{xtick={1}, xticklabels={$1$}}

    \addplot[white, line width=1.5pt] coordinates {(0.5,0) (0.5,0.15)};

    \path[fill=\ColorCp, opacity=\OpacityPlot]
    (0.2,0) -- (1/2,0) -- (0.8,0.3) -- cycle;   
    \path[fill=\ColorCp, opacity=\OpacityPlot]
    (0.8,1/2) -- (1/2,1/2) -- (0.2,0.2) -- cycle;   

    \fill (0.5,0.15) circle (0.01) node[above left] {$A'$};

    \fill (0.2,0) circle (0.01);
    \fill (0.8,0.3) circle (0.01);

    \fill (0.5,0) circle (0.01) node[below] {$h(A_1)$};
    
    \addplot[black, very thick] coordinates {(1/2,0) (1/2,0.15)};
    
    \fill (0.42,0.05) circle (0) node {$\cc_\rg$};
    \fill (0.58,0.14) circle (0) node {$\cc_\rd$};

  }
  {
    \path[fill=\ColorCm, opacity=\OpacityPlot]
    (0,0) -- (0.5,0.15) -- (1,0.5) -- cycle;
    
    \path[fill=\ColorOut, opacity=\OpacityPlot]
    (0,0) -- (0.2,0) -- (1/2,0.15) -- cycle;

    \path[fill=\ColorOut, opacity=\OpacityPlot]
    (1/2,0.15) -- (0.8,0.3) -- (1,1/2) -- cycle;

  }
\end{RegionPlotBernoulli}
\caption{Case: $q=p'=1/2$ ($r' = 0.15$). The vertical bold line starting from $h(A_1)$ corresponds to $\cc_q$.}
\label{fig:ber_cm_q=1/2_p=1/2}
\end{subfigure}

\vspace{0.5cm}

\begin{subfigure}[c]{0.49\textwidth}
\centering
\begin{RegionPlotBernoulli}{0.44}{0.06}{0.3}
  {  
    \pgfplotsset{xtick={1}, xticklabels={$1$}}

   \fill (0.44,0.06) circle (0.01) node[above left] {$A'$};

    \fill (0.7,0) circle (0.01) node[below] {$h(A_1)$};
    
    \fill (0.55,0.03) circle (0) node {$\cc_q$};
    \fill (0.75,0.09) circle (0) node {$\cc_\rd$};
    \fill (0.27,0) circle (0) node[below] {$\cc_\rg$};

    \fill (0.7,0.09545454545454544) circle (0.01); 
    \fill (0.7843808085167396, 0.152986914897777) circle (0.01);
    \fill (0.88, 0.88+0.3-1) circle (0.01);

    \fill (0.3,0) circle (0.01);
    \fill (0.258080465972004, 0) circle (0.01);
  }
  {
    \path[fill=\ColorCm, opacity=\OpacityPlot]
    (0,0) -- (0.44,0.06) -- (1,0.3) -- cycle;
  
    \path[fill=\ColorOut, opacity=\OpacityPlot]
    (0.44,0.06) -- (0.7, 0.09545454545454544) -- plot[domain=0.7843808085167396:0.88, smooth, samples = \nsamples](\x,{0.3*\x - 0.1008*\x/(\x + 0.176)})   -- (1,0.3) 
    -- cycle;

    \begin{scope}
      \clip (0,0) -- (0.3,0) -- (0.44,0.06) -- cycle;
      \fill[\ColorOut, opacity=\OpacityPlot, even odd rule]
        (0,0) -- (0.3,0) -- (0.44,0.06) -- cycle
        [rounded corners=0.75pt]
        ([xshift=-0.75pt]0.258080465972004,0) --
        ([xshift=0.75pt]0.4,0) --
        ([xshift=0.75pt,yshift=1.5pt]0.4,0) --
        ([xshift=-0.75pt,yshift=1.5pt]0.258080465972004,0) --
        cycle;
    \end{scope}

    \begin{scope}[transparency group, opacity=\OpacityPlot]
      \path[fill=\ColorCp]
      (0.3,0) -- (0.44,0.06) -- (0.7, 0.09545454545454544) -- plot[domain=0.7843808085167396:0.88, smooth, samples = \nsamples](\x,{0.3*\x - 0.1008*\x/(\x + 0.176)}) -- (0.7,0) -- cycle;

      \draw[\ColorCp, line width=1.5pt, yshift=1.5/2pt, line cap=round] (0.258080465972004,0) -- (0.4,0);
    \end{scope}
  }
\end{RegionPlotBernoulli}
\caption{Case: $\cc_\rg$ being a segment ($q=0.3$, $p' = 0.44$ and $r'=0.06$).}
\label{fig:ber_cm_rg_segment}
\end{subfigure}
\hfill
\begin{subfigure}[c]{0.49\textwidth}
\centering
\begin{RegionPlotBernoulli}{0.5}{0.05}{0.3}
  {
    \pgfplotsset{xtick={1}, xticklabels={$1$}}
    \fill (0.5,0.05) circle (0.01) node[above left] {$A'$};

    \fill (0.7,0) circle (0.01) node[below] {$h(A_1)$};
    
    \fill (0.6,0.03) circle (0) node {$\cc_q$};
    \fill (0.75,0.07) circle (0) node {$\cc_\rd$};

    \fill (0.4,0) circle (0.01);
    \fill (0.7,0.07) circle (0.01); 
    \fill (0.76234753829798, 0.11779977936178464) circle (0.01);
    \fill (0.838516480713450, 0.838516480713450+0.3-1) circle (0.01);
  }
  {
    \path[fill=\ColorCm, opacity=\OpacityPlot]
    (0,0) -- (0.5,0.05) -- (1,0.3) -- cycle;
        
    \path[fill=\ColorOut, opacity=\OpacityPlot]
    (0,0) -- (0.4,0) -- (0.5,0.05) -- cycle;

    \path[fill=\ColorOut, opacity=\OpacityPlot]
    (1/2,0.05) -- (0.7, 0.07) -- plot[domain=0.76234753829798:0.838516480713450, smooth, samples = \nsamples](\x,{0.3*\x - 0.14*\x/(\x + 0.2)})   -- (1,0.3) 
    -- cycle;

    \path[fill=\ColorCp, opacity=\OpacityPlot]
    (0.4,0) -- (1/2,0.05) -- (0.7, 0.07) -- plot[domain=0.76234753829798:0.838516480713450, smooth, samples = \nsamples](\x,{0.3*\x - 0.14*\x/(\x + 0.2)}) -- (0.7,0) -- cycle;   
  }  
\end{RegionPlotBernoulli}
\caption{Case: $\cc_\rg$ being empty ($q=0.3$, $p' = 0.5$ and $r'=0.05$).}
\label{fig:ber_cm_rg_empty}
\end{subfigure}

\caption{Description of $C^+_{q, A'}\cap D_q^-$ (in blue) with bullets on its boundary delimiting the change of active constraint.}
\label{fig:description_cm}
\end{figure}

Notice in particular that the set $C_{q,A'}^+$ is closed. 
We give some instance of the set $C_{q,A'}^+$:
in Fig.~\ref{fig:ber_cm_p_in_Rg} for  $A' \in R_\rg$; in Fig.~\ref{fig:ber_cm_p_in_Rd} for  $A' \in R_\rd$; and in Fig.~\ref{fig:ber_cm_p_in_Iq} for  $A' \in R_\m$. In particular, the case $q=1/2$ is illustrated in Fig.~\ref{fig:ber_cm_q=1/2} and the case $q = p' = 1/2$ in Fig.~\ref{fig:ber_cm_q=1/2_p=1/2}.

The set $C^+_{q,A'}$ is described by a number of constraints and a natural question is whether they are all always necessary. Indeed, in some cases, the geometry of the set is simpler and fewer constraints are necessary to describe this set. 
In section~\ref{sec:saturated_constraints} we identify the active constraints. Note that in the case $q=1/2$, the function $\psi_\rg$ and $\psi_\rd$ are equal and affine, see Fig.~\ref{fig:description_cp-cas-1demi}.

Next, we give some particular geometric properties of $\cc_i$. 

\begin{figure}[htbp]
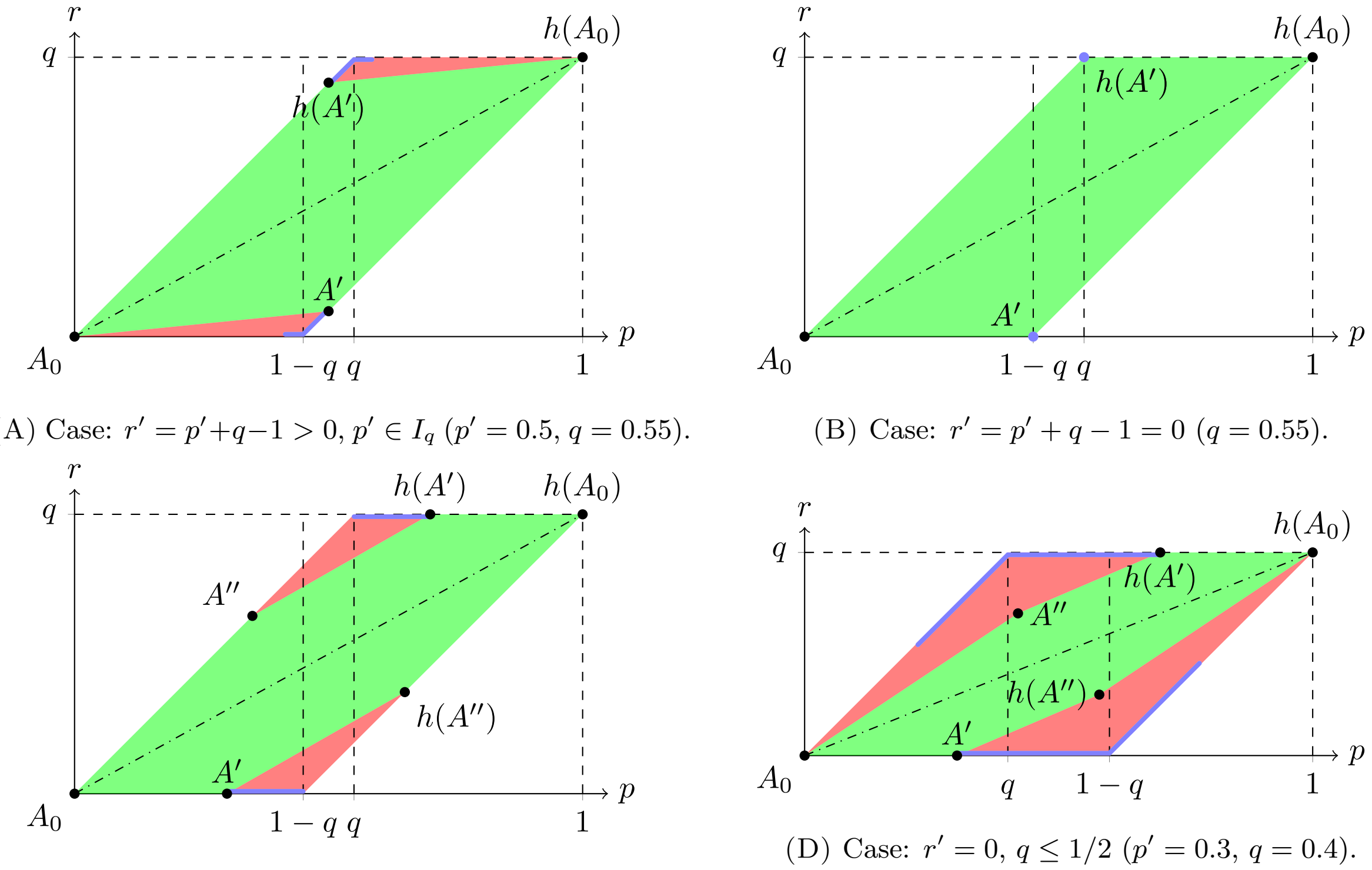

\begin{subfigure}[c]{0.49\textwidth}
\centering

\begin{RegionPlotBernoulli}{0.5}{0.05}{0.55}
  {
    \draw[white, opacity = 1, line width=1.5pt, yshift=1.5/2pt, line cap=round] (0.4145898033750315,0) -- (0.45,0) -- (0.5,0.05);
    \draw[white, opacity = 1, line width=1.5pt, yshift=-1.5/2pt, line cap=round] (1-0.4145898033750315,0.55) -- (1-0.45,0.55) -- (1-0.5,0.55-0.05);

    \draw[\ColorCp, opacity = \OpacityPlot, line width=1.5pt, yshift=1.5/2pt, line cap=round] (0.4145898033750315,0) -- (0.45,0) -- (0.5,0.05);
    \draw[\ColorCp, opacity = \OpacityPlot, line width=1.5pt, yshift=-1.5/2pt, line cap=round] (1-0.4145898033750315,0.55) -- (1-0.45,0.55) -- (1-0.5,0.55-0.05);

    \fill (0.5,0.05) circle (0.01) node[above] {$A'$};
    \fill (0.5,0.5) circle (0.01) node[below] {$h(A')$};
  }
  {
    \path[fill=\ColorCm, opacity=\OpacityPlot]
    (0,0) -- (0.5,0.05) -- (1,0.55) -- cycle;
    \path[fill=\ColorOut, opacity=\OpacityPlot]
    (0,0) -- (0.45,0) -- (0.5,0.05)-- cycle;
  }
\end{RegionPlotBernoulli}
\caption{Case: $r'=p'+q-1 >0$, $p' \in I_q$ ($p'=0.5$, $q=0.55$). }
\label{fig:r_prime_eg_p_prime+q-1>0}
\end{subfigure}
\hfill
%\hfill
\begin{subfigure}[c]{0.49\textwidth}
\centering

\begin{RegionPlotBernoulli}{0.45}{0}{0.55}
  {
    \fill[white, opacity = 1] (0.45,0) circle (0.01) node[black, above left] {$A'$};
    \fill[white, opacity = 1] (0.55,0.55) circle (0.01) node[black, below right] {$h(A')$};
    \fill[\ColorCp, opacity = \OpacityPlot] (0.45,0) circle (0.01) node[black, above left] {$A'$};
    \fill[\ColorCp, opacity = \OpacityPlot] (0.55,0.55) circle (0.01) node[black, below right] {$h(A')$};
  }
  {
    \path[fill=\ColorCm, opacity=\OpacityPlot] (0,0) -- (0.45,0) -- (1,0.55) -- cycle;
  }
\end{RegionPlotBernoulli}
\caption{Case: $r'=p'+q-1=0$ ($q=0.55$).}
\label{fig:p_prime_eg_1-q_r_prime_0}
\label{fig:X=Y}
\end{subfigure}

\begin{subfigure}[c]{0.49\textwidth}
\centering

\begin{RegionPlotBernoulli}{0.3}{0}{0.55}
  {
    \draw[white, opacity = 1, line width=1.5pt, yshift=1.5/2pt, line cap=round] (0.3,0) -- (0.45,0);
    \draw[white, opacity = 1, line width=1.5pt, yshift=-1.5/2pt, line cap=round] (0.7,0.55) -- (0.55,0.55);
    \draw[\ColorCp, opacity = \OpacityPlot, line width=1.5pt, yshift=1.5/2pt, line cap=round] (0.3,0) -- (0.45,0);
    \draw[\ColorCp, opacity = \OpacityPlot, line width=1.5pt, yshift=-1.5/2pt, line cap=round] (0.7,0.55) -- (0.55,0.55);

    \fill (0.3,0) circle (0.01) node[above] {$A'$};
    \fill (0.7,0.55) circle (0.01) node[above] {$h(A')$};

    \fill (0.65,0.2) circle (0.01) node[below right] {$h(A'')$};
    \fill (0.35,0.35) circle (0.01) node[above left] {$A''$};
  }
  {
    \path[fill=\ColorCm, opacity=\OpacityPlot]
    (0,0) -- (0.3,0) -- (0.65,0.2) -- (1,0.55) -- cycle;
    \path[fill=\ColorOut, opacity=\OpacityPlot]
    (0.3,0) -- (0.45,0) -- (0.65,0.2) -- cycle;
  }
\end{RegionPlotBernoulli}
\caption{Case: $r'=0$, $p' \leq 1-q$ ($p'=0.3$, $q=0.55$).}
\label{fig:r_prime=0_p_prime<=1-q}
\end{subfigure}
\hfill
\begin{subfigure}[c]{0.49\textwidth}
\centering
\begin{RegionPlotBernoulli}{0.3}{0}{0.4}
  {
    \draw[white, opacity = 1, line width=1.5pt, yshift=1.5/2pt, line cap=round] (0.3,0) -- (0.6,0) -- (0.7772001872658766,0.17720018726587664);
    \draw[white, opacity = 1, line width=1.5pt, yshift=-1.5/2pt, line cap=round] (1-0.3,0.4) -- (1-0.6,0.4) -- (1-0.7772001872658766,0.4-0.17720018726587664);
    \draw[\ColorCp, opacity = \OpacityPlot, line width=1.5pt, yshift=1.5/2pt, line cap=round] (0.3,0) -- (0.6,0) -- (0.7772001872658766,0.17720018726587664);
    \draw[\ColorCp, opacity = \OpacityPlot, line width=1.5pt, yshift=-1.5/2pt, line cap=round] (1-0.3,0.4) -- (1-0.6,0.4) -- (1-0.7772001872658766,0.4-0.17720018726587664);

    \fill (0.3,0) circle (0.01) node[black, above] {$A'$};
    \fill (0.7,0.4) circle (0.01) node[black, below] {$h(A')$};

    \fill (0.42,0.28) circle (0.01) node[black, right] {$A''$};
    \fill (0.58,0.12) circle (0.01) node[black, left] {$h(A'')$};

  }
  {
    \path[fill=\ColorCm, opacity=\OpacityPlot]
    (0,0) -- (0.3,0) -- (0.58,0.12) -- (1,0.4) -- cycle;
    \path[fill=\ColorOut, opacity=\OpacityPlot]
    (0.3,0) -- (0.6,0) -- (1,0.4)-- (0.58,0.12) -- cycle;
  }
\end{RegionPlotBernoulli}
\caption{Case: $r'= 0$, $q \leq 1/2$ ($p' = 0.3$, $q=0.4$).}
\label{fig:r_prime=0,q<=1/2}
\end{subfigure}
%\vspace{0.5cm}
\caption{Description of $C_{q,A'}^-$ (in green) and $C^+_{q, A'}$ (in blue) with $A'$ belonging to the boundary of $D_q$. See Fig.~\ref{fig:X=Y} for the degenerated case $X=Y$ a.s..}

\label{fig:A'-on-boundary}
\end{figure}

\begin{rem}[When is $\cc_i$ empty?]
Notice that $\cc_q$ is never empty as it contains  the point $A_1$ (or its symmetric $h(A_1)$). This point also belongs to  (the boundary of) $\cc_\rg$ whenever $1 - q \le q$, and to  (the boundary of) $\cc_\rd$ whenever $q \leq 1 - q$. 
The sets $\cc_\rg$ and $\cc_\rd$ may be empty, see 
Fig.~\ref{fig:p_prime_eg_1-q_r_prime_0} (with the sufficient conditions $A'\not \in \Delta_q$ and $A'\in E_q$), 
Fig.~\ref{fig:r_prime_eg_p_prime+q-1>0} and~\ref{fig:r_prime=0_p_prime<=1-q} (with the sufficient conditions $A'\not \in \Delta_q$,  $A'\in \partial D_q$ and $p'\in \mathring{I_q}$), and Fig.~\ref{fig:ber_cm_rg_empty} which is an instance where $A'\not \in \Delta_q\cup   \partial D_q$ and $p'\in \mathring{I_q}$. 
\end{rem}

\subsection{Proof of Proposition~\ref{prop:bernoulli_cm}}
\label{sec:proof-bernoulli_cm}

Set $A'=(p', r')\in D_q$.
We deduce from Proposition~\ref{prop:C+-cup-C-}, that  $C^-_{q, A'}=D_q$ if $A'\in \Delta_q$, and that
$C^+_{q, A'}=E_q$ if $A'\in E_q$.

We shall prove Proposition~\ref{prop:bernoulli_cm}~\ref{it:C-3}.
So we assume that $A'\not \in \Delta_q\cup E_q$, and thus $q\in (0, 1)$. 
We also deduce from Propositions~\ref{prop:C+-cap-C-} and~\ref{prop:C+-cup-C-} that $A_0$ and  $h_q(A_0)$ do not belong to $C^+_{q, A'}$, so that $A=(p,r)\in C^+_{q, A'}$ implies $p\in (0,1)$.

We derive from~\eqref{eq:C-=invertC+} the following description of $C^{+-}_{q,A'}=C^+_{q, A'}\cap D_q^-$:
\begin{equation}
    \label{eq:C-=C+}
C^{+-}_{q,A'}=\left\{A= (p,r) \in D_q^-\, \colon\, A'\in C^-_{q,A} \right\}.
\end{equation}
Recall $\beta=q \wedge (1-q) $.

\subsubsection{The case $A' \in D_q^- \setminus \Delta_q$ and $A=(p,r) \in R_\rg$ with $p \in (0, \beta]$.} Notice that $q\in (0, 1)$. 
Recall $h_q$ is an involution. Let $A^*=(p^*, r^*)$ be such that  $h_q(A^*)$ is defined as $A''$ with $(p',r')$ replaced by $(p,r)$ in~\eqref{eq:bernoulli_cp''}. 
We consider the affine functions $y=\Delta_i^{A}(x)$, for $i=1,2,3$, which are the equation of the line going through respectively  $A_0$ and $A$, $A$ and $A^*$, and  $A^*$ and $h_q(A_0)$:
\[
\Delta_1^{A}(x) = \frac{r}{p} x,\quad 
\Delta_2^{A}(x) = \frac{r^* - r}{p^* -p}(x-p)+r \quad \text{and} \quad
\Delta_3^{A}(x) = \frac{q - r^*}{1-p^*}(x-1)+q.
\]
We deduce from~\eqref{eq:C-=C+} that $A=(p,r)\in D_q^-$ belongs to $C^{+-}_{q,A'}$ if and only if $r' \geq \max_{i=1,2,3}\Delta_i^{A}(p')$. Elementary computation gives, that for $A=(p,r)\in D_q^-$ with $p \in (0,\beta]$:
\[
A = (p,r) \in C^{+-}_{q,A'} 
\iff
r' \geq \max_{i=1,2,3}\Delta_i^{A}(p')
\iff
r \leq \min\left( \Lambda_\rg(p), \psi_\rg(p), \Lambda_\rg'(p) \right),
\]
with:
\begin{equation}
\label{eq:bernoulli_proof_cm_def_lambdag}
\Lambda_\rg(x) = \frac{r'}{p'} x \quad \text{and} \quad
\Lambda_\rg'(x) = x\left(q-  \frac{(qp'-r')(1-\beta)}{\beta(1-p')}\right) = xs(p',r'),
\end{equation}
and the maps $s$ and $\psi_\rg$ defined in~\eqref{eq:def-s} and~\eqref{eq:bernoulli_cm_psi}. 
By observing that:
\[
r\leq \Lambda_\rg(p) \iff A \in  R_\rg \cap \ch_\rg; 
\quad
r\leq \Lambda_\rg'(p) \iff A \in  R_\rg \cap \ch_\rg'; 
\quad \text{and} \quad
r\leq \psi_\rg(p) \iff A \in  R_\rg \cap \ce_\rg,
\]
we obtain that $C_{q, A'}^{+-}\cap R_\rg=\cc_\rg$, with $\cc_\rg=R_\rg \cap \ch_\rg \cap \ch'_\rg \cap \ce_\cg$.

\medskip 

We will occasionally refer to irredundant intersections, so let us recall the definition.
We say, that an intersection of sets $ \ca = \cap_{j\in J} C_j$ is irredundant if no $C_j$ can be omitted from this intersection: for all $J'\subset J$ and $J'\neq J$, we have $\cap_{j\in J'} C_j\neq \cap_{j\in J} C_j $.

\begin{rem}[$\cc_\rg$ as an irreduntant intersections]
    \label{rem:sature-Cg}
We now discuss  representations of $\cc_\rg$ as an irredundant intersection, which is given in Section~\ref{sec:saturated_constraints}. 
The difference of the slopes of $\Lambda'_\rg$ and $\Lambda_\rg$ is equal to:
\[
q-  \frac{(qp'-r')(1-\beta)}{\beta(1-p')} - \frac{r'}{p'} = \frac{(qp'-r')(\beta- p')}{\beta p'(1-p')},
\]
which implies that the constraint $\ch'_\rg$ (resp.\ $\ch_\rg$) is not saturated for $p'\in (0, \beta]$ (resp.\ for $p'\in [\beta, 1)$), that is:
\begin{equation}
    \label{eq:Cg-minimal}
\cc_\rg = 
\begin{cases}
    R_\rg \cap \ch_\rg \cap \ce_\rg & \text{if $p' \in (0, \beta]$,} \\
    R_\rg \cap \ch_\rg' \cap \ce_\rg & \text{if $p'\in [\beta, 1)$.}
\end{cases} 
\end{equation}
The constraint that $ R_\rg$ is essential, as we assume throughout this part that $A \in R_\rg$. Removing the constraint $\ce_\rg$ would imply that $A_0 \in \cc_\rg$, which is a contradiction. Furthermore, we observe that:
\[
\psi_\rg(\beta) - \Lambda_\rg(\beta) = \frac{(\beta-p')(2-3\beta)(qp'-r')}{p' ((1-\beta)+(1-p')(1-2\beta))},
\]
and:
\[
\psi_\rg(\beta) - \Lambda_\rg'(\beta) = \frac{(p'-\beta)(1-\beta)(p'q-r')}{(1-p')((1-\beta)+(1-p')(1-2\beta))}\cdot
\]
We have $\psi_\rg(\beta) > \Lambda_\rg(\beta)$ if and only if $\beta > p'$, and $\psi_\rg(\beta) > \Lambda_\rg'(\beta)$ if and only if $\beta < p'$. In particular, if $p' \in (0,\beta)$, then the point $(\beta, \psi_\rg(\beta))$ does not belong to $\cc_\rg$, so the constraint $\ch_\rg$ is necessary. Similarly, if $p' \in (\beta,1)$, then $(\beta, \psi_\rg(\beta)) \notin \cc_\rg$, and the constraint $\ch'_\rg$ is necessary. Thus, Equation~\eqref{eq:Cg-minimal} is irredundant if and only if $p' \neq \beta$.
For $p' = \beta$, the constraint $\ch_\rg$ and $\ch'_\rg$ are not saturated, and we get the following irredundant intersection:
\[
\cc_\rg =  R_\rg \cap \ce_\rg.
\]
\end{rem}

\subsubsection{The case  $A' \in D_q^- \setminus \Delta_q$ and $A=(p,r) \in R_\m$.} Notice that $q\in (0, 1)$. We consider the affine functions $y=\Delta_i^{A}(x)$, for $i=1,2$, which are the equation of the line going through respectively  $A_0$ and $A$, $A$ and $h_q(A_0)$:
\[
\Delta_1^{A}(x) = \frac{r}{p} x\quad \text{and}\quad
\Delta_2^{A}(x) = \frac{q-r}{1-p} (x-p) + r.
\]
We deduce from~\eqref{eq:C-=C+} that  $A=(p,r)\in D_q^-$ belongs to $C^{--}_{q,A'}$ if and only if $r' \geq \max_{i=1,2}\Delta_i^{A}(p')$. Elementary computation gives, that for $A=(p,r)\in D_q^-$ with $p\in I_q$:
\[
A = (p,r) \in C^{+-}_{q,A'} 
\iff
r' \geq \max_{i=1,2}\Delta_i^{A}(p')
\iff
r \leq \min_{i=\rg,\rd} \Lambda_i(p),
\]
with $\Lambda_\rg$ defined in~\eqref{eq:bernoulli_proof_cm_def_lambdag} and:
\begin{equation}
    \label{eq:bernoulli_proof_cm_def_lambdad}
    \Lambda_\rd(x) = q - (1-x) \frac{q-r'}{1-p'}\cdot
\end{equation}
By observing that:
\[
r\leq \Lambda_\rg(p) \iff A \in  R_\m \cap \ch_\rg; 
\quad \text{and} \quad
r\leq \Lambda_\rd(p) \iff A \in  R_\m \cap \ch_\rd, 
\]
we obtain that $C_{q, A'}^{+-}\cap R_\m=\cc_q$, with  $\cc_q =  R_\m \cap \ch_\rg \cap \ch_\rd$.

\begin{rem}[$\cc_q$ as an irreduntant intersections]
    \label{rem:sature-Cq}
We now discuss  representations of $\cc_q$ as an irredundant intersection, which is given in Section~\ref{sec:saturated_constraints}. 
The functions $\Lambda_\rg$ and $\Lambda_\rd$ are the lines passing through $A_0$ and $A'$, and through $A'$ and $h_q(A_0)$, respectively. In particular, we have:
\[
\Lambda_\rg(x) \geq \Lambda_\rd(x) \iff x \leq p' \quad \text{and}\quad
\Lambda_\rg(x) = \Lambda_\rd(x) \iff x = p',
\]
which implies also that:
\begin{equation}
    \label{eq:Cq-minimal}
    \cc_q = 
\begin{cases}
    R_\m \cap \ch_\rg & \text{if $p'\in (0, \beta]$,} \\
    R_\m \cap \ch_\rg \cap \ch_\rd & \text{if $p' \in (\beta, 1- \beta)$,} \\
    R_\m \cap \ch_\rd & \text{if $p'\in [1- \beta, 1)$.}
\end{cases}
\end{equation}
The constraint  $R_\m$ is essential, as we assume throughout this part that $A \in R_\m$. For $p' \in (0, \beta]$, removing the constraint $\ch_\rg$ would imply $\cc_q=R_\m$, which is false because $\Delta_q \cap R_\m \not \subset \cc_q$ since $qp' -r' >0$. Similarly, for $p' \in [1-\beta, 1)$, removing the constraint $\ch_\rd$ would imply $\cc_q=R_\m$, which is false because $\Delta_q \cap R_\m \not \subset \cc_q$ since $qp' -r' >0$. For $p' \in (\beta, 1-\beta)$, we have $\Lambda_\rg(\beta) > \Lambda_\rd(\beta)$ and $\Lambda_\rg(1-\beta) < \Lambda_\rd(1-\beta)$ then both $\ch_\rg$ and $\ch_\rd$ are necessary. Thus, the intersections in~\eqref{eq:Cq-minimal} are irredundant.
\end{rem}

\subsubsection{The case $A' \in D_q^- \setminus \Delta_q$ and $A=(p,r) \in R_\rd$ with $p \in [1- \beta, 1)$.} Notice that $q\in (0, 1)$. 
Let $A^*=(p^*,r^*)$ be defined as $A''$
in~\eqref{eq:bernoulli_cp''} with  $(p',r')$ replaced by $(p,r)$. 
We consider the affine functions $y=\Delta_i^{A}(x)$, for $i=1,2,3$, which are the equation of the line going through respectively  $A_0$ and $A^*$, $A^*$ and $A$, $A$ and $h_q(A_0)$:
\[
\Delta_1^{A}(x) = \frac{r^*}{p^*} x,\quad 
\Delta_2^{A}(x) = \frac{r^* - r}{p^* -p}(x-p)+r \quad \text{and} \quad
\Delta_3^{A}(x) = \frac{q - r}{1-p}(x-1)+q.
\]
We deduce from~\eqref{eq:C-=C+} that $A=(p,r)\in D_q^-$ belongs to $C^{+-}_{q,A'}$ if and only if $r' \geq \max_{i=1,2,3}\Delta_i^{A}(p')$. Elementary computation gives, that for $A=(p,r)\in D_q^-$ with $p \in [1-\beta, 1)$:
\[
A = (p,r) \in C^{+-}_{q,A'} 
\iff
r' \geq \max_{i=1,2,3}\Delta_i^{A}(p')
\iff
r \leq \min\left( \Lambda_\rd'(p), \psi_\rd(p), \Lambda_\rd(p) \right),
\]
with $\psi_\rd$ defined in~\eqref{eq:bernoulli_cm_psi},  $\Lambda_\rd$ defined in~\eqref{eq:bernoulli_proof_cm_def_lambdad}, and with $s$ defined   in~\eqref{eq:def-s}:
\[
\Lambda_\rd'(x) =\left( \frac{r'}{p'} + \frac{qp'-r'}{\beta p'}\right) ( x-1) + q =s(1-p',q-r')(x-1) + q.
\]
By observing that:
\[
r\leq \Lambda_\rd(p) \iff A \in  R_\rd \cap \ch_\rd; 
\quad
r\leq \Lambda_\rd'(p) \iff A \in  R_\rd \cap \ch_\rd'; 
\quad \text{and} \quad
r\leq \psi_\rd(p) \iff A \in  R_\rd \cap \ce_\rd,
\]
we obtain that  $C_{q, A'}^{+-}\cap R_\rd=\cc_\rd$, with $\cc_\rd=R_\rd \cap \ch_\rd \cap \ch'_\rd \cap \ce_\cd$.

\medskip 

\begin{rem}[$\cc_\rd$ as an irreduntant intersections]
    \label{rem:sature-Cd}
We now discuss  representations of $\cc_\rd$ as an irredundant intersection, which is given in Section~\ref{sec:saturated_constraints}. The difference of the slopes of 
  $\Lambda_\rd$ and $\Lambda'_\rd$ is equal to:
\[
\frac{q-r'}{1-p'}- \left( \frac{r'}{p'} + \frac{qp'-r'}{\beta p'}  \right) = \frac{(qp'-r')(p'- (1-\beta))}{\beta p'(1-p')},
\]
which implies that the constraint $\ch_\rd$ (resp.\ $\ch'_\rd$) is not saturated for $p'\in (0, 1-\beta]$ (resp.\ for $p'\in [1-\beta, 1)$), that is:
\begin{equation}
    \label{eq:Cd-minimal}
 \cc_\rd = 
\begin{cases}
    R_\rd \cap \ch_\rd' \cap \ce_\rd & \text{if $p' \in (0, 1-\beta]$,} \\
    R_\rd \cap \ch_\rd \cap \ce_\rg & \text{if $p' \in [1-\beta, 1)$.}
\end{cases} 
\end{equation}
The constraint  $ R_\rd$ is essential, as we assume throughout this part that $A \in R_\rd$. Removing the constraint $\ce_\rd$ would imply that $h_q(A_0) \in \cc_\rd$, which is a contradiction. Furthermore, we observe that:
\[
\psi_\rd(\beta) - \Lambda_\rd(\beta) = \frac{(p'-(1-\beta))(2-3\beta)(qp'-r')}{(1-p')((1-\beta) + p'(1-2\beta))},
\]
and:
\[
\psi_\rd(\beta) - \Lambda_\rd'(\beta) = \frac{((1-\beta)-p')(1-\beta)(qp'-r')}{p'((1-\beta) + p'(1-2\beta))} \cdot
\]
We have $\psi_\rd(\beta) > \Lambda_\rd(\beta)$ if and only if $p' > 1-\beta$, and $\psi_\rd(\beta) > \Lambda_\rd'(\beta)$ if and only if $1-\beta > p'$. In particular, if $p' \in (0,1-\beta)$, then the point $(1-\beta, \psi_\rd(1-\beta))$ does not belong to $\cc_\rd$, so the constraint $\ch_\rd'$ is necessary. Similarly, if $p' \in (1-\beta,1)$, then $(1-\beta, \psi_\rd(1-\beta)) \notin \cc_\rd$, and the constraint $\ch_\rg$ is necessary. Thus, Equation~\eqref{eq:Cd-minimal} is irredundant if and only if $p' \neq 1-\beta$.
For $p' = 1-\beta$, the constraint $\ch_\rd$ and $\ch'_\rd$ are not saturated, and we get the following irredundant intersection:
\[
\cc_\rd =  R_\rd \cap \ce_\rd.
\]
\end{rem}

\subsubsection{Conclusion}
We deduce from the previous sections that  if $A'\in D_q^-\setminus \Delta_q$, then 
$C^{+-}_{q, A'}= \cc_\rg \cup \cc_q\cup \cc_\rd$. 
Since 
$C^-_{q, A'}=C^{+-}_{q, A'} \cup h_q(C^{+-}_{q, A'})$, we deduce 
that Proposition~\ref{prop:bernoulli_cm}~\ref{it:C-3}  holds when furthermore $A'\in D_q^-$. Next, use that 
$C^{+-}_{q, A'}=C^{+-}_{q, h_q(A')}$ by symmetry to get 
that Proposition~\ref{prop:bernoulli_cm}~\ref{it:C-3}  holds in general.

\subsubsection{Saturated constraints}
\label{sec:saturated_constraints}

Let us identify the active constraints.  

\textbf{The general case $q \in (0,1)$.}
Let $A' \in D_q \setminus (\Delta_q \cup E_q)$. 
By comparing the slopes of the boundaries of $\ch_\rg$ and $\ch_\rg'$, as well as those of $\ch_\rd$ and $\ch_\rd'$, we obtain for $i\in \{\rg, \rd\}$, and $j$ such that $\{i, j\}=\{\rg, \rd\}$:
\[
\cc_i = 
\begin{cases}
    R_i \cap \ch_i \cap \ce_i & \text{if $A' \in R_i$,} \\
    R_i \cap \ch_i' \cap \ce_i & \text{if $A' \in R_\m \cup R_j$,}
\end{cases} 
\quad \text{and} \quad
\cc_q = 
\begin{cases}
    R_\m \cap \ch_i & \text{if $A' \in R_i$,} \\
    R_\m \cap \ch_\rg \cap \ch_\rd & \text{if $A' \in R_\m$ with $p'\in (\beta, 1- \beta)$.} 
\end{cases} 
\]  
We get from Remark~\ref{rem:sature-Cg} for $\cc_\rg$
from Remark~\ref{rem:sature-Cd} for $\cc_\rd$, and 
from Remark~\ref{rem:sature-Cq} for $\cc_q$
 that   all those  intersections are irredundant, unless $i=\rg$ and $p'=\beta$,   or  $i=\rd$ and $p'=1-\beta$. 
 We also have the irredundant representation  $\cc_\rg=R_\rg \cap \ce_\rg$ for  $p'=\beta$, and 
$\cc_\rd=R_\rd \cap \ce_\rd$
for $p'=1-\beta$. 

We deduce that for  $A'\not\in \Delta_q \cup E_q$
the boundary of $C_{q,A'}^+$ restricted to the interior of $D_q$ contains a segment with positive length unless $q=p'=1/2$. We shall use this result in the proof of Proposition~\ref{prop:Cp_subset_Fp}.

\textbf{The particular case $q = 1/2$.} 
When $q = 1/2$, we have  $\psi_\rg=\psi_\rd$ and this function is  linear, see~\eqref{eq:bernoulli_psi_1/2}, so the set  $\ce_\rg=\ce_\rd$ is an  half-plane, see Fig.~\ref{fig:q1/2_p_prime_ne_1/2}. 
Furthermore, if $p' = 1/2$, we have $\ch_i=\ch_i'$ for $i = \rg, \rd$ and  we get $\cc = D_q \cap  \ce_\rg = D_q \cap  \ce_\rd$, which are both irredundant intersections,  and $\cc$ is the convex hull of $3$ points in convex position. An instance  is provided in Fig.~\ref{fig:q_1/2_p_prime_1/2}, Fig.~\ref{fig:ber_cm_q=1/2} and Fig.~\ref{fig:ber_cm_q=1/2_p=1/2}.

We shall use in the proof of Proposition~\ref{prop:Cp_subset_Fp} that for $q=p'=1/2 $ and $A'\not\in \Delta_q \cup E_q$, the boundary of $C_{q,A'}^-$ restricted to the interior of $D_q$ contains a segment with positive length.

\subsection{Proof of Lemma~\ref{lem:Cm_boundaries}}
\label{sec:proof:Cm_boundaries}
Let $q\in (0,1)$ and $A' =(p',r') \in D_q \setminus (\Delta_q \cup E_q)$. We denote by $B_1$ and $B_2$ the two boundaries of $D_q$ belonging to $D_q^-$, that is:
\[
B_1 = \left\{ (p,r) \in D_q \colon r = 0 \right\}
\quad \text{and} \quad
B_2 = \left\{ (p,r) \in D_q \colon r = p+q-1 \right\}.
\]
Since $C_{q,A'}^+ = C_{q,h_q(A')}^+$, we may assume that $A' \in D_q^-$ without loss of generality. Suppose $A' \in \partial D_q$, we have, $r' = 0$ or $r' = p' + q - 1$. If $r' = 0$, then we have we have $A'\in B_1$ and:
\[
s(h_q(A')) = \begin{cases}
1 & \text{if $q \leq 1/2$,} \\
\frac{q}{1-q} & \text{if $q \geq 1/2$}.
\end{cases}
\]
Hence we get $\cc_\rg = R_\rg \cap \ch_\rg \subset B_1$, $\cc_q = R_\m \cap \ch_\rg \subset B_1$ and $\cc_\rd = R_\rd \cap \ch_\rd' \subset B_2$. Symilarly, if $r' = p' + q - 1$, then we have $A'\in B_2$ and:
\[
s(A') = \begin{cases}
2- \frac{1}{q} & \text{if $q \leq 1/2$,} \\
0 & \text{if $q \geq 1/2$}.
\end{cases}
\] 
Hence we get $\cc_\rg = R_\rg \cap \ch_\rg' \subset B_1$, $\cc_q = R_\m \cap \ch_\rd \subset B_2$ and $\cc_\rd = R_\rd \cap \ch_\rd \subset B_2$. In both cases, we get the following $
\cc_\rg \subset B_1$, $ \cc_q \subset B_1 \cup B_2$ and $\cc_\rd \subset B_2$.
This gives  that $C_{q,A'}^+$ is a union of segments belonging to $\partial D_q$. 

Conversely, since $C_{q,A'}^+$ is a union of segments contained in $\partial D_q$ and $A' \in C_{q,A'}^+$ we get $A' \in \partial D_q$.

\subsection{Proof of Proposition~\ref{prop:C+-cap-C-}}
\label{sec:C+-cap-C-}

Set $A=(p,r)$. First notice, as   $A'$ belongs to $\Delta_q$ if and only if $Z'=1$,  that if $A'\in \Delta_q$, then we have 
$A\in C^+_{q, A'}\cap C^-_{q, A'}$ if and only if $A\in \Delta_q$. 

We now assume that $A'$ and $A$ do not belong to $\Delta_q$, which also implies that $q$, $p'$ and $p$ belong to $(0, 1)$. Replacing $A'$ by $h_q(A')$ and/or $A$ by $h_q(A)$ if necessary, we can assume without loss of generality that  $r<qp$ and $r'< qp'$. By Remark~\ref{rem:stat-ordre}, we deduce that $\max(a_1, a_2)< 1 < \min (a_3, a_4)$ and  and that $a_1=a_2$ (resp. $a_3=a_4$) if and only if $p=1-q$ (resp. $p=q$), and similarly for the $a'_i$. 
The density ratios  $Z'$ and $Z$ have the same distribution (that is,  $A$ belongs to $C^+_{q, A'}\cap C^-_{q, A'}$)   if and only if  one of the following case holds:
\begin{enumerate}[(i)]
    \item $a_1\neq a_2$, $\rho_1=\rho'_1$ and $a_1=a'_1$, that is, $A=A'$ and $p'=p\neq 1-q$;
    \item $a_3\neq a_4$, $\rho_3=\rho'_3$ and $a_3=a'_3$, that is, $A=A'$ and $p'=p\neq q$;
    \item $a_1= a_2=a'_1=a'_2$ and $a_3=a_4=a'_3=a'_4$, that is, $A=A'$ and $q=p=p'=1/2$;
      \item $a_1=a'_2\neq a_2=a'_1$, $a_3=a'_4\neq  a_4=a'_3$, $\rho_1=\rho'_2$ and $\rho_3=\rho'_4$, that is, $q=1/2$ and 
      $A=A''$ if $1-p=p'>1/2$
      or $A= h_q(A'')$ if $1-p=p'<1/2$, 
\end{enumerate}
and, to complete the description, we mention that the case $a_1\neq a_2$, $\rho_1=\rho'_2$ and $a_3=a_4=a'_3=a'_4$, as well as the (symmetric) case
 $a_3\neq a_4$, $\rho_3=\rho'_4$ and $a_1=a_2=a'_1=a'_2$, are not possible.

\subsection{On the union $C_{q,A'}^+\cup C_{q,A'}^-$}
\label{sec:C+-cup-C-}

In the following proposition, we give the expression of $C_{q,A'}^+\cup C_{q,A'}^-$.

\begin{prop}[On $C_{q,A'}^+\cup C_{q,A'}^-$]
\label{prop:C+-cup-C-}
Let $q\in [0, 1]$ and $A'=(p',r')\in D_q$. We have:
\begin{equation}
  \label{eq:equiv-Dq}
  C^-_{q, A'}=\Delta_q
  \iff 
  C^+_{q, A'}=D_q
  \iff
  A'\in \Delta_q,
 \end{equation}
and this corresponds to  $Y$ and $X'$ being independent. 
If $q\in (0, 1)$ (and thus $E_q \not \subset \Delta_q$), then we have:
\begin{equation}
  \label{eq:equiv-Eq}
  C^-_{q, A'}=D_q
  \iff 
  C^+_{q, A'}= E_q
  \iff
  A'\in E_q,
\end{equation}
and this corresponds to  $X'$ equal a.s.\ to $Y$ or a.s.\ to $1-Y$. 
Furthermore, we have:
\begin{equation}
  \label{eq:equiv-=Dq}
  C^+_{q, A'} \cup C^-_{q, A'} = D_q
  \iff 
  A' \in \Delta_q  \cup E_q .
 \end{equation}
\end{prop}

\begin{proof}

Let $q \in [0,1]$ and let $A' = (p',r') \in D_q$. We simply write $C^\varepsilon$ for $C^\varepsilon_{q, A'}$.  If $q \in \{0,1\}$, then, thanks to~\eqref{eq:Sf=0iff}, 
we have $E_q\subset \Delta_q=D_q=C^+=C^-$ and the proof of~\eqref{eq:equiv-Dq} and~\eqref{eq:equiv-=Dq} are immediate. From now on, we assume that  $q \in (0,1)$. 

Equation~\eqref{eq:equiv-Dq} is a direct consequence of~\eqref{eq:Sf=0iff} and 
the fact that the Csiszár indices are nonnegative. 

We then  prove~\eqref{eq:equiv-Eq}. Recall the set $E_q$ has two distinct points and that $E_q=h_q(E_q)$. 
Since $D_q$ is the convex hull of $E_q\cup\{A_0, h_q(A_0)\}$
and that this set has four distinct points which are in convex position (as $q\not \in \{0, 1\}$), we deduce from Proposition~\ref{prop:bernoulli_cp} that the equality $C^-=D_q$ holds if and only if $A'\in E_q$. Then use Proposition~\ref{prop:C+-cap-C-}~\ref{it:C+capC-=2} to deduce that if $A'\in E_q$ and $C^-=D_q$ then $C^+=E_q $. Eventually, if $C^+=E_q $, since $A'$ belongs to $C^+$, we deduce that $A'\in E_q$. This conclude the proof of~\eqref{eq:equiv-Eq}.

The implication in~\eqref{eq:equiv-=Dq} is a consequence of~\eqref{eq:equiv-Dq} and~\eqref{eq:equiv-Eq}.  
We now prove the converse. 
Assume that $A'\not \in E_q$. Recall
that $q\in (0, 1)$ and thus $E_q$ is not a subset of $\Delta_q$. 
Using Proposition~\ref{prop:bernoulli_cp}, we deduce from the fact that the points in $\ca$ are in convex position, that  $E_q\subset D_q \setminus C^-$. 
Since $C^-$ is closed (as the convex hull of a finite number of points), 
$\Delta_q \subset C^-=\Cv(\ca)$ and $A_1\not \in C^-$, then 
we deduce that all the segments from $A_1$ to $A \in \Delta_q$ have a unique point belonging to the boundary of $C^-$ in $D_q$.
 Denote by $\cc$ the set of such points. Notice that 
 Lemma~\ref{lem:equivalence_convex_order_int_cdf} implies that $C^+$ is closed (see also Proposition~\ref{prop:bernoulli_cm} below). 
 We deduce that if  $C^+ \cup C^-=D_q$ then we have  $\cc\subset C^+$ and thus $\cc\subset C^+ \cap C^-$. Since $\cc$ is infinite, we deduce from 
Proposition~\ref{prop:C+-cap-C-} that $A'\in \Delta_q$. 
This concludes the proof of~\eqref{eq:equiv-=Dq}. 
\end{proof}

\subsection{Additional results on the Csiszár order and correlation}
\label{sec:Cp_subset_Fp}

In this section, we give additional results on the Csiszár order and correlation.

For $q \in [0,1]$ and $A' \in D_q$, we consider the regions corresponding the squared correlation order:
\begin{align*}
F_{q,A'}^+ = \left\{ A = (p,r) \in D_q \,\colon \, \corr^2(X,Y) \geq \corr^2(X',Y) \right\}, \\
F_{q,A'}^- = \left\{ A = (p,r) \in D_q \,\colon \, \corr^2(X,Y)  \leq  \corr^2(X',Y)  \right\}. 
\end{align*}
These sets are analogous to $C_{q,A'}^+$ and $C_{q,A'}^-$, with the Csiszár order replaced by the ordering induced by the squared correlations.

We recall that the set $C^+_{q,A'} \cap C^-_{q,A'}$ corresponds to distributions where $X \cseq X'$ w.r.t.\ $Y$ and is characterized in Proposition~\ref{prop:C+-cap-C-}. 
The next proposition is a direct consequence of the expression of the correlation in terms of the Csiszár index. Indeed, the squared correlation of $(X,Y)$, a pair of Bernoulli random variables, is a particular Csisz\'ar index associated to the Pearson $\chi^2$ divergence defined with $f_\mathrm{P}$: 
\begin{equation}
  \label{eq:SfP=corr2}
  S_{\mathrm{P}}(X,Y) = \corr^2(X,Y).
\end{equation}

\begin{prop}[Comparison of the uniform Csisz\'ar index order and the  squared correlation order]
    \label{prop:Cp_subset_Fp}
    Let $q \in [0,1]$ and $A'\in D_q$. Then, the Csisz\`ar order is stronger than the squared correlation order:
    \begin{equation}
      \label{eq:inclusion_C+_F+}
      C^+_{q,A'} \subset F^+_{q,A'} \quad \text{and} \quad C^-_{q,A'} \subset F^-_{q,A'};
    \end{equation}
     it is strictly stronger unless $|\corr(X', Y)|\in \{0, 1\}$:
       \begin{equation}
      \label{eq:equality_C+F+}
      C^+_{q,A'} = F^+_{q,A'} \iff  C^-_{q,A'} = F^-_{q,A'} \iff A' \in \Delta_q \cup E_q;
    \end{equation}
    and the two orders cannot be strictly inverse:
    \begin{equation}
      \label{eq:intersection_C+_F-}
      C^-_{q,A'} \cap F^+_{q,A'} = C^+_{q,A'} \cap F^-_{q,A'} = C^+_{q,A'} \cap C^-_{q,A'}.
    \end{equation}
 \end{prop}

\begin{proof}
Let $q\in [0, 1]$ and $A'\in D_q$ be fixed. Unless there is an ambiguity, for $\varepsilon\in \{-1, +1\}$, we shall write for simplicity $C^\varepsilon$ and $F^\varepsilon$ for $C_{q, A'}^\varepsilon$ and $F_{q, A'}^ \varepsilon$.  

We first prove the first inclusion in~\eqref{eq:inclusion_C+_F+}. Let  $A=(p,r) \in C^-$, that is, $S_f(X',Y) \geq S_f(X,Y)$ for all $f\in \F$. In particular, this holds for $f=f_\mathrm{P}$, and thus $\corr(X',Y)^2 \geq \corr(X,Y)^2$,  that is, $A\in F^-$. The proof of the second inclusion is similar; thus  Equation~\eqref{eq:inclusion_C+_F+} holds.

\medskip

We now prove~\eqref{eq:equality_C+F+}.  If $A'\in \Delta_q$ (resp.\ $A'\in E_q$), then we have $C^-= \Delta_q$ (resp.\ $C^-=D_q$) by Proposition~\ref{prop:C+-cup-C-} and $  F^-= \Delta_q$ (resp.\ $F^-=D_q$) thanks to~\eqref{eq:corr=01}. This gives $A'\in \Delta_q\cup E_q \implies C^-=F^-$, and using similar arguments, we have  $A'\in \Delta_q\cup E_q \implies C^+=F^+$. 
We prove the converse. 
  
Let $A'\not \in \Delta_q\cup E_q$. We define the functions
 $\psi_+$ and $\psi_-$ on $[0, 1]$ by:
  \[
  \psi_\varepsilon(p) = qp + \varepsilon \sqrt{p(1-p)} \, \frac{\vert qp' - r'\vert}{\sqrt{p'(1-p')}} \quad \text{for $\varepsilon \in \{-1, 1\}$}.
  \]
Notice that the function $\psi_+$ is strictly concave and the function $\psi_-$ is strictly convex, and $\psi_\varepsilon(0)=0$,  $\psi_\varepsilon(1)=q$ and $|\psi'_\varepsilon(0)|=|\psi'_\varepsilon(1)|=+\infty$. 
We have:
 \begin{equation}
     \label{eq:rep.F+}
  F^- = \left\{ (p,r) \in D_q\,  \colon \, \psi_+(p) \geq r \geq \psi_-(p) \right\},
 \end{equation}
and denoting by $\cb$ the common boundary of $F^+$ and $F^-$, we get:
\[
\cb=F^+ \cap F^-=\partial F^+ \cap \partial F^-=\bigcup_{\varepsilon\in \{-1, +1\}}\{(p,r)\in D_q\, \colon\, r=\psi_\varepsilon(p)\quad\text{and}\quad p\in (0, 1)\}.
\]
 Using that $|\corr(X',Y)|<1$, we deduce that $E_q$ 
belongs to the complementary of $F^-$, which implies that $\cb$ is not  empty; as $|\psi'_\varepsilon(0)|=|\psi'_\varepsilon(1)|=+\infty$, we get that $\cb $ is closed; it is also easy to see that $\cb$ is the union of two disjoint closed (connected) paths not reduced to a point. 
Furthermore, thanks to the concavity property of $\varepsilon \psi_\varepsilon$, we deduce that $\cb$ contains no  segment with positive length. 
 
Thanks to the properties of $\cb$, the boundary of  $F^-$,  and thus also  of $F^+$, restricted to the interior of $D_q$  does not contain segment of positive length.
According to Proposition~\ref{prop:bernoulli_cp}, the set $C^-\neq D_q$ is a convex hull of $4$ or $6$ points in convex order, and thus the boundary of $C^-$ restricted to the interior of $D_q$ contains a segment with positive length. 
Similarly, according to Proposition~\ref{prop:bernoulli_cm}
and Section~\ref{sec:saturated_constraints},
the boundary of $C^+$ restricted to the interior of $D_q$ contains a segment with positive length. 
We deduce that  $C^-$ is not equal to $F^-$
and that  $C^+$ is not equal to $F^+$, that is, 
$A'\not \in \Delta_q\cup E_q \implies C^-\neq F^- $ and $C^+\neq F^+$ .
This concludes the proof of~\eqref{eq:equality_C+F+}.

\medskip
  
We now prove~\eqref{eq:intersection_C+_F-}.
Notice that~\eqref{eq:intersection_C+_F-} is obvious for $A'\in \Delta_q\cup E_q$ by~\eqref{eq:equality_C+F+}. We now assume that 
$A'\not \in \Delta_q\cup E_q$.
Notice that $C^-\cap F^+=C^-\cap \cb= \partial C^- \cap \cb$. 
Recall that $C^-$ is the convex hull of 4 or 6 points, say $\ca=\ca_0\cup h_q(\ca_0)$, which are in convex order. 
Using that $\cb$ is the union of  (the restriction of ) the graph of $\psi_+$, which is strictly concave, and of $\psi_-$, which is strictly convex, and that $C^-$ is below the graph of $\psi_+$ and above the graph of $\psi_-$, we deduce that   $C^-\cap \cb \subset \ca$.  By definition of $\cb$, we have that $A_0=(0,0)$ and $h_q(A_0)=(1, q)$ do not belong to $\cb$.  
By construction, we have that $A'$, as well as $h_q(A')$, belong to 
$C^-\cap \cb$. 
So we only have to consider the case where $\ca$ has $6$ points, that is, $p'\not \in I_q$, and check  whether $A''=(p'', r'')$ (and thus also $h_q(A'')$)
from Proposition~\ref{prop:bernoulli_cp}
belongs to $\cb$ or not. The former case corresponds to $r''\in \{\psi_+(p''), \psi_-(p'')\}$, where $(p'', r'')$ is given in~\eqref{eq:bernoulli_cp''}.

We first consider the case 
$q \in [1/2,1)$,  $r' <q p'$ and  $p' < 1-q$. This gives:
\[
p'' = \frac{p'(1-p')}{2(1-q)-p'} \quad \text{and} \quad r'' = p'' q + p'' (qp'- r') \frac{1-q}{qp} \cdot
\]
Since $p'' q > r''$, we have $r'' > \psi_-(p'')$. Elementary computations gives that 
$r''= \psi_+(p'') $ if and only if $q=1/2$. 
This gives that $A''$ belongs to $\cb$ 
if and only if 
$q=1/2$, and thus $A''$ belongs to $C^-\cap F^+$ if and only if 
$q=1/2$. We deduce from Proposition~\ref{prop:C+-cap-C-} that, in this case, $C^+ \cap C^-$ is equal to $\{A', h_q(A')\}$ if $q>1/2$ and to $\{A', A'', h_q(A'), h_q(A'')\}$ if $q=1/2$. In conclusion, we get that
 $C^+ \cap C^-=F^+ \cap C^-$ if $q\in [1/2, 1)$, $r'<qp'$ and $p'<1-q$. 
The proof is similar for $p'>q$, and 
using symmetries arguments as in Section~\ref{sec:proof_bernoulli_sym_expression_finale} on the proof of Proposition~\ref{prop:bernoulli_cp}, we also get the result for $q\in (0, 1/2]$ 
and $p'\not \in I_q$. In conclusion we have that 
 $C^+ \cap C^-=F^+\cap C^-$ for $A'\not \in \Delta_q \cup E_q$, and thanks to the argument above, we deduce that for $A'\in D_q$:
 \begin{equation}
     \label{eq:C+F-_proof}
     F^+_{q,A'} \cap C^-_{q,A'}=C^+_{q,A'} \cap C^-_{q,A'} 
 \end{equation}
 
To conclude, we need to check the first equality in~\eqref{eq:intersection_C+_F-}. This equality is obvious as, using~\eqref{eq:C+F-_proof} with $A'$ replaced by $A$ for the second equivalence, we get for $A\in D_q$:
\[
A \in C^+_{q,A'} \cap F^-_{q,A'} \iff A' \in C^-_{q,A} \cap F^+_{q,A} \iff A' \in C^+_{q,A} \cap C^-_{q,A} \iff A \in C^-_{q,A'} \cap C^+_{q,A'}.
\]

\end{proof}

\section{Proofs of the results in Section~\ref{sec:ThresholdedGaussian}}

In Section~\ref{sec:Thresholded_Model_Results}, we introduce the model under study and the results needed to prove Theorem~\ref{thm:thresholded}. Some technical results are then collected in Section~\ref{sec:Tecnical_results_copula_thresholded}. We prove Proposition~\ref{prop:copule_thresholded} and Corollary~\ref{cor:reverse_implication_thresholded} in Sections~\ref{sec:proof_copule_thresholded} and~\ref{sec:proof_reverse_implication_thresholded}, respectively, and finally we prove Theorem~\ref{thm:thresholded} in Section~\ref{sec:proof_thresholded}.

\subsection{Model and results}
\label{sec:Thresholded_Model_Results}

Let $(X, Y)$ be a centered 2-dimensional Gaussian vector with unit variances and correlation between $X$ and $Y$ given by  $\rho\in [-1, 1]$. 
Let $\Phi$ denotes the cdf of the real standard Gaussian 
distribution, $\Phi_\rho$  the cdf of $(X,Y)$, and $C_\rho$ the corresponding copula, that is, for $\rho\in [-1, 1]$ and $p,q\in [0, 1]$:
\[
C_\rho(p,q) = \P(\Phi(X) \leq p, \Phi(Y) \leq q) = \Phi_\rho \left(\Phi^{-1}(p), \Phi^{-1}(q)\right).
\]
We recover the lower Fréchet-Hoeffding copula for $\rho=-1$, the independence copula for $\rho = 0$, and the upper Fréchet-Hoeffding copula for  $\rho = 1$:
\begin{equation}
  \label{eq:GaussianCopulaRho-101}
  C_{-1}(p,q) = (p+q-1)_+, \quad C_0(p,q) = pq \quad \text{and}\quad C_1(p,q) = \min(p, q).
\end{equation}

Let $\phi=\partial_x \Phi$ denote the density of the standard Gaussian distribution and for $\rho\in (-1, 1)$ let  $\phi_\rho= \partial^2_{x,y} \Phi_\rho$ denote the density of the distribution of $(X, Y)$.
It is elementary to check that for $\rho\in (-1, 1)$ and $q \in (0, 1)$ the function $p \mapsto C_\rho(p,q)$ is continuously differentiable on $[0,1]$ and twice continuously differentiable on $(0,1)$. Moreover, for $\rho\in (-1, 1)$, $q \in (0, 1)$ and $p \in [0,1]$, we have: 
\begin{equation}
  \label{eq:deriv_p_copule_gauss}
  \partial_p C_\rho(p,q) = \Phi\left(\frac{\Phi^{-1}(q) - \rho  \Phi^{-1}(p)}{\sqrt{1- \rho^2}}\right),
\end{equation}
Since $\partial_p C_\rho(p,q)$ is positive for $p \in (0,1)$, we get that the function $p \mapsto C_\rho(p,q)$ is increasing on $[0,1]$ for $\rho\in (-1, 1)$ and $q\in (0, 1)$. 
Similarly, for $p,q \in [0,1]$, the function $\rho \mapsto C_{\rho}(p,q)$ is continuous on $[-1,1]$ and for $\rho \in (-1,1)$ and $p,q \in [0,1]$, we have:
\begin{equation}
  \label{eq:deriv_rho_copule_gauss}
  \partial_\rho C_\rho(p,q) = \phi_\rho(\Phi^{-1}(p), \Phi^{-1}(q)) .
\end{equation}
Moreover, for $p,q \in (0,1)$ and $\rho \in (-1,1)$, the quantity $\partial_{\rho} C_{\rho}(p,q)$ is positive. Hence, the function $\rho \mapsto C_{\rho}(p,q)$ is increasing on $[-1,1]$ for all $p,q \in (0,1)$.

Let $(\alpha,\beta) \in [-\infty, \infty]^2$. We consider the Bernoulli random variables given by the thresholded Gaussian random variable $(X,Y)$:
\begin{equation}
  \label{eq:DefUV}
  U= \indicatrice{X \leq \alpha} \quad \text{and} \quad V = \indicatrice{Y \leq \beta},
\end{equation}
with  the associated parameter $(p,  q, r)$:
\[
p = \P(U = 1) = \Phi(\alpha), \quad \quad q = \P(V = 1)=\Phi(\beta)
\quad\text{and}\quad 
r = \P(U = 1, V = 1) =\Phi_\rho(\alpha, \beta)=C_\rho(q,p),
\]
with $q \in [0,1]$ and $A = (p,r) \in D_q$ defined in~\eqref{eq:def-Dq}.
Recall that $A_0 = (0,0)$, $h_q(p,r) = (1-p, q-r)$, $D_q^\varepsilon = \{(p,r)\in D_q\, \colon\, \varepsilon (r-qp)\geq 0\}$ for $\varepsilon \in \{-,+\}$ and $\Delta_q = D_q^- \cap D_q^+$, see Section~\ref{sec:Bernoulli}. We shall be interested in the possible values of such $(p,r)$ as $p$ runs through the set $[0,1]$ and $(0,1)$:
\[
\Upsilon_{q, \rho}=\{A=(p,r)\, \colon\, p\in [0, 1]
\quad\text{and}\quad r=C_\rho(q,p)\}
\quad \text{and} \quad
\Upsilon_{q, \rho}^* = \Upsilon_{q, \rho} \setminus \{A_0, h_q(A_0)\}.
\]
Notice that  $\Upsilon_{q, \rho}\subset D_q$,
$\Upsilon_{q, 0}=\Delta_q$ (case $\rho=0$), $\Upsilon_{q, 1}=\partial D_q \cap D_q^+$ (case $\rho=1$), and that $h_q(\Upsilon_{q, \rho})=\Upsilon_{q,-\rho}$. %, where $h_q$ is defined by~\eqref{eq:def-hq}. 
For $B, B' \subset \R^2$, we write $B \leq B'$ (resp.\ $B < B'$) if for all $x,y,y' \in \R$ such that $(x,y) \in B$ and $(x,y') \in B'$, we have $y \leq y'$ (resp.\ $y < y'$).
We deduce from the monotonicity of the map $\rho \mapsto C_\rho(p,q)$ that:
\begin{equation}
  \label{eq:monotony_upsilon_rho}
  \Delta_q < \Upsilon_{q, \rho}^* < \Upsilon_{q, \rho'}^* \quad \text{for $0 < \rho < \rho' \leq 1$}.
\end{equation}
In particular, we deduce that $U$ and $V$ are independent if and only if $X$ and $Y$ are independent or $\alpha\in \{-\infty, + \infty\}$ or $\beta\in \{-\infty, + \infty\}$.

Let $X'$ be a standard Gaussian random variable and $\alpha'\in [-\infty, \infty]$. Assume that $(X',Y)$ is a Gaussian vector and define $\rho'$, $U'$,  $p'$ and $r'$ as $\rho$, $U$, $p$ and $r$ with $X$ and $\alpha$ replaced by $X'$ and $\alpha'$. 
Recall that $A''$ is defined in~\eqref{eq:bernoulli_cp''} when $p' \notin I_q$
For $q = 1/2$ and $p' \neq q$, we also define: 
\begin{equation}
  \label{eq:def_A''+}
  A''_+ = (1-p', r''_+) \quad \text{and} \quad r''_+ = C_\rho(1-p', q),
\end{equation}
which is equal to $A''$ if $p' > 1/2$ or to $h_q(A'')$ if $p' < 1/2$, with $A''$ defined in~\eqref{eq:bernoulli_cp''}. Notice that $A''_+ \in D_q^+$. The proof of the following proposition is given in Section~\ref{sec:proof_copule_thresholded}.

\begin{figure}[htbp]
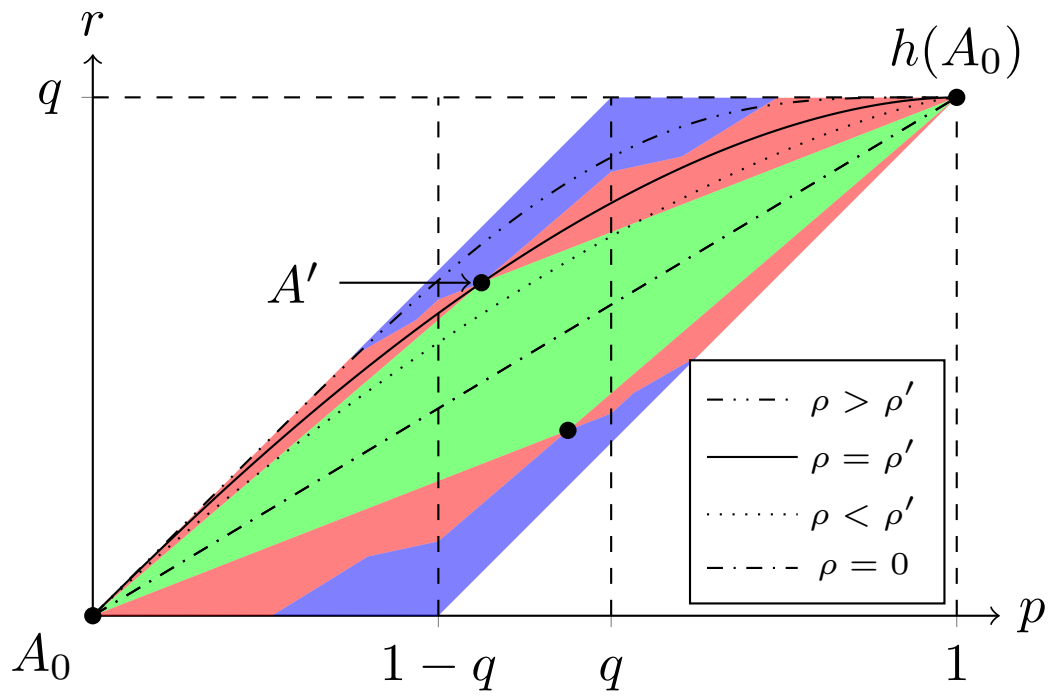
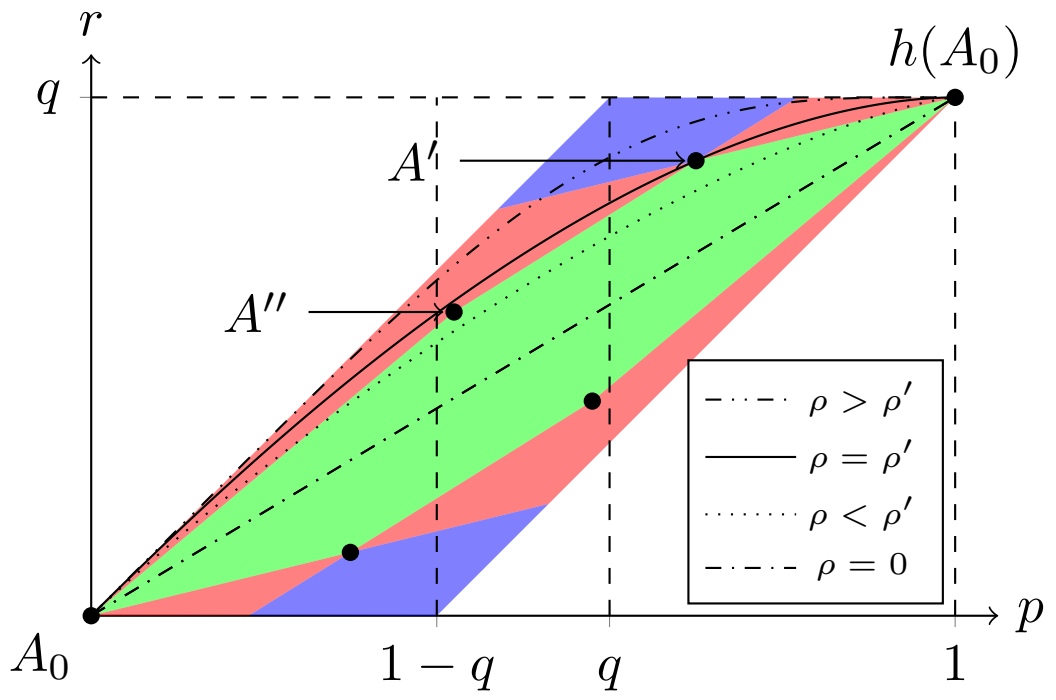
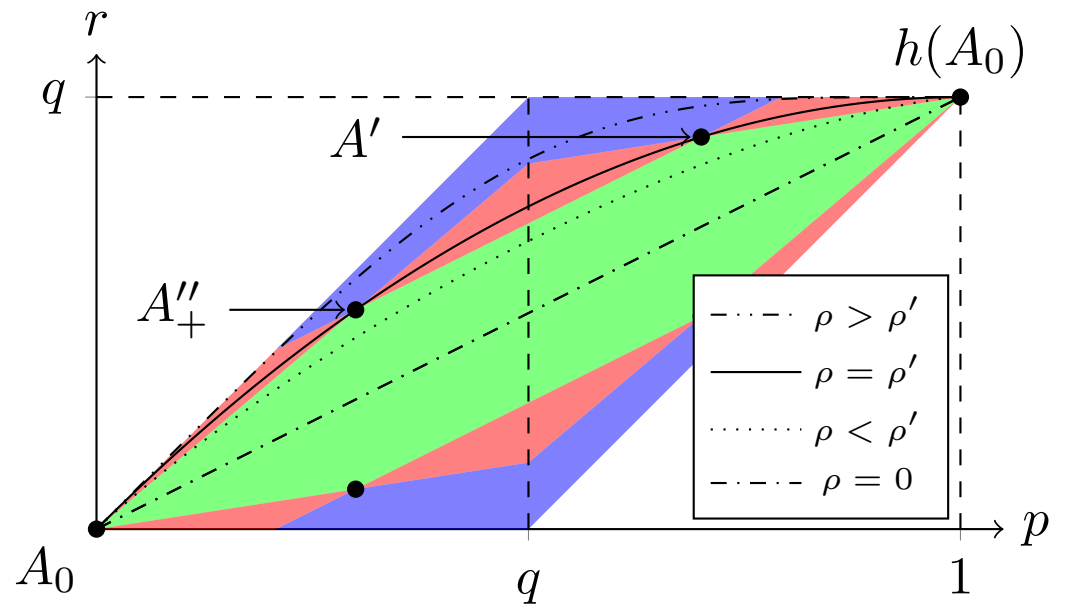

  \begin{subfigure}[c]{0.49\textwidth}
  \centering
  \begin{RegionPlotBernoulli}{0.45}{0.38544160967578955}{0.6}
    {
      \pgfplotsset{legend style={font=\tiny, at={(0.94,0.02)}, anchor=south east}}

      \addplot[black, dashdotdotted] table {Plots/GaussianCopula/copula_q0.6_rho0.9.dat};
      \addlegendentry{$\rho > \rho'$}
      
      \addplot[black, solid] table {Plots/GaussianCopula/copula_q0.6_rho0.7.dat};
      \addlegendentry{$\rho = \rho'$}

      \addplot[black, dotted] table {Plots/GaussianCopula/copula_q0.6_rho0.5.dat};
      \addlegendentry{$\rho < \rho'$}

      \addplot[black, dashdotted] coordinates {(0,0) (1,0.6)};
      \addlegendentry{$\rho = 0$}

      \fill (0.45,0.38544160967578955) circle (0.01) node[below] {}; % {$A'$};
      \fill (0.55,0.21455839032421042) circle (0.01) node[below] {}; %{$h(A')$};

      \draw[\ArrowStyle, thin, shorten >=2pt] (0.38544160967578955-0.1,0.38544160967578955) -- (0.45,0.38544160967578955);
      \fill (0.38544160967578955-0.1,0.38544160967578955) circle (0) node[left] {$A'$};
    }
    {
    \path[fill=\ColorCm, opacity=\OpacityPlot]
    (0,0) -- (0.45,0.38544160967578955) -- (1,0.6) -- cycle;
    
    \path[fill=\ColorOut, opacity=\OpacityPlot]
    plot[domain=0.681970514902493:0.792658710826822, smooth, samples = \nsamples](\x,{0.6*\x + 0.138529931610947*\x/(\x + 0.09)})
    -- (1,0.6)
    -- (0.45, 0.38544160967578955)
    -- (0.6, 0.513922146234386)
    -- cycle;

    \path[fill=\ColorOut, opacity=\OpacityPlot]
    plot[domain=0.299331826888419:0.374219511452964, smooth, samples = \nsamples](\x,{0.6*\x - (0.138529931610947*\x - 0.138529931610947)/(1.11 - \x)})
    -- (0.4, 0.365936301464498)
    -- (0.45, 0.38544160967578955)
    -- (0,0)
    -- cycle;

    \path[fill=\ColorCp, opacity=\OpacityPlot]
    plot[domain=0.299331826888419:0.374219511452964, smooth, samples = \nsamples](\x,{0.6*\x - (0.138529931610947*\x - 0.138529931610947)/(1.11 - \x)})
    -- (0.4, 0.365936301464498)
    -- (0.45,0.38544160967578955)
    -- (0.6, 0.513922146234386)
    -- plot[domain=0.681970514902493:0.792658710826822, smooth, samples = \nsamples](\x,{0.6*\x + 0.138529931610947*\x/(\x + 0.09)})
    -- (0.6,0.6)
    -- cycle;
    }
  \end{RegionPlotBernoulli}
  \caption{Case: $q = 0.6 $ and $p'=0.45 \in I_q$.}
\label{fig:upsilon_p0.45_q0.6_rho0.7}
\end{subfigure}
\hfill
  \begin{subfigure}[c]{0.49\textwidth}
  \centering
  \begin{RegionPlotBernoulli}{0.7}{0.526669584339166}{0.6}
    {
    \pgfplotsset{legend style={font=\tiny, at={(0.94,0.02)}, anchor=south east}}

    \addplot[black, dashdotdotted] table {Plots/GaussianCopula/copula_q0.6_rho0.9.dat};
    \addlegendentry{$\rho > \rho'$}
    
    \addplot[black, solid] table {Plots/GaussianCopula/copula_q0.6_rho0.7.dat};
    \addlegendentry{$\rho = \rho'$}

    \addplot[black, dotted] table {Plots/GaussianCopula/copula_q0.6_rho0.5.dat};
    \addlegendentry{$\rho < \rho'$}

    \addplot[black, dashdotted] coordinates {(0,0) (1,0.6)};
    \addlegendentry{$\rho = 0$}

    \fill (0.7,0.526669584339166) circle (0.01) node[below] {}; % {$A'$};
    \fill (0.42,0.351558278716555) circle (0.01) node[below] {}; % {$A''$};

    \fill (0.3,0.07333041566083398) circle (0.01) node[above] {}; %{$h(A')$};
    \fill (0.58,0.248441721283445) circle (0.01) node[above] {}; %{$h(A'')$};

    \draw[\ArrowStyle, thin, shorten >=2pt] (0.526669584339166-0.1,0.526669584339166) -- (0.7,0.526669584339166);
    \fill (0.526669584339166-0.1,0.526669584339166) circle (0) node[left] {$A'$};

    \draw[\ArrowStyle, thin, shorten >=2pt] (0.351558278716555 - 0.1,0.351558278716555) -- (0.42,0.351558278716555);
    \fill (0.351558278716555 - 0.1,0.351558278716555) circle (0) node[left] {$A''$};
    }
    {
    \path[fill=\ColorCm, opacity=\OpacityPlot]
    (0,0) -- (0.42,0.351558278716555) -- (0.7, 0.526669584339166) -- (1,0.6) -- cycle;
    
    \path[fill=\ColorOut, opacity=\OpacityPlot]
    (0,0)
    -- (0.42, 0.351558278716555)
    -- (0.7, 0.526669584339166)
    -- (0.470595049839400, 0.470595049839400)
    -- cycle;

    \path[fill=\ColorOut, opacity=\OpacityPlot]
    plot[domain=0.7:0.817842860086056, smooth, samples = \nsamples](\x,{0.6*\x + 0.128003501206999*\x/(\x + 0.14)})
    -- (1, 0.6)
    -- cycle;

    \path[fill=\ColorCp, opacity=\OpacityPlot]
    plot[domain=0.7:0.817842860086056, smooth, samples = \nsamples](\x,{0.6*\x + 0.128003501206999*\x/(\x + 0.14)})
    -- (0.6, 0.6)
    -- (0.470595049839400, 0.470595049839400)
    -- cycle;
    }
  \end{RegionPlotBernoulli}
  \caption{Case: $q = 0.6$ and $p'=0.7 \notin I_q$.}
\label{fig:upsilon_p0.7_q0.6_rho0.7}
\end{subfigure}
\vspace{0.5cm}
\begin{subfigure}[c]{0.49\textwidth}
  \centering
  \begin{RegionPlotBernoulli}{0.7}{0.45378773084189233}{0.5}
    {
      \pgfplotsset{legend style={font=\tiny, at={(0.94,0.02)}, anchor=south east}}

      \addplot[black, dashdotdotted] table {Plots/GaussianCopula/copula_q0.5_rho0.9.dat};
      \addlegendentry{$\rho > \rho'$}
      
      \addplot[black, solid] table {Plots/GaussianCopula/copula_q0.5_rho0.7.dat};
      \addlegendentry{$\rho = \rho'$}

      \addplot[black, dotted] table {Plots/GaussianCopula/copula_q0.5_rho0.5.dat};
      \addlegendentry{$\rho < \rho'$}

      \addplot[black, dashdotted] coordinates {(0,0) (1,0.5)};
      \addlegendentry{$\rho = 0$}

      \fill (0.7,0.45378773084189233) circle (0.01) node[below] {}; % {$A'$};
      \fill (0.3,1/2-0.45378773084189233) circle (0.01) node[above] {}; %{$h(A')$};

      \fill (0.3,0.2537877308418924) circle (0.01) node[above] {}; %{$A''$};
      \fill (0.7,1/2-0.2537877308418924) circle (0.01) node[above] {}; %{$h(A'')$};

      \draw[\ArrowStyle, thin, shorten >=2pt] (0.2537877308418924-0.1,0.2537877308418924) -- (0.3,0.2537877308418924);
      \fill (0.2537877308418924-0.1,0.2537877308418924) circle (0) node[left] {$A''_+$};

      \draw[\ArrowStyle, thin, shorten >=2pt] (0.45378773084189233-0.1,0.45378773084189233) -- (0.7,0.45378773084189233);
      \fill (0.45378773084189233-0.1,0.45378773084189233) circle (0) node[left] {$A'$};
    }
    {
    \path[fill=\ColorCm, opacity=\OpacityPlot]
    (0,0) -- (0.3,0.2537877308418924) -- (0.7,0.45378773084189233) -- (1,0.5) -- cycle;
    
    \path[fill=\ColorOut, opacity=\OpacityPlot]
    (0,0) -- (0.3,0.2537877308418924) -- (0.2075754616837847,0.2075754616837847) -- cycle;

    \path[fill=\ColorOut, opacity=\OpacityPlot]
    (0.3,0.2537877308418924) -- (0.7,0.45378773084189233) -- (0.5,0.4229795514031539) -- cycle;

    \path[fill=\ColorOut, opacity=\OpacityPlot]
    (0.7924245383162153,0.5) -- (0.7,0.45378773084189233) -- (1,1/2) -- cycle;
    
    \path[fill=\ColorCp, opacity=\OpacityPlot]
    (0.2075754616837847,0.2075754616837847) -- (0.3,0.2537877308418924) -- (0.5,0.4229795514031539) -- (0.7,0.45378773084189233) -- (0.7924245383162153,0.5) -- (0.5,0.5) -- cycle;  
    }
  \end{RegionPlotBernoulli}
\caption{Case: $q=1/2$ and $p'=0.7 \neq q$.}
\label{fig:upsilon_p0.7_q0.5_rho0.7}
\end{subfigure}
\caption{Plot of the curves $\Upsilon_{q,\rho}$ for $\rho \in \{0, 0.5,\rho' ,0.9\}$, with $\rho' = 0.7$.}
\label{fig:upsilon_rho0.7}
\end{figure}

\begin{prop}[No inversion of the uniform Csiszár order]
\label{prop:copule_thresholded}
Let $q\in (0, 1)$ (then $\mathring{D}_q^+ \neq \emptyset$) and $\rho' \in (0, 1)$. Let $A'=(p', r')$ with $p'\in (0, 1)$ and $r'=C_\rho(q, p')$, then $A' \in \mathring{D}_q^+$. We have:
\begin{enumerate}[(i)]
  \item \label{it:Up-0} If $\rho=0$, then $\Upsilon_{q,\rho}= \Delta_q \subset  C^+_{q,A'}$.
  \item \label{it:Up-0rhop} If $\rho\in (0, \rho')$, then  $\Upsilon_{q,\rho}^* \cap C^+_{q,A'}\neq \emptyset$, $\Upsilon_{q,\rho}^* \cap (C^+_{q,A'})^c\neq \emptyset$ and  $\Upsilon_{q,\rho} \cap C^-_{q,A'}= \emptyset$.
  \item \label{it:Up-q} If $\rho=\rho'$ and $q\neq 1/2$ or $q=p'=1/2$, then:
  \[
    \Upsilon_{q,\rho}^* \cap C^+_{q,A'} = \Upsilon_{q,\rho} \cap C^-_{q,A'} = \{A'\}.
  \]
  \item \label{it:Up-1/2} If $\rho=\rho'$, $q=1/2$ and $p'\neq q$, then with $A''_+$ defined in~\eqref{eq:def_A''+}, we have:
  % , with $A''_+=(1-p', r'')$ and $r''_+=C_\rho(q,1-p')$ (and $A''_+$ is also equal to $A''$ if $p'>1/2$ or to $h_q(A'')$ if $p'<1/2$):
  \[
    \Upsilon_{q,\rho}^* \cap C^+_{q,A'} = \Upsilon_{q,\rho} \cap C^-_{q,A'} = \{A', A''_+\}.
  \]
  \item \label{it:Up-rhop1} If $\rho\in (\rho', 1)$, then  $\Upsilon_{q,\rho}^* \cap C^+_{q,A'}= \emptyset$, $\Upsilon_{q,\rho}^* \cap C^-_{q,A'}\neq \emptyset$ and $\Upsilon_{q,\rho}^* \cap (C^-_{q,A'})^c\neq \emptyset$.
  \item \label{it:Up-1} If $\rho=1$, then $\Upsilon_{q,\rho}= \partial D_q \cap D_q^+$.
\end{enumerate}
\end{prop}
Using the symmetry $h_q$, we can trivially extend this results to all values of  $\rho'\in (-1, 1)\setminus\{0\}$ and $\rho\in [-1, 1]$. If $q \in \{0, 1\}$ or $p' \in \{0, 1\}$ or $\rho' = 0$, then $\Upsilon_{q,\rho'}$ is reduced to $\Delta_q$ (and $U'$ and $V$ are independent). 
Morevoer, if $\vert \rho' \vert = 1$ (that is $|X|=|Y|$ a.s.), then $\Upsilon_{q,\rho'} \subset \partial D_q$.  The Points~\ref{it:Up-q} and~\ref{it:Up-1/2} are to be related to the intersection between $C^+_{q, A'}$ and $C^-_{q, A'}$ given in Proposition~\ref{prop:C+-cap-C-}. Several instances of $\Upsilon_{q,\rho}$ are represented in Fig.~\ref{fig:upsilon_rho0.7}: Figs.~\ref{fig:upsilon_p0.45_q0.6_rho0.7} and~\ref{fig:upsilon_p0.7_q0.6_rho0.7} illustrate the Point~\ref{it:Up-q} of the previous proposition, while Fig.~\ref{fig:upsilon_p0.7_q0.5_rho0.7} illustrates the Point~\ref{it:Up-1/2}. 
From this proposition, we deduce the Corollary~\ref{cor:reverse_implication_thresholded} which gives the implication of the Theorem~\ref{thm:thresholded}. Its proof is given in Section~\ref{sec:proof_reverse_implication_thresholded}.

\begin{cor}[Reverse implication for thresholded Gaussian random variables]
\label{cor:reverse_implication_thresholded}
For $(X,X',Y) \in \R^3$ a Gaussian random vector, we have:
% such that $\Var(X)\Var(X')\Var(Y) > 0$
\begin{equation}
    \label{eq:corr_reverseimp_eq}
    \exists \, \alpha, \alpha', \beta \in \R \quad\text{such that}\quad \ind_{\{X \leq \alpha\}} \cseq \ind_{\{X' \leq \alpha'\}} \quad \text{w.r.t.}\quad  \ind_{\{Y \leq \beta\}} \implies X \cseq X' \quad \text{w.r.t.}\quad Y,
\end{equation}
\begin{equation}
    \label{eq:corr_reverseimp_ineg}
    \exists \, \alpha, \alpha', \beta \in \R \quad\text{such that}\quad \ind_{\{X \leq \alpha\}} \csls \ind_{\{X' \leq \alpha'\}} \quad \text{w.r.t.}\quad  \ind_{\{Y \leq \beta\}} \implies X \csls X' \quad \text{w.r.t.}\quad Y.
\end{equation}
\end{cor}

\subsection{Technical results}
\label{sec:Tecnical_results_copula_thresholded}
The proof of Proposition~\ref{prop:copule_thresholded} relies on the following two technical lemmas.
Recall the following convention on the sign function:  for $x\in \R$, $\sgn(x)=0$ if and only if $x=0$. 

\begin{lem} 
  \label{lem:copulegaussienne_geq_r-second}
  Let  $\mathrm{f}$ be the function defined on $\ddom_0 = \{(p,q,\rho) \in [0,1]^2 \times[-1,1] \, \colon \,q \in [1/2, 1), p \in [q, 1] \}$ by:
  \begin{equation}
  \label{eq:difference_copule_rsecond}
    \rmf(p,q,\rho) = (1-p) \left(C_\rho(p',q) - p'q \right)  - \frac{p'(1-q)}{q} \left(C_\rho(p,q) - pq \right), \quad \text{with}\quad p' = \frac{p(1-p)}{p-2q+1} \cdot
  \end{equation}
  % We have $\rmf(p,q,\rho) = 0$ if $q = 1/2$ or $p \in \{q, 1\}$ or $\rho = 1$; and $\sgn(\rmf(p,q,\rho)) = \sgn(\rho)$ otherwise.
     We have $\sgn(\rmf(p,q,\rho)) = \sgn(\rho)$, unless
     $q = 1/2$ or $p \in \{q, 1\}$ or $\rho = 1$ for which      
     $\rmf(p,q,\rho) = 0$:
  \[
  \sgn(\rmf(p,q,\rho)) = \sgn(\rho) \, \ind_{\{ q \neq 1/2, \, p\in (q, 1), \, \rho<1\}}.
    \]   
\end{lem}

\begin{proof}[Proof of Lemma~\ref{lem:copulegaussienne_geq_r-second}]

Let $\ddom_1 = \{(p,x) \in [0,1]^2 \, \colon \, 0 \leq x < p < 1 \}$. Notice that for any $(p,x) \in \ddom_1$ we have $p-x \in (0,1)$. We define the function $\rmf_1$ on $\ddom_1$ as follows: 
\[
\rmf_1(p,x) = \Phi^{-1}(p-x)^2  - \Phi^{-1}(p)^2 - 2 \log \left( \frac{p(1-p)}{(p+x)(1-p+x)} \right).
\]
The function is differentiable w.r.t.\ $x\in [0, p)$ and we have:
\[
\frac{1}{2} \, \partial_x \rmf_1(p,x) = \frac{1}{p+x} + \frac{1}{1-p+x} - \frac{\Phi^{-1}(p-x)}{\phi(\Phi^{-1}(p-x))}\cdot
\]
If $p-x \in (0,1/2]$, then $\Phi^{-1}(p-x)$ is nonpositive and $\partial_x \rmf_1(p,x) >0$. If $p-x \in (1/2,1)$, noting that for $y \in (0,+\infty)$ we have:
\[
1-\Phi(y) - \frac{\phi(y)}{y} = - \int_y^{+\infty} \frac{\phi(t)}{t^2} dt < 0,
\]
we deduce with $y = \Phi^{-1}(p-x)$ that:
\[
\frac{1}{2} \, \partial_x \rmf_1(p,x) =\frac{1}{p+x} +  \frac{1}{1-p+x} - \frac{\Phi^{-1}(p-x)}{\phi(\Phi^{-1}(p-x))} > \frac{1}{p+x}>0.
\]
Therefore, the function $x \mapsto \rmf_1(p,x)$ is increasing on $[0,p)$ and since $\rmf_1(p,0) = 0$, we deduce that  $\rmf_1(p,x) > 0$ for $(p,x) \in \mathring{\ddom}_1$. Using $\Phi^{-1}(p-x) = -\Phi^{-1}(1-p+x)$, elementary computations give:
\begin{equation}
    \label{eq:proof_eq_px_imp_ineg}
    (p,x) \in \mathring{\ddom}_1 \implies \frac{\phi(\Phi^{-1}(p))}{\phi(\Phi^{-1}(1-p+x))} > \frac{p(1-p)}{(p+x)(1-p+x)} \cdot
\end{equation}

We define the function $\rmf_2$ on $\ddom_2 = \{(p,q) \in [0,1]^2 \, \colon \, 1/2 < q \leq p < 1 \}$ as follow:
\[
\rmf_2(p,q) = 2 \Phi^{-1}(q) - \Phi^{-1}(p) - \Phi^{-1}(p'),
\]
with $p'$ defined in~\eqref{eq:difference_copule_rsecond}. Notice that since $p \in [q,1)$, we have $p' \in (0,q]$. The function $\rmf_2$ is differentiable w.r.t.\ $p$ on $[q,1)$ and we have:
\[
\partial_p \rmf_2(p,q) = - \frac{1}{\phi(\Phi^{-1}(p))} + \frac{p'+2p-1}{p-2q+1} \frac{1}{\phi(\Phi^{-1}(p'))} \cdot
\]
We deduce from~\eqref{eq:proof_eq_px_imp_ineg} with $x=p'+p-1$ (and thus $0<x<p<1$) that :
\[
\frac{\phi(\Phi^{-1}(p))}{\phi(\Phi^{-1}(p'))} > \frac{p-2q+1}{p'+2p-1} \quad \text{on $\ddom_2$,}
\]
and thus $\partial_p \rmf_2(p,q) > 0$ on $\ddom_2$. Since $\rmf_2(q,q) = 0$, we deduce that:
\begin{equation}
    \label{eq:proof_eq_pq_imp_ineg}
    (p,q) \in \mathring{\ddom}_2 \implies \rmf_2(p,q) >0.
\end{equation}

We define the function $\rmf_3$ on $\ddom_3 =  \mathring{\ddom}_2 \times [-1,1]$ by:
\[
\rmf_3(p,q,\rho) = (\rho^2-1) \log\left( \frac{p'(1-q)}{(1-p)q} \right) + 2 \left(\Phi^{-1}(p) - \Phi^{-1}(p') \right) \left(\Phi^{-1}(p) + \Phi^{-1}(p') -2 \rho \Phi^{-1}(q)  \right).
\]
For $(p,q) \in \mathring{\ddom}_2$, we have $(2q-1)(p-q) > 0$ and thus $p'(1-q)/((1-p)q) < 1$, so its logarithm is negative and the function $\rho \in [-1,1] \mapsto \rmf_3(p,q,\rho)$ is a quadratic polynomial with a negative leading coefficient. Since $q > 1/2$, we have $\Phi^{-1}(q) > 0$. Moreover, since $\Phi^{-1}$ is increasing and $1-p < p' < p$, it follows that $\Phi^{-1}(1-p) < \Phi^{-1}(p') < \Phi^{-1}(p)$ and thus $\Phi^{-1}(p) + \Phi^{-1}(p') + 2 \Phi^{-1}(q) > 0$. Furthermore, according to~\eqref{eq:proof_eq_pq_imp_ineg}, we have $\Phi^{-1}(p) + \Phi^{-1}(p') - 2\Phi^{-1}(q) < 0$. Therefore, we have:
\begin{align*}
  &\rmf_3(p,q,-1) = 2 \left(\Phi^{-1}(p) - \Phi^{-1}(p') \right) \left(\Phi^{-1}(p) + \Phi^{-1}(p') + 2 \Phi^{-1}(q) \right) > 0, \\
  %&h_3(p,q,0) = - \log\left( \frac{p'(1-q)}{(1-p)q} \right) + 2 \left(\Phi^{-1}(p) - \Phi^{-1}(p') \right) \left(\Phi^{-1}(p) + \Phi^{-1}(p') \right) > 0, \\
  &\rmf_3(p,q,1) = 2 \left(\Phi^{-1}(p) - \Phi^{-1}(p') \right) \left(\Phi^{-1}(p) + \Phi^{-1}(p') -2 \Phi^{-1}(q) \right) < 0.
\end{align*}
The function $\rmf_3$ being a quadratic polynomial in $\rho$, there exists $\rho_0 \in (-1,1)$ such that $\sgn(\rmf_3(p,q,\rho)) = \sgn(\rho_0-\rho)$ for $(p,q,\rho) \in \ddom_3$.

Recall the function $\rmf$ is defined on $\ddom_0 = \ddom_2 \times [-1,1]$ by \eqref{eq:difference_copule_rsecond}. Let $(p,q) \in \mathring{\ddom}_2$. The function $\rho \mapsto \rmf(p,q,\rho)$ is differentiable on $(-1,1)$ and using~\eqref{eq:deriv_rho_copule_gauss} we have:
\[
\partial_\rho \rmf(p,q,\rho) = (1-p) \phi_\rho(\Phi^{-1}(p'),\Phi^{-1}(q)) - \frac{p'(1-q)}{(1-p)} \phi_\rho(\Phi^{-1}(p),\Phi^{-1}(q)) \cdot
\]
Using the expression of the Gaussian density function, we get $\sgn(\partial_\rho \rmf(p,q,\rho)) = \sgn(\rmf_3(p,q,\rho)) =  \sgn(\rho_0-\rho)$. The function $\rho \mapsto \rmf(p,q,\rho)$ is therefore increasing on $[-1,\rho_0]$ and decreasing on $[\rho_0,1]$ for $(p,q) \in \mathring{\ddom}_2$. Since $\rmf(p,q,0)=\rmf(p,q,1) = 0$, we obtain $\sgn(\rmf(p,q,\rho))=\sgn(\rho)$ for $\rho \in [-1,1)$ and $(p,q) \in \mathring{\ddom}_2$.

To conclude, it is easy to check that if $\rho = 1$ or $q = 1/2$ or $p \in \{q, 1\}$ (the two latter conditions correspond to $(p,q) \in \ddom_2 \setminus \mathring{\ddom}_2$), then we have $\rmf(p,q,\rho) = 0$. This ends the proof. 
\end{proof}

The following lemma provides a way to identify points that are comparable to $A$ for the same correlation $\rho$, namely points belonging to $\Upsilon_{q,\rho} \cap (C_{q,A}^+ \cup C_{q,A}^-)$.

\begin{lem}
\label{lem:copule_notin_CpCm}
Let $q, p \in (0,1)$, $\rho \in (0, 1)$ and $A = (p,r)$ with $r = C_\rho(p,q)$. We have:
\begin{equation}
  \label{eq:copule_notin_CpCm}
\Upsilon_{q,\rho}^* \cap \left(C_{q,A}^+ \cup C_{q,A}^-\right) 
= \Upsilon_{q,\rho} \cap C_{q,A}^+ \cap C_{q,A}^- 
= D_q^+\cap C_{q,A}^+ \cap C_{q,A}^-.
\end{equation}
\end{lem}

\begin{rem}
\label{rem:descriptionCpCapCmThresholded}
For $q,p \in (0,1)$ and $\rho \in [-1,1] \setminus \{0\}$, we have $r = C_\rho(p,q) \neq pq$, that is $A = (p,r) \notin \Delta_q$. According to Proposition~\ref{prop:C+-cap-C-}, we have:
\begin{enumerate}[(i)]
    \item \label{it:RemdescriptionCpCapCmThresholded1}
    If $q\neq 1/2$ or $p=q=1/2$, then $D_q^+\cap C^+_{q, A}\cap C^-_{q, A}$ is reduced to $\{A\}$ if  $\rho\in (0,1]$ and to $\{h_q(A)\}$ if $\rho \in [-1, 0)$ and $A \in \Upsilon_{q,\rho}$.
    \item \label{it:RemdescriptionCpCapCmThresholded2}
    If $p\neq q= 1/2$, then $D_q^+\cap C^+_{q, A}\cap C^-_{q, A}$ is reduced to $\{A, A''_+\}$ if $\rho \in (0,1]$ and to $\{h_q(A), h_q(A''_+)\}$ if $\rho \in [-1,0)$ and $\{A, A''_+\} \subset \Upsilon_{q,\rho}$ (where $A''_+$ is defined in~\eqref{eq:def_A''+}).
\end{enumerate}
\end{rem}

\begin{proof}[Proof of Lemma~\ref{lem:copule_notin_CpCm}]
Let $p,q \in (0,1)$ and $\rho \in (-1,1) \setminus \{0\}$.
We first note some properties of the function $t \mapsto C_{\rho}(t,q)$. According to~\eqref{eq:deriv_p_copule_gauss} the function $t \in [0,1] \mapsto C_\rho(t,q)$ is increasing for $\rho \in (-1,1)$. Moreover, for $t \in (0,1)$, we have:
\[
\partial_t^2 C_\rho(t,q) = -\frac{\rho}{1-\rho^2} \frac{\phi_\rho(\Phi^{-1}(t), \Phi^{-1}(q))}{\phi( \Phi^{-1}(t))^2} \cdot
\]
Thus, the function $t \in [0,1] \mapsto C_\rho(t,q)$ is strictly convex for $\rho \in (-1,0)$ and strictly concave for $\rho \in (0,1)$. 

\medskip
Let $\rho \in (-1,1) \setminus \{0\}$, $r = C_{\rho}(p,q)$ and define $A = (p,r)$. Let us recall that $I_q = [\beta, 1-\beta]$, with $\beta = \min(q, 1-q)$. 
We first show that:
\begin{equation}
  \label{eq:copule_notin_Cp} 
  \Upsilon_{q,\rho}^* \cap C_{q,A}^+ = D_q^{\sgn(\rho)} \cap C_{q,A}^+ \cap C_{q,A}^-. 
\end{equation}
It suffices to show that $\Upsilon_{q,\rho}^* \cap C_{q,A}^+ \subset C_{q,A}^+ \cap C_{q,A}^-$, since the converse inclusion is immediate by Remark~\ref{rem:descriptionCpCapCmThresholded}. 
We recall Proposition~\ref{prop:bernoulli_cp} which states that if $A \notin \Delta_q$ (that is $\rho \neq 0$) then $C^+_{q,A}$ is the convex hull of either $4$ or $6$ points, corresponding respectively to the cases $p \in I_q$ and $p \notin I_q$. 

\medskip 
We first consider $\rho \in (0,1)$ and thus $\Delta_q \leq \Upsilon_{q,\rho}$ and $A \in D_q^+$. 
We distinguish the following two cases: $q \in (0,1)$ with $p \in I_q$; and $q \in [1/2,1)$ with $p \in (q,1)$. The remaining cases are treated by symmetry.

\medskip

\textbf{Case 1: $\rho \in (0,1)$, $q\in (0,1)$ and $p \in I_q$.} According to Proposition~\ref{prop:bernoulli_cp}, the set $C^+_{q,A} \cap D_q^+$ is the convex hull of $\{A_0, A, h_q(A_0)\}$. The point $(t, C_{\rho}(t,q))$ coincides with $A_0$, $A$, and $h_q(A_0)$ for $t=0$, $t=p$, and $t=1$, respectively.  
% For $\rho \in (0,1)$, we have for all $t \in (0,1)$ that $C_{\rho}(t,q) > tq$. Thus, the points $(t, C_{\rho}(t,q))$ and $h_q(A)$ lie on opposite sides of the line passing through $A_0$ and $h_q(A_0)$. 
Since the function $t \mapsto C_{\rho}(t,q)$ is strictly concave on $[0,1] $, it follows that $\Upsilon_{q,\rho} \cap C_{q,A}^+$ is reduced to $\{A_0, A, h_q(A_0)\}$, thus Equation~\eqref{eq:copule_notin_Cp} holds in this case.

\medskip

\textbf{Case 2: $\rho \in (0,1)$, $q \in [1/2,1)$ and $p \in (q,1)$.} 
According to Proposition~\ref{prop:bernoulli_cp}, the set $C^+_{q,A} \cap D_q^+$ is the convex hull of  $\{A_0, A, A', h_q(A_0)\}$, where $A' = (p',r')$ with:
\[
p' = \frac{p(1-p)}{p-2q+1} \quad \text{and} \quad r' = qp' + p' (r-qp) \frac{1-q}{(1-p)q} \cdot
\]
The point $(t, C_{\rho}(t,q))$ coincides with $A_0$, $A$, and $h_q(A_0)$ for $t=0$, $t=p$, and $t=1$, respectively. According to Lemma~\ref{lem:copulegaussienne_geq_r-second}, for $q > 1/2$, we have:
\begin{equation}
    \label{eq:proof_copule_notin_CpCm_diff}
    (1-p) \left(C_\rho(p',q) - qp' \right) - \frac{p'(1-q)}{q} \left(r - qp \right) > 0,
\end{equation}
that is $C_\rho(p',q) > r'$ which is equivalent to $A' < \Upsilon_{q,\rho}$. Since the function $t \mapsto C_{\rho}(t,q)$ is strictly concave, it follows that $\Upsilon_{q,\rho} \cap C_{q,A}^+$ is reduced to $\{A_0, A, h_q(A_0)\}$ and thus Equation~\eqref{eq:copule_notin_Cp} holds for $q > 1/2$ thanks to Remark~\ref{rem:descriptionCpCapCmThresholded}~\ref{it:RemdescriptionCpCapCmThresholded1}. 

Moreover, for $q = 1/2$, we have $p' = 1-p \neq q$ and according to Lemma~\ref{lem:copulegaussienne_geq_r-second}, Equation~\eqref{eq:proof_copule_notin_CpCm_diff} is an equality.  Since the function $t \mapsto C_{\rho}(t,q)$ is strictly concave, it follows that $\Upsilon_{q,\rho} \cap C_{q,A}^+$ is reduced to $\{A_0, A, A', h_q(A_0)\}$ and thus Equation~\eqref{eq:copule_notin_Cp} holds for $q = 1/2$ thanks to Remark~\ref{rem:descriptionCpCapCmThresholded}~\ref{it:RemdescriptionCpCapCmThresholded2}.

\medskip 

\textbf{Proof of~\eqref{eq:copule_notin_Cp} in the general case.} Cases $1$ and $2$ can be treated exactly as above when $\rho \in (-1,0)$. Therefore, Equation~\eqref{eq:copule_notin_Cp} holds in Cases $1$ and $2$ for all $\rho \in (-1,1)\setminus\{0\}$.

We now consider the case where $q \geq 1/2$ and $p < 1-q$. Using the previous subcase, and the fact that $C^+_{q,A} =  C^+_{q,h_q(A)} = h_q(C_{q,A}^+)$ and $h_q(\Upsilon_{q,\rho}) = \Upsilon_{q,-\rho}$, Equation~\eqref{eq:copule_notin_Cp} holds for $\rho \in (-1,1) \setminus \{0\}$, $q \geq 1/2$ and $p < 1-q$ and thus for $p \in (0,1)$.

We now suppose that $q < 1/2$, $p \in (0,1)$ and $\rho \in (-1,1) \setminus \{0\}$. We introduce a transformation which reduce to the case $q > 1/2$. Let $U$ and $V$ be two Bernoulli random variables defined as follows:
\[
\P(U = 1) = p, \quad \P(V = 1) = q \quad \text{and} \quad \P(U = 1, V = 1) = C_\rho(p,q) .
\]
Let $t \in (0,1) \setminus \{p\}$ and let $U'$ be  defined in the same way as $U$, with $p$ replaced by $t$. 
Let $\tilde V = 1-V$, $A = (p,r)$ and $\tilde A = (p, p-r)$ with 
$r=C_\rho(p,q) $. 
We have $\tilde q = \P(\tilde V = 1) = 1-q > 1/2$. Since $t \neq p$ and $q \neq 1/2$ and using Propostion~\ref{prop:C+-cap-C-} we get $(t, C_{-\rho}(t,\tilde q)) \notin D_{\tilde q}^{\sgn(-\rho)} \cap C_{\tilde q, \tilde A}^+ \cap C_{\tilde q, \tilde A}^-$. Since~\eqref{eq:copule_notin_Cp} holds for $q > 1/2$, $p \in (0,1)$ and $\rho \in (-1,1) \setminus \{0\}$, we use it with $(A, q, \rho)$ replaced by $(\tilde A,\tilde q,-\rho)$ and get that $(t, C_{-\rho}(t,\tilde q)) \notin \Upsilon_{\tilde q, \tilde A} \cap C^+_{\tilde q,\tilde A}$. Since $(t, C_{-\rho}(t,\tilde q)) \in \Upsilon_{\tilde q, - \rho}$, we get that $(t, C_{-\rho}(t,\tilde q)) \notin C^+_{\tilde q,\tilde A}$ and there exists a function $f$ belonging to $\F$ such that:
\[
S_{f}(U, \tilde V) < S_{f}(U', \tilde V).
\]
Using Proposition~\ref{prop:invariance_csiszar_index}, we have $S_f(U, V) = S_f(U, \tilde V)$ and $S_f(U', V) = S_f(U', \tilde V)$. We get that $(t, C_{\rho}(t,q)) \notin C^+_{q,A}$. We get that Equation~\eqref{eq:copule_notin_Cp} holds for $q < 1/2$, $p \in (0,1)$ and $\rho \in (-1,1) \setminus \{0\}$. 
We get finally get that Equation~\eqref{eq:copule_notin_Cp} holds for $p,q \in (0,1)$ and $\rho \in (-1,1) \setminus \{0\}$.

\medskip

\textbf{Proof of~\eqref{eq:copule_notin_CpCm}.} Let $p,q \in (0,1)$, $\rho \in (-1,1) \setminus \{0\}$ and define $A = (p, r)$ with $r = C_\rho(p,q)$. Let $A^* = (p^*, r^*)$ with $r^* = C_\rho(p^*, q))$ such that $A^* \in \Upsilon_{q,\rho}^* \cap (C_{q,A}^+ \cup C_{q,A}^-)$. 
If $A^*$ belongs to $C_{q,A}^+$, using~\eqref{eq:copule_notin_Cp} we deduce that $A^* \in D_q^{\sgn(\rho)} \cap C_{q,A}^+ \cap C_{q,A}^-$. 
If $A^*$ belongs to $C_{q,A}^-$, then $A$ belongs to $C_{q,A^*}^+$. Using~\eqref{eq:copule_notin_Cp} with $A$ replaced by $A^*$, we get that $A \in D_q^{\sgn(\rho)} \cap C_{q,A^*}^+ \cap C_{q,A^*}^-$, and $A^* \in D_q^{\sgn(\rho)} \cap C_{q,A}^+ \cap C_{q,A}^-$ since we have $A \in C_{q,A^*}^+ \cap C_{q,A^*}^-$ if and only if $A^* \in C_{q,A}^+ \cap C_{q,A}^-$. 
% $A = A^*$ if $q \neq 1/2$ and $A \in \{A^*, (1-p^*, C_\rho(1-p^*, q)) \}$ otherwise. If $A = A^*$, we get that $A^* \in D_q^{\sgn(\rho)} \cap C_{q,A}^+ \cap C_{q,A}^-$. If $A = (1-p^*, C_\rho(1-p^*, q))$, $q = 1/2$. Using the radial symmetry of the Gaussian copula (that is $C_\rho(1-u, 1-v) = u + v - 1 + C_\rho(u,v)$ for $u,v \in [0,1]$) we get that $A^* = (1-p, C_\rho(1-p, q))$ and $A^* \in D_q^{\sgn(\rho)} \cap C_{q,A}^+ \cap C_{q,A}^-$. 
Thus we get that $\Upsilon_{q,\rho}^* \cap (C_{q,A}^+ \cup C_{q,A}^-) \subset D_q^{\sgn(\rho)} \cap C_{q,A}^+ \cap C_{q,A}^-$. The reverse inclusion is immediate by Remark~\ref{rem:descriptionCpCapCmThresholded}, and therefore we obtain that:
\[
\Upsilon_{q,\rho}^* \cap \left(C_{q,A}^+ \cup C_{q,A}^-\right) 
= D_q^{\sgn(\rho)} \cap C_{q,A}^+ \cap C_{q,A}^- .
\]
Moreover, using Remark~\ref{rem:descriptionCpCapCmThresholded}, we get that $D_q^{\sgn(\rho)} \cap C_{q,A}^+ \cap C_{q,A}^- = \Upsilon_{q,\rho} \cap C_{q,A}^+ \cap C_{q,A}^-$ and thus~\eqref{eq:copule_notin_CpCm} holds. 
\end{proof}

\subsection{Proof of Proposition~\ref{prop:copule_thresholded}}
\label{sec:proof_copule_thresholded}

If $\rho = 0$, then $C_{\rho}$ reduces to the independent copula. In this case, we have $\Upsilon_{\rho,q} = \Delta_q$. Moreover, according to Proposition~\ref{prop:bernoulli_cp}, the set $\Delta_q$ is always contained in $C_{q,A'}^+$. This proves Point~\ref{it:Up-0}.

\medskip
We now establish Points~\ref{it:Up-q} and~\ref{it:Up-1/2}. Assume that $\rho = \rho'$. Since $\rho' = \rho > 0$, using Proposition~\ref{prop:C+-cap-C-} or Remark~\ref{rem:descriptionCpCapCmThresholded}
we deduce from Lemma~\ref{lem:copule_notin_CpCm} with $A$ replaced by $A'$, that $\Upsilon^*_{q,\rho} \cap C^+_{q,A'} = \Upsilon_{q,\rho} \cap C^-_{q,A'}$ is reduced to $\{A'\}$ if $q \neq 1/2$ or $q=p'=1/2$ and to $\{ A', A_+''\}$ if $q=1/2$ and $p'\neq q$.

\medskip

We now prove Point~\ref{it:Up-0rhop}. Let $\rho \in (0,\rho')$. We first show that $\Upsilon_{q,\rho}^* \cap C^+_{q,A'} \neq \emptyset$. According to Proposition~\ref{prop:bernoulli_cp}, it follows that for all $r \in [qp', r']$, the point $(p',r)$ belongs to $C^+_{q,A'}$. In particular, this holds for $r = C_{\rho}(p',q)$ since $\Upsilon_{q,\rho}^* < \Upsilon_{q,\rho'}^*$ by~\eqref{eq:monotony_upsilon_rho}. Hence, $\Upsilon_{q,\rho}^* \cap C^+_{q,A'} \neq \emptyset$. 

We now prove $\Upsilon_{q,\rho}^* \cap  (C^+_{q,A'})^c\neq \emptyset$.
According to Proposition~\ref{prop:bernoulli_cp}, the set $C^+_{q,A'}$ is the convex hull of either four or six points. Among its boundary segments, one connects $A'$ either to $A_0$ or to $h_q(A_0)$. 
We first consider the case where $A'$ is connected to $A_0$. We denote by $\Delta'$ the line passing through $A_0$ and $A'$ so that $C_{q,A'}^+ \leq \Delta'$. Since $\rho' \in (0,1)$, we have $A' \notin \partial D_q$, and therefore the slope of $\Delta'$ is strictly smaller than $1$. According to~\eqref{eq:deriv_p_copule_gauss}, 
% the tangent at $0$ to 
the function $p \in [0,1] \mapsto C_{\rho}(p,q)$ has slope equal to $1$ at $p=0$. 
We deduce that $\Upsilon_{q,\rho} \cap (C^+_{q,A'})^c \neq \emptyset$.
In the second case, we denote by $\Delta'$ the line passing through $A'$ and $h_q(A_0)$ so that $C_{q,A'}^+ \leq \Delta'$. Since $A' \notin \partial D_q$, the slope of this line is positive. According to~\eqref{eq:deriv_p_copule_gauss}, the function $p \in [0,1] \mapsto C_{\rho}(p,q)$ has slope equal to $0$ at $p=1$. 
We deduce that $\Upsilon_{q,\rho} \cap (C^+_{q,A'})^c \neq \emptyset$.
In both cases, we conclude that $\Upsilon_{q,\rho} \cap (C^+_{q,A'})^c \neq \emptyset$.

We now prove $\Upsilon_{q,\rho} \cap C^-_{q,A'}= \emptyset$. Since $\rho \in (0,1)$, 
we have $\Delta_q \leq \Upsilon_{q,\rho}$. So we only need to check that $\Upsilon_{q,\rho} \cap C^-_{q,A'} \cap D_q^+ \neq \emptyset$. 
% it is easy to see that $C^-_{q,A'} \cap D_q^- < \Delta_q < C^-_{q,A'} \cap D_q^+$ and $\Delta_q \leq \Upsilon_{q,\rho}$ which gives that $\Upsilon_{q,\rho} \cap C^-_{q,A'} \cap D_q^- = \emptyset$. 
By Points~\ref{it:Up-q} and~\ref{it:Up-1/2}, since $q \in (0,1)$, we have $\Upsilon_{q,\rho'} \leq C^-_{q,A'} \cap D_q^+$. Moreover, since $\rho < \rho'$, we obtain $\Upsilon_{q,\rho}^* < \Upsilon_{q,\rho'}^*$ by~\eqref{eq:monotony_upsilon_rho}. This implies $\Upsilon_{q,\rho}^* \cap C^-_{q,A'} \cap D_q^+ = \emptyset$. Since $A_0$ and $h_q(A_0)$ don't belong to $C^-_{q,A'}$, we get $\Upsilon_{q,\rho} \cap C^-_{q,A'} \cap D_q^+ = \emptyset$ and finally $\Upsilon_{q,\rho} \cap C^-_{q,A'}= \emptyset$. 

\medskip

We now prove Point~\ref{it:Up-rhop1}. Let $\rho \in (\rho', 1)$. We first prove $\Upsilon_{q,\rho}^* \cap C^+_{q,A'} = \emptyset$. By Points~\ref{it:Up-q} and~\ref{it:Up-1/2}, and since $p,q \in (0,1)$, we have $C^+_{q,A'} \leq \Upsilon_{q,\rho'}$. Moreover, since $\rho' < \rho$, we obtain $\Upsilon_{q,\rho'}^* < \Upsilon_{q,\rho}^*$ by~\eqref{eq:monotony_upsilon_rho} and thus $\Upsilon_{q,\rho}^* \cap C^+_{q,A'}= \emptyset$. 

We now prove $\Upsilon_{q,\rho} \cap C^-_{q,A'} \neq \emptyset$. Using Proposition~\ref{prop:bernoulli_cm}, all points in $D_q$ lying above $A'$ belong to $C^-_{q,A'}$, that is, for all $r \in [r',\min(p',q)]$, we have $(p',r) \in C^-_{q,A'}$. 
Since $\Upsilon_{q,\rho'}^* < \Upsilon_{q,\rho}^*$, we deduce that the point $(p',r)$ with $r = C_{\rho}(p',q)$ belongs to $C^-_{q,A'}$ and thus $\Upsilon_{q,\rho}^* \cap C^-_{q,A'} \neq \emptyset$.

We now prove $\Upsilon_{q,\rho}^* \cap  (C^-_{q,A'})^c\neq \emptyset$. According to Proposition~\ref{prop:C+-cup-C-}, since $A' \notin \Delta_q$ we have $A_0 \notin C_{q,A'}^-$. Moreover, since $C_{q,A'}^-$ is a closed set, there exists a closed ball $B_0$, centered at $A_0$ with positive radius, such that $B_0 \cap C_{q,A'}^- = \emptyset$. By continuity of the function $p \in [0,1] \mapsto C_{\rho}(p,q)$, the path $\Upsilon_{q,\rho}$ is connected and since it contains $A_0$, we deduce that $\Upsilon_{q,\rho} \cap  (C^-_{q,A'})^c\neq \emptyset$. Since $A_0$ and $h_q(A_0)$ don't belong to $\Upsilon_{q,\rho}$, we get $\Upsilon_{q,\rho} \cap (C^-_{q,A'})^c\neq \emptyset$.

\medskip

The proof of Point~\ref{it:Up-1} is immediate.
This conclude the proof of Proposition~\ref{prop:copule_thresholded}.

\subsection{Proof of Corollary~\ref{cor:reverse_implication_thresholded}}
\label{sec:proof_reverse_implication_thresholded}

We start the proof by considering the cases where one of the random variables $X$, $X'$ or $Y$ is constant, and then we consider the case where all three random variables have positive variance.

If $\Var(Y) = 0$, then $Y$ is constant; and $X$ and $X'$ are independent of $Y$. Then, for all $\alpha, \alpha', \beta \in \R$, $\ind_{\{X \leq \alpha\}}$ and $\ind_{\{X' \leq \alpha'\}}$ are independent of $\ind_{\{Y \leq \beta\}}$, and thus~\eqref{eq:corr_reverseimp_eq} and~\eqref{eq:corr_reverseimp_ineg} hold. 

If $\Var(X) = 0$, then $X$ is constant; and $X$ and $Y$ are independent. If there exists $\alpha, \alpha', \beta \in \R$ such that $\ind_{\{X \leq \alpha\}} \cseq \ind_{\{X' \leq \alpha'\}}$ w.r.t.\ $\ind_{\{Y \leq \beta\}}$, then $\ind_{\{X' \leq \alpha'\}}$ is independent of $\ind_{\{Y \leq \beta\}}$, and thus $X'$ and $Y$ are independent. This gives~\eqref{eq:corr_reverseimp_eq}. Similarly, if there exists $\alpha, \alpha', \beta \in \R$ such that $\ind_{\{X \leq \alpha\}} \csls \ind_{\{X' \leq \alpha'\}}$ w.r.t.\ $\ind_{\{Y \leq \beta\}}$, then $\ind_{\{X' \leq \alpha'\}}$ is not independent of $\ind_{\{Y \leq \beta\}}$, and thus $X'$ and $Y$ are not independent. This gives~\eqref{eq:corr_reverseimp_ineg}. Same arguments hold if $\Var(X') = 0$.

We now suppose that $\Var(X) \Var(X')\Var(Y) > 0$ and we denote by $\rho$ (resp.\ $\rho'$) the correlation between $X$ (resp.\ $X'$) and $Y$.
According to Proposition~\ref{prop:invariance_csiszar_index}, any bijective transformation of a random variable leaves the Csiszár index unchanged. Therefore, without loss of generality, we can assume that  $(X,X',Y)$ is a centered Gaussian vector with unit variances, $\rho \geq 0$ and $\rho' \geq 0$.
Suppose there exist $\alpha, \alpha', \beta \in \R$ such that $\ind_{\{X \leq \alpha\}} \cseq \ind_{\{X' \leq \alpha'\}}$ w.r.t.\ $\ind_{\{Y \leq \beta\}}$ and let:
\[
p = \P(X \leq \alpha), \quad q = \P(Y \leq \beta), \quad \text{and} \quad r = \P(X \leq \alpha, Y \leq \beta),
\] 
and $A' = (p', r')$ defined as $A = (p,r)$ with $X$ replaced by $X'$.
This implies that $A$ belongs to $\Upsilon_{q,\rho'}^*$ by Lemma~\ref{lem:copule_notin_CpCm} for $\rho' >0$ and the fact that $\Upsilon_{q,\rho'}=\Delta_q$ for $\rho'=0$, and thus $\Upsilon_{q,\rho}^*= \Upsilon_{q,\rho'}^*$ as the sequence $t\mapsto \Upsilon_{q,t}^*$ is increasing. 
This gives $\rho=\rho'$, that is $X \cseq X'$ w.r.t.\ $Y$ and Equation~\eqref{eq:corr_reverseimp_eq} is proved.

Similarly, suppose there exist $\alpha, \alpha', \beta \in \R$ such that $\ind_{\{X \leq \alpha\}} \csls \ind_{\{X' \leq \alpha'\}}$ w.r.t.\ $\ind_{\{Y \leq \beta\}}$ and let $A$, $A'$ and $q$ be defined as above. 
This implies that $\{A\} < \Upsilon_{q,\rho'}^*$ by Lemma~\ref{lem:copule_notin_CpCm} for $\rho' >0$ and the fact that $\Upsilon_{q,\rho'}=\Delta_q$ for $\rho'=0$, and thus $\Upsilon_{q,\rho}^* < \Upsilon_{q,\rho'}^*$ as the sequence $t\mapsto \Upsilon_{q,t}^*$ is increasing. 
This gives $\rho < \rho'$, that is $X \csls X'$ w.r.t.\ $Y$ and Equation~\eqref{eq:corr_reverseimp_ineg} is proved.

\subsection{Proof of Theorem~\ref{thm:thresholded}}
\label{sec:proof_thresholded}

By Corollary~\ref{cor:reverse_implication_thresholded}, it suffices to prove the implication of~\eqref{eq:thm_thresholded_eq} and ~\eqref{eq:thm_thresholded_ineg}. To this end, it is enough to exhibit suitable values of $\alpha$, $\alpha'$, and $\beta$. Such choices have already been discussed in Remark~\ref{rem:expressions_aab_thresholded}.

Let $(X,X',Y)\in \R^3$ be a Gaussian random vector. Let $\rho$ (resp.\ $\rho'$) be the correlation between $X$ (resp.\ $X'$) and $Y$ with the convention that $\rho = 0$ (resp.\ $\rho' = 0$) if $\Var(X) \Var(Y) = 0$ (resp.\ $\Var(X') \Var(Y) = 0$). According to Proposition~\ref{prop:invariance_csiszar_index}, any bijective transformation of a random variable leaves the Csiszár index unchanged. Therefore, without loss of generality, we can assume that  $(X,X',Y)$ is a centered Gaussian vector with zero or unit variances, $\rho \geq 0$ and $\rho' \geq 0$. 

We divide the proof into three cases. The first case is when $\rho = 0$ or $\rho' = 0$ (that also includes the case where $X$, $X'$, or $Y$ has zero variance). The second case is when $\rho = 1$ or $\rho' = 1$. The third case is when $\rho$ and $\rho'$ belong to $(0,1)$.

\medskip

\textbf{Case $\rho=0$ or $\rho'=0$:} Note that is case also includes the case where $X$, $X'$, or $Y$ has zero variance.
It is easy to see that any choice of real $\alpha$, $\alpha'$, and $\beta$ is valid. Thus, using Corollary~\ref{cor:reverse_implication_thresholded}, the equivalences in~\eqref{eq:thm_thresholded_eq} and~\eqref{eq:thm_thresholded_ineg} hold.

\medskip

\textbf{Case $\rho=1$ or $\rho'=1$:} 
We first prove the implication in~\eqref{eq:thm_thresholded_eq}. It is easy to see that any choice of real $\alpha$, $\alpha'$, and $\beta$ such that $\alpha = \alpha'$ is valid to get the implication in~\eqref{eq:thm_thresholded_eq}. 
We now prove the implication in~\eqref{eq:thm_thresholded_ineg}. 
If $\rho = 1$, then $X = Y$ a.s.\ and there no $X'$ such that $X \csls X'$ w.r.t.\ $Y$. If $\rho' = 1$, then $X' = Y$ a.s.. Moreover, if $X \csls X'$ w.r.t.\ $Y$, then any choice of real $\alpha$, $\alpha'$, and $\beta$ is valid to get the implication in~\eqref{eq:thm_thresholded_ineg}.
Thus, using Corollary~\ref{cor:reverse_implication_thresholded}, the equivalences in~\eqref{eq:thm_thresholded_eq} and~\eqref{eq:thm_thresholded_ineg} hold.

\medskip  

\textbf{Case $\rho, \rho' \in (0, 1)$.} We first prove the implication in~\eqref{eq:thm_thresholded_eq}. Suppose that $X \cseq X'$ w.r.t.\ $Y$. According to Corollary~\ref{cor:gaussian_case_dim_1}, we get that $\rho = \rho'$. Taking any real $\alpha,\alpha',\beta$ such that $\alpha = \alpha'$, we get the implication in~\eqref{eq:thm_thresholded_eq}.

We now prove the implication in~\eqref{eq:thm_thresholded_ineg}. Suppose that $X \csls X'$ w.r.t.\ $Y$. According to Corollary~\ref{cor:gaussian_case_dim_1}, we get that $\rho < \rho'$. Let $\alpha', \beta \in \R$ such that $\beta \neq 0$, and $A' = (p',r')$ and $q$ defined as:
\[
p' = \P(X' \leq \alpha'), \quad q = \P(Y \leq \beta), \quad r' = \P(X' \leq \alpha', Y \leq \beta).
\]
According to the Point~\ref{it:Up-0rhop} of Proposition~\ref{prop:copule_thresholded}, we get that $\Upsilon_{q,\rho}^* \cap C_{q,A'}^+ \neq \emptyset$. Thus, there exists $A = (p,r) \in C_{q,A'}^+$ such that $r = C_\rho(p,q)$, that is $\ind_{\{X \leq \alpha\}} \csleq \ind_{\{X' \leq \alpha'\}}$ w.r.t.\ $\ind_{\{Y \leq \beta\}}$ with $\alpha = \Phi^{-1}(p)$. 
Since $\rho < \rho'$ we get that $A \neq A'$. Using that $\beta \neq 0$ (and thus $q \neq 1/2$) and Proposition~\ref{prop:C+-cap-C-}, we get that $\ind_{\{X \leq \alpha\}}$ and $\ind_{\{X' \leq \alpha'\}}$ are not equal in the Csiszár order w.r.t.\ $\ind_{\{Y \leq \beta\}}$. Thus, we get the implication in~\eqref{eq:thm_thresholded_eq}. 

Thus, using Corollary~\ref{cor:reverse_implication_thresholded}, the equivalences in~\eqref{eq:thm_thresholded_eq} and~\eqref{eq:thm_thresholded_ineg} hold.

\end{document}